\documentclass[a4paper,10pt]{amsart}
\usepackage{amsmath,amsfonts,amssymb, amsthm, xypic}
\usepackage{epsfig,psfrag}
\usepackage[all]{xy}

\def\dim{\mathrm{dim}}

\def\Z{\mathbb{Z}}
\def\Q{\mathbb{Q}}

\def\H{\mathbf{H}}
\def\U{\mathbf{U}}
\def\qed{$\hfill \checkmark$}
\def\endexample{$\hfill \triangle$}
\def\N{\mathbb{N}}
\def\tto{\twoheadrightarrow}
\def\lto{\longrightarrow}
\def\a{\alpha}
\def\b{\beta}
\def\ga{\gamma}

\def\C{\mathbb{C}}
\def\CC{\mathbf{C}}
\def\O{\mathcal{O}}
\def\B{\mathbf{B}}
\def\x{\mathbf{x}}
\def\y{\mathbf{y}}
\def\z{\mathbf{z}}
\def\p{\mathbf{p}}

\def\A{\mathcal{A}}
\def\g{\mathfrak{g}}
\def\bo{\mathfrak{b}}
\def\h{\mathfrak{h}}
\def\n{\mathfrak{n}}
\def\LLambda{\boldsymbol{\Lambda}}
\def\v{\nu}
\def\xpl{\mathbb{X}_{\mathbf{p},\boldsymbol{\lambda}}}
\def\llambda{\boldsymbol{\lambda}}

\newtheorem{theo}{\bf{Theorem}}[section]
\newtheorem{lem}[theo]{Lemma}
\newtheorem{cor}[theo]{Corollary}
\newtheorem{prop}[theo]{Proposition}

\newtheorem{conj}[theo]{Conjecture}

\numberwithin{equation}{section}

\title{Lectures on Hall algebras}
\author{Olivier Schiffmann}

\begin{document}
\maketitle
\tableofcontents


\centerline{\large{\textbf{Introduction}}}
\addcontentsline{toc}{section}{Introduction}

\vspace{.1in}

These notes represent the written, expanded and improved version of a series of lectures given at the winter school ``Representation theory and related topics'' held at the ICTP in Trieste in January 2006, and at the summer school ''Geometric methods in representation theory'' held at Grenoble in June~2008. The topic for the lectures was ``Hall algebras'' and I have tried to give a survey of what I believe are the most fundamental results and examples in this area. The material was divided into five sections, each of which initially formed the content of (roughly) one lecture. These are, in order of appearance on the blackboard~:

\vspace{.05in}

\begin{enumerate}
\item[]  $\bullet$ Lecture 1. Definition and first properties of (Ringel-)Hall algebras,

\vspace{.05in}

\item[] $\bullet$ Lecture 2. The Jordan quiver and the classical Hall algebra,

\vspace{.05in}

\item[] $\bullet$ Lecture 3. Hall algebras of quivers and quantum groups,

\vspace{.05in}

\item[] $\bullet$ Lecture 4. Hall algebras of curves and quantum loop groups,

\vspace{.05in}

\item[] $\bullet$ Lecture 5. The Drinfeld double and Hall algebras in the derived setting.
\end{enumerate}

\vspace{.05in}

By lack of time, chalk, (and yes, competence !), I was not able to survey with the proper due respect several important results (notably Peng and Xiao's Hall Lie algebra associated to a $2$-periodic derived category \cite{PX}, Kapranov and To\"en's versions of Hall algebras for derived categories, see \cite{Kap2}, \cite{Toen}, or the recent theory of Hall algebras of cluster categories, see \cite{CalKel}, \cite{CalKeller}, or the recent use of Hall algebra techniques in counting invariants such as in Donaldson-Thomas theory, see \cite{Joyce}, \cite{KoSoi}, \cite{ReinekeDT},...). These are thus largely absent from these notes. Also missing is the whole geometric theory of Hall algebras, initiated by Lusztig \cite{Lu1}~: although crucial for some important applications of Hall algebras (such as the theory of crystal or canonical bases in quantum groups), this theory requires a rather different array of techniques (from algebraic geometry and topology) and I chose not to include it here, but in the companion survey \cite{SLectures2}.
More generally, I apologize to all those whose work deserves to appear in any reasonable survey on the topic, but is unfortunately not to be found in this one. Luckily, other texts are available, such as \cite{Ringelsurvey1}, \cite{Ringelsurvey2}, \cite{Hubery1}.  There are essentially no new results in this text.

\vspace{.1in} 

Let me now describe in a few words the subject of these notes as well as the content of the various lectures. 

\vspace{.1in}

Roughly speaking, the Hall, or Ringel-Hall algebra $\H_{\mathcal{A}}$ of a (small) abelian category $\mathcal{A}$ encodes the structure of the space of \textit{extensions} between objects in $\mathcal{A}$. In slightly more precise terms, $\H_{\mathcal{A}}$ is defined to be the $\C$-vector space with a basis consisting of symbols $\{[M]\}$, where $M$ runs through the set of isomorphism classes of objects in $\mathcal{A}$; the multiplication between two basis elements $[M]$ and $[N]$ is a linear combination of elements $[P]$, where $P$ runs through the set of extensions of $M$ by $N$ (i.e. middle terms of short exact sequences $0 \to N \to P \to M \to 0$), and the coefficient of $[P]$ in this product is obtained by counting the number of ways in which $P$ may be realized as an extension of $M$ by $N$ (see Lecture~1 for details). Of course, for this counting procedure to make sense $\mathcal{A}$ has to satisfy certain strong finiteness conditions (which are coined under the term \textit{finitary}), but there are still plenty of such abelian categories around. Another fruitful, slightly different (although equivalent) way of thinking about the Hall algebra $\H_{\mathcal{A}}$ is to consider it as the algebra of finitely supported functions on the ``moduli space'' $\mathcal{M}_{\mathcal{A}}$ of objects of $\mathcal{A}$ (which is nothing but the set of isoclasses of objects of $\mathcal{A}$, equipped with the discrete topology), endowed with a natural convolution algebra structure (this is the point of view that leads to some more geometric versions of Hall algebras, as in \cite{Lu1}, \cite{Lau}, \cite{Sinv}).

\vspace{.1in}

Thus, whether one likes to think about it in more algebraic or more geometric terms, Hall algebras provide rather subtle invariants of finitary abelian categories. Note that it is somehow the ``first order'' homological properties of the category $\mathcal{A}$ (i.e. the structure of the groups ${Ext}^1(M,N)$) which directly enter the definition of $\H_{\mathcal{A}}$, but $\mathcal{A}$ may a priori be of arbitrary (even infinite) homological dimension. However, as discovered by Green \cite{Green}, when $\mathcal{A}$ is \textit{hereditary} , i.e. of homological dimension one or less, it is possible to define a comultiplication $\Delta: \H_{\mathcal{A}} \to \H_{\mathcal{A}} \otimes \H_{\mathcal{A}}$ and, as was later realized by Xiao \cite{Xiao1}, an antipode $S: \H_{\mathcal{A}} \to \H_{\mathcal{A}}$. These three operations are all compatible and endow (after a suitable and harmless twist which we prefer to ignore in this introduction) $\H_{\mathcal{A}}$ with the structure of a Hopf algebra.  All these constructions are discussed in details in Lecture 1.

\vspace{.1in}

As the reader can well imagine, the above formalism was invented only after some motivating examples were discovered. In fact, the above construction appears in various (dis)guises in domains such as modular or $p$-adic representation theory (in the form of the functors of parabolic induction/restriction), number theory and automorphic forms (Eisenstein series for function fields), and in the theory of symmetric functions. The first occurence of the concept of a Hall algebra can probably be traced back to the early days of the twentieth century in the work of E. Steinitz (a few years before P. Hall was born) which, in modern language, deals with the case of the category $\mathcal{A}$ of abelian $p$-groups for $p$ a fixed prime number. This last example, the so-called \textit{classical} Hall algebra is of particular interest due to its close relation to several fundamental objects in mathematics such as symmetric functions (see \cite{Mac}), flag varieties and nilpotent cones. After studying in some details Steintiz's classical Hall algebra we briefly state some of the other occurences of (examples of) Hall algebras in Lecture~2. 

\vspace{.1in}

The interest for Hall algebras suddenly exploded after C. Ringel's groundbreaking discovery (\cite{Ri}) in the early 1990s that the Hall algebra $\H_{{Rep}\;\vec{Q}}$ of the category of $\mathbb{F}_q$-representations of a Dynkin quiver $\vec{Q}$ (equiped with an arbitrary orientation) provides a realization of the positive part $\U(\bo)$ of the enveloping algebra $\U(\mathfrak{g})$ of the simple complex Lie algebra $\mathfrak{g}$ associated to the same Dynkin diagram (to be more precise, one gets a quantized enveloping algebra $\U_v (\mathfrak{g})$, where the deformation parameter $v$ is related to the order $q$ of the finite field $\mathbb{F}_q$).
 
 \vspace{.1in}
 
It is also at that time that the notion of a Hall algebra associated to a finitary category was formalized (see \cite{R2}). These results were subsequently extended to arbitrary quivers in which case one gets (usually infinite-dimensional) Kac-Moody algebras, and were later completed by Green. The existence of a close relationship between the representation theory of quivers on one hand, and the structure of simple or Kac-Moody Lie algebras on the other hand was well-known since the seminal work of Gabriel, Kac and others on the classification of indecomposable representations of quivers (see \cite{Gabriel}, \cite{Kac}). Hall algebras thus provide a concrete, beautiful (and useful !) realization of this correspondence. After recalling the forerunning results of Gabriel and Kac, we state and prove Ringel's and Green's fundamental theorems in the third Lecture.

\vspace{.1in}

Apart from the categories of $\mathbb{F}_q$-representations of quivers, a large source of finitary categories of global dimension one is provided by the categories $Coh(X)$ of coherent sheaves on some smooth projective curve $X$ defined over a finite field $\mathbb{F}_q$. As pointed out by Kapranov in \cite{Kap1}, 
the Hall algebra $\H_{Coh(X)}$ may be interpreted in the context of automorphic forms over the function field of $X$. Using this interpretation, he wrote down a set of relations satisfied by $\H_{Coh(X)}$ for an arbitrary $X$ (these relations involve as a main component the zeta function of $X$). These relations turn out to determine completely $\H_{Coh(X)}$ when $X \simeq \mathbb{P}^1$ but this is most likely not true in higher genus (see \cite{SV3}, however,  for a combinatorial approach).

\vspace{.1in}

In another direction, H. Lenzing discovered in the mid-80's some important generalizations $Coh(\mathbb{X}_{p,\lambda})$ of the category $Coh(\mathbb{P}^1)$-- the so-called weighted projective lines-- which depend on the choice of points $\lambda_1, \ldots, \lambda_r \in \mathbb{P}^1$ and multiplicities $p_1, \ldots, p_r \in \N$ associated to each point (\cite{Lenzing1}).
The category $Coh(\mathbb{X}_{p,\lambda})$ is hereditary and shares many properties with the categories $Coh(X)$ of coherent sheaves on curves (not necessarily of genus zero). In fact, in good characteristics, $Coh(\mathbb{X}_{p,\lambda})$ is equivalent to the category of \textit{$G$-equivariant} coherent sheaves on some curve $Y$ acted upon by a finite group $G$, for which $Y/G \simeq \mathbb{P}^1$. The Hall algebras $\H_{Coh(\mathbb{X}_{p,\lambda})}$ are studied in \cite{SDuke}, where it is shown that they provide a realization of the positive part of quantized enveloping algebras of \textit{loop algebras} of Kac-Moody algebras. Note that these algebras are in general not Kac-Moody algebras~: for instance when $\mathbb{X}_{p,\lambda}$ is of ``genus one'' one gets the double affine, or elliptic Lie algebras $\mathcal{E}_{\g}=\mathfrak{g}[t^{\pm 1}, s^{\pm 1}] \oplus \mathbb{K}$ for a Lie algebra $\mathfrak{g}$ of type $D_4, E_6, E_7$ or $E_8$. Simultaneously, Crawley-Boevey was led in his beautiful work on the Deligne-Simpson problem \cite{CBSimpson} to study the classes of indecomposable sheaves in $Coh(\mathbb{X}_{p,\lambda})$ and found them to be related to loop algebras of Kac-Moody algebras as well (see \cite{CB2}).
The above results concerning Hall algebras of coherent sheaves on curves form the content of Lecture 4, and should be viewed as analogues, in the context of curves, of Gabriel's, Kac's and Ringel's theorems for quivers.

\vspace{.1in}

Finally in the last lecture, we state various results and conjectures regarding the behavior of Hall algebras under \textit{derived} equivalences. Recall that taking the Drinfeld double is a process which turns a Hopf algebra $\H$ into another one $\mathbf{DH}$ which is twice as big as $\H$ and which is self-dual; for instance the Drinfeld double of the positive part $\U_v(\mathfrak{b})$ of a quantized enveloping algebra is isomorphic to the whole quantized enveloping algebra $\U_v(\mathfrak{g})$. The guiding heuristic principle --which has recently been established in a wide class of cases by T. Cramer \cite{Cramer}-- is that although the Hall algebras $\H_{\mathcal{A}}$ and $\H_{\mathcal{B}}$ of two derived equivalent finitary \textit{hereditary} categories need not be isomorphic, their Drinfeld doubles $\mathbf{D}\H_{\mathcal{A}}$ and $\mathbf{D}\H_{\mathcal{B}}$ should be. More generally, any fully faithful triangulated functor  $F:D^b(\mathcal{A}) \to D^b(\mathcal{B})$ between derived categories should give rise to a homomorphism of algebras ${F}_{\star}:\mathbf{D}\H_{\mathcal{A}} \to \mathbf{D}\H_{\mathcal{B}}$. In particular, the group of autoequivalences of the derived category $D^b(\mathcal{A})$ is expected to act on $\mathbf{DH}_{\mathcal{A}}$ by algebra automorphisms. As supporting example and motivation for the above principle, we show how the group $Aut(D^b(Coh(X)))$ for an elliptic curve $X$ acts on $\mathbf{DH}_{Coh(X)}$. This action turns out to be the key point in understanding the structure of the algebra $\mathbf{DH}_{Coh(X)}$ (the \textit{elliptic Hall algebra} studied in \cite{BS}).

\vspace{.1in}

A recent theorem of Happel \cite{Happel} states that any (connected) hereditary category which is linear over an algebraically closed field $k$ and which possesses a tilting object (see Lecture 5.) is derived equivalent to either ${Rep}_{k}\vec{Q}$ for some quiver $\vec{Q}$ or $Coh(\mathbb{X}_{p,\lambda})$ for some weighted projective line $\mathbb{X}_{p,\lambda}$.  Although the case of categories which are linear over a finite field $k$ is slightly more complicated (see \cite{ReitenHappel}, and also \cite{ReitenVdB}), if one believes the above heuristic principle then the results of Lectures 3 and 4 essentially describe the Hall algebra of any finitary hereditary category which possesses a tilting object. Of course the case of finitary hereditary categories which \textit{do not} possess a tilting object (this corresponds to curves of higher genus) is still very mysterious, as is the case of categories of higher global dimension (this corresponds to higher-dimensional varieties) for which virtually nothing is known.

\vspace{.3in}

A final word concerning the style of these Lecture notes. They follow a leisurely pace and many examples are included and worked out in details. Nevertheless, because they are mostly (though not only !) aimed at people interested in representation theory of finite-dimensional algebras, I have decided to assume some basic homological algebra and, starting from Lecture 3, a little familiarity with quivers. On the other hand, I assume nothing from Lie algebras and quantum groups theory. Hence I have included in a long appendix a ``crash course'' on simple and Kac-Moody Lie algebras, loop algebras, and the corresponding quantum groups. 

The first four Lectures follow each other in a logical order, but a reader allergic to examples could well jump to Lecture~5 directly after Lecture~1.

\newpage

\centerline{\large{\textbf{Lecture~1.}}}
\addcontentsline{toc}{section}{\tocsection {} {} {Lecture~1.}}

\setcounter{section}{1}

\vspace{.2in}

The aim of this first Lecture is to introduce in as much generality as possible the notion of the Hall algebra of a finitray abelian category, and to describe in details all the extra structures (coproduct, antipode,...) which have been discovered over the time and which one can put on such an algebra. A final paragraph briefly discusses some functoriality properties of this construction. Examples of Hall algebras abound in Lectures~2,~3 and~4, and the reader is invited to have a look at them as he proceeds through this first Lecture. 

\vspace{.15in}

\centerline{\textbf{1.1. Finitary categories}}
\addcontentsline{toc}{subsection}{\tocsubsection {}{}{\; 1.1. Finitary categories}}

\vspace{.15in}

\paragraph{} A small abelian category $\mathcal{A}$ is called \textit{finitary} if the following two conditions are satisfied~:
\begin{enumerate}
\item[i)] For any two objects $M,N \in \text{Ob}(\mathcal{A})$ we have $|{Hom}(M,N)| < \infty$,
\item[ii)]For any two objects $M,N \in \text{Ob}(\mathcal{A})$ we have $|{Ext}^1(M,N)| < \infty$.
\end{enumerate}

\vspace{.05in}

In most, if not all examples of finitary categories which we will be considering in these notes, $\mathcal{A}$ is linear over some finite field $\mathbb{F}_q$, and we have 
\begin{equation}\label{E:finitary}
{dim}\;{Hom}(M,N) < \infty, \qquad {dim}\;{Ext}^1(M,N) < \infty
\end{equation}
for any pair of objects $M,N \in {Ob}(\mathcal{A})$. Examples of such categories are provided by the categories ${Rep}_{\mathbb{F}_q}\vec{Q}$ of (finite dimensional) $\mathbb{F}_q$-representations of a quiver, or more generally by the categories ${Mod}\,A$ of finite-dimensional representations of a finite-dimensional $\mathbb{F}_q$-algebra $A$. For another class of examples of a more geometric flavor, one may consider the categories $Coh(X)$ of coherent sheaves on some projective scheme defined over $\mathbb{F}_q$ (the finiteness property (\ref{E:finitary}) holds by a famous theorem of Serre, see e.g. \cite{Hart}).

We denote by $K(\mathcal{A})$ the Grothendieck group (over $\Z$) of an abelian category $\mathcal{A}$. In most situations of interest for us, this will be a free $\Z$-module.
If $\A$ is a finite length category (i.e. if any object of $\mathcal{A}$ has a finite composition sequence with simple factors) then $K(\mathcal{A})$ is freely generated by the classes of the simple objects.

\vspace{.2in}

\centerline{\textbf{1.2. Euler form and symmetric Euler form.}}
\addcontentsline{toc}{subsection}{\tocsubsection {}{}{\; 1.2. Euler form and symmetric Euler form.}}

\vspace{.15in}

\paragraph{} Let $\mathcal{A}$ be a finitary category, and let us make the additional assumptions that ${gldim}(\mathcal{A})<\infty$ and that property ii) above is satisfied for the groups $Ext^i(M,N)$ for all $i$\footnote{Here we implicitly assume that the groups $Ext^i$ are well-defined. This is the case for all examples discussed above (modules over a finite-dimensional algebra, coherent sheaves over smooth projective varieties, ...).}. For any two objects $M,N$ of $\mathcal{A}$ we put
\begin{equation}\label{E:Eulerform1}
\langle M, N \rangle_m=\left(\prod_{i=0}^\infty (\# {Ext}^i(M,N))^{(-1)^i}\right)^{\frac{1}{2}}.
\end{equation} 
Since $\mathcal{A}$ is of finite global dimension, ${Ext}^i(M,N)=\{0\}$ for $i \gg 0$ and the product is finite. 
Note that the definition of $\langle M,N\rangle$ implicitly involves a choice of a square root.
An easy application of the long exact sequences in homology associated to the functor $\text{Hom}$ shows that $\langle M, N \rangle_m$ only depend on the classes of $M$ and $N$ in the Grothendieck group and (\ref{E:Eulerform1}) thus defines a form $\langle\;,\; \rangle_m~:K(\mathcal{A}) \times K(\mathcal{A}) \to \mathbb{C}$ which is called the (square root of the) multiplicative \textit{Euler form}. It is also useful to introduce the multiplicative \textit{symmetric Euler form} $(M,N)_m=\langle M,N \rangle_m \cdot \langle N,M\rangle_m$.

\vspace{.1in}

When $\mathcal{A}$ is $k$-linear then one usually considers additive versions of the Euler forms instead, which are defined by $\langle M,N \rangle_{a}=\sum_i (-1)^i {dim}\;{Ext}^i(M,N)$ and $(M,N)_{a}=\langle M,N\rangle_{a} + \langle N,M \rangle_a$. In this very simple way, we have associated to any finitary $k$-linear category a lattice $(K(\mathcal{A}), (\;,\;)_a)$, that is a 
(usually free and finite rank) $\Z$-module equipped with a $\Z$-valued symmetric bilinear form. As we will see, this seemingly rather coarse invariant already carries a lot of information regarding $\mathcal{A}$. 

\vspace{.2in}

\centerline{\textbf{1.3. The name of the game.}} 
\addcontentsline{toc}{subsection}{\tocsubsection {}{}{\; 1.3. The name of the game.}}

\vspace{.15in}

\paragraph{}Let $\mathcal{A}$ be a finitary category. We now introduce the main character of these notes, namely the \textit{Hall algebra} $\H_{\mathcal{A}}$ of $\A$. Let $\mathcal{X}={Ob}(\mathcal{A})/\sim$ be the set of isomorphism classes of objects in $\mathcal{A}$. Consider a vector space 
$$\H_{\mathcal{A}}:=\bigoplus_{M \in \mathcal{X}} \C[M]$$
linearly spanned by symbols $[M]$, where $M$ runs through $\mathcal{X}$. We will now define a multiplication on $\H_{\mathcal{A}}$. Given any three objects $M,N,R$, let $\mathcal{P}_{M,N}^R$ denote the set of short exact sequences $0 \to N \to R \to M \to 0$, and put $\mathbf{P}_{M,N}^R=|\mathcal{P}_{M,N}^R|$. Observe that $\mathcal{P}_{M,N}^R$ is indeed finite since by assumption
${Hom}(N,R)$ and ${Hom}(R,M)$ are finite. For any object $P$, we put $a_{P}=|{Aut}(P)|$.

\begin{prop}[Ringel, \cite{R2}]\label{P:multiplicationHall} The following defines on $\H_{\mathcal{A}}$ the structure of an associative algebra ~:
\begin{equation}\label{E:multiplicationHall}
[M] \cdot [N]=\langle M,N\rangle_m \sum_R \frac{1}{a_M a_N}\mathbf{P}_{M,N}^R [R].
\end{equation}
The unit $i: \C \to \H_{\A}$ is given by $i(c)=c[0]$, where $0$ is the zero object of $\mathcal{A}$.
\end{prop}

The proof of this result will be quite easy and natural once we have reinterpreted the above definition of $\H_{\mathcal{A}}$ from a slightly more geometric perspective. We view $\mathcal{X}$ as some kind of ``moduli space of objects in $\mathcal{A}$'', and $\mathbf{H}_{\mathcal{A}}$ as the set of finitely supported functions on $\mathcal{X}$
$$\mathbf{H}_{\mathcal{A}}=\big\{f: \mathcal{X} \to \C\;|\; \text{supp}(f) \;\text{is\;finite}\;\big\}$$
by identifying the symbol $[M]$ with the characteristic function $1_M$. We claim that the product 
(\ref{E:multiplicationHall}) can be rewritten as follows~:
\begin{equation}\label{E:multiplicationHall2}
(f \cdot g) (R)=\sum_{Q \subseteq R} \langle R/Q,Q \rangle_m f(R/Q)g(Q)
\end{equation}
Indeed, by bilinearity it is enough to check that (\ref{E:multiplicationHall2}) coincides with (\ref{E:multiplicationHall}) when $f=1_M$ and $g=1_N$. This is in turn a consequence of the following Lemma~:

\begin{lem} For any three objects $M,N,R$ of $\mathcal{A}$ we have
\begin{equation}\label{E:PMNR}
\frac{1}{a_Ma_N}\mathbf{P}^R_{M,N}=\left| \big\{ L \subset R\;|\;L \simeq N \;\text{and}\; R/L \simeq M\big\}\right|.
\end{equation}
\end{lem}
\noindent
\textit{Proof.} The group ${Aut}(M)\times {Aut}(N)$ acts freely on $\mathcal{P}_{M,N}^R$, and  the quotient is canonically identified with the right-hand side of (\ref{E:PMNR}). \qed

\vspace{.2in}

\noindent
\textit{Proof of Proposition~\ref{P:multiplicationHall}.} It will be more convenient to formulate the proof using (\ref{E:multiplicationHall2}) rather than (\ref{E:multiplicationHall}). Assume that $f=1_M$ and $g=1_N$. Then the right-hand side of (\ref{E:multiplicationHall2}) is a function supported on the set $\{R\}$ of extensions of $M$ by $N$. As ${Ext}^1(M,N)$ is finite, there are only finitely many such extensions, hence $1_M\cdot 1_N$ does indeed belong to $\H_{\mathcal{A}}$. By bilinearity it follows that $f \cdot g \in \H_{\mathcal{A}}$ for any $f,g \in \H_{\mathcal{A}}$. We now prove the associativity.
Let $f,g,h \in \H_{\mathcal{A}}$ and let $M$ be an object of $\mathcal{A}$. We compute 
\begin{equation}\label{E:equalityone}
\begin{split}
\big(f \cdot (g \cdot h)\big)(M)&= \sum_{N \subseteq M} \langle M/N,N\rangle_m f(M/N) (g\cdot h)(N)\\
&=\sum_{N \subseteq M, L \subseteq N} \langle M/N,N\rangle_m\langle N/L,L\rangle_m f(M/N) g(N/L) h(L)\\
&=\sum_{L \subseteq N \subseteq M} \langle M/N,L\rangle_m\langle M/N,N/L\rangle_m\langle N/L,L\rangle_m f(M/N) g(N/L) h(L)
\end{split}
\end{equation}
where we have used the bilinearity of the multiplicative Euler form. Similarly,
\begin{equation}\label{E:equalitytwo}
\begin{split}
\big((f \cdot g) \cdot h\big)(M)&= \sum_{R \subseteq M} \langle M/R,R\rangle_m (f\cdot g)(M/R) h(R)\\
&=\sum_{R \subseteq M, S \subseteq M/R} \langle M/R,R\rangle_m\langle (M/R)/S,S\rangle_m f((M/R)/S) g(S) h(R)
\end{split}
\end{equation}
Observe that there is a natural bijection $S \mapsto S'$ between the set $\{S\;|\;S \subseteq M/R\}$ and the set $\{S'\;|\; R \subseteq S' \subseteq M\}$ satisfying $S \simeq S'/R$ and $(M/R)/S \simeq M/S'$. Using this, we may rewrite (\ref{E:equalitytwo}) as
$$\sum_{R \subseteq S' \subseteq M} \langle M/S',R\rangle_m\langle S'/R,R\rangle_m \langle M/S',S'/R\rangle_m f(M/S') g(S'/R) h(R)$$
which is none other than (\ref{E:equalityone}). This shows that the product is associative. The statement concerning the unit is obvious.\qed

\vspace{.15in}

\addtocounter{theo}{1}
\noindent \textbf{Remarks \thetheo .} i) As pointed out in the course of the proof of Proposition~\ref{P:multiplicationHall}, the product $[M] \cdot [N]$ is a linear combination of elements $[R]$, where $R$ is an \textit{extension} of $M$ by $N$. Hence, $\H_{\mathcal{A}}$ encodes the structure of the set of short exact sequences in $\mathcal{A}$. In fact, as observed by A. Hubery \cite{Hubery1} it is possible to define a Hall algebra for any \textit{exact} category satisfying the finiteness conditions of Section~1.1.

\vspace{.05in}

\noindent
ii) The algebra $\H_{\mathcal{A}}$ is naturally graded by the Grothendieck group $K(\mathcal{A})$ of $\mathcal{A}$~: we have 
$$\H_{\mathcal{A}}=\bigoplus_{\a \in K(\mathcal{A})} \H_{\A}[\a], \qquad \text{where}\qquad
 \H_{\A}[\a]=\bigoplus_{\overline{M}=\a} \C[M].$$
iii) Assume that $\A$ is a \textit{semisimple} category and let $S=\{S_i\}_{i \in I}$ be the set of simple objects. Then if $i \neq j$ one easily checks that
$$[S_i] \cdot [S_j]=[S_i \oplus S_j]=[S_j] \cdot [S_i]$$
$$[S_i] \cdot [S_i]=|End(S_i)|^{1/2}(|{End}(S_i)| +1) [S_i \oplus S_i]$$
and in fact $\H_{\A}$ is a free commutative polynomial algebra in the generators $\{[S_i]\}_{i \in I}$.

\vspace{.05in}

\noindent
iv) In the case of a semisimple category $\mathcal{A}$ as above the algebra $\H_{\A}$ is commutative. However, this is a rather rare phenomenon~:  in general the set of extensions of $M$ by $N$ and the set of extensions of $N$ by $M$ differ and the algebra $\H_{\mathcal{A}}$ is \textit{not} commutative.

\vspace{.05in}

\noindent
v) The product of more than two elements also has an interpretation in terms of number of filtrations~: if $M_1, \ldots, M_r$ are objects in $\A$ then
$$[M_1] \cdots [M_r]=\sum_{R}\;\big( \prod_{i<j} \langle M_i,M_j \rangle_m \left|\big\{ L_r \subset \cdots \subset L_1=R\;|\; L_i/L_{i+1} \simeq M_i\big\}\right| \big)[R].$$

\vspace{.1in}

The integers $\mathbf{P}_{M,N}^R$ or $\frac{1}{a_Ma_N}\mathbf{P}_{M,N}^R$ are usually called \textit{Hall numbers}. We will give many examples of computations of Hall numbers and Hall algebras in the second, third and fourth Lectures.

\vspace{.2in}

\centerline{\textbf{1.4. Green's coproduct.}}
\addcontentsline{toc}{subsection}{\tocsubsection {}{}{\; 1.4. Green's coproduct.}}

\vspace{.15in}

\paragraph{} Let $\A$ denote again a fixed finitary category. 
The multiplication in $\H_{\A}$ encodes, essentially, all the ways of putting an object $M$ on top of an object $N$. It seems natural to try to define a dual operation, which \textit{breaks up} a given object in all possible ways. This is exactly what Green's coproduct achieves. There is, however a subtle point here. Although there are (in finitary categories) only finitely many possible extensions between any two given object, there are in general \textit{infinitely} many ways of splitting any given object into two pieces. This explains why we are forced to consider certain completions of $\H_{\mathcal{A}}$ and $\H_{\mathcal{A}} \otimes \H_{\mathcal{A}}$.

We will again define the comultiplication in two (equivalent) ways, one purely algebraic and a second more geometric. For $\a,\b \in K(\A)$, set 
$$\H_{\A}[\a] \widehat{\otimes} \H_{\A}[\b]=
\prod_{\overline{M}=\a, \overline{M}=\b} \C [M] \otimes \C[N],$$
$$\H_{\A} \widehat{\otimes} \H_{\A}=\prod_{\a,\b} \H_{\A}[\a] \widehat{\otimes} \H_{\A}[\b].$$
In other words, $\H_{\A} \widehat{\otimes} \H_{\A}$ is simply the space of all formal (infinite) linear combinations $\sum_{M,N} c_{M,N} [M] \otimes [N]$. 

\vspace{.1in}

\begin{prop}[Green, \cite{Green}]\label{P:coproductHall}  The following defines on $\H_{\mathcal{A}}$ the structure of a topological coassociative coproduct~:
\begin{equation}\label{E:coproductHall1}
\Delta([R])=\sum_{M,N} \langle M,N\rangle_m \frac{1}{a_R} \mathbf{P}^R_{M,N} [M] \otimes [N],
\end{equation} with counit $\epsilon: \H_{\A} \to \C$ defined by $\epsilon([M])=\delta_{M,0}$. 
\end{prop}

Before giving a proof of this result, let us spell out the meaning of the word \textit{topological} here. Formula (\ref{E:coproductHall1}) only defines a map $\Delta: \H_{\A} \to \H_{\A} \widehat{\otimes} \H_{\A}$, and \textit{not} $\Delta: \H_{\A} \to \H_{\A} {\otimes} \H_{\A}$, as in a genuine coalgebra. Moreover, for the coassociativity to even make sense, we must check that the two maps $(\Delta \otimes 1) \circ \Delta,
(1 \otimes \Delta) \circ \Delta : \H_{\A} \to \H_{\A} \widehat{\otimes} \H_{\A} \widehat{\otimes} \H_{\A}$ are well-defined (this is not obvious since we are composing functions with values in spaces consisting of infinite sums). Luckily, this follows from the fact that $\A$ is finitary~: indeed, the only terms in $\H_\A \widehat{\otimes} \H_{\A}$ which may contribute to $[M_1] \otimes [M_2] \otimes [M_3]$ in $(\Delta \otimes 1)\circ \Delta$ (resp. in $(1 \otimes \Delta) \circ \Delta$) are of the form $[N] \otimes [M_3]$ for some extension $N$ of $M_1$ by $M_2$ (resp. of the form $[M_1] \otimes N$ for some extension $N$ of $M_2$ by $M_1$), and there are only finitely many such terms.

\vspace{.1in}

\noindent
\textit{Proof of Proposition~\ref{P:coproductHall}.} The coassociativity is easily seen to be equivalent to the following set of relations, for any quadruple of objects $(M, N, Q, R)$~:
$$\sum_S \langle S,Q \rangle_m \langle M,N\rangle_m \frac{1}{a_Sa_R} \mathbf{P}^S_{M,N} \mathbf{P}^R_{S,Q}=\sum_T \langle N,Q \rangle_m \langle M,T\rangle_m \frac{1}{a_Ta_R} \mathbf{P}^T_{N,Q} \mathbf{P}^R_{M,T},$$
which can be rewritten as
$$\sum_S  \frac{1}{a_S} \mathbf{P}^S_{M,N} \mathbf{P}^R_{S,Q}=\sum_T  \frac{1}{a_T} \mathbf{P}^T_{N,Q} \mathbf{P}^R_{M,T}.$$
This last equality expresses in fact the associativity of the multiplication in $\H_{\A}$: up to multiplication by the factor $\langle M,N \rangle_m \langle M, Q \rangle_m \langle N, Q \rangle_m \frac{1}{a_Ma_Na_Q}$, the left-hand side is the coefficient of $[R]$ in $([M] \cdot [N]) \cdot [Q]$ while the
right-hand side is the coefficient of $[R]$ in $[M] \cdot ([N] \cdot [Q])$. \qed

\vspace{.2in}

As for the product, there is a more ``geometric'' interpretation of the map $\Delta$, at least when $\A$ is of global dimension at most one. Observe that, in the same way as $\H_{\A}$ was identified with the set of finitely supported functions on $\mathcal{X}$, we may identify $\H_{\A} \widehat{\otimes} \H_{\A}$ with the set of (arbitrary) functions on $\mathcal{X} \times \mathcal{X}$.

\begin{prop}\label{P:coprophage} Assume that $gldim(\A) \leq 1$. For any $f: \mathcal{X} \to \C$ belonging to $\H_{\A}$ and for any objects $M,N$ we have
$$\Delta(f)(M,N)=\frac{\langle M,N\rangle_m^{-1}}{|{Ext}^1(M,N)|} \sum_{\xi \in {Ext}^1(M,N)} f(X_{\xi}),$$
where $X_\xi$ is the middle term of the extension of $M$ by $N$ which is associated to $\xi$.\end{prop}
\noindent
\textit{Proof.} To prove the Proposition, we have to show that for any fixed objects $M,N,R$,
\begin{equation}\label{E:ext1}
\frac{1}{| {Ext}^1(M,N)|} \left| \big\{ \xi \in {Ext}^1(M,N)\;|\; X_{\xi} \simeq R\big\}\right|=
\langle M,N\rangle_m^2 \frac{1}{a_R} \mathbf{P}^R_{M,N}.
\end{equation}
By definition, $\left|\big\{ \xi \in {Ext}^1(M,N)\;|\; X_{\xi} \simeq R\big\}\right|$ is equal to the number of equivalence classes of short exact sequences 
\begin{equation}\label{E:ab}
\xymatrix{ 0 \ar[r]&  N \ar[r]^-{a} & R \ar[r]^-{b} & M \ar[r] & 0}.
\end{equation}
The group ${Aut}(R)$ acts on the set $\mathcal{P}^R_{M,N}$ of all exact sequences as above by
$\phi \cdot (a,b)=(\phi \circ a,b \circ \phi^{-1})$ and two exact sequences $(a,b)$ and $(a',b')$
are equivalent if and only if there exists $\phi \in {Aut}(R)$ such that $(a',b')=\phi\cdot(a,b)$. We claim that the stabilizer of any short exact sequence $(a,b)$ under the above action of ${Aut}(R)$ is isomorphic to ${Hom}(M,N)$. Indeed,  $\phi\cdot (a,b)=(a,b)$ if and only if we have $\phi_{|{Im}(a)}=Id: {Im}(a) \to R$ and if the induced map $\phi': R/{Im}(a) \to R/{Im}(a)$ is the identity. This holds exactly when
$\phi \in Id \oplus {Hom}(R/{Im}(a), {Im}(a)) \simeq Id \oplus {Hom}(M,N)$. 
Therefore, we get $\left| \big\{ \xi \in {Ext}^1(M,N)\;|\; X_{\xi} \simeq R\big\}\right|=| \mathcal{P}^R_{M,N}| \frac{| {Hom}(M,N)|}{a_R}$ which (for hereditary $\mathcal{A}$) immediately yields (\ref{E:ext1}).\qed 

\vspace{.2in}

\addtocounter{theo}{1}
\noindent \textbf{Remarks \thetheo .}  Properties of the coproduct map $\Delta: \H_{\A} \to \H_{\A} \widehat{\otimes} \H_\A$ are essentially dual to those of the multiplication map $m: \H_{\A} \otimes \H_{\A} \to \H_{\A}$. Namely,

\vspace{.05in}

\noindent
i) The map $\Delta$ respects the grading by $K(\A)$, that is 
$$\Delta (\H_{\A}[\ga]) \subset \prod_{\a+\b=\ga} \H_{\A}[\a] \widehat{\otimes} \H_{\A}[\b].$$

\vspace{.05in}

\noindent
ii) For $S$ a simple object we have $\Delta([S])=[S] \otimes 1 + 1 \otimes [S]$, i.e. $[S]$ is a primitive element of $\H_{\A}$. The converse is \textit{not} true (game~: find some counterexamples in Lecture~2.).

\vspace{.05in}

\noindent
iii) As for the product, there is absolutely no reason for the coproduct to be cocommutative in general.

\vspace{.05in}

\noindent
iv) As $\Delta$ is coassociative, we may consider an iterated comultiplication $\Delta^r: \H_{\A} \to 
\H_{\A} \widehat{\otimes} \cdots \widehat{\otimes} \H_{\A}$, which has the following interpretation~:
\begin{equation*}
\begin{split}
\Delta^r([R])=\sum_{M_1, \ldots, M_r}  \big( \prod_{i<j} \langle M_i, M_j \rangle_m  \frac{a_{M_1} \cdots a_{M_r}}{a_R} \times &\\
\times \big| \big\{ L_r \subset \cdots \subset L_1=R\;|\; L_i/L_{i+1}& \simeq M_i\big\} \big| \big) [M_1] \otimes \cdots \otimes [M_n].
\end{split}
\end{equation*}

\vspace{.05in}

\noindent
v) The coproduct $\Delta$ takes values in $\H_{\A} \otimes \H_{\A}$ (and not in the completion $\H_{\A} \widehat{\otimes} \H_{\A}$) if and only if the following condition is satisfied~: 
\begin{equation}\label{E:FS}
 \textit{Any \;fixed \;object \;$R$ of $\mathcal{A}$\;has\;only\;finitely\;many\;subobjects\;$N \subset R$}.
 \end{equation}
We will call this condition the \textit{finite subobjects condition}. It holds for categories of representations of quivers (see Lecture~2, Lecture~3), but not for the categories of coherent sheaves on curves (Lecture~4).

\vspace{.15in}

Let us finish this Section with a rather useful and completely general result. For $\gamma \in K(\mathcal{A})$, let 
$$\mathbf{1}_{\gamma}=\sum_{\overline{M}=\gamma} [M]$$
be the sum of all objects in $\A$ of class $\gamma$ (this sum may be infinite for some categories, so strictly speaking $\mathbf{1}_{\gamma}$ belongs to the formal completion of $\H_{\A}$). 

\vspace{.1in}

\begin{lem}\label{L:coprodun} Assume that $gldim(\A) \leq 1$. Then we have
\begin{equation}\label{E:coprodun}
{\Delta}(\mathbf{1}_{\gamma})= \sum_{\a+\beta=\gamma} \langle \a,\beta \rangle_m^{-1} 
\mathbf{1}_{\a}
\otimes \mathbf{1}_{\b}. 
\end{equation}
\end{lem}

\noindent
\textit{Proof.} By definition, for any objects $M$ and 
$N$ such that $\overline{M} = \alpha$ and $\overline{N} = \beta$ 
the coefficient of the element $[M]  \otimes [N]$ in
${\Delta}(\mathbf{1}_{\gamma})$ is equal to $\langle \a,\b \rangle_m 
\sum_{L} \frac{1}{a_{L}} 
\mathbf{P}^{L}_{M,N}$. From 
Yoneda's description of ${Ext}^1(M,N)$ we derive 
$$|{Ext}^1(M,N)|=\sum_{L} \big|\big( \{ 0 \to N 
\stackrel{a}{\to} L \stackrel{b}{\to} 
M \to 0\}/{Aut}(L)\big)\big|,$$
where the action of $Aut(L)$ on the exact sequence $\{ 0 \to N 
\stackrel{a}{\to} L \stackrel{b}{\to} 
N \to 0\}$ is given by the rule  $g \cdot (a,b)=(ga,bg^{-1})$.
It remains to observe that for any $L$ the stabilizer of any such sequence  
is isomorphic to ${Hom}(M,N)$, and thus
$$\langle \a,\b \rangle_m \sum_{L} 
\frac{1}{a_{L}} 
\mathbf{P}^{L}_{M,N}=\langle \a,\b \rangle_m 
\frac{| {Ext}^1(M,N)|}{|{Hom}(M,N)|}=
\langle\a,\b\rangle_m^{-1},$$
proving (\ref{E:coprodun}).
\qed

\vspace{.2in}

\centerline{\textbf{1.5. The Hall bialgebra and Green's theorem.}}
\addcontentsline{toc}{subsection}{\tocsubsection {}{}{\; 1.5. The Hall bialgebra and Green's theorem.}}

\vspace{.15in}

\paragraph{} As we have seen, any finitary category $\A$ gives rise to a $\C$-vector space $\H_{\A}$ which is both an algebra and a coalgebra (the latter in a topological sense). It is natural to ask whether these two operations are compatible and endow $\H_{\A}$ with the structure of a \textit{bialgebra}. This turns out to be false in general, but essentially true for \textit{hereditary} categories (categories of global dimension at most one), as shown by Green \cite{Green}. In this paragraph we explain the (difficult) proof of this important result. \textit{We make the assumption throughout that $\A$ is hereditary}.

\vspace{.1in}

First of all, it is necessary to slightly twist the mutiplication in $\H_{\A} \otimes \H_{\A}$~: if $x, y, z, w$ are homogeneous elements in $\H_{\A}$ of respective weight $wt(x), wt(y), wt(z),$ $wt(w) \in K(\A)$, we define a new multiplication by
$$(x \otimes y) \cdot (z \otimes w)=( wt(y), wt(z))_m (xz \otimes yw)$$
(i.e. we introduce an extra factor $(wt(y), wt(z))_m$ of the symmetric Euler form when
$y$ ``jumps over'' $z$).  

Next, we need to face the fact that $\Delta$ takes values in the completion $\H_{\A} \widehat{\otimes} \H_{\A}$ rather than $\H_{\A} \otimes \H_{\A}$. The problem here is that the multiplication map
$$(\H_{\A} \otimes \H_{\A}) \otimes (\H_{\A} \otimes \H_{\A}) \to \H_{\A} \otimes \H_{\A}$$
clearly doesn't extend to the completion. We will say that a product 
$$(\sum_i a_i \otimes b_i)\cdot (\sum_j c_j \otimes d_j)$$
 of elements of $\H_{\A} \widehat{\otimes} \H_{\A}$ is \textit{convergent} if for any pair of objects $R,S$ in $\A$ the coefficient of $[R] \otimes [S]$ in
$(a_i \otimes b_i) \cdot (c_j \otimes d_j)$ is nonzero for only finitely many values of $(i,j)$. In this case of course the product is a well-defined element in $\H_{\A} \widehat{\otimes} \H_{\A}$.
Luckily for us, the product of any two elements in the \textit{image} of $\Delta$ converges, as the following Lemma shows~:

\begin{lem}\label{L:deltaextends} Let $M, N$ be objects in $\A$. Then the product $\Delta([M]) \cdot \Delta([N])$ converges.\end{lem}
\noindent
\textit{Proof.} The coefficient of $[T_1] \otimes [L_1]$ is nonzero in $\Delta([M])$ only if there exists a short exact sequence $0 \to L_1 \to M \to T_1 \to 0$. Similarly, the term $[T_2] \otimes [L_2]$ appears in $\Delta([N])$ only if there exists a sequence $0 \to L_2 \to N \to T_2\to 0$. Now let $R,S$ be objects in $\A$. The coefficient of $[S]$ in $[T_1]\cdot [T_2]$ is nonzero only when there exists a sequence
$0 \to T_2 \to S \to T_1 \to 0$, and $[R]$ appears in $[L_1] \cdot [L_2]$ only when there exists a sequence $0 \to L_2 \to R \to L_1 \to 0$. But if all these conditions are satisfied then $L_1$ is isomorphic to the image of a morphism $R \to M$ and $T_2$ is isomorphic to the image of a morphism $N \to S$. By the finitary condition on $\A$, the sets ${Hom}(R,M)$ and ${Hom}(N,S)$ are finite and in particular there are only finitely many choices for $L_1$ and $T_2$. But then there are also only finitely many possibilities for $T_2$ and $L_2$. The Lemma is proved.\qed

\vspace{.1in}

Slightly abusing notions (since, as explained above, $\H_{\A} \widehat{\otimes} \H_{\A}$ is \textit{not} an algebra) we may state Green's fundamental theorem as follows~:

\begin{theo}[Green, \cite{Green}] The map $\Delta: \H_{\A} \to \H_{\A} \widehat{\otimes} \H_{\A}$ is a morphism of algebras, i.e. for any $x, y \in \H_{\A}$ we have $\Delta(x\cdot y)=\Delta(x) \cdot \Delta(y)$.\end{theo}
\noindent
\textit{Proof.}  There is a very detailed proof in \cite{RingelGreen}. For the reader's convenience, we sketch it here. By bilinearity of the product, it is enough to show that for any pair of objects $M,N$ in $\A$ we have 
\begin{equation}\label{E:Greenproof1}
\Delta([M] \cdot [N])=\Delta([M]) \cdot \Delta([N]).
\end{equation}
The left hand side of (\ref{E:Greenproof1}) may be rewritten as follows~:
\begin{equation*}
\begin{split}
\Delta([M][N])&= \langle M, N \rangle \sum_J \frac{1}{a_M a_N} \mathbf{P}_{M,N}^J \Delta([J])\\
&=\sum_{K,L} \langle M,N \rangle \langle K,L \rangle \frac{1}{a_Ma_N} \sum_J \frac{1}{a_J} \mathbf{P}^J_{M,N} \mathbf{P}^J_{K,L} [K] \otimes [L].
\end{split}
\end{equation*}
As for the right hand side of (\ref{E:Greenproof1}), it may be expressed as

\begin{equation*}
\begin{split}
&\Delta([M])\Delta([N])\\
&=\left( \sum_{K_2,L_2} \langle K_2,L_2 \rangle \frac{1}{a_M} \mathbf{P}^M_{K_2,L_2} [K_2] \otimes [L_2] \right) \left(  \sum_{K_1,L_1} \langle K_1,L_1 \rangle \frac{1}{a_N} \mathbf{P}^N_{K_1,L_1} [K_1] \otimes [L_1] \right)\\
&=\sum_{K_1,L_1,K_2,L_2} \frac{\langle K_1,L_1 \rangle \langle K_2,L_2\rangle \langle K_1,L_2 \rangle \langle L_2,K_1\rangle}{a_Ma_N} \mathbf{P}^{M}_{K_2,L_2} \mathbf{P}^N_{K_1,L_1} [K_2][K_1] \otimes [L_2][L_1]
\end{split}
\end{equation*}
\begin{equation*}
\begin{split}
&= \sum_{K,L}\sum_{K_1,L_1,K_2,L_2} \frac{\langle K_1,L_1 \rangle \langle K_2,L_2\rangle \langle K_1,L_2 \rangle \langle L_2,K_1\rangle\langle K_2,K_1\rangle \langle L_2,L_1\rangle }{a_Ma_Na_{K_1}a_{K_2}a_{L_1}a_{L_2}} \times \\
& \qquad \qquad \qquad \qquad \qquad \qquad \qquad \qquad \qquad \times\mathbf{P}^{K}_{K_2,K_1}\mathbf{P}^{L}_{L_2,L_1}\mathbf{P}^{M}_{K_2,L_2} \mathbf{P}^N_{K_1,L_1}
[K] \otimes [L]\\
\\
&=\sum_{K,L} \frac{\langle M,N \rangle \langle K,L\rangle}{a_Ma_N} \sum_{K_1, K_2, L_1, L_2} \frac{| {Ext}(K_2,L_1)|}{|{Hom}(K_2,L_1)|}\frac{\mathbf{P}^{K}_{K_2,K_1}\mathbf{P}^{L}_{L_2,L_1}\mathbf{P}^{M}_{K_2,L_2} \mathbf{P}^N_{K_1,L_1}}{a_{K_1}a_{K_2}a_{L_1}a_{L_2}}
[K] \otimes [L]
\end{split}
\end{equation*}
where we have used the bilinearity of the Euler form. 
Simplifying both sides, we have to prove that, for each pair of objects $K,L$, the number
\begin{equation}\label{E:Greenproof3}
\sum_{K_1,L_1,K_2,L_2}\frac{|{Ext}(K_2,L_1)|} {|{Hom}(K_2,L_1)|}\frac{\mathbf{P}^{K}_{K_2,K_1}\mathbf{P}^{L}_{L_2,L_1}\mathbf{P}^{M}_{K_2,L_2} \mathbf{P}^N_{K_1,L_1}}{a_{K_1}a_{K_2}a_{L_1}a_{L_2}}
\end{equation}
which counts (with a certain weight) the set of ``squares''
\begin{equation}\label{E:gurz}
\xymatrix{
& 0 \ar[d] & & 0\ar[d] & \\
0 \ar[r] & L_1 \ar[r]^-{u} \ar[d]^-{u'} & L\ar[r]^-{v} & L_2 \ar[r] \ar[d]^-{x}& 0\\
& N \ar[d]^-{v'} & & M \ar[d]^-{y} &\\
0\ar[r] & \ar[r] K_1 \ar[d] \ar[r]^-{x'} & K \ar[r]^-{y'} & K_2 \ar[r] \ar[d] & 0\\
& 0 & & 0 &}
\end{equation}
\noindent is equal to the number
\begin{equation}\label{E:Greenproof2}
\sum_J \frac{1}{a_J} \mathbf{P}_{M,N}^J \mathbf{P}^J_{K,L}
\end{equation} which counts (with a certain weight) the set of ``crosses''
\begin{equation}\label{E:crosses}
\xymatrix{
& & 0 \ar[d] & &\\
& & L \ar[d]^-{a'} & & \\
0 \ar[r] & N \ar[r]^-{a} & J \ar[r]^{b} \ar[d]^-{b'} & M \ar[r] & 0\\
& & K \ar[d] & &\\
& & 0 & &}
\end{equation}

We call two crosses as in (\ref{E:crosses}), with respective middle term $J$ and $J'$ equivalent if there exists an isomorphism $\phi: J \simeq J'$ making all relevant diagrams commute, and we let $C_J$ (resp. $\widetilde{C}_J$) stand for the set of all crosses with middle term $J$ (resp. the set of all such crosses, up to equivalence). In the same fashion we define, for a triple $(K_1, K_2, L_1, L_2)$, the set $D_{(K_1, K_2, L_1, L_2)}$ of all squares as in (\ref{E:gurz}) with vertices $K_1, K_2, L_2, L_2$ and the set $\widetilde{D}_{(K_1,K_2,L_1,L_2)}$ of all such squares, up to equivalence.\\

Note that for a given quadruple $(K_1, K_2, L_1, L_2)$, the group ${Aut}(K_1) \times {Aut}(K_2) \times {Aut}(L_1) \times {Aut}(L_2)$ acts freely on $D_{(K_1, K_2, L_1, L_2)}$ and hence
\begin{equation*}
\begin{split}
| \widetilde{D}_{(K_1,K_2, L_1, L_2)}|&=\frac{| D_{(K_1, K_2, L_1, L_2)}|}{a_{K_1}a_{K_2}a_{L_1}a_{L_2}}\\
&=\frac{\mathbf{P}^{K}_{K_2,K_1}\mathbf{P}^{L}_{L_2,L_1}\mathbf{P}^{M}_{K_2,L_2} \mathbf{P}^N_{K_1,L_1}}{a_{K_1}a_{K_2}a_{L_1}a_{L_2}}
\end{split}
\end{equation*}
so that the expression (\ref{E:Greenproof3}) is equal to
\begin{equation}\label{E:Greenproof4}
\sum_{K_1, K_2, L_1, L_2}\frac{| {Ext}(K_2, L_1)|}{| {Hom}(K_2, L_1)|}  |\widetilde{D}_{(K_1, K_2, L_1, L_2)}|.
\end{equation}
Similarly, the stabilizer in ${Aut}(J)$ of a cross $C=(J,a,a',b,b')$ in $C_J$ is the set of $\phi$ satisfying
$\phi_{|{Ker} \;b}=Id,\; \phi_{|{Ker} \;b'}=Id$, and for which the induced maps  satisfy $\phi_{|{Coker} \;a}=Id,\; \phi_{|{Coker} \;a'}=Id$. A simple diagram chase shows that such automorphisms are in bijection with elements in ${Hom}({Coker}\;b'a, {Ker}\;b'a)$. It follows that the expression in (\ref{E:Greenproof2}) is equal to
\begin{equation}\label{E:Greenproof5}
\sum_J \frac{1}{a_J} | C_J|=\sum_J \sum_{(J,a,a',b,b') \in \widetilde{C}_J} \frac{1}{|{Hom}({Coker}\; b'a, {Ker}\;b'a)|}.
\end{equation}

\vspace{.1in}

We will show that (\ref{E:Greenproof4}) and (\ref{E:Greenproof5}) are equal by using a canonical map $\widetilde{\Phi}: \bigsqcup \widetilde{C}_J \to \bigsqcup \widetilde{D}_{(K_1, K_2, L_1, L_2)}$, constructed as follows~: given a cross $C=(J,a,a',b,b')$ we may complete it into a big commutative diagram by setting
$$L_1={Ker}\; b'a \simeq {Ker}\; ba', \qquad L_2={Im}\;ba', $$
$$ K_1={Im}\;b'a, \qquad K_2={Coker}\; ba' \simeq {Coker}\;b'a.$$
The maps $u,v,u',v',x,y,x',y'$ are all tautological, and it is an easy exercise to check that the resulting square $S=\Phi(C)$, obtained by deleting the central object $J$, does indeed belong to
${D}_{(K_1,K_2,L_1,L_2)}$. For instance, when $\mathcal{A}$ is the category of modules over some ring, we have $L_1 \simeq L\cap N, \;L_2 \simeq L/(L\cap N), \;K_1\simeq N/(L \cap N)$ and $K_2 \simeq J/(L \oplus N)$. The map $\Phi$ defined above descends to the desired map $\widetilde{\Phi}$ between equivalence classe of crosses and squares.

Combining (\ref{E:Greenproof4}) and (\ref{E:Greenproof5}), we see that the proof of Green's theorem boils down to the following fact~:

\vspace{.05in}

\noindent
\textbf{Claim.} The fiber of $\widetilde{\Phi}$ over any point of $\widetilde{D}_{(K_1, K_2, L_1, L_2)}$ is of cardinality $|{Ext}(K_2,L_1)|$.

\vspace{.05in}

\noindent
\textit{Proof of Claim.} Let us fix a square $S \in D_{(K_1, K_2, L_1, L_2)}$. The sequences
\begin{align*}
&\xymatrix{ 0 \ar[r] & K_1 \ar[r] & K \ar[r] & K_2 \ar[r] & 0},\\
&\xymatrix{ 0 \ar[r] & L_2 \ar[r] & M \ar[r] & K_2 \ar[r] & 0}
\end{align*}
determine elements $\sigma_1 \in {Ext}(K_2,K_1), \sigma_2 \in {Ext}(K_2,L_2)$ respectively, which give rise to an element $\eta=\sigma_1\oplus\sigma_2 \in {Ext}(K_2, K_1 \oplus L_2)$ corresponding to an extension
$$\xymatrix{ 0 \ar[r] &K_1 \oplus L_2 \ar[r]^-{x''} & X \ar[r]^-{y''}& K_2 \ar[r] & 0}$$
where $X=K \times_{K_2} M =\{(k \oplus m) \in K \oplus M\;|\; y'(k)=y(m)\}$. Similarly, there exists a canonical exact sequence
$$\xymatrix{ 0 \ar[r] &L_1 \ar[r]^-{u''} & Y \ar[r]^-{v''}& K_1 \oplus L_2 \ar[r] & 0}$$
where $Y=L \sqcup_{\small{L_1}} \hspace{-.05in}N =L \oplus N / \{u(l)\oplus u'(l)\;|\;l \in L_1\}$. 

Assume given a cross $C=(a,a',b,b',J)$ such that $\Phi(C)\simeq S$. As is easily seen, we may fit $J$ into an exact, commutative diagram
\begin{equation}\label{E:Cross1}
\xymatrix{
& 0&0& &\\
0 \ar[r] & K_1 \oplus L_2 \ar[r]^-{x''} \ar[u] & X \ar[u] \ar[r]^-{y''}& K_2 \ar[r] & 0\\
0 \ar[r] & Y \ar[u]^-{v''} \ar[r]^-{a''} & J \ar[u]^-{b''} \ar[r]^-{t} & K_2 \ar[r] \ar@{=}[u] & 0\\
& L_1 \ar[u]^-{u''} \ar@{=}[r] & L_1 \ar[u]^-{s} & & \\
& 0 \ar[u] & 0 \ar[u] & &}
\end{equation}

Conversely, one shows (see \cite{Green}, Prop. (2.6d)) that, once $J$ is fixed, the set $\Phi_J^{-1}(S)$ of crosses in the fiber $\Phi^{-1}(S)$ with central object $J$ is in \textit{bijection} with the set of exact, commutative diagrams

\begin{equation}\label{E:Cross2}
\xymatrix{
0 \ar[r] & K_1 \oplus L_2 \ar[r]^-{x''}  & X  \ar[r]^-{y''}& K_2 \ar[r] & 0\\
0 \ar[r] & Y \ar[u]^-{v''} \ar[r]^-{a''} & J \ar[u]^-{b''} \ar[r]^-{t} & K_2 \ar[r] \ar@{=}[u] & 0}
\end{equation}
(where $v'',x'',y''$ are fixed and $a'',b''$ are allowed to vary). Of course, instead of (\ref{E:Cross2}), we could just as well have chosen the first two columns of (\ref{E:Cross1}). Our task will now be to compute, for a given $J$, the number of diagrams of the form (\ref{E:Cross2}).

\vspace{.05in}

First of all, for a given extension
\begin{equation}\label{E:Cross3}
\xymatrix{
0 \ar[r] & Y \ar[r]^-{a''} & J \ar[r]^-{t} & K_2 \ar[r]  & 0}
\end{equation}
there exists a map $b''$ making (\ref{E:Cross2}) commute if and only if the element $\xi \in {Ext}(K_2,Y)$ associated to (\ref{E:Cross3}) satisfies $\phi(\xi)=\eta \in {Ext}(K_2, K_1 \oplus L_2)$, where $\eta$ is defined previously and where $\phi: {Ext}(K_2,Y) \to {Ext}(K_2, K_1 \oplus L_2)$ fits in the long exact sequence

\begin{equation}\label{E:Cross4}
\xymatrix{
0 \ar[r] & {Hom}(K_2,L_1) \ar[r] & {Hom}(K_2,Y) \ar[r]& {Hom}(K_2,K_1 \oplus L_2) \ar[r] &}
\end{equation}
$$\qquad\xymatrix{
 \ar[r] & {Ext}(K_2,L_1) \ar[r] & {Ext}(K_2,Y) \ar[r]^-{\phi} & {Ext}(K_2,K_1 \oplus L_2) \ar[r] & 0}
$$ 
Note that the sequence stops after two lines since by assumption ${gldim}(\mathcal{A}) \leq 1$. Let ${Ext}_{\eta,J}(K_2,Y) \subset {Ext}(K_2,Y)$ be the subset of $\phi^{-1}(\eta)$ whose associated middle term is isomorphic to $J$. We will use the following two facts, whose (easy) proofs we leave to the reader (see Prop.~\ref{P:coprophage} for statement i))~:

\vspace{.05in}

i) the number of sequences (\ref{E:Cross3}) corresponding to any fixed extension class $\xi \in {Ext}(K_2,Y)$ is equal to $a_{J} / |{Hom}(K_2,Y)|$,

\vspace{.05in}

ii) for any such sequence (\ref{E:Cross3}) associated to an element in ${Ext}_{\eta,J}(K_2,Y)$, the number of choices for the map $b''$ in (\ref{E:Cross2}) (i.e. making (\ref{E:Cross2}) commute) is equal to $|{Hom}(K_2, K_1 \oplus L_2)|$.

\vspace{.1in}

Hence, we see that, all together, there are 
$$|{Ext}_{\eta,J}(K_2,Y)| a_J |{Hom}(K_2,K_1 \oplus L_2)| /  |{Hom}(K_2,Y)|$$
 sequences as in (\ref{E:Cross2}) involving $J$, and hence as many crosses belonging to 
 $\Phi^{-1}(S)$ with $J$ as a central object. Remembering that the stabilizer in ${Aut}(J)$ of any such cross is isomorphic to ${Hom}(K_2,L_1)$, we see that there are
\begin{equation*}
 |{Ext}_{\eta,J}(K_2,Y)| \; |{Hom}(K_2,L_1)| \; |{Hom}(K_2,K_1 \oplus L_2)| /  |{Hom}(K_2,Y)|
\end{equation*}
equivalence classes of crosses in $\widetilde{\Phi}^{-1}(S)$ with central object isomorphic to $J$.

Summing up over all $J$ (up to isomorphism), we get
\begin{equation*}
\begin{split}
 |\widetilde{\Phi}^{-1}(S)|&=
 \sum_J |{Ext}_{\eta,J}(K_2,Y)| \; |{Hom}(K_2,L_1)| \; |{Hom}(K_2,K_1 \oplus L_2)| /  |{Hom}(K_2,Y)|\\
&=| \phi^{-1}(\eta)|\; |{Hom}(K_2,L_1)| \; |{Hom}(K_2,K_1 \oplus L_2)| /  |{Hom}(K_2,Y)|.
\end{split}
\end{equation*}
But from the long exact sequence (\ref{E:Cross4}), we deduce that 
$$ |\phi^{-1}(\eta)|=|{Ext}(K_2,L_1)| \;|{Hom}(K_2,Y)|/  |{Hom}(K_2, K_1 \oplus L_2)| \; |{Hom}(K_2,L_1)|,$$
from which we finally obtain $|\widetilde{\Phi}^{-1}(S)|=|{Ext}(K_2,L_1)|$ as wanted. This concludes the proof of the claim, and thus the proof of Green's theorem.
\qed

\vspace{.1in}

It is quite useful to restate Green's theorem in a different way, by adding a piece of ``degree zero'' to $\H_{\A}$. This way we will avoid twisting the multiplication in $\H_{\A} \otimes \H_{\A}$ as above.
 This may seem rather artificial at first glance, but it is in fact very natural given the analogy with quantum groups (see Lecture~3 ). Let $\mathbf{K}=\C[K(\A)]$ be the group algebra of the Grothendieck group $K(\A)$. To avoid any confusion we denote by $\mathbf{k}_{{M}}$ (resp. $\mathbf{k}_{\a}$) the element of $\mathbf{K}$ corresponding to the class of an object $M$ (resp. to the class $\a$). We also write $\overline{N}$ for the class of an object $N$ in $K(\A)$. We equip the vector space $\widetilde{\H}_{\A} := \H_{\A} \otimes_{\C} \mathbf{K}$ with the structure of an algebra (containing $\H_{\A}$ and $\mathbf{K}$ as subalgebras) by imposing the relations 
$$\mathbf{k}_{\a} [M] \mathbf{k}_{\a}^{-1}=(\a, M )_m[M].$$ 
The algebra $\widetilde{\H}_{\A}$ is still $K(\A)$-graded, where $deg(\mathbf{k}_{\a})=0$ for any $\a$.
We also extend the comultiplication to a map ${\Delta}: \widetilde{\H}_{\A} \to \widetilde{\H}_{\A} \widehat{\otimes} \widetilde{\H}_{\A}$ by setting
\begin{equation}\label{E:hallcoprodext1}
{\Delta}(\mathbf{k}_{\a})=\mathbf{k}_{\a} \otimes \mathbf{k}_{\a},
\end{equation}
\begin{equation}\label{E:hallcoprodext2}
{\Delta}([R]\mathbf{k}_{\a})=\sum_{M,N} \langle M, N \rangle_m \frac{1}{a_R} \mathbf{P}^R_{M,N} [M]\mathbf{k}_{\overline{N}+\a} \otimes [N]\mathbf{k}_{\a}.
\end{equation}
Finally, we equip $\widetilde{\H}_{\A} \otimes \widetilde{\H}_{\A}$ with the \textit{standard} multiplication, i.e. $(x \otimes y) (z \otimes w)=xz \otimes yw$. We call $\widetilde{\H}_{\A}$ the \textit{extended} Hall algebra of $\A$. To avoid confusion, we will sometimes write $\Delta'$ for the \textit{old} comultiplication, defined on $\H_{\A}$.

\vspace{.1in}

Green's theorem may now be expressed in the following fashion.

\begin{cor} The map ${\Delta}: \widetilde{\H}_{\A} \to \widetilde{\H}_{\A} \widehat{\otimes} \widetilde{\H}_{\A}$ is a morphism of algebras.\end{cor}

\vspace{.1in}

Of course, concerning the multiplication in the completion $\widetilde{\H}_{\A} \widehat{\otimes} \widetilde{\H}_{\A}$, the reservations voiced at the beginning of this section still apply.

To turn $\widetilde{\H}_{\A}$ into a bialgebra, it only remains to introduce a counit map. Define a $\C$-linear morphism $\epsilon: \widetilde{\H}_{\A} \to \C$ by
$$
\epsilon([M]\mathbf{k}_{\a})=
\begin{cases} 0  \qquad & \text{if}\; M \neq 0,\\
1 \qquad & \text{if}\; M=0.
\end{cases}
$$
The reader will easily check that $(\widetilde{\H}_{\A}, i, m, \epsilon, {\Delta})$ is a (topological) bialgebra.

\vspace{.2in}

\addtocounter{theo}{1}
\noindent \textbf{Remarks \thetheo .} i) Of course, when $\A$ satisfies the finite subobjects condition (\ref{E:FS}), there is no need to consider any completion at all, and $\widetilde{\H}_{\A}$ is a genuine bialgebra. 

\vspace{.05in}

\noindent
ii) When the symmetrized Euler form is trivial, the twisted multiplication in the tensor product $\H_{\A} \otimes \H_{\A}$ coincides with the usual multiplication, and there is \textit{a priori} no need to introduce the extension $\widetilde{\H}_{\A}$.

\vspace{.05in}

\noindent
iii) Although it is possible to define a comultiplication in $\H_{\A}$ (or $\widetilde{\H}_{\A}$) for exact categories, Green's theorem only holds for \textit{abelian} categories.

\vspace{.2in}

\centerline{\textbf{1.6. Green's scalar product.}}
\addcontentsline{toc}{subsection}{\tocsubsection {}{}{\; 1.6. Green's scalar product.}}

\vspace{.15in}

\paragraph{} The bialgebras which arise as Hall algebras of some abelian category $\A$ have an additional important feature~: they are \textit{self-dual}. This means that the dual space $\widetilde{\H}_{\A}^*$, equipped with the multiplication which is dual to $\Delta$ and the comultiplication dual to $m$, is isomorphic to $\widetilde{\H}_{\A}$ as a bialgebra. The best way to state this property is to use a natural nondegenerate bilinear form on $\widetilde{\H}_{\A}$, introduced by Green in \cite{Green}. In this paragraph as in the previous one, $\A$ denotes a \textit{hereditary} finitary abelian category. 

\begin{prop}[\cite{Green}] The nondegenerate scalar product $(\,,\,): {\H}_{\A} \otimes {\H}_{\A} \to \C$ defined by
\begin{equation}\label{E:Greenscalar}
([M],[N])=\frac{\delta_{M,N}}{a_M}
\end{equation}
is a Hopf pairing, that is for any triple $x,y,z$ of elements in $\H_{\A}$ we have $(xy,z)=(x \otimes y, \Delta'(z))$.
\end{prop}
\noindent
\textit{Proof.} By bilinearity, it is enough to check this when $x=[M], y=[N]$ and $z=[P]$ for some objects $M,N,P$ of $\mathcal{A}$. Then $xy=\langle M,N\rangle_m \sum_R \frac{1}{a_Ma_N} \mathbf{P}^R_{M,N}[R]$ hence
$$(xy,z)=\langle M,N \rangle_m \frac{1}{a_Ma_Na_P}\mathbf{P}^P_{M,N}.$$
On the other hand, we have $\Delta'(z)=\sum_{R,S} \langle R,S \rangle_m \frac{1}{a_P} \mathbf{P}^P_{R,S} [R] \otimes [S]$ therefore
$$(x \otimes y,\Delta'(z))=\langle M,N \rangle_m \frac{1}{a_Ma_Na_P}\mathbf{P}^P_{M,N}.$$
We are done. \qed

\vspace{.15in}

Recall that $\mathbf{H}_{\mathcal{A}}$ is in general not a bialgebra. However, it is an easy task to extend Green's scalar product to $\widetilde{\mathbf{H}}_{\mathcal{A}}$.

\begin{cor} The scalar product $\widetilde{\mathbf{H}}_{\mathcal{A}} \otimes \widetilde{\mathbf{H}}_{\mathcal{A}} \to \mathbb{C}$ defined by
\begin{equation}\label{E:Greenscalar2}
([M]\mathbf{k}_{\a},[N]\mathbf{k}_{\beta})=\frac{\delta_{M,N}}{a_M}(\a,\beta)_m
\end{equation}
is a Hopf pairing, i.e. for any triple $x,y,z$ of elements in $\H_{\A}$ we have $(xy,z)=(x \otimes y, \Delta(z))$.
\end{cor}

\vspace{.12in}

As the examples of Lecture~3 and Lecture~4 will attest, the existence of this nondegenerate scalar product is a very strong property of Hall algebras, and despite its simple form, it ``encodes''  a surprising amount of information concerning $\H_{\A}$.

\newpage

\centerline{\textbf{1.7. Xiao's antipode and the Hall Hopf algebra.}} 
\addcontentsline{toc}{subsection}{\tocsubsection {}{}{\; 1.7. Xiao's antipode and the Hall Hopf algebra.}}

\vspace{.15in}

\paragraph{}The results of paragraphs~1.4. and 1.5. show that the Hall algebra $\widetilde{\H}_{\A}$ of a hereditary, finitary abelian category $\mathcal{A}$ may be equipped with the structure of a bialgebra, at least in a topological sense.  At this point, it is natural to ask whether $\widetilde{\H}_{\A}$ may actually be upgraded to a Hopf algebra, i.e. whether one can define a morphism $S: \widetilde{\H}_{\A} \to \widetilde{\H}_{\A}$ satisfying the axioms of an antipode, namely
\begin{equation}\label{E:antipode1}
m \circ (1 \otimes S) \circ {\Delta} = i \circ \epsilon, \qquad
m \circ (S \otimes 1) \circ {\Delta} = i \circ \epsilon,
\end{equation}
\begin{equation}\label{E:antipode2}
S(xy)=S(y)S(x), \qquad {\Delta} \circ S = (S \otimes S) \circ {\Delta}^{op},
\end{equation}
\begin{equation}\label{E:antipode3}
S \circ i=i, \qquad \epsilon \circ S=\epsilon,
\end{equation}

When $\mathcal{A}$ is the category of representations of a quiver, Xiao discovered such an antipode map (\cite{Xiao1}). His construction can be directly extended to any category satisfying the finite subobjects condition (\ref{E:FS}). Let $\mathcal{A}$ be such an abelian finitary category. 
For any object $M$ and integer $r$ let $S_{M,r}$ denote the set of \textit{strict} $r$-step filtrations
$$0 \neq L_r \subsetneq \cdots \subsetneq L_2 \subsetneq L_1=M$$
and let us set
\begin{equation}\label{E:antipodef}
\begin{split}
S([M])=&\mathbf{k}_M^{-1} \bigg( \sum_{r \geq 1} (-1)^r \sum_{\underline{L}_{\bullet} \in S_{M,r}} \prod_{i=1}^r \langle L_i/L_{i+1}, L_{i+1} \rangle_m \frac{a_{L_r} a_{L_{r-1}/L_r} \cdots a_{M/L_2}}{a_M} \times \\
&\qquad \qquad \qquad \qquad \qquad \qquad \qquad \; \;  [M/L_2]  \cdot [L_2/L_3] \cdots [L_{r-1}/L_r] \cdot [L_r] \bigg)\\
=& -\mathbf{k}_{M}^{-1}[M] + \sum_{0 \neq N \subsetneq M} \langle M/N,N \rangle_m \frac{a_{M/N} a_N}{a_M} \mathbf{k}_M^{-1}[M/N] \cdot [N] + \cdots 
\end{split}
\end{equation}
In plain words, we consider all strict filtrations of $M$, take the corresponding successive subquotients, multiply them together in the Hall algebra \textit{in the same order}, and sprinkle a few signs here and there. Of course, this sum is finite by our assumption on $\mathcal{A}$. We extend the map $S$ to the whole of $\widetilde{\H}_{\A}$ by bilinearity and by setting $S([M]\mathbf{k}_{\a})=\mathbf{k}_{\a}^{-1} S([M])$. This map respects the grading of $\widetilde{\H}_{\A}$.

\begin{theo}[Xiao, \cite{Xiao1}] The map $S: \widetilde{\H}_{\A} \to \widetilde{\H}_{\A} $ is an antipode for the Hall bialgebra $\widetilde{\H}_{\A}$, i.e. the identities (\ref{E:antipode1}), (\ref{E:antipode2}) and (\ref{E:antipode3}) are all satisfied.\end{theo}

\noindent
\textit{Proof.} The relations (\ref{E:antipode3}) follow directly from the definitions. In addition, it is well-known in the theory of Hopf algebras that relations (\ref{E:antipode2}) are consequences of (\ref{E:antipode1}) and (\ref{E:antipode3}). Hence we only need to prove relations (\ref{E:antipode1}). We will deal with the first equality, and leave the second one to the reader. It is clearly enough to show that for any non zero object $M$ in $\mathcal{A}$, we have $m \circ (1 \otimes S) \circ {\Delta} ([M])=0$. We have, using the definitions (\ref{E:antipodef}) and (\ref{E:hallcoprodext2})

\begin{equation}\label{E:proofantipode}
\begin{split}
m& \circ (1 \otimes S) \circ {\Delta}([M])\\
&= [M] + \mathbf{k}_MS([M]) +\sum_{0 \neq N \subsetneq M} \langle M/N,N \rangle_m \frac{a_{M/N}a_N}{a_M} [M/N] \mathbf{k}_N S([N])
\end{split}
\end{equation}
\begin{equation*}
\begin{split}
&=[M] + \sum_{r \geq 1} (-1)^r \sum_{\underline{L}_{\bullet} \in S_{M,r}} \prod_{i=1}^r \langle L_i/L_{i+1}, L_{i+1} \rangle_m \frac{\prod_{i=1}^r a_{L_i/L_{i+1}}}{a_M} [M/L_2] \cdots [L_{r-1}/L_r] \cdot [L_r] \\
& \; + \sum_{0 \neq N \subsetneq M} \langle M/N,N \rangle_m \frac{a_{M/N}a_N}{a_M} [M/N] \bigg\{
\sum_{s \geq 1} (-1)^s \sum_{\underline{K}_{\bullet} \in S_{N,s}} \prod_i \langle K_i/K_{i+1},K_{i+1} \rangle_m \\
&\qquad \qquad \qquad \qquad \qquad \qquad \qquad \; \;\qquad \qquad \qquad  \frac{\prod_i a_{K_i/K_{i+1}}}{a_N} [N/K_2] \cdots [K_{s-1}/K_s] \cdot [K_s] \bigg\}. 
\end{split}
\end{equation*}

The last term in the above expression may be rewritten as

\begin{equation*}
\begin{split}
&\sum_{s \geq 1} (-1)^s \sum_{0 \neq N \subsetneq M}\sum_{\underline{K}_{\bullet} \in S_{N,s}}
\langle M/N,N\rangle_m \langle N/K_2, K_2 \rangle_m \cdots \langle K_{s-1}/K_s,K_s \rangle_m\\
& \qquad \qquad \qquad \qquad \qquad \qquad \frac{a_{M/N} a_{N/K_2} \cdots a_{K_s}}{a_M} [M/N] \cdot [N/K_2] \cdots [K_{s-1}/K_s] \cdot [K_s]\\
&= -\sum_{r \geq 2} (-1)^r \sum_{\underline{L}_{\bullet} \in S_{M,r}} \prod_{i=1}^r \langle L_i/L_{i+1}, L_{i+1} \rangle_m \frac{\prod_{i=1}^r a_{L_i/L_{i+1}}}{a_M} [M/L_2] \cdots [L_{r-1}/L_r] \cdot [L_r].
\end{split}
\end{equation*}

It follows that the whole quantity in (\ref{E:proofantipode}) vanishes, as desired. The Theorem is proved. \qed

\vspace{.2in}

\addtocounter{theo}{1}
\noindent \textbf{Remark \thetheo .}  What happens when $\mathcal{A}$ does not satisfy the finite subobject condition (for example, when $\mathcal{A}$ is a category of coherent sheaves on a curve) ?
In this case, the sum in (\ref{E:antipodef}) does not make sense in general (indeed, it may be an infinite sum of positive integers !). Depending on the situation, there may be two ways around this difficulty~:

\vspace{.05in}

\noindent
i) The category $\mathcal{A}$ may be $\mathbb{F}_q$-linear, and all the structure constants for the antipode may be expressed as power series in $q$.

\vspace{.05in}

\noindent
ii) Instead of an antipode map $S$, one may try to define the inverse map $S^{-1}$. When $\mathcal{A}$ satisfies the finite subobjects condition, such a map is given by the expression 
\begin{equation}\label{E:inverseantipode}
\begin{split}
S^{-1}([M])=&\bigg( \sum_{r \geq 1} (-1)^r \sum_{\underline{L}_{\bullet} \in S_{M,r}} \prod_{i=1}^r \langle L_i/L_{i+1}, L_{i+1} \rangle_m \frac{a_{L_r} a_{L_{r-1}/L_r} \cdots a_{M/L_2}}{a_M} \times \\
&\qquad \qquad \qquad \qquad \qquad \qquad \; \;  [L_r] \cdot [L_{r-1}/L_r]  \cdots [L_2/L_3] \cdot [M/L_2] \bigg)\mathbf{k}_{M}.
\end{split}
\end{equation}
In this sum, the coefficient of $[N]$ counts, with a certain weight, the number of pairs of filtrations
$ \big( L_r \subsetneq \cdots \subsetneq L_2 \subsetneq L_1=M; \;\; L'_r \subsetneq \cdots \subsetneq L'_2 \subsetneq L'_1=N \big)$ satisfying $L_i/L_{i+1}=L'_{r+1-i}/L'_{r+2-i}$ for $i=1, \ldots, r$. If $\mathcal{A}$ is such that for any $M,N \in Ob(\mathcal{A})$ there are finitely many such pairs of filtrations then the expression in (\ref{E:inverseantipode}) converges, and we may define an inverse antipode map $S^{-1}: \widetilde{\H}_{\A} \to \widetilde{\H}_{\A}^{c}$. Here $\widetilde{\H}_{\A}^{c}$ is a certain formal completion of $\widetilde{\H}_{\A}$.

The map $S^{-1}$, which serves the same basic purposes as $S$, satisfies the following relations~:
\begin{equation*}
m \circ (1 \otimes S^{-1}) \circ {\Delta}^{op} = i \circ \epsilon, \qquad
m \circ (S^{-1} \otimes 1) \circ {\Delta}^{op} = i \circ \epsilon,
\end{equation*}
\begin{equation*}
S^{-1}(xy)=S^{-1}(y)S^{-1}(x), \qquad {\Delta}^{op} \circ S^{-1} = (S^{-1} \otimes S^{-1}) \circ {\Delta},
\end{equation*}
\begin{equation*}
S^{-1} \circ i=i, \qquad \epsilon \circ S^{-1}=\epsilon.
\end{equation*}

\newpage

\centerline{\textbf{1.8. Functorial properties.}} 
\addcontentsline{toc}{subsection}{\tocsubsection {}{}{\; 1.8. Functorial properties.}}

\vspace{.15in}

\paragraph{} Let $\mathcal{A}$ and $\mathcal{B}$ be two finitary categories, and let $F: \mathcal{A} \to \mathcal{B}$ be a functor. What does the existence of $F$ imply about the Hall algebras $\H_{\A}$ and $\H_{\mathcal{B}}$ ? The functor $F$ gives rise to a map between 
sets of isomorphism classes of objects $f:\;\mathcal{X}_{\A} \to \mathcal{X}_{\mathcal{B}}$, and we may define two natural maps \textit{of vector spaces}
\begin{align*}
f^*~: \H_{\mathcal{B}} &\to \H_{\mathcal{A}}^c,\\
[M] &\mapsto \sum_{R \in \mathcal{X}_{\A}, F(R) \simeq M} [R] 
\end{align*}
and
\begin{align*}
f_*~: \H_{\mathcal{A}} &\to \H_{\mathcal{B}},\\
 [M] & \mapsto [F(M)]
\end{align*}

\noindent
(in the above, $\H_{\A}^c=\prod_{M \in \mathcal{X}_A} \C [M]$ is the formal completion of $\H_{\A}$).
Without any further assumptions on the functor $F$, these maps are in general neither morphisms of algebras nor morphisms of coalgebras. Let us briefly examine under which circumstances they are.

\vspace{.1in}

Let us start with $f^*$. By definition, if $g$ and $h$ belong to $\H_{\mathcal{B}}$ and $R \in \mathcal{X}_{\A}$ we have
$$f^* (g \cdot h) (R)=g \cdot h (F(R))=\sum_{S \subset F(R)} \langle F(R)/S, S \rangle_m g(F(R)/S) h(S)$$
while
$$f^*g \cdot f^* h (R)= \sum_{T \subset R} \langle R/T, T \rangle_m g(F(R/T)) h(F(T)).$$
A sufficient condition for these two quantitites to be equal, and hence for $f^*$ to be a morphism of algebras, is that~: $F$ preserves Euler forms (i.e. $\langle M, N \rangle_m = \langle F(M), F(N) \rangle_m$ for any $M,N \in \mathcal{X}_{\A}$) and that for any $R \in \mathcal{X}_{\A}$, the functor $F$ sets up a bijection between the subobjects of $R$ and $F(R)$.

Suppose now that $\A$ and $\mathcal{B}$ are hereditary and let us consider the compatibility of $f^*$ with the coproduct. We have for any function $g \in \H_{\mathcal{B}}$ and objects $M,N \in \mathcal{X}_{\A}$ (see Proposition~\ref{P:coprophage}),
$$\Delta (f^* g) (M,N)=\frac{1}{| {Ext}^1(M,N)|}  \sum_{\xi \in {Ext}^1(M,N)} g(F(R_{\xi}))$$
and
$$f^* \otimes f^* (\Delta(g))(M,N)= \frac{1}{|{Ext}^1(F(M),F(N))|} \sum_{\mu \in {Ext}^1(F(M),F(N))}g(R_{\mu}).$$
Therefore, a natural sufficient condition for $f^*$ to be a morphism of coalgebras is that $F$ is exact and that the associated map $F: {Ext}^1(M,N) \to {Ext}^1(F(M),F(N))$ is an isomorphism for any $M,N \in \mathcal{X}_{\A}$. This last condition amounts to saying that $F$ sets up a bijection between the sets of short exact sequences $\mathcal{P}_{M,N}^R$ and $\mathcal{P}^{F(R)}_{F(M),F(N)}$, and that furthermore any extension of $F(N)$ by $F(M)$ can be obtained in this way.

\vspace{.1in}

The properties of the map $f_*$ are also easily worked out. For $f_*$ to be a morphism of algebras, it is sufficient that $F$ is exact and that the associated map $F: {Ext}^1(M,N) \to {Ext}^1(F(M),F(N))$ is an isomorphism for any $M,N \in \mathcal{X}_{\A}$. For $f_*$ to be a morphism of coalgebras, it is sufficient that $F$ preserves Euler forms and that for any $R \in \mathcal{X}_{\A}$, the functor $F$ sets up a bijection between the subobjects of $R$ and $F(R)$. Thus the conditions for $f_*$ are in a certain sense dual to those for $f^*$.

\vspace{.15in}

All in all, we see that the conditions to impose on a functor $F$ for it to induce a (co)algebra morphism between Hall algebras are rather restrictive. We summarize the above discussion in the following result which is sufficient for many purposes. We will call an exact functor $F: \A \to \mathcal{B}$ \textit{extremely faithful} if it defines isomorphisms $Ext^i(M,N) \stackrel{\sim}{\to} Ext^i(F(M),F(N))$ for any $M,N \in \A$ and any $i \geq 0$.

\begin{cor}\label{C:cathall} Let $\mathcal{A}$ and $\mathcal{B}$ be finitary categories, and let $F: \A \to \mathcal{B}$ be an extremely faithful exact functor . Then $f_*$ is an embedding of algebras and if $\A, \mathcal{B}$ are hereditary then $f^*$ is a morphism of coalgebras. If in addition $F(\mathcal{A})$ is essentially stable under taking subobjects ($\mathrm{in}\;\mathcal{B}$) then $f_*$ is a morphism of coalgebras and $f^*$ is a morphism of algebras. \end{cor}

\vspace{.2in}

Though we have not explained the definition of the Hall algebra of an exact category, let us give a last (essentialy obvious) functoriality property of Hall algebras.

\begin{cor}\label{C:catexacthall} Let $\mathcal{B}$ be a finitary category and let $\mathcal{A} \subset \mathcal{B}$ be an exact full subcategory stable under extensions. Let $f: \mathcal{A} \to \mathcal{B}$ be the embedding. Then $f_*: \H_{\A} \to \H_{\mathcal{B}}$ is an embedding of algebras.
\end{cor}

As Lecture~4 will demonstrate, the above Corollary is particularly useful in the context of Hall algebras of categories of coherent sheaves on curves (where abelian subcategories are few and far between while exact subcategories abound).

\newpage

\centerline{\large{\textbf{Lecture~2.}}}
\addcontentsline{toc}{section}{\tocsection {}{}{Lecture~2.}}

\setcounter{section}{2}
\setcounter{equation}{0}
\setcounter{theo}{0}

\vspace{.2in}

After having introduced all the basic concepts of Hall algebras in Lecture~1, it is now high time for us to provide the reader with some concrete examples. This is what this and the next two Lectures are devoted to. In doing so, we hope to illustrate the following ``abstract nonsense'' principle, which we learned from Yves Benoist~:

\vspace{.1in}

\textit{
``There are many more theories than fundamental objects in mathematics.''}

\vspace{.1in}

A direct corollary of this principle is that behind many \textit{distinct} interesting theories lie in fact the \textit{same} fundamental objects. Of course, each theory sheds its own light on these objects, and combining the various perspectives is likely to be very fruitful.
 
As the examples of Lectures~2,3 and 4 will show, the theory of Hall algebras turns out to be intimately related to the structure theory of semisimple Lie algebras or Kac-Moody algebras, or to the combinatorial theory of symmetric functions. This interaction has proved to be extremely useful for all parties involved.

Incidentally, we will take the above principle as an excuse for not answering the question of \textit{why} such a deep relation exists and concentrate instead on the question of \textit{how} this relation exists.

\vspace{.2in}

The aim of the present lecture is to provide the first and perhaps most fundamental example of a Hall algebra --the so-called \textit{classical} Hall algebra, introduced by Steinitz at the turn of the twentieth century. This will allow the reader to see all the notions appearing in Lecture~1 in action. We will finish by briefly describing other occurences of Hall algebras in various mathematical contexts.

\vspace{.2in}

\centerline{\textbf{2.1. The Jordan quiver.}} 
\addcontentsline{toc}{subsection}{\tocsubsection {}{}{\; 2.1. The Jordan quiver.}}

\vspace{.15in}

\paragraph{}Let $k$ be any field. The Jordan quiver is the oriented graph $\vec{Q}_0$ with a single vertex $i$ and a single loop $h: i \to i$. 

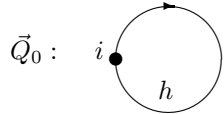
\vspace{.25in}

\centerline{
\begin{picture}(300, 10)
\put(80,0){$\vec{Q}_0:$}
\put(120,0){\circle*{5}}
\put(112,0){$i$}
\put(140,0){\circle{40}}
\put(140,20){\vector(1,0){2}}
\put(136,-15){$h$}
\end{picture}}

\vspace{.3in}

By definition, a representation of $\vec{Q}_0$ over $k$ is a pair $(V,x)$ consisting of a finite-dimensional $k$-vector space $V$ and a $k$-linear map $x: V \to V$. A morphism between two representations $(V,x)$ and $(V',x')$ of $\vec{Q}_0$ is simply a linear map $f: V \to V'$ making the following diagram commute
\begin{equation}\label{E:quiverdiag}
\xymatrix{
V \ar[r]^-{x} \ar[d]_-{f} & V \ar[d]_-{f}\\
V' \ar[r]^-{x'} & V'}
\end{equation}
The collection of all representations of $\vec{Q}_0$ over $k$ thus acquires the structure of a $k$-linear category $Rep_k\vec{Q}_0$ which is easily seen to be abelian. In fact, two seconds of thought will be amply enough for the reader to convince himself that this category is equivalent to the category of finite-dimensional modules over the polynomial ring $k[x]$. In particular, it is of global dimension one.

We will be interested here in the full subcategory consisting of those representations $(V,x)$ for which $x$ is a nilpotent endomorphism of $V$. This is the only one we will be considering here, we denote it by $Rep^{nil}_k\vec{Q}_0$.
The structure of $Rep^{nil}_k\vec{Q}_0$ has been essentially well-known since the end of the nineteenth century, due to the work of Jordan and Kronecker (see \cite{Brechi} for a fascinating discussion of the famous controversy, surrounding this result, which opposed these two great mathematicians). In modern language, it reads~:

\begin{theo}[Jordan and/or Kronecker]\label{T:Jordan} The following hold~:
\begin{enumerate}
\item[i)] The object $S=(k,0)$ is the only simple object of $Rep^{nil}_k\vec{Q}_0$,
\item[ii)] We have ${Hom}(S,S)=k$ and ${Ext}^1(S,S)=k$,
\item[iii)] For any $n \in \N$ there exists a unique indecomposable object $I_n$ of length $n$ in $Rep^{nil}_k\vec{Q}_0$. It is given by the endomorphism $x \in {End}(k^n)$ with matrix
$$x=\begin{pmatrix} 0 & 1 & 0 & \cdots & 0\\ 0 & 0 & 1 & \cdots & 0\\ 0 & 0 & 0 & \ddots & \vdots\\ 0 & 0 & 0 & \cdots & 1\\ 0 & 0 & 0 & \cdots & 0 \end{pmatrix}.$$ 
\end{enumerate}
\end{theo}

\vspace{.15in}

From i) it follows that $K(Rep^{nil}_k\vec{Q}_0) \simeq \Z$ and from ii) we deduce that the Euler form $\langle\;,\;\rangle_a$ vanishes. In addition, from iii) and the Krull-Schmidt theorem we see that any object $(V,x)$ in $Rep^{nil}_k\vec{Q}_0$ is isomorphic to a direct sum $I_1^{\oplus l_1} \oplus \cdots \oplus I_r^{\oplus l_r}$ for certain integers $r, l_1, \ldots, l_r$. Thus $\mathcal{X}_{Rep^{nil}_k\vec{Q}_0}$ is canonically isomorphic to the set $\Pi$ of all partitions, via the assignement
$${\lambda}=(\lambda_1, \lambda_2, \ldots, \lambda_r) \mapsto I_{\lambda}=I_{\lambda_1} \oplus \cdots \oplus I_{\lambda_r}.$$

\vspace{.2in}

\centerline{\textbf{2.2. Computation of some Hall numbers.}}
\addcontentsline{toc}{subsection}{\tocsubsection {}{}{\; 2.2. Computation of some Hall numbers.}}

\vspace{.15in}

\paragraph{}Let us assume from now on that $k=\mathbb{F}_q$ is a finite field with $q$ elements, so that $Rep_k\vec{Q}_0$ and $Rep^{nil}_k\vec{Q}_0$ are finitary hereditary categories. Our aim will be to describe the Hall algebra of $Rep^{nil}_k\vec{Q}_0$ in some details. For this, we need to understand the structure constants $P^{\nu}_{\mu,\lambda}:=\frac{1}{a_{I_{\mu}}a_{I_{\lambda}}}\mathbf{P}_{I_{\mu},I_{\lambda}}^{I_{\nu}}$. Obviously, $P^{\nu}_{\mu,\lambda}=0$ unless $|\nu|=|\lambda|+|\mu|$.

\vspace{.2in}

\addtocounter{theo}{1}
\paragraph{\textbf{Example \thetheo.}} Let us compute explicitly the first few values of $P^{\nu}_{\mu,\lambda}$.\\
$\bullet$ $P_{(1),(1)}^{(1^2)}$ counts the number of submodules $R \subset S \oplus S$ which are isomorphic to $S$, and such that $(S \oplus S) /R $ is isomorphic to $S$. Any one-dimensional submodule $R$ will do, and there are $| \mathbb{P}^1(k)|=q+1$ such submodules.\\
$\bullet$ $P^{(2)}_{(1),(1)}$ counts the number of submodules $R \subset I_2$ which are isomorphic to $S$ and such that $I_2/R$ is isomorphic to $S$. Again, any one-dimensional submodule $R$ will do but this time there is only one such submodule (the kernel of the map $x$). Hence
$P^{(2)}_{(1),(1)}=1$.\\
$\bullet$ $P^{(2,1)}_{(2),(1)}$ counts the number of submodules $R \subset M =I_2 \oplus S$ isomorphic to $S$ and such that $M /R$ is isomorphic to $I_2$. We have ${Im}\;x \subsetneq {Ker}\;x \subsetneq M$. A one-dimensional subspace $R$ will be isomorphic to $S$ if $R \subset {Ker}\;x$ while for $M/R$ to be isomorphic to $I_2$ we need $R \neq {Im}\;x$. Thus there are $|\mathbb{P}^1 (k)| -1 =q$ allowed choices for $R$.\\
$\bullet$ $P^{(2,1)}_{(1),(2)}$ counts the number of submodules of $R \subset M =I_2 \oplus S$ isomorphic to $I_2$ and such that $M /R$ is isomorphic to $S$. We have ${Im}\;x \subsetneq {Ker}\;x \subsetneq M$. A subspace $R$ of dimension two will be isomorphic to $I_2$ if $R \supset {Im}\;x$ while  $R \neq {Ker}\;x$. Thus there are again $|\mathbb{P}^1 (k)| -1 =q$ valid choices for $R$. 
\endexample

\vspace{.15in}

\addtocounter{theo}{1}
\paragraph{\textbf{Example \thetheo.}} The above very simple examples already afford some interesting generalisations~:\\
$\bullet$ $P^{(1^t)}_{(1^{t-r}),(1^r)}$ counts the number of $r$-dimensional submodules $R$ of the trivial module $M=S^{\oplus t}$ which are trivial and such that $M/R$ is trivial. We may pick any $r$-dimensional submodule, hence 
$P^{(1^t)}_{(1^{t-r}),(1^r)}$ is equal to the number of points over $\mathbb{F}_q$ of the Grassmanian $Gr(r,t)$ of $r$-planes in $t$-space. This number is equal to a $q$-binomial coefficient 
$$|Gr(r,t)|=\begin{bmatrix} t\\ r \end{bmatrix}_+=\frac{[t]_+ \cdot [t-1]_+ \cdots [t-r+1]_+}{[2]_+ \cdots [r-1]_+[r]_+} $$
where for any integer $n$ we set $[n]_+=1+q+ \cdots + q^{n-1}$. It is known that this $q$-binomial coefficient in fact belongs to $\N[q]$.\\
$\bullet$ $P^{(t)}_{(t-r),(r)}$ counts the number of indecomposable submodules $R$ of dimension $r$ of $I_t$  for which $I_t/R$ is again indecomposable. There is a unique $r$-dimensional submodule of $I_t$ and it satisfies the required conditions. Hence $P^{(t)}_{(t-r),(r)}=1$.
\endexample

\vspace{.15in}

\addtocounter{theo}{1}
\paragraph{\textbf{Example \thetheo.}}  As a final, more complicated example let us compute
$P^{\nu}_{\mu,(1^r)}$. Thus, given a representation $I_{\nu}$ we want to count the number of submodules $R$ isomorphic to $S^{\oplus r}$ such that $I_{\nu}/M$ is of a given type. Let us set
$$K={Ker}\;x \subset I_{\nu}, \qquad K_i=K \cap {Im}\;x^i,$$
so that we have $K=K_0 \supset K_1 \supset K_2 \supset \cdots$. Now let $R \subset K$ be an $r$-dimensional submodule and let us set $r_i={dim}(R \cap K_i)$, so that $r=r_0 \geq r_1 \geq r_2 \cdots$.  We claim that the type of $I_{\nu}/R$ is completely determined by the integers $r_i$. In fact, we claim that $I_{\nu}/R$ is isomorphic to $I_{\mu}$ where if $\nu=(1^{l_1}, 2^{l_2}, \ldots, n^{l_n})$ then 
\begin{equation}\label{E:hallpolform}
\mu=(1^{l_1+2r_1-r_0-r_2}, 2^{l_2+2r_2-r_1-r_3}, \ldots, n^{l_n+r_n-r_{n-1}}).
\end{equation}
To see this, observe that a representation $(W,y)$ is isomorphic to $I_{\lambda}$ where $\lambda=(1^{m_1}, 2^{m_2}, \ldots,n^{m_n} )$ if and only if ${dim}({Ker}\;y^{i+1})-{dim}({Ker}\;y^i)=m_{i+1} + \cdots + m_n$ for all $i \geq 0$. Note that if $y \in {End}(I_{\nu}/R)$ denotes the operator induced by $x$ then
\begin{align*}
{dim}({Ker}\; y)&={dim}({Ker}\;x)-{dim}(R)+{dim}(R \cap {Im}\;x)=l_1+\ldots l_n -r_0+r_1\\
{dim}({Ker}\; y^2)&={dim}({Ker}\;x^2)-{dim}(R)+{dim}(R \cap {Im}\;x^2)\\
&=l_1+2(l_2+\ldots l_n) -r_0+r_2
\end{align*}
and in general
\begin{equation*}
\begin{split}
{dim}({Ker}\; y^i)&={dim}({Ker}\;x^i)-{dim}(R)+{dim}(R \cap {Im}\;x^i)\\
&=l_1+2l_2 + \cdots + (i-1)l_{i-1} + i(l_i+\ldots l_n) -r_0+r_i.
\end{split}
\end{equation*}
Formula (\ref{E:hallpolform}) now follows by a simple computation.

To sum up, we have shown that for a given $\mu$ there are as many possible choices of a valid submodule $R$ as there are $r$-dimensional subspaces $R$ in $K$ which intersect the flag of subspaces $K_1 \supset K_2 \supset \cdots \supset K_{n-1} \supset K_n=\{0\}$ in some certain \textit{specified} dimensions ${dim}(R \cap K_i)=r_i$ (which can be obtained from $\mu$ using (\ref{E:hallpolform})). The cardinality of this set can be computed as follows. Let ${O} \subset
Gr(r,\sum l_i)$ be the above subset of $r$-planes in $K$. There is a natural projection map
$$\pi: O \to Gr(r_{n-1},l_n) \times Gr(r_{n-2}-r_{n-1},l_{n-1}) \times \cdots \times Gr(r_0-r_1,l_1)$$
which associates to $R$ the sequence of subspaces 
$$R \cap K_{n-1}, (R \cap K_{n-2})/(R \cap K_{n-1}), \ldots,  (R \cap K_{1})/(R \cap K_2), R/(R \cap K_1)$$ 
in $K_{n-1}, K_{n-2}/K_{n-1}, \ldots, K_1/K_2, K/K_1$. We claim that $\pi$ is an affine fibration of rank
\begin{equation}\label{E:rankpi}
\begin{split}
t=&(r_{n-2}-r_{n-1})(l_n-r_{n-1}) + (r_{n-3}-r_{n-2})(l_{n-1}+l_n-r_{n-2}) + \cdots \\
&+ (r_0-r_1)(l_2+ \cdots +l_n-r_1)
\end{split}
\end{equation}
To see why this is true, let us pick arbitrary subspaces $S_{n-1}, \ldots, S_1, S_0$ in the above product of Grassmanians and look at the fiber of $\pi$ at that point. For $R$ to belong to it, we need to have in particular $R \cap K_{n-1}=S_{n-1}$ and $(R \cap K_{n-2})/(R \cap K_{n-1})=S_{n-2} \in K_{n-2}/K_{n-1}$. Fixing a basis of $S_{n-1}$ and $S_{n-2}$ we see that $R \cap K_{n-2}$ is determined up to a choice of a linear map from $S_{n-2}$ to $K_{n-1}/(R \cap K_{n-1})$. But then once $(R \cap K_{n-2})$ and $S_{n-3}$ are fixed, $(R \cap K_{n-3})$ is determined up to a linear map from $S_{n-3}$ to $K_{n-2} / (R \cap K_{n-2})$, and so on. Continuing in this manner, we arrived at the desired result.

Finally, since the number of points over $\mathbb{F}_q$ of the Grassmaninan $Gr(r,t)$ is $\begin{bmatrix} t \\ r \end{bmatrix}_+$, we obtain \textit{in fine}
$$P^{\nu}_{\mu, (1^r)}=\begin{bmatrix} l_n \\ r_{n-1} \end{bmatrix}_+ \cdot \begin{bmatrix} l_{n-1} \\ r_{n-2}-r_{n-1} \end{bmatrix}_+ \cdots \begin{bmatrix} l_1 \\ r_0-r_1 \end{bmatrix}_+ q^t$$
where $\nu=(1^{l_1}, \ldots, n^{l_n})$, the integers $r_i$ are determined by (\ref{E:hallpolform}) and $t$ is given by (\ref{E:rankpi}). 
\endexample

\vspace{.15in}

All the above examples suggest that $P^{\nu}_{\mu,\lambda}$ is in fact given by the evaluation at $t=q$ of some ``universal'' polynomial $P_{\mu,\lambda}^{\nu}(t) \in \Z[t]$. Here ``universal'' means that $P_{\mu,\lambda}^{\nu}(t) \in \Z[t]$ is independent of the choice of the ground (finite) field $k$. This is indeed the case, as we will see at the end of the next section.

\vspace{.2in}

\centerline{\textbf{2.3. Steinitz's classical Hall algebra.}}
\addcontentsline{toc}{subsection}{\tocsubsection {}{}{\; 2.3. Steinitz's classical Hall algebra.}}

\vspace{.15in}

\paragraph{} Let $\H_{cl}=\H_{Rep^{nil}_k\vec{Q}_0}$ be the Hall algebra\footnote{the index `cl' stands for `classical'.} of $Rep^{nil}_k\vec{Q}_0$. By definition, it has a basis $\{[I_{\lambda}]\;|;\; \lambda \in \Pi\}$ and relations
$$[I_{\mu}] \cdot [I_{\lambda}]=\sum_{\nu} P^{\nu}_{\mu,\lambda} [I_{\nu}]$$
(recall that the Euler form vanishes). It is $\N$-graded (by the dimension of the representation, which is also the class in the Grothendieck group). As the Euler form vanishes, the Hall algebra $\H_{cl}$ is already a bialgebra as it stands, and there is no need to consider its extension $\widetilde{\H}_{cl}$. Moreover, the finite subobjects condition being verified, $\H_{cl}$ is in fact a genuine (and not only a topological) Hopf algebra. The first few values of the (co)multiplication and the antipode are given in the following examples~:

\vspace{.15in}

\addtocounter{theo}{1}
\paragraph{\textbf{Example \thetheo.}} Using the computations of the previous section, we have~:
$$[I_{(1)}]  \cdot [I_{(1)}]=(q+1) [I_{(1^2)}] + [I_{(2)}],$$
$$[I_{(1)}] \cdot [I_{(2)}]=q [I_{(2,1)}] + [I_{(3)}],$$
$$ [I_{(2)}] \cdot [I_{(1)}]=q [I_{(2,1)}] + [I_{(3)}]. $$
For the coproduct, we have
$$\Delta ([I_{(1)}])=1 \otimes [I_{(1)}] + [I_{(1)}] \otimes 1,$$
$$\Delta([I_{(2)}])=1 \otimes [I_{(2)}] + (1-q^{-1})[I_{(1)}] \otimes [I_{(1)}] +[I_{(2)}] \otimes 1,$$
$$\Delta([I_{(1^2)}])=1 \otimes [I_{(1^2)}] + q^{-1}[I_{(1)}] \otimes [I_{(1)}] +[I_{(1^2)}] \otimes 1,$$
\begin{equation*}
\begin{split}
\Delta([I_{(2,1)}])=& 1 \otimes [I_{(2,1)}]+ [I_{(2,1)}] \otimes 1 + (1-q^{-2}) \big( [I_{(1)}] \otimes [I_{(1^2)}] + [I_{(1^2)}] \otimes [I_{(1)}] \big) +\\
&+ q^{-1} \big( [I_{(1)}] \otimes [I_{(2)}] + [I_{(2)}] \otimes [I_{(1)}]\big).
\end{split}
\end{equation*}
In the above computations, we have used the following formulas
$$a_{I_{(1)}}=(q-1), \quad a_{I_{(1^2)}}=(q^2-1)(q^2-q), \quad a_{I_{(2)}}=(q^2-q), \quad a_{I_{(2,1)}}=(q^3-q^2)(q^2-q).$$
Finally, for the antipode map, we have
$$S([I_{(1)}])=-[I_{(1)}],$$
$$S([I_{(1^2)}])=-[I_{(1^2)}] + (q+1) \frac{(q-1)^2}{(q^2-1)(q^2-q)} [I_{(1)}]^2=q^{-1} [I_{(2)}] + q^{-1}[I_{(1^2)}],$$
$$S([I_{(2)}])=-[I_{(2)}] + \frac{(q-1)^2}{(q^2-q)} [I_{(1)}]^2=-q^{-1} [I_{(2)}] +(q- q^{-1})[I_{(1^2)}].$$
\endexample

\vspace{.15in}

Despite these seemingly complicated formulas, the (abstract) structure of $\H_{cl}$ is startingly simple~:

\begin{theo}[Steinitz, Hall, Macdonald]\label{T:Steinitz} The following hold~:
\begin{enumerate}
\item[i)] $\H_{cl}$ is both commutative and cocommutative,
\item[ii)] $\H_{cl} \simeq \C[[I_{(1)}], [I_{(1^2)}], \ldots ]$ is a free polynomial algebra in the infinitely many generators $[I_{(1)}], [I_{(1^2)}], \ldots$.
\item[iii)] For any integer $n$ we have 
\begin{equation}\label{E:coupcoup}
\Delta ([I_{(1^n)}])=\sum_{r=0}^n q^{-r(n-r)} [I_{(1^r)}] \otimes [I_{(1^{n-r})}].
\end{equation}
\end{enumerate}
\end{theo}
\noindent
\textit{Proof.} To prove point i), we have to show that for any triple $\nu, \mu, \lambda$, we have
$P^{\nu}_{\mu,\lambda}=P^{\nu}_{\lambda,\mu}$. For this we will use the natural (\textit{contravariant}) duality functor in the category $Rep^{nil}_k\vec{Q}_0$, which associates to a representation $(V,x)$ its transpose
$(V^*, x^t)$. Note that any object of $Rep^{nil}_k\vec{Q}_0$ is isomorphic to its transpose. On the other hand, this duality $M \mapsto M^*$ gives rise to a bijection
\begin{align*}
\{R \subset I_{\nu}\;|\; R \simeq I_{\lambda},\; I_{\nu}/R \simeq I_{\mu} \}&\stackrel{1:1}{\longrightarrow}
\{S \subset I_{\nu}^* \simeq I_{\nu}\;|\; I_{\nu}^*/S \simeq I_{\lambda}^* \simeq I_{\lambda} \simeq,\; S \simeq I_{\mu}^* \simeq I_{\mu}\}\\
R & \mapsto S =R^{\perp}
\end{align*}
It now suffices to recall that $P^{\nu}_{\mu,\lambda}$ is the cardinality of the left-hand side while $P^{\nu}_{\lambda,\mu}$ is the cardinality of the right-hand side.\\
We turn to the second statement. Define a partial order on the set of partitions $\Pi$ as follows~:
$\lambda=(1^{l_1}, 2^{l_2}, \ldots ) \succeq \mu=(1^{m_1}, 2^{m_2}, \ldots)$
if for any $i \geq 1$ we have 
\begin{equation*}
\begin{split}
l_1+2l_2 + \cdots + (i-1)l_{i-1}&+ i(l_i + l_{i+1} + \cdots)\\
& \leq m_1+2m_2 + \cdots + (i-1)m_{i-1}+ i(m_i + m_{i+1} + \cdots)
\end{split}
\end{equation*}
(this is the \textit{transpose} of the usual partial order on partitions). Let us fix a partition $\lambda=(1^{l_1}, 2^{l_2}, \ldots, n^{l_n})$ and let us consider the product 
$$X=[I_{(1^{l_n})}] \cdot [I_{(1^{l_{n-1}+l_n})}] \cdots [I_{(1^{l_1+ \cdots + l_n})}].$$
By construction, if $I_{\nu}$ appears in $X$ with a nonzero coefficient and if $(V,x)$ represents $I_{\nu}$ then there exists a filtration $\{0\}=V_0 \subset V_1 \subset \cdots \subset V_n=V$ such that ${dim}(V_i/V_{i-1})=l_i+ \cdots + l_n$ and $x(V_i) \subset V_{i-1}$. In that situation, ${dim}({Ker}\;x^i) \geq {dim}\;V_i=l_1 + 2l_2 + \cdots + i(l_i + l_{i+1} + \cdots)$. But writing $\nu=(1^{m_1},2^{m_2}, \ldots)$ we have ${dim}({Ker}\;x^i)= m_1 + 2m_2 + \cdots + i(m_i + m_{i+1} + \cdots)$, hence $\nu \preceq \lambda$. In addition, it is easy to see that if $\nu=\lambda$ then there exists a unique admissible filtration $(V_i)$ as above. Therefore we have
$$X \in [I_{\lambda}] + \oplus_{\nu \prec \lambda} \C [I_{\nu}].$$
We deduce from this that the collection $\mathcal{I}$ of products $[I_{(1^{n_1})}] \cdot [I_{(1^{n_2})}] \cdots [I_{(1^{n_t})}]$ (ordered so that $n_1 \leq n_2 \leq \cdots \leq n_t$) is obtained from the basis $\{[I_{\lambda}]\}$ by acting by a matrix which is upper triangular with respect to $\preceq$ and has $1$'s along the diagonal. Therefore $\mathcal{I}$ forms a basis of $\H_{cl}$ and point ii) follows.\\
Finally, iii) is shown by a direct calculation. Since any submodule or quotient of $I_{(1^n)}$ is again of the form $I_{(1^r)}$ for some $r \leq n$, we have
$$\Delta([I_{(1^n)}])=\sum_{r =0}^n \frac{a_{I^{(1^{n-r})}} a_{I_{(1^r)}}}{a_{I_{(1^n)}}} P^{(1^n)}_{(1^{n-r}),(1^r)} [I_{(1^{n-r})}] \otimes [I_{(1^r)}].$$
The formulas
$$a_{I_{(1^u)}}=| GL(u,\mathbb{F}_q)|=(q^u-1)(q^u-q) \cdots (q^u-q^{u-1}) = (q-1)^u q^{\frac{u (u-1)}{2}} [u] [u-1] \cdots [1],$$
$$P^{(1^n)}_{(1^{n-r}),(1^r)}=\begin{bmatrix} n\\r \end{bmatrix}_+=\frac{[n]_+ [n-1]_+ \cdots [1]_+}{([r]_+ [r-1]_+ \cdots [1]_+) ([n-r]_+ [n-r-1]_+ \cdots [1]_+)}$$
yield the desired result.\qed

\vspace{.15in}

The above Theorem is enough to prove the existence of the so-called \textit{Hall polynomials} $P^{\nu}_{\mu,\lambda}(t)$ alluded to at the end of Section~2.2.
 
\begin{prop}\label{P:Hallpols} For any triple of partitions $\lambda,\mu,\nu$ satisfying $|\lambda|+|\mu|=|\nu|$ there exists a unique polynomial $P_{\mu,\lambda}^{\nu}(t) \in \Z[t]$ such that for any prime power $q$ and finite field $k$ with $q$ elements we have $P_{\mu,\lambda}^{\nu}=P_{\mu,\lambda}^{\nu}(q)$.\end{prop}

\noindent
\textit{Proof.} Note that the unicity is obvious since any two polynomials taking the same values at an infinite set of points are equal. Hence it is enough to show existence. We start by slightly refining the proof of point ii) of Theorem~\ref{T:Steinitz} above. By Example~2.4. we know that for any $\mu \in \Pi$ and $r \geq 1$ there exists polynomials $P^{\nu}_{\mu,(1^r)}(t) \in \Z[t]$ such that
$$[I_{\mu}] \cdot [I_{(1^r)}] =\sum_{\nu} P^{\nu}_{\mu,(1^r)}(q) [I_{\nu}].$$
By iteration we deduce that for any $r_1, \ldots, r_t$ there exists polynomials $P^{\nu}_{(1^{r_t}), \ldots (1^{r_1})}(t)$ such that
\begin{equation}\label{E:proofhall1}
[I_{(1^{r_t})}] \cdots [I_{(1^{r_1})}] =\sum_{\nu} P^{\nu}_{(1^{r_t}), \ldots,(1^{r_1})}(q) [I_{\nu}].
\end{equation}
In particular, if $\lambda=(1^{l_1}, 2^{l_2} ,\ldots, n^{l_n})$ and $r_i=l_i + \cdots + l_n$ then
$$[I_{(1^{r_t})}] \cdots [I_{(1^{r_1})}] =[I_{\lambda}] +\sum_{\nu \prec \lambda} P^{\nu}_{(1^{r_t}), \ldots,(1^{r_1})}(q) [I_{\nu}].$$
Inverting the upper triangular matrix \textit{of polynomials} $(P^{\nu}_{(1^{r_t}), \ldots,(1^{r_1})}(t))$ we see that conversely for any $\nu$ there exists polynomials $Q^{\nu}_{(1^{r_t}), \ldots,(1^{r_1})}(t)$ such that
\begin{equation}\label{E:proofhall2}
[I_{\nu}]=\sum_{(r_i)} Q^{\nu}_{(1^{r_t}), \ldots,(1^{r_1})}(q) [I_{(1^{r_t})}] \cdots [I_{(1^{r_1})}] ,
\end{equation}
where the sum ranges over all possible sets of ordered indices $(r_1 \geq r_2 \cdots \geq r_t)$. But now, using (\ref{E:proofhall2}) we may rewrite any product $[I_{\mu}]\cdot [I_{\lambda}]$ as a linear combination with \textit{polynomial} coefficients of terms of the form $[I_{(1^{r_t})}] \cdots [I_{(1^{r_1})}]$. Then we use (\ref{E:proofhall1}) to rewrite each such element as a linear combination of $[I_{\nu}]$'s, again with polynomial coefficient.  We are done. \qed

\vspace{.15in}

An important corollary of the existence of Hall polynomials is that we may now consider a ``universal'', or ``generic'' version $\underline{\H}_{cl}$ of the Hall algebra $\H_{cl}$, which is defined over the ring $\C[t,t^{-1}]$.

\vspace{.1in}

\noindent
We now define the \textit{generic classical Hall algebra} to be the algebra 
$$\underline{\H}_{cl}=\bigoplus_{\lambda} \C[t,t^{-1}] [I_{\lambda}]$$
where the multiplication is defined by
$$[I_{\mu}] \cdot [I_{\lambda}]=\sum_{\nu} P^{\nu}_{\mu,\lambda}(t) [I_{\nu}].$$

\vspace{.15in}

The algebra $\underline{\H}_{cl}$ is clearly commutative, and the proof of Proposition~\ref{P:Hallpols} shows that it is a free polynomial ring over $\C[t,t^{-1}]$ in the generators $[I_{(1)}], [I_{(1^2)}], \ldots$. We may define a coproduct $\Delta$ on $\underline{\H}_{cl}$ by (\ref{E:coupcoup}). Indeed, this satisfies all the required properties (coassociativity, compatibility with the multiplication, $\ldots$ ) over $\C[t,t^{-1}]$ since the same is true for an infinite number of specializations $t=q$. Finally, we leave it as an exercise to the reader to check the formula 
$$S([I_{(1^n)}])=(-1)^n q^{-n(n-1)/2}\mathbf{1}_n,$$
where $\mathbf{1}_n$ is the characteristic function of the set of all representations of dimension $n$; this allows us to 
also define an antipode map $S$ for $\underline{\H}_{cl}$. We have obtained in this manner a genuine Hopf algebra $\underline{\H}_{cl}$ over $\C[t,t^{-1}]$. Note that the product is
defined at $t=0$, but not the coproduct (or the antipode).

\vspace{.15in}

There is one last piece of the structure of Hall algebras which we didn't mention here so far~: Green's scalar product. Recall that by definition, it is given by the formula
$$([I_{\lambda}],[I_{\mu}])=\delta_{\lambda,\mu} \frac{1}{a_{I_{\lambda}}}.$$

\vspace{.1in}

\begin{lem} The number $a_{I_{\lambda}}$ of automorphism of $I_{\lambda}$ is given by the formula
\begin{equation}\label{E:HLnorm}
a_{I_{\lambda}}=q^{|\lambda| + 2 n(\lambda)}\prod_i (1-q^{-1}) (1-q^{-2}) \cdots (1-q^{-l_i})
\end{equation}
where $\lambda=(1^{l_1}, 2^{l_2}, \ldots , m^{l_m})$ and $n(\lambda)=\sum_i (i-1) \lambda_i$.
\end{lem}

\noindent
\textit{Proof.} We have $I_{\lambda}=I_1^{\oplus l_1} \oplus I_2^{\oplus l_2} \cdots$. Let $x^{(i)}_1, x^{(i)}_2, \ldots , x^{(i)}_{l_i}$ be fixed generators of $I_i^{\oplus l_i}$. An endomorphism $f$ of $I_{\lambda}$ is completely determined by the image under it of the elements $x^{(i)}_j$. For $f$ to be an automorphism, the restriction of $f$ to $I_i^{\oplus l_i}$ composed with the projection to $I_i^{\oplus l_i}$ should be an automorphism. In addition, for $f$ to be well-defined, we should have $f(x_j^{(i)}) \in {Ker}\; x^i$. Thus there are
\begin{equation}\label{E:automo1}
\begin{split}
 &|GL(l_i, \mathbb{F}_q)| \cdot \bigg| \left( \bigg\{(\text{Ker}\; x^i) \cap \bigoplus_{j \neq i} I_j^{\oplus l_j}) \oplus (\text{Ker}\; x^{i-1}) \cap I_i^{\oplus l_i})\bigg\} \otimes k^{l_i} \right)\bigg| \\
&=|GL(l_i, \mathbb{F}_q)| q^{l_i(l_1 + 2 l_2 + \cdots + (i-1)l_{i-1} + (i-1) l_i + i(l_{i+1} + \cdots))}
\end{split}
\end{equation}
possible choices for $(f(x_1^{(i)}), \ldots , f(x_{l_i}^{(i)}))$. Hence the number of automorphism is obtained by multiplying together (\ref{E:automo1}) for all values of $i$. A little arithmetic on partitions
using the identities
$$|\lambda|=\sum_i i l_i, \qquad n(\lambda)=\frac{1}{2}\sum_i i l_i (l_i-1) + \sum_{i<j} i l_i l_j$$ 
brings this expression to the form (\ref{E:HLnorm}).\qed

\vspace{.2in}

\centerline{\textbf{2.4. Link with the ring of symmetric functions.}}
\addcontentsline{toc}{subsection}{\tocsubsection {}{}{\; 2.4. Link with the ring of symmetric functions.}}

\vspace{.15in}

\paragraph{}As we have seen, the classical Hall algebra $\H_{cl}$ , or its generic version $\underline{\H}_{cl}$, are $\N$-graded polynomial rings with one generator in every degree. There is another famous such, which plays a very important role in much of mathematics, and particularly in combinatorics and representation theory~: Macdonald's ring of symmetric functions (\cite{Mac}). Let us briefly recall its definition.

Consider, for $n \geq 1$, the ring of symmetric polynomials in $n$ variables~: $\Lambda_n=\C[x_1, \ldots, x_n]^{\mathfrak{S}_n}$. These form a projective system via the maps $\Lambda_{n+1} \to \Lambda_n$ obtained by setting the last variable $x_{n+1}$ to zero. The projective limit (in the category of graded rings) $\LLambda=\varprojlim \Lambda_n$ can thus be considered as the ring ``$\C[x_1, x_2, \ldots ]^{\mathfrak{S}_{\infty}}$'' of symmetric functions in infinitely many variables. 

The ring $\LLambda$ is also equipped with a canonical coproduct, which was formally introduced by Zelevinsky \cite{Zelevinsky} ~: for $n \geq 1$ consider the map $\Delta_{n}: \Lambda_{2n} \to \Lambda_n \otimes \Lambda_n$ induced by the embedding 
\begin{equation*}
\begin{split}
\C[x_1, \ldots, x_{2n}]^{\mathfrak{S}_{2n}} &\hookrightarrow \C[x_1, \ldots, x_{2n}]^{\mathfrak{S}_{n}\times \mathfrak{S}_n}\\
&= \C[x_1, \ldots, x_{n}]^{\mathfrak{S}_n} \otimes \C[x_{n+1}, \ldots ,x_{2n}]^{\mathfrak{S}_n},
\end{split}
\end{equation*}
where in the second term the first copy of $\mathfrak{S}_n$ permutes together the variables $x_1, \ldots, x_n$ while the second copy of $\mathfrak{S}_n$ permutes together the variables $x_{n+1}, \ldots, x_{2n}$. In the projective limit, the maps $\Delta_n$ give rise to a coproduct $\Delta: \LLambda \to \LLambda \otimes \LLambda$.

\vspace{.1in}

There are many bases of $\LLambda$, and it is often the matrix relating one such basis to another which carries the interesting combinatorial information. The simplest of all these bases is probably the basis of elementary symmetric functions~: for $r \in \N$ set
$$e_r=\sum_{i_1 < i_2 < \cdots < i_r} x_{i_1} x_{i_2} \cdots x_{i_r} \in \LLambda$$
and for a partition $\lambda=(\lambda_1, \lambda_2, \ldots)$ put $e_{\lambda}=e_{\lambda_1} e_{\lambda_2} \cdots e_{\lambda_n}$.

\vspace{.1in}

\begin{theo}[Macdonald, \cite{Mac}] The set $\{e_{\lambda}\;|\; \lambda \in \Pi\}$ forms a basis of $\LLambda$, i.e. $\LLambda \simeq \C[e_1, e_2, \ldots]$.\end{theo}

\vspace{.1in}

By the above Theorem, we may construct an algebra isomorphism $\Phi_{(q)}: \H_{cl} \to \LLambda$ by imposing
$\Phi_{(q)}([I_{(1^r)}])=q^{-\frac{r (r-1)}{2}}e_r$.  From the definition, it is easy to check that
$\Delta(e_r)=\sum_{r=0}^n e_{n-r} \otimes e_r$. From this and from Theorem~\ref{T:Steinitz} iii) it follows that $\Phi$ is a morphism of bialgebras as well. Hence $\H_{cl}$ provides a new model or realization of $\LLambda$. Moreover, it comes almost for free with a very canonical one-parameter deformation $\underline{\H}_{cl}$ defined over $\C[t,t^{-1}]$ and there is an isomorphism $\Phi: \underline{\H}_{cl} \to \LLambda \otimes \C[t,t^{-1}]$.

\vspace{.1in}

Let us mention two immediate applications~: we can use $\Phi$ to transport on $\LLambda \otimes \C[t,t^{-1}]$ the nondegenerate Hopf pairing $(\,,\,)$ of Green on $\underline{\H}_{cl}$; and the set $\{\Phi([I_{\lambda}])\;|\; \lambda \in \Pi\}$ is an orthogonal basis for this scalar product. These turn out to be very interesting~: $\Phi_{\star}((\,,\,))$ is the Hall-Littlewood scalar product, which is the scalar product uniquely determined by the conditions
$$\{x,yz\}=\{\Delta(x) ,y\otimes z\},$$
\begin{equation}\label{E:HLscalarprod}
\{p_r,p_s\}=\delta_{r,s}\frac{r}{q^r-1}
\end{equation}
where $p_r=\sum_i x_i^r$ stands for the power sum symmetric function. Therefore the basis $\{\Phi([I_{\lambda}])\;|\; \lambda \in \Pi\}$ is (up to a harmless renormalization) the basis of Hall-Littlewood polynomials. In particular, the norm of these polynomials (with respect to the Hall-Littlewood scalar product) is given by (\ref{E:HLnorm}).

 There are many more applications of the classical Hall algebra to the theory of symmetric functions. As discussing these would distract us too much from our topic, we prefer to refer the interested reader directly to the Scriptures \cite{Mac}.

\vspace{.2in}

\addtocounter{theo}{1}
\noindent \textbf{Remark \thetheo .}  One may wonder why the maps $\Phi_{(q)}$ and $\Phi$ are defined in this way and not in another. The first answer is that this is the easiest way to obtain a morphism of bialgebras. A better answer is that there exists a ``geometric lift'' of the map $\Phi_{(q)}$ which makes it indeed very canonical; this involves (affine) flag varieties and (spherical) Hecke algebras --see the survey \cite{SToronto} for details.

\vspace{.2in}

\centerline{\textbf{2.5. Other occurences of Hall algebras.}}
\addcontentsline{toc}{subsection}{\tocsubsection {}{}{\; 2.5. Other occurences of Hall algebras.}}

\vspace{.15in}

\paragraph{} In this final section, which is independent of the rest of these notes, we quickly describe a few \textit{classical} mathematical contexts in which Hall algebras arise naturally.

\vspace{.2in}

\noindent
\textit{Steinitz's formulation of the classical Hall algebra.} Needless to say, neither quivers nor abelian categories were around in the early days of the twentieth century. Instead, Steinitz considered the set of all abelian $p$-groups; of course, these share many properties with nilpotent representations of the Jordan quiver (over $\mathbb{F}_p$), namely there is a unique simple object and a unique indecomposable object of any given length. In fact, it is easy to see that the Hall numbers computed from both categories coincide, and hence Sections~2.2. and 2.3. are essentially a direct reformulation of Steinitz's work. A third alternative formulation, which is used by Macdonald in \cite{Mac}, is to consider the category of finite length modules over some discrete valuation ring $R$ whose residue field $R/\mathfrak{m}$ is finite.

\vspace{.2in}

\noindent
\textit{Parabolic induction for $GL(n)$.} For $n \geq 1$, put $G_n=GL(n,\mathbb{F}_q)$, and let $R_n$ be the character ring of $G_n$, i.e. the ring of class functions on $G_n$. The vector space $R=\bigoplus_n R_n$ has a ring structure given by the parabolic induction~: if $f \in R_n$ and $g \in R_m$ then we put
$$f \circ h=ind_{P_{n,m}}^{G_{n+m}} (f \otimes h),$$
where $P_{n,m} \subset G_{n+m}$ is the maximal parabolic subgroup of type $(n,m)$, and $f \otimes h$ is pulled back to $P_{n,m}$ by the projection $P_{n,m} \tto G_n \times G_m$. By definition, the value of $f \circ h$ a at class $\mu \in G_{n+m}$ is equal to $\sum_{t \in G_{n+m}/P_{n,m}} f \otimes h(t \mu t^{-1})$.

For any $x \in G_l$, let $V_x$ be the $\mathbb{F}_q[t]$-module structure on $\mathbb{F}_q^l$ obtained by having $t$ act as $x$. Clearly, $V_x \simeq V_y$ if and only if $x$ and $y$ are conjugate. Hence there is a well-defined isomorphism class of $\mathbb{F}_q[t]$-module associated to any conjugacy class $\mu$ in $G_l$, for any $l$. Using this notation, we may write
$$f \circ h (\mu)=\sum_{\mu_1,\mu_2} g_{\mu_1,\mu_2}^\mu f(\mu_1) \cdot h(\mu_2)$$
where $g_{\mu_1,\mu_2}^{\mu}$ is the number of $\mathbb{F}_q[t]$-submodules $W \subset V_{\mu}$ such that $W \simeq V_{\mu_2}$ and $V_{\mu}/W \simeq V_{\mu_1}$. Thus, after a small reformulation, parabolic induction for the groups $GL(l,\mathbb{F}_q)$ may be considered as a simple, special case of Hall algebra multiplication. There is a similar interpretation for parabolic induction for the groups $GL(l,\mathbb{Q}_p)$.

\vspace{.2in}

\noindent
\textit{Automorphic forms for function fields.}  This may be viewed as a global analog of the previous example. Let $X$ be a smooth projective curve defined over a finite field $\mathbb{F}_q$ and let $Bun_{r}(X)$ denote the set of isomorphism classes of vector bundles on $X$ of rank $r$. A complex-values function $f$ on $Bun_r(X)$ may be interpreted as an automorphic form for the group $GL(r)$ over the function field $\mathbb{F}_q(X)$ of $X$. In this setting, the induction product for automorphic forms coincides with the product in the Hall algebra $\H_{Vec(X)}$ of the exact category $Vec(X)$ of vector bundles on $X$. Many properties of automorphic forms may be stated in a simple manner using the language of Hall algebras. For instance, an automorphic form $f$ is a cusp form if $\Delta(f) \in \H_{Vec(X)} \otimes \H_{Tor(X)} \oplus \H_{Tor(X)} \otimes \H_{Vec(X)}$, where $Tor(X)$ denotes the category of torsion sheaves on $X$; the Hecke operators on the space of automorphic forms are given by the adjoint action by certain elements of $\H_{Tor(X)}$, etc. We refer to \cite{Kap1} for much more in this direction. 

\vspace{.15in}

\newpage

\centerline{\large{\textbf{Lecture~3.}}}
\addcontentsline{toc}{section}{\tocsection {}{}{Lecture~3.}}

\setcounter{section}{3}
\setcounter{theo}{0}
\setcounter{equation}{0}

\vspace{.15in}

In this Lecture, we describe the main examples of Hall algebras~: those associated to the categories of representations of quivers over finite fields. This finds its source in Gabriel's insight that the set of indecomposable representations of finite type quivers can be parametrized by root systems of simple Lie algebras \cite{Gabriel}, in Kac's generalization of this result to arbitrary quivers \cite{Kac}, and culminates in Ringel's discovery \cite{Ri} (later completed by Green \cite{Green}) that the Hall algebra of a quiver $\vec{Q}$ contains a copy of the quantized enveloping algebra of the Kac-Moody algebra associated to the graph $Q$ underlying $\vec{Q}$. After recalling Gabriel's and Kac's theorems and some representation theory of quivers, we present the fundamental results of Ringel and Green. The last three Sections are devoted to some more advanced topics concerning tame quivers (which we mostly state without proofs). These are important for Lecture~4.
For the reader's convenience,  just enough of the structure theory of Kac-Moody Lie algebras and quantum groups is reviewed in Appendices~A.1 through A.4. As for the representation theory of quivers, all of the results used below may be found in \cite{Barot} or \cite{CrawleyNotes}.

\vspace{.2in}

\centerline{\textbf{3.1. Quivers.}} 
\addcontentsline{toc}{subsection}{\tocsubsection {}{}{\; 3.1. Quivers.}}

\vspace{.15in}

\paragraph{}Let $\vec{Q}$ be a quiver with vertex set $I$ and edge set $H$. We allow $\vec{Q}$ to have multiple edges and cycles, but no loops. The target, resp. source of an edge $h$ will be denoted $t(h)$, resp. $s(h)$. By a representation of $\vec{Q}$ over a field $k$ we mean a pair $(V,\underline{x})$ where $V=\bigoplus_i V_i$ is a finite-dimensional $I$-graded vector space and $\underline{x}=(x_h)_{h \in H} \in \bigoplus_h {Hom}(V_{s(h)},V_{t(h)})$. A representation $(V,\underline{x})$ is \textit{nilpotent} if there exists $N \gg 0$ such that for any $n >N$, $x_{h_n}x_{h_{n-1}} \cdots x_{h_1}=0$ for any path $h_1 \cdots h_n$ of length $n$ in $\vec{Q}$. We let $Rep^{nil}_k\vec{Q}$ stand for the category of nilpotent representation of $\vec{Q}$ over $k$. It is an abelian, $k$-linear category satisfying the Krull-Schmidt property. Of course, when $\vec{Q}$ has no oriented cycles then any representation is nilpotent.

\vspace{.1in}

To any vertex is attached a simple representation $S_i$ such that $V_i =k$, $V_j=\{0\}$ for $j\neq i$ and $\underline{x}=0$.

\vspace{.1in}

\begin{prop} The collection $\{S_i\;|\; i \in I\}$ is a complete set of simple objects of $Rep^{nil}_k\vec{Q}$. \end{prop}
\noindent
\textit{Proof.} Let $(V,\underline{x})$ be simple and let $i_0$ be any vertex such that $V_i \neq \{0\}$.
Let $v \in V_{i_0}$ be a nonzero vector. By the nilpotency condition, there exists a path $h_1 \cdots h_n$ starting at $i_0$ such that $x_{h_n} \cdots x_{h_1} \cdot v \neq 0$ but $x_{h_{n+1}}x_{h_n} \cdots x_{h_1} \cdot v = 0$ for all edges $h_{n+1}$ leaving the terminal vertex, say $i_1$ of $h_n$. This means that $S_{i_1} \subset (V,\underline{x})$ . But since $V$ is simple we deduce that $V=S_{i_1}$.\qed

\vspace{.1in}

\begin{cor} The following hold~:
\begin{enumerate}
\item[i)] Any representation $(V,\underline{x})$ in $Rep^{nil}_k\vec{Q}$ admits a finite composition series consisting of $S_i$'s,
\item[ii)] We have $K(Rep^{nil}_k\vec{Q}) \simeq \Z^I$.
\end{enumerate}
\end{cor}

\vspace{.1in}

The class of a representation $(V,\underline{x})$ in $K(Rep^{nil}_k\vec{Q})$ is simply its dimension vector $\underline{{dim}}\; V=({dim}\;V_i)_{i \in I}$.

\vspace{.1in}

It is well-known that the category $Rep^{nil}_k\vec{Q}$ is hereditary.
As a consequence, we may easily compute the Euler form on $K(Rep^{nil}_k\;\vec{Q})$. Let us assume from now on that $k=\mathbb{F}_q$ is a finite field. Let 
$$c_{ij}=\#\{h\in H\;|\; s(h)=i, t(h)=j\}$$
be the number of oriented edges in $\vec{Q}$ going from $i$ to $j$. 

\begin{prop} We have $\langle S_i, S_j \rangle_a=\delta_{ij}-c_{ij}$ and
$$\langle S_i,S_j \rangle_m=q^{\frac{1}{2}(\delta_{ij}-c_{ij})}.$$
\end{prop}
\noindent
\textit{Proof.} This is a consequence of the facts that $Rep^{nil}_k\vec{Q}$ is hereditary and that 
${dim}( {Ext}^1(S_i,S_j))=c_{ij}$. \qed

\vspace{.15in}

The matrix of the additive Euler form is $A=(a_{ij})_{i,j \in I}$ with $a_{ij}=2\delta_{i,j} -c_{ij}-c_{ji}$.
In particular, it is a symmetric integral matrix satisfying
$$a_{ii}=2, \qquad a_{ij} \leq 0 \qquad \text{if\;} i \neq j.$$
Therefore $A$ is a generalized (symmetric) Cartan matrix and we may associate to it a Kac-Moody algebra $\g$ (see Appendix~A.1, A.2.). The Dynkin diagram of $\g$ is simply the unoriented graph underlying $\vec{Q}$.

Let $\{\a_i\;|\; i \in I\}$ be the set of simple roots, let $Q=\bigoplus_i \Z \a_i$ be the root lattice and let $(\;,\;)$ be the restriction of the Cartan-Killing form to $Q$.

\vspace{.15in}

\begin{cor}\label{C:KrepQ} The map 
\begin{align*}
\rho~: K(Rep^{nil}_k\vec{Q}) &\to Q\\
\overline{V} &\mapsto  \sum_i \text{dim}(V_i) \a_i
\end{align*}
is an isomorphism of $\Z$-modules. It maps the symmetrized Euler form $(\;,\;)_a$ to the Cartan-Killing form $(\;,\;)$.
\end{cor}

\vspace{.1in}

This is a first, purely numerical, indication of a link between categories of representations of quivers on the one hand and Kac-Moody algebra on the other. This link will be made much stronger in the coming paragraphs.

\vspace{.15in}

For completeness, we add the following theorem, which explains the importance of quivers in the representation theory of finite dimensional algebras.

\vspace{.1in}

\begin{theo}[Gabriel]\label{T:Gabriel} Let $k$ be an algebraically closed field. Any finite dimensional $k$-algebra of global dimension at most one is Morita equivalent to the category $Rep^{nil}_k\vec{Q}$ for some (uniquely determined) quiver $\vec{Q}$.\end{theo}

\vspace{.1in}

\addtocounter{theo}{1}
\noindent \textbf{Remark \thetheo .} There is a version of the above theorem when $k$ is a finite field~: one then has to consider not only quivers, but \textit{species}, or equivalently quivers equipped with automorphisms (see \cite{RingelTame}).

\vspace{.2in}

\centerline{\textbf{3.2. Gabriel's and Kac's Theorems.}}
\addcontentsline{toc}{subsection}{\tocsubsection {}{}{\; 3.2. Gabriel's and Kac's theorems.}}

\vspace{.15in}

\paragraph{} Since the category $Rep^{nil}_k\vec{Q}$ is Krull-Schmidt, the next important invariant (after the Grothendieck group and the Euler form) is the set of indecomposable objects. Call a quiver \textit{of finite type} if it has only finitely many indecomposables, and call it \textit{tame} if the indecomposables lying in each class in $K(Rep^{nil}_k\vec{Q})$ can be arranged in finitely many one-parameter families. A quiver not falling in these two sets is called \textit{wild}.

\vspace{.1in}

\begin{theo}[Gabriel] The category $Rep^{nil}_k\vec{Q}$ is of finite type if and only if $A$ is positive definite, i.e. if and only if $\g$ is a simple Lie algebra. Moreover, in this case the map $M=(V,\underline{x}) \mapsto \sum_i \text{dim}(V_i) \a_i \in Q$ establishes a bijection between the set of indecomposable objects and the set $\Delta^+$ of positive roots of $\g$.\end{theo}

\vspace{.1in}

In other words, once we have identified $K(Rep^{nil}_k\vec{Q})$ and $Q$ using Corollary~\ref{C:KrepQ}, the positive root system $\Delta^+$ pinpoints exactly the classes of indecomposables and there is a unique indecomposable belonging to each such class.

\vspace{.15in}

\addtocounter{theo}{1}

\paragraph{\textbf{Example \thetheo.}} Suppose that $\vec{Q}$ is given by 

\vspace{.15in}

\centerline{
\begin{picture}(250, 10)
\put(60,0){\circle*{5}}
\put(100,0){\circle*{5}}
\put(190,0){\circle*{5}}
\put(75,0){\vector(1,0){5}}
\put(115,0){\vector(1,0){5}}
\put(175,0){\vector(1,0){5}}
\put(57,-10){$1$}
\put(97,-10){$2$}
\put(187,-10){$n$}
\put(60,0){\line(1,0){68}}
\put(135,0){\line(1,0){3}}
\put(145,0){\line(1,0){3}}
\put(155,0){\line(1,0){3}}
\put(165,0){\line(1,0){25}}
\end{picture}}

\vspace{.3in}

The associated Lie algebra is $\g=\mathfrak{sl}_{n+1}(\C)$. If $\{\a_1, \ldots, \a_n\}$ are the simple roots then (see Appendix~A.1, Example~A.5.)
$$\Delta^+=\{\a_{i} + \a_{i+1} + \cdots + \a_{j}\;|\; i \leq j\}$$
These correspond to the indecomposable representations

\vspace{.35in}

\centerline{
\begin{picture}(300, 10)
\put(10,0){$I_{i \cdots j}=$}
\put(60,0){\circle*{2}}
\put(85,0){\circle*{2}}
\put(130,0){\circle*{2}}
\put(160,0){\circle*{2}}
\put(240,0){\circle*{2}}
\put(280,0){\circle*{2}}
\put(138,8){$\simeq$}
\put(168,8){$\simeq$}
\put(218,8){$\simeq$}
\put(135,7){\vector(1,0){20}}
\put(165,7){\vector(1,0){20}}
\put(215,7){\vector(1,0){20}}
\put(57,-10){$1$}
\put(82,-10){$2$}
\put(92,0){\line(1,0){4}}
\put(104,0){\line(1,0){4}}
\put(116,0){\line(1,0){4}}
\put(185,0){\line(1,0){4}}
\put(197,0){\line(1,0){4}}
\put(209,0){\line(1,0){4}}
\put(245,0){\line(1,0){4}}
\put(257,0){\line(1,0){4}}
\put(269,0){\line(1,0){4}}
\put(127,-10){$i$}
\put(157,-10){$i+1$}
\put(237,-10){$j$}
\put(127,5){$k$}
\put(157,5){$k$}
\put(237,5){$k$}
\put(277,-10){$n$}
\end{picture}}

\vspace{.3in}

The case of a quiver with the same underlying graph but different orientation is entirely similar.\endexample

\vspace{.15in}

\addtocounter{theo}{1}

\paragraph{\textbf{Example \thetheo.}} Now let us suppose that $\vec{Q}$ is 

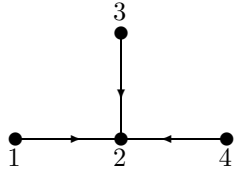
\vspace{.5in}

\centerline{
\begin{picture}(200, 10)
\put(60,-15){\circle*{5}}
\put(100,-15){\circle*{5}}
\put(140,-15){\circle*{5}}
\put(100,25){\circle*{5}}
\put(57,-25){$1$}
\put(97,-25){$2$}
\put(137,-25){$4$}
\put(97,30){$3$}
\put(60,-15){\line(1,0){80}}
\put(100,25){\line(0,-1){40}}
\put(80,-15){\vector(1,0){5}}
\put(120,-15){\vector(-1,0){5}}
\put(100,5){\vector(0,-1){5}}
\end{picture}}

\vspace{.55in}

The associated Lie algebra is now $\mathfrak{so}_8(\C)$. If $\a_1, \ldots, \a_4$ are the simple roots then the positive roots are (see Appendix~A.1, Example~A.6.)
\begin{equation*}
\begin{split}
\Delta^+=& \{\a_i\;|\; i=1, \ldots, 4\}  \cup \{\a_1+\a_2, \;\a_2+\a_3, \;\a_2+\a_4\}\\
&\cup \{\a_1+\a_2+\a_3,\;\a_1+\a_2+\a_4,\; \a_2+\a_3+\a_4\}\\
&\cup  \{\a_1 + \a_2+\a_3+\a_4\} \cup \{\a_1+2\a_2+\a_3+\a_4\}.
\end{split}
\end{equation*}

Of all these roots only the last two do not have support in a subdiagram of type $A$. The corresponding indecomposable representations are

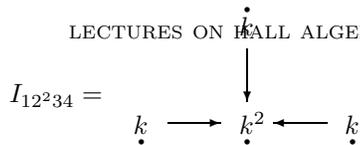
\vspace{.5in}

\centerline{
\begin{picture}(200, 10)
\put(10,0){$I_{1234}=$}
\put(60,-15){\circle*{2}}
\put(100,-15){\circle*{2}}
\put(140,-15){\circle*{2}}
\put(100,35){\circle*{2}}
\put(57,-12){$k$}
\put(97,-12){$k$}
\put(137,-12){$k$}
\put(97,25){$k$}
\put(100,20){\vector(0,-1){20}}
\put(70,-8){\vector(1,0){20}}
\put(130,-8){\vector(-1,0){20}}
\put(102,12){$\simeq$}
\put(70,-6){$\simeq$}
\put(115,-6){$\simeq$}
\end{picture}}

\vspace{.35in}

and

\vspace{.5in}

\centerline{
\begin{picture}(200, 10)
\put(10,0){$I_{12^234}=$}
\put(60,-15){\circle*{2}}
\put(100,-15){\circle*{2}}
\put(140,-15){\circle*{2}}
\put(100,35){\circle*{2}}
\put(57,-12){$k$}
\put(97,-12){$k^2$}
\put(137,-12){$k$}
\put(97,25){$k$}
\put(100,20){\vector(0,-1){20}}
\put(70,-8){\vector(1,0){20}}
\put(130,-8){\vector(-1,0){20}}
\end{picture}}

\vspace{.35in}

In this last case, the three maps $k \to k^2$ have to be injective, and thus define three lines in $k^2$, or three points $(\lambda_1, \lambda_2, \lambda_3)$ in $\mathbb{P}^1(k)$. The representation is indecomposable if and only if these points are distinct. Moreover the action by conjugation at the vertex $2$ of $GL(2,k)$ corresponds to the action on $(\lambda_1, \lambda_2, \lambda_3)$ by an element of ${Aut}(\mathbb{P}^1)$. It is well-known that any triple of distinct points in $\mathbb{P}^1$ may be brought to $(0,1,\infty)$ by an automorphism, hence $I_{12^234}$ is indeed unique up to isomorphism.
 
 \endexample
 
\vspace{.2in}

By Cartan's classification (see Appendix~A.1.) $A$ is positive definite if and only if the quiver $\vec{Q}$ is of type $A_n, D_n$ or $E_l$ with $l=6,7,8$. The above examples treat the (easy) cases of $A_n$ and $D_n$. The indecomposables for the exceptional  quivers are worked out in \cite{RingelTame}.

\vspace{.2in}

Concerning tame quivers, we have the following result, established independently by Nazarova \cite{Nazarova} and Donovan-Freislich \cite{DF}~:

\vspace{.1in}

\begin{theo}[Nazarova, Donovan-Freislich]\label{T:Godknows} The category $Rep^{nil}_k\vec{Q}$ is tame if and only if the matrix $A$ is positive semi-definite and has corank equal to $1$, i.e. if and only if $\g$ is an affine Lie algebra. Moreover,
\begin{enumerate}
\item[i)] For $\a \in Q$, there exists an indecomposable representation $M=(V,\underline{x})$ with
$\sum_i {dim}(V_i)\a_i$ if and only if $\a \in \Delta^+$,
\item[ii)] If $\a \in \Delta^+$ is a real root then there exists a unique such indecomposable; if $\a \in \Delta^+$ is imaginary then, unless $\vec{Q}$ is a cyclic quiver, there exists a one-parameter family of such indecomposables.
\end{enumerate}
\end{theo}

\vspace{.15in}

The new phenomenon here, as opposed to the finite type case, is the existence of dimensions for which an infinite number of indecomposables exist (when the ground field is itself infinite, of course). The case of a cyclic quiver is special, and is discussed separately in Section~3.5. As the root system of affine Lie algebras contains a single line $\Z\delta$ of imaginary roots, all dimension vectors for which there exists several indecomposables are multiple of the \textit{indivisible imaginary root} $\delta$.

\vspace{.15in}

\addtocounter{theo}{1}

\paragraph{\textbf{Example \thetheo.}} The simplest tame quiver is the Kronecker quiver~:

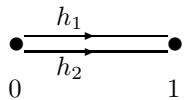
\vspace{.15in}

\centerline{
\begin{picture}(200, 10)
\put(60,0){\circle*{5}}
\put(120,0){\circle*{5}}
\put(85,-3){\vector(1,0){5}}
\put(85,3){\vector(1,0){5}}
\put(63,3){\line(1,0){54}}
\put(63,-3){\line(1,0){54}}
\put(57,-20){$0$}
\put(117,-20){$1$}
\put(75,7){$h_1$}
\put(75,-12){$h_2$}
\end{picture}}

\vspace{.4in}

The matrix of the Euler form is $A=\begin{pmatrix} 2 & -2 \\ -2 & 2 \end{pmatrix}$ and the relevant Kac-Moody
algebra is $\widehat{\mathfrak{sl}}_2(\C)$. If $\{a_0, \a_1\}$ are the simple roots then the indivisible imaginary root is $\delta=\a_0+\a_1$ and
$$\Delta^+=\{\a_1+n \delta\;| \; n \geq 0\} \sqcup \{-\a_1 + n\delta\;|\; n >0\} \sqcup \{n\delta\;|\; n >0\}.$$
Here are some indecomposables associated to real roots (all maps below are isomorphisms)~:

\vspace{.5in}

\centerline{
\begin{picture}(300, 10)
\put(10,0){$I_{\a_1+\delta}=$}
\put(57,0){\circle*{2}}
\put(105,0){\circle*{2}}
\put(75,-2){$h_1$}
\put(75,22){$h_2$}
\put(63,7){\vector(2,1){34}}
\put(63,7){\vector(1,0){34}}
\put(54,-20){$0$}
\put(102,-20){$1$}
\put(54,5){$k$}
\put(102,5){$k$}
\put(102,15){$\oplus$}
\put(102,25){$k$}
\put(160,0){$I_{-\a_1+\delta}=$}
\put(217,0){\circle*{2}}
\put(265,0){\circle*{2}}
\put(235,-2){$h_1$}
\put(235,20){$h_2$}
\put(223,24){\vector(2,-1){34}}
\put(223,7){\vector(1,0){34}}
\put(214,-20){$0$}
\put(262,-20){$1$}
\put(262,5){$k$}
\put(214,5){$k$}
\put(214,15){$\oplus$}
\put(214,25){$k$}
\end{picture}}

\vspace{.35in}

Let us now consider dimensions corresponding to the imaginary root $\delta$~:

\vspace{.35in}

\centerline{
\begin{picture}(200, 10)
\put(10,0){$I_{\delta}^{(\lambda,\mu)}=$}
\put(60,-7){\circle*{2}}
\put(120,-7){\circle*{2}}
\put(67,8){\vector(1,0){46}}
\put(67,2){\vector(1,0){46}}
\put(57,-20){$0$}
\put(117,-20){$1$}
\put(75,12){$ \cdot \lambda$}
\put(75,-7){$\cdot \mu$}
\put(57,0){$k$}
\put(117,0){$k$}
\put(160,0){$(\lambda,\mu) \neq (0,0)$}
\end{picture}}

\vspace{.5in}

It is clear that $I_{\delta}^{(\lambda,\mu)}$ is isomorphic to $I_{\delta}^{(\lambda',\mu')}$ if and only if there exists $t \in k$ such that $(\lambda',\mu')=(t\lambda,t\mu)$. Hence the indecomposable representations of class $\delta$ are parametrized by points of $\mathbb{P}^1(k)$. \endexample

\vspace{.35in}

\addtocounter{theo}{1}

\paragraph{\textbf{Example \thetheo.}} Let us now choose for $\vec{Q}$ the following quiver of type $D_4^{(1)}$~:

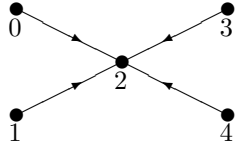
\vspace{.5in}

\centerline{
\begin{picture}(200, 10)
\put(60,-20){\circle*{5}}
\put(100,0){\circle*{5}}
\put(140,-20){\circle*{5}}
\put(140,20){\circle*{5}}
\put(60,20){\circle*{5}}
\put(57,-30){$1$}
\put(57,10){$0$}
\put(97,-10){$2$}
\put(137,-30){$4$}
\put(137,10){$3$}
\put(60,-20){\line(2,1){40}}
\put(60,20){\line(2,-1){40}}
\put(100,0){\line(2,-1){40}}
\put(100,0){\line(2,1){40}}
\put(80,-10){\vector(2,1){5}}
\put(80,10){\vector(2,-1){5}}
\put(120,-10){\vector(-2,1){5}}
\put(120,10){\vector(-2,-1){5}}
\end{picture}}

\vspace{.6in}

The Kac-Moody Lie algebra here is $\widehat{\mathfrak{so}}_8(\C)$. The indivisible imaginary root is $\delta=\a_0+\a_1+2\a_2+\a_3+\a_4$. Indecomposable representations of dimension $\delta$ are given by a quadruple of injective maps $k \hookrightarrow k^2$

\vspace{.5in}

\centerline{
\begin{picture}(200, 10)
\put(-40,0){$I_{\delta}^{(\lambda_0, \lambda_1,\lambda_3,\lambda_4)}=$}
\put(60,-20){\circle*{2}}
\put(100,0){\circle*{2}}
\put(140,-20){\circle*{2}}
\put(140,20){\circle*{2}}
\put(60,20){\circle*{2}}
\put(57,-30){$1$}
\put(57,10){$0$}
\put(97,-10){$2$}
\put(137,-30){$4$}
\put(137,10){$3$}
\put(57,-17){$k$}
\put(57,23){$k$}
\put(97,3){$k^2$}
\put(137,-17){$k$}
\put(137,23){$k$}
\put(68,-11){\vector(2,1){24}}
\put(68,23){\vector(2,-1){24}}
\put(132,-11){\vector(-2,1){24}}
\put(132,23){\vector(-2,-1){24}}
\end{picture}}

\vspace{.6in}

\noindent
landing in at least three different lines in $k^2$. Each such line corresponds to a point $\lambda_l$ in $\mathbb{P}^1(k)$. Hence the set of indecomposables is in bijection with the set of ordered quadruples $(\lambda_0, \lambda_1,\lambda_3,\lambda_4)$ of points in $\mathbb{P}^1(k)$, with at least three of the $\lambda_i$'s distinct, and all up to an automorphism of $\mathbb{P}^1(k)$. This set is a little bit more tricky than one might first think~:

\vspace{.05in}

\noindent
\textit{Case a).} $\lambda_i \neq \lambda_j$ for $i=1,3,4$. Then there is a \textit{unique} $\phi \in {Aut}(\mathbb{P}^1(k))$ sending $(\lambda_1,\lambda_3,\lambda_4)$ to $(0,1,\infty)$. It sends $\lambda_0$ to some point $x \in \mathbb{P}^1(k)$.

\vspace{.1in}

\noindent
\textit{Case b).} Two among $\lambda_1,\lambda_3,\lambda_4$ are equal; then, up to the action of $Aut(\mathbb{P}^1(k))$ we may reduce $(\lambda_0,\lambda_1,\lambda_3,\lambda_4)$ to exactly one of~:
$$(0,1,1,\infty), \qquad (0,1,\infty,\infty), \qquad (0,1,\infty,1).$$
Hence, geometrically, the set of indecomposables ressembles a projective line with three points being ``doubled''~:

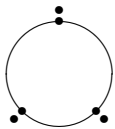
\vspace{.15in}

\begin{equation}\label{E:ponethreepoints}
\centerline{
\begin{picture}(300, 10)
\put(150,0){\circle{40}}
\put(150,20){\circle*{3}}
\put(150,24){\circle*{3}}
\put(136,-14){\circle*{3}}
\put(133,-17){\circle*{3}}
\put(164,-14){\circle*{3}}
\put(167,-17){\circle*{3}}
\end{picture}}
\end{equation}

\endexample

\vspace{.3in}

In a remarkable work \cite{Kac}, Kac managed to extend the theorem of Gabriel to the case of an \textit{arbitrary} quiver.

\vspace{.1in}

\begin{theo}[Kac]\label{T:Kac} Let $\vec{Q}$ be an arbitrary quiver.
\begin{enumerate}
\item[i)] For $\a \in Q$, there exists an indecomposable representation $M=(V,\underline{x})$ with
$\sum_i {dim}(V_i)\a_i$ if and only if $\a \in \Delta^+$,
\item[ii)] If $\a \in \Delta^+$ is a real root then there exists a unique such indecomposable; if $\a \in \Delta^+$ is imaginary then there exists many (i.e. more than one) such indecomposables.
\end{enumerate}
\end{theo}

\vspace{.2in}

\centerline{\textbf{3.3. Hall algebras of quivers.}}
\addcontentsline{toc}{subsection}{\tocsubsection {}{}{\; 3.3. Hall algebras of quivers.}}

\vspace{.15in}

\paragraph{} We are now ready to compute the extended Hall algebra $\widetilde{\mathbf{H}}_{\vec{Q}}$ of $Rep^{nil}_k\vec{Q}$. We assume in this Section that $k=\mathbb{F}_q$ so that $Rep^{nil}_k\vec{Q}$ is a finitary, hereditary category. Moreover, as this category satisfies the finite subobject condition (see Remark~1.6) its Hall algebra $\widetilde{\mathbf{H}}_{\vec{Q}}$ is an honest Hopf algebra. It will be convenient to introduce a new notation $\nu=q^{\frac{1}{2}}$ (note that there is a choice of
the square root involved here).

\vspace{.15in}

Let us start with some simple calculations of Hall numbers.

\vspace{.1in}

\addtocounter{theo}{1}
\paragraph{\textbf{Example \thetheo.}} The simplest quiver in the galaxy~: $\vec{Q}=\bullet$ (one vertex, zero arrows). There is but a single simple object $S$, and any object is isomorphic to $S^{\oplus j}$ for some $j$. We have ${Hom}(S,S)=k$ and ${Ext}^1(S,S)=\{0\}$, so that $\langle S,S \rangle_m=q^{\frac{1}{2}}=\v$. We have
$$[S^{\oplus i}] \cdot [S^{\oplus j}]=\v^{ij} \# Gr(j,i+j) [S^{\oplus i+j}]=\v^{ij} \begin{bmatrix} i+j \\ j \end{bmatrix}_+ [S^{\oplus i+j}]$$
since the Hall number in this case simply counts the number of $i$-dimensional subspaces in an $i+j$-dimensional vector space.  In particular,
\begin{equation}\label{E:powerE}
[S]^n=\v^{\frac{n(n-1)}{2}} [n]!_+ [S^{\oplus n}]=\v^{n(n-1)}[n]! [S^{\oplus n}].
\end{equation}
(see Example~2.3 and Appendix~A.4. for the definitions of $\v$-binomial numbers).

\endexample

\vspace{.15in}

\addtocounter{theo}{1}
\paragraph{\textbf{Examples \thetheo.}} Assume that $\vec{Q}$ has two vertices $1$ and $2$, with $c_{12}$ arrows going from $1$ to $2$, and $c_{21}$ arrows going from $2$ to $1$. The simple modules associated to the vertices are denoted $S_1$ and $S_2$; Of course, ${Hom}(S_1,S_2)={Hom}(S_2,S_1)=\{0\}$.

\vspace{.05in}

\noindent
$\bullet$ If $c_{12}=c_{21}=0$ (i.e. the two vertices are not connected) then
there are no nontrivial extensions between $S_1$ and $S_2$ and $\langle S_1,S_2 \rangle_m=\langle S_2,S_1 \rangle_m=1$. Hence  
$$[S_1] \cdot [S_2]=[S_1 \oplus S_2]=[S_2] \cdot [S_1].$$

\vspace{.05in}

\noindent
$\bullet$ Assume now that $c_{12}=1, c_{21}=0$, i.e. that vertices $1$ and $2$ are connected by a single edge going from, say, $1$ to $2$. We have ${Ext}^1(S_1,S_2)=\{0\}$ and ${Ext}^1(S_2,S_1)=k$. Then
$\langle S_1,S_2 \rangle_m=q^{-\frac{1}{2}}=\v^{-1}$ and there is a unique nontrivial extension $I_{12}$ of $S_2$ by $S_1$. As ${Hom}(S_2,I_{12})=k$, it is easy to see that
$$[S_1] \cdot [S_2]=\v^{-1} \big( [S_1 \oplus S_2] + [I_{12}]\big).$$
On the other hand, there are no nontrivial extensions of $S_1$ by $S_2$, thus $\langle S_2,S_1\rangle_m=1$ and
$$[S_2] \cdot [S_1] = [S_1 \oplus S_2].$$
This simple example already shows that $\widetilde{\mathbf{H}}_{\vec{Q}}$ and $\H_{\vec{Q}}$ are not commutative in general.

\vspace{.05in}

\noindent
$\bullet$ Let us keep the same quiver as above. As there are still no extensions of $S_2$ by $S_1$, we have 
\begin{equation}\label{E:EEex1}
[S_2] \cdot [S_1]^2 =\v(\v^2+1)[S_2] \cdot [S_1^{\oplus 2}]= \v(\v^2+1)[S_1^{\oplus 2} \oplus S_2].
\end{equation}
In the other direction, 
\begin{equation}\label{E:EEex2}
\begin{split}
[S_1]^2 \cdot [S_2]&= \v(\v^2+1)[S_1^{\oplus 2}] \cdot [S_2]\\
&=\v^{-1}(\v^2+1)\big( [S_1^{\oplus 2} \oplus S_2] + [S_1 \oplus I_{12}] \big)
\end{split}
\end{equation}
since ${Hom}(S_2,S_1^{\oplus 2} \oplus S_2)={Hom}(S_2, S_1 \oplus I_{12})=k$.

Finally, to compute $[S_1] \cdot [S_2] \cdot [S_1]=[S_1] \cdot [S_1 \oplus S_2]$, note that there are $q+1$ submodules of $S_1^{\oplus 2} \oplus S_2$ isomorphic to $S_1 \oplus S_2$, but only one submodule of $I_{12} \oplus S_2$ isomorphic to $S_1 \oplus S_2$. We deduce that
\begin{equation}\label{E:EEex3}
[S_1] \cdot [S_2] \cdot [S_1]= (\v^2+1) [S_1^{\oplus 2} \oplus S_2] + [S_1 \oplus I_{12}].
\end{equation}

Observe that combining (\ref{E:EEex1}), (\ref{E:EEex2}) and (\ref{E:EEex3}) yields the following equation satisfied by $[S_1]$ and $[S_2]$~:
\begin{equation}\label{E:rk21}
[S_1]^2 \cdot [S_2] - (\v+\v^{-1})[S_1]\cdot [S_2] \cdot [S_1] + [S_2] \cdot [S_1]^2=0.
\end{equation}
Similar computations give the dual relation~:
\begin{equation}\label{E:rk22}
[S_2]^2 \cdot [S_1] - (\v+\v^{-1})[S_2]\cdot [S_1] \cdot [S_2] + [S_1] \cdot [S_2]^2=0.
\end{equation}
In fact, as will follow from Ringel's Theorem~\ref{T:RingelHall}, these are the \textit{only} relations satisfied by $[S_1]$ and $[S_2]$.

\endexample

\vspace{.2in}

Recall that $\g$ is the Kac-Moody algebra associated to the Euler matrix $A=(a_{ij})_{i,j \in I}$ of $\vec{Q}$. Let $\U_\nu(\bo'_+)$ be the positive half of the quantum group $\U_\nu(\g')$ associated to $\g$ (more precisely, to the derived algebra $\g'$ of $\g$)--see Appendix~A.4. Here we consider the version of the quantum group $\U_v(\bo'_+)$ which is \textit{specialized at} $v=\nu$. 

The Hopf algebra  $\U_\nu(\bo'_+)$ is generated by elements $E_i, K_i^{\pm 1},$ for $i \in I$ subject to the following set of relations~:

\vspace{.05in}

\begin{equation}\label{E:KKetEK}
K_i K_j=K_j K_i, \qquad
K_i E_j K_i^{-1}=\v^{a_{ij}} E_j, \qquad
\forall\; i,j \in I, 
\end{equation}
\begin{equation}\label{E:EE}
\sum_{l=0}^{1-a_{ij}} (-1)^l \begin{bmatrix} 1-a_{ij} \\ l \end{bmatrix} E_i^l E_j E_i^{1-a_{ij}-l}=0
\end{equation}

\vspace{.05in}

As for the coproduct and antipode, they are given by
\begin{equation}\label{E:copKE}
\Delta(K_i)=K_i \otimes K_i, \quad \Delta(E_i)=E_i \otimes 1 + K_i \otimes E_i, 
\end{equation}
\begin{equation}\label{E:antipKE}
S(K_i)=K_i^{-1}, \qquad S(E_i)=-K_i^{-1}E_i.
\end{equation}

\vspace{.1in}

\begin{theo}[Ringel, Green]\label{T:RingelHall} The assignement $E_i \mapsto [S_i],\; K_i \mapsto \mathbf{k}_{S_i}$ for $i \in I$ extends to an embedding of Hopf algebras 
$$\Psi: \U_{\nu}(\bo'_+) \to \widetilde{\mathbf{H}}_{\vec{Q}}.$$
The map $\Psi$ is an isomorphism if and only if $\vec{Q}$ is of finite type, (i.e. if and only if $\g$ is a simple Lie algebra).
\end{theo}

\noindent
\textit{Proof.} There are two steps in the proof. First, we need to check that relations (\ref{E:KKetEK}),(\ref{E:EE}), (\ref{E:copKE}) and (\ref{E:antipKE}) hold in $\widetilde{\mathbf{H}}_{\vec{Q}}$; this will prove the existence of the map $\Psi$. Secondly, we need to show that $\Psi$ is injective.

\vspace{.1in}

The first step is essentially one big (but straightforward) computation. Relations (\ref{E:KKetEK}), (\ref{E:copKE}) and (\ref{E:antipKE}) all hold trivially by definition. The only difficulty lies in checking (\ref{E:EE}). Note that this relation only involves two vertices $i,j \in I$. When $|a_{ij}| \leq 1$, this is treated in Example~3.15. We deal here with the general case.
In order to unburden the notation, let us set $r=c_{ij}, s=c_{ji}$ and $t=r+s=-a_{ij}$. By Example~3.14, (\ref{E:powerE}), we have
$[S_i]^{(l)}:=\frac{[S_i]^l}{[l]!}=\v^{l(l-1)} [S_i^{\oplus l}]$. As $\langle S_i^{\oplus l},S_j\rangle_m=\v^{-lr}$, we have
$$[S_i]^{(l)}[S_{j}]=\v^{l(l-1)-lr}\sum_{M \in
\mathcal{I}_1}[M],$$ where $\mathcal{I}_1=\{M\;|\;\exists \;N \subset M
\;s.t.\; N \simeq S_{j},  \; M/N \simeq S_i^{\oplus l}\}$. 

Next, for a
representation $L$ of $\vec{Q}$ of dimension $(r+1)\epsilon_i +
\epsilon_j$, we define
$$U_L=\bigcap_{i \stackrel{h}{\to} j} \mathrm{Ker}\;x_{h},
\qquad  V_L=\sum_{j \stackrel{h}{\to} i} \mathrm{Im}\;x_{h},$$
and set $u_L=\mathrm{dim}(U_L), w_L=\mathrm{dim}(V_L)$.
A direct
computation shows that 
$$[S_i]^{(l)}
[S_{j}][S_i]^{(n)}=\v^{-ns-lr+nl+l(l-1)+n(n-1)}\sum_{[M]}
\mathbf{p}_{M,n}[M],$$
where $\mathbf{p}_{M,n}=0$ unless $V_M \subset U_M$, in which case we have
$$\mathbf{p}_{M,n}=\#Gr(n-w_M,\;u_M-w_M)=\v^{(u_M-n)(n-w_M)}
\begin{bmatrix} u_M-w_M\\ n-w_M\end{bmatrix}.$$

\vspace{.1in}

Setting $n=t+1-l$ and summing up, we obtain
$$\sum_{l=0}^{t+1} (-1)^l [S_i]^{(l)}[S_{j}][S_i]^{(n)}=
\sum_{[M]\;s.t.V_M \subset U_M} \mathbf{p}_M[M],$$
where
\begin{equation*}
\begin{split}
\mathbf{p}_M &=\sum_{l=0}^{t+1} (-1)^l
\v^{-ns-lr+nl+l(l-1)+n(n-1)+(u_M-n)(n-w_M)}
\begin{bmatrix} u_M-w_M\\ n-w_M\end{bmatrix}\\
&=\v^{(t+1)s -u_Mw_M}
\sum_{n=0}^{t+1}(-1)^{t+1-n}\v^{-(2s+1-u_M-w_M)n} \begin{bmatrix}
u_M-w_M\\ n-w_M\end{bmatrix}.
\end{split}
\end{equation*}

\vspace{.1in}

Clearly, $u_M \geq s+1 >w_M$ for any $M$. We deduce that
$1-u_M-w_M \leq 2s+1-u_M-w_M \leq u_M +w_M -1$.  Now we use the following identity
(see, for example, \cite[(3.2.8)]{Kas})~:

\vspace{.05in}

\paragraph{\textbf{Lemma.}} \textit{Let $m \geq 1$ and let $1-m \leq d \leq m-1$ with $d \equiv
m-1  \;(\mathrm{mod}\;2)$. Then}
$$\sum_{n=0}^m (-1)^n \v^{-dn}
\left[\begin{matrix} m\\n\end{matrix}\right]=0.$$

\vspace{.05in}

As a consequence, we obtain

$$\sum_{l=0}^{t+1} (-1)^l[S_i]^{(l)} [S_j] [S_i]^{(t+1-l)}=0$$
which is nothing else than the $\v$-Serre relation (\ref{E:EE}). This shows that $\Psi$ is well-defined. The fact that $\Psi$ is compatible with the coproduct and the antipode is straightforward given the definitions.

\vspace{.05in}

Let $\{\,,\,\}: \U_\nu(\bo'_+) \otimes \U_{\nu}(\bo'_+) \to \C$ be the pullback under $\Psi$ of Green's scalar product on $\widetilde{\H}_{\vec{Q}}$. As $\Psi$ is a morphism of Hopf algebras and Green's scalar product is a Hopf pairing, $\{\,,\,\}$ is also a Hopf pairing. Moreover, by construction it satisfies
$$ \{E_i,E_j\}=\frac{\delta_{ij}}{\v^2-1}, \qquad (K_i,K_j)=\v^{a_{ij}}, \qquad (E_i, K_j)=0.$$
By Theorem~A.18 this completely determines $\{\,,\,\}$, which  coincides with the restriction to $\U_\v(\bo'_+)$ of Drinfeld's scalar product on $\U_v(\bo_+)$. In particular, the restriction of $\{\,,\,\}$ to $\U_\v(\n_+)$ is nondegenerate. Hence 
$$({Ker}\;\Psi) \cap \U_\v(\n_+) \subset ({Ker}\;\{\,,\,\}) \cap \U_\v(\n_+) =\{0\}.$$ 
The injectivity of $\Psi$ easily follows from the fact that $\U_\v(\bo'_+) = \C[K_i^{\pm 1}]_{i \in I} \otimes \U_\v(\n_+)$ while $\widetilde{H}_{\vec{Q}}=\C[\mathbf{k}^{\pm 1}_{S_i}]_{i \in I} \otimes \H_{\vec{Q}}$. Of course, the restriction of $\Psi$ gives an embedding $\U_\v(\n_+) \to \H_{\vec{Q}}$. This concludes the proof of the first statement. 

The second statement, concerning finite type quivers, follows from Kac's theorem~\ref{T:Kac}; indeed, by a quantum version of the PBW theorem (see Theorem~A.16) for $\U_{\v}(\n_+)$ we have for any weight $\beta$
$$dim\;\U_{\v}(\n_+)[\beta]=|\{ (n_\a) \in \N^{\Delta^+}\;|\; \sum n_{\a} \a=\beta\}|.$$
On the other hand, by construction
\begin{equation*}
\begin{split}
dim\;\H_{\vec{Q}}[\beta]=&|\{M \in \mathcal{X}_{\vec{Q}}\;|\; \underline{dim}\;M=\beta\}|\\
=&|\{(m_\a) \in \N^{Indec\;\vec{Q}}\;|\; \sum m_I \underline{dim}\;I=\beta\}|
\end{split}
\end{equation*}
where $Indec\;\vec{Q}$ stands for the set of all indecomposable representations of $\vec{Q}$. These two expressions for graded dimensions are equal if and only if there is
precisely one indecomposable of dimension $\a$ for all positive roots $\a \in \Delta^+$. This happens if and only if $\vec{Q}$ is of finite type.
We are done. \qed

\vspace{.15in}

This seems to be the right place to summarize the correspondence between the category of representations of a quiver $\vec{Q}$ and the associated Kac-Moody algebra~:

\vspace{.15in}

$$
\begin{tabular}{|c|c|}
\hline
&\\
Category $\mathcal{C}=Rep_k \vec{Q}$ & Kac-Moody Lie algebra $\g$\\
&\\
\hline
&\\
Grothendieck group $K(\mathcal{C})$ & Root lattice $Q=\bigoplus_i \Z \a_i$\\
&\\
\hline
&\\
Symmetrized additive Euler form $(\;,\;)$ & Cartan-Killing form $(\;,\;)$\\
&\\
\hline
&\\
(classes of) simple objects $\{S_i\}$ & Simple roots $\{\a_i\}$\\
&\\
\hline
&\\
(classes of) indecomposable objects & Positive root system $\Delta^+$\\
&\\
\hline
&\\
Hall algebra $\mathbf{H}_{\vec{Q}}$ & Quantum group $\U_v(\n_+)$\\
&\\
\hline
&\\
Group algebra of $K(\mathcal{C})$ & Cartan  $\U_v(\h)$\\ 
&\\
\hline
&\\
Extended Hall algebra $\widetilde{\H}_{\vec{Q}}$ & Quantum group $\U_v(\bo'_+)$\\
&\\
\hline
&\\
Finite type, tame type, wild type & Simple Lie algebra, affine Lie algebra, \\
&\qquad \qquad \qquad the rest  (!)\\
&\\
\hline
\end{tabular}
$$

\vspace{.3in}

\addtocounter{theo}{1}
\paragraph{\textbf{Remarks \thetheo.}} i) The existence of the map $\Psi$ as a morphism of algebras is due to Ringel \cite{Ri}, where he also shows that $\Psi$ is an isomorphism for finite type quivers. The fact that $\Psi$ is also compatible with a comultiplication and a scalar product, and its rather easy corollary that $\Psi$ is injective, are due to Green \cite{Green}.\\
ii) The image of $\Psi$, i.e. the subalgebra of $\mathbf{H}_{\vec{Q}}$ generated by the simple objects, is the \textit{composition subalgebra} $\mathbf{C}_{\vec{Q}}$. Lusztig's theory of canonical bases provides a geometric characterization of $\mathbf{C}_{\vec{Q}}$ as a subalgebra of $\mathbf{H}_{\vec{Q}}$ but as this relies on intersection cohomology methods it is highly unlikely that there exists any elementary description of the elements of $\mathbf{C}_{\vec{Q}}$ for a general quiver. However, such a description \textit{does} exist when $\vec{Q}$ is tame-- see Section~3.7.\\
iii) There is a version of this result for quivers in which loops are allowed. The corresponding Lie algebras are then the so-called \textit{generalized}, or \textit{Borcherds Kac-Moody algebras} --see \cite{KangS}. There is also a version of this result for species (or quivers with automorphisms). This provides a construction of quantum groups for non necessarily simply laced Kac-Moody algebras.

\vspace{.2in}

\centerline{\textbf{3.4. PBW bases (finite type).}}
\addcontentsline{toc}{subsection}{\tocsubsection {}{}{\; 3.4. PBW bases (finite type).}}

\vspace{.15in}

\paragraph{} We assume here that $\vec{Q}$ is a quiver of finite type, and hence that $\g$ is a simple Lie algebra. As in the previous Section, $k=\mathbb{F}_q$ and we set $\nu=q^{\frac{1}{2}}$. By Theorem~\ref{T:RingelHall}, there is an isomorphism
$$\Psi_{\nu} :\U_{\nu}(\bo'_+) \stackrel{\sim}{\to} \widetilde{\H}_{\vec{Q}}$$
which restricts to an isomorphism
$$\Psi_{\nu} :\U_{\nu}(\n_+) \stackrel{\sim}{\to} {\H}_{\vec{Q}}.$$
We have added the index $\nu$ to $\Psi$ because we will soon vary the field $k$ and it will be important for us to remember the value of $q$. The pullback under $\Psi_{\nu}$ of the natural basis
$\{[M]\;|\; M \in Obj(Rep_{k}\vec{Q})\}$ yields a basis $\{\mathbf{f}_{\nu,M}\;|\; M \in Obj(Rep_{k}\vec{Q})\}$ of the quantum group $\U_{\nu}(\n_+)$ \textit{specialized at} $v=\nu$. This basis, usually called the \textit{PBW basis}\footnote{strictly speaking, the PBW basis is obtained from $\{\mathbf{f}_M\}$ by a slight renormalization which we choose to ignore here.} is a very useful tool in quantum groups theory (see e.g. \cite{LusztigPBW}). What is crucial here, and a priori not obvious at all, is that the bases $\{\mathbf{f}_{\nu,M}\}$ for different values of $q$ all come from a \textit{same} basis of the integral form $\U_v^{res}(\n_+)$ (see Appendix A.4.). This, as we will see, is more or less equivalent to the existence of Hall polynomials.

\vspace{.1in}

Before we proceed, let us make the important remark that \textit{when $\vec{Q}$ is of finite type}, the set of objects of $Rep^{nil}_k\vec{Q}$ is independent of the field $k$. Indeed, by Gabriel's Theorem~\ref{T:Gabriel}, it can be canonically identified with the set of maps $\Delta^+ \to \N$ by giving the multiplicity of each indecomposable. We denote this set of objects by $\mathcal{O}_{\vec{Q}}$.

\vspace{.1in}

\begin{prop}[Ringel]\label{P:Hallpolquiver} Let $\vec{Q}$ be a quiver of finite type. Then
\begin{enumerate}
\item[i)] For each object $L \in \mathcal{O}_{\vec{Q}}$ there exists a unique element $\mathbf{f}_{L} \in \U_v^{res}(\n_+)$ such that, for any finite field $k=\mathbb{F}_q$ we have $(\mathbf{f}_{L})_{|v=\nu}=\Psi_{\nu}^{-1}([L])$.
\item[ii)] Let $M,N,R$ be objects of $\mathcal{O}_{\vec{Q}}$. There exists a unique polynomial $P_{M,N}^R(t) \in \Q[t]$ ($\mathrm{independent}$ $\mathrm{of}$ $\mathrm{the}$ $\mathrm{field}$ $k$), such that for any choice of finite field $k=\mathbb{F}_q$ we have
$$\frac{1}{a_Ma_N}\mathbf{P}^R_{M,N}=P_{M,N}^R(q).$$
\end{enumerate}
\end{prop}

\vspace{.15in}

The proof  hinges on the following Lemma, which is a byproduct of Auslander-Reiten theory (see \cite{CrawleyNotes}).

\vspace{.1in}

\begin{lem}\label{L:orderfinite} Let $\vec{Q}$ be a finite quiver. There exists a total ordering $\preceq$ of the set of indecomposables (independent of the ground field $k$) such that
$$M \prec N \Rightarrow {Hom}(N,M) = {Ext}^1(M,N)=\{0\}.$$
Moreover, for any indecomposable $M$, we have ${Ext}^1(M,M)=\{0\}$ and ${End}(M)=k$.
\end{lem}

\vspace{.15in}

\addtocounter{theo}{1}

\paragraph{\textbf{Example \thetheo.}} Suppose that $\vec{Q}$ is the equioriented quiver of type $A_n$ of Example~3.8. Then
$${Hom}(I_{i, \ldots, j}, I_{i', \ldots, j'})=\{0\}\qquad \text{if\;} i<i'\;\text{or}\;i=i'\;\text{and}\;j<j',$$
$${Ext}^1(I_{i, \ldots, j}, I_{i', \ldots, j'})=\{0\}\qquad \text{if\;} i>i'\;\text{or}\;i=i'\;\text{and}\;j\geq j'.$$
Hence we may set 
$I_{i,\ldots,j} \preceq I_{i',\ldots, j'}$ if $i>i'$ or $i=i'$ and $j \geq j'$.
\endexample

\vspace{.15in}

\addtocounter{theo}{1}

\paragraph{\textbf{Example \thetheo.}} Let us now consider the quiver of type $D_4$ of Example~3.9. We leave it to the reader to check that we may take for $\preceq$ any total ordering refining the following partial order~:

\vspace{.2in}

\centerline{
\begin{picture}(220,20)
\put(0,0){$I_2$}
\put(20,0){$\prec$}
\put(40,-15){$I_{24}$}
\put(40,0){$I_{23}$}
\put(40,15){$I_{12}$}
\put(60,0){$\prec$}
\put(75,0){$I_{12^234}$}
\put(100,0){$\prec$}
\put(120,-15){$I_{234}$}
\put(120,0){$I_{124}$}
\put(120,15){$I_{123}$}
\put(140,0){$\prec$}
\put(160,0){$I_{1234}$}
\put(180,0){$\prec$}
\put(200,-15){$I_{234}$}
\put(200,0){$I_{124}$}
\put(200,15){$I_{123}$}
\end{picture}}

\vspace{.3in}

\noindent
(the subscript indicates the support of the indecomposable).
\endexample

\vspace{.2in}

\noindent
\textit{Proof of Proposition~\ref{P:Hallpolquiver}~:} We will first derive statement i). Clearly, if $i \in I$ and $n \in \N$ then by (\ref{E:powerE}) we have 
$$\mathbf{f}_{[S_i^{\oplus n}]}=v^{-n(n-1)}E_i^{(n)}=v^{-n(n-1)}\frac{E_i^n}{[n]!}.$$
Let us now argue by induction on the dimension vector $\underline{\text{dim}}\;L$. Assume that
we have proved the existence of $\mathbf{f}_{M}$ for any $M$ such that $\underline{\text{dim}}\;M < \underline{\text{dim}}\;L$. We distinguish two cases~:

\vspace{.05in}

\noindent
\textit{Case a)~: $L$ has at least two nonisomorphic indecomposable summands.} Thus we may write $L=I_1^{\oplus n_1} \oplus \cdots \oplus I_r^{\oplus n_r}$ where $r >1$ and $I_l$ are indecomposables which we may as well order in such a way that $I_j \prec I_h$ if
$j <h$. In particular, we have ${Ext}^1(I_j,I_h)=\{0\}$ if $j < h$. We claim that, in $\H_{\vec{Q}}$, 
\begin{equation}\label{E:prodhouih}
[I_1^{\oplus n_1}] \cdots [I_r^{\oplus n_r}]=\nu^{\underline{d}}[I_1^{\oplus n_1} \oplus \cdots \oplus I_r^{\oplus n_r}]=\nu^{\underline{d}}[L]
\end{equation}
where $\underline{d}=\sum_{j <h} n_jn_h {dim}\;{Hom}(I_j,I_h)$. Indeed, any extension of the
$I_j^{\oplus n_j}$ in that order is necessarily trivial so that we only have to compute the Hall number
$P^L_{I_1^{\oplus n_1}, \ldots, I_r^{\oplus n_r}}$. But as ${Hom}(I_r,I_j)=\{0\}$ for any $j <r$, there is only one submodule of $L$ isomorphic to $I_r^{\oplus n_r}$, and then only one submodule of
$L/I_r^{\oplus n_r}$ which is isomorphic to $I_{r-1}^{\oplus n_{r-1}}$, and so on. Hence $P^L_{I_1^{\oplus n_1}, \ldots, I_r^{\oplus n_r}}=1$ and (\ref{E:prodhouih}) is proved.
Observe that $\underline{d}$ is independent of $k$. We may therefore set $\mathbf{f}_L=v^{-\underline{d}} \mathbf{f}_{I_1^{\oplus n_1}} \cdots  \mathbf{f}_{I_r^{\oplus n_r}}$.
Note that as $r>1$ we have ${\underline{dim}}\;I_j^{\oplus n_j} < {\underline{dim}}\;L$ for all $j$.

\vspace{.05in}

\noindent
\textit{Case b)~: $L=I^{\oplus n}$ for some indecomposable representation $I$.} Restricting $\prec$ to the set of simple objects we get a total ordering such that $S_i \prec S_j$ if there exists an arrow $j \to i$ in $\vec{Q}$. Let us write all simple objects in this order as $\{S_{i_1}, S_{i_2}, \ldots, S_{i_m}\}$.
Let $l_i={dim}(L_i)$ so that ${\underline{dim}}(L)=\sum_i l_i \epsilon_i$. We claim that
\begin{equation}\label{E:hallprodhuuh}
[S_{i_m}^{\oplus l_{i_m}}] \cdots [S_{i_1}^{\oplus l_{i_1}}]=\nu^{\underline{d}}\sum_N [N]
\end{equation}
where $\underline{d}=\sum_{j<h} l_{i_j}l_{i_h} \langle S_{i_j},S_{i_h} \rangle_a$ and where $N$ ranges over all representations of $\vec{Q}$ of dimension $d_L=\sum_i l_i \epsilon_i$. Indeed, since there are no arrows in $\vec{Q}$ leaving the vertex $i_1$, any representation $N$ of dimension $d_L$ has a unique submodule isomorphic to $S_{i_1}^{\oplus l_{i_1}}$. Similarly, there is a unique submodule of $N/S_{i_1}^{\oplus l_{i_1}}$ isomorphic to $S_{i_2}^{\oplus l_{i_2}}$, and so on. 
Hence 
$$P^N_{S_{i_m}^{\oplus l_{i_m}}, \ldots, S_{i_1}^{\oplus l_{i_1}}}=1.$$
This proves (\ref{E:hallprodhuuh}).
Because there is at most one indecomposable representation of any given dimension,
we may rewrite (\ref{E:hallprodhuuh}) as
$$[M]=\nu^{-\underline{d}}[S_{i_m}^{\oplus l_{i_m}}] \cdots [S_{i_1}^{\oplus l_{i_1}}] -\sum_N [N]$$
where $N$ ranges over all representations of dimension $d_L$ which are different from $M$. Note that these all have at least two nonisomorphic indecomposable summands~: indeed,
the dimension of any indecomposable is a positive root and it is known that no two positive roots span the same line over $\Q$. Thus, by case a) above, there exists an element 
$\mathbf{f}_{N}$ with  the required properties for each of these $N$. Hence we may set
$$\mathbf{f}_{M}=v^{-\underline{d}}\mathbf{f}_{S_{i_m}^{\oplus l_{i_m}}} \cdots \mathbf{f}_{S_{i_1}^{\oplus l_{i_1}}}-\sum_N \mathbf{f}_{N}.$$
This finishes the induction step and concludes the proof of statement i). 

\vspace{.1in}

We turn to statement ii). Since $\{\mathbf{f}_{\nu,M}\}$ is a basis of
$\U_{\nu}(\n_+)$, the set $\{\mathbf{f}_{M}\}$ is linearly independent. By graded dimensions considerations, we deduce that $\{\mathbf{f}_M\}$ is a basis of $\U_v(\n_+)$ (over $\C(v)$). In particular, for any triple $M,N,R$ of objects the coefficient of $\mathbf{f}_R$ in the product $\mathbf{f}_M \cdot \mathbf{f}_N$ is a rational function $T_{M,N}^R(v) \in \C(v)$. We set $P_{M,N}^R(t)=t^{-\frac{1}{2}\langle M,N\rangle_a}T_{M,N}^R(t^{\frac{1}{2}})$. A priori, this is only a rational function in $t^{\frac{1}{2}}$, but note that by construction we have $P_{M,N}^R(q) \in \N$ for any prime power $q$. We leave it to the reader to check that this last property in fact forces $P_{M,N}^R(q)$ to belong to $\Q[t]$.\qed

\vspace{.15in}

As a consequence of the existence of Hall polynomials, we may, just like in Lecture~2, define a 
 ``universal'', or ``generic'' version $\underline{\H}_{\vec{Q}}$ of the Hall algebra $\H_{\vec{Q}}$, which is defined over the ring $\C[t^{\frac{1}{2}},t^{-\frac{1}{2}}]$.

\vspace{.1in}

\noindent
\textbf{Definition.} The generic Hall algebra of $\vec{Q}$ is the algebra 
$$\underline{\H}=\bigoplus_{M \in \mathcal{O}_{\vec{Q}}} \C[t^{\frac{1}{2}},t^{-\frac{1}{2}}] \mathbf{f}_{M}$$ in which the multiplication is defined by
$$\mathbf{f}_M \cdot \mathbf{f}_{N}=\sum_{R} t^{\frac{1}{2}\langle M,N \rangle_a}P^{R}_{M,N}(t) \mathbf{f}_R.$$

\vspace{.1in}

\begin{cor} There is an isomorphism of algebras $\underline{\Psi}: \U_{v}(\n_+) \stackrel{\sim}{\to} \underline{\H}_{\vec{Q}}$, where $v=t^{\frac{1}{2}}$. 
\end{cor}

\noindent
\textit{Proof.} By Proposition~\ref{P:Hallpolquiver} there is a canonical embedding of $\C[v,v^{-1}]$-algebras $i:\underline{\H}_{\vec{Q}} \hookrightarrow \U_v^{res}(\n_+)$. Conversely, by definition $\U_v^{res}(\n_+)$ is generated by the elements $E_i^{(n)}$ for $i \in I$. As these all belong to $i(\underline{\H}_{\vec{Q}})$ we deduce that $i$ is in fact an isomorphism and we may take for $\underline{\Psi}$ its inverse.\qed

\vspace{.1in}

The generic extended Hall algebra is defined to be the tensor product $$\widetilde{\underline{\H}}_{\vec{Q}}=\underline{\H}_{\vec{Q}} \otimes \C[\mathbf{k}^{\pm 1}_i]_{i \in I}$$
with relations
$$[\mathbf{k}_i,\mathbf{k}_j]=0,$$
$$\mathbf{k}_i \mathbf{f}_M \mathbf{k}_i^{-1}=t^{\frac{1}{2}(M,S_i)_a}\mathbf{f}_M.$$

Of course, a Corollary to the Corollary says that there is an isomorphism $\underline{\Psi}: \U_v(\bo_+) \stackrel{\sim}{\to} \underline{\widetilde{\H}}_{\vec{Q}}$. Thus $\underline{\widetilde{\H}}_{\vec{Q}}$ is a Hopf algebra and there exists ``Hall polynomials'' for the comultiplication $\Delta$ or the antipode $S$. The basis $\{\mathbf{f}_{M}\;|\; M \in \mathcal{O}_{\vec{Q}}\}$ is called the \textit{generic PBW basis} of $\U_v^{res}(\n_+)$.

\vspace{.15in}

\addtocounter{theo}{1}
\noindent \textbf{Remarks \thetheo .} i) With a little more work, one can show that in fact the Hall polynomials $P_{M,N}^R(t)$ belong to $\Z[t]$. See \cite{RingelHallpols} for a thorough treatment of these polynomials (some interesting applications are outlined at the very end of that paper).\\
ii) As one can well guess from Example~3.20. for instance, the ``homological'' ordering $\preceq$ on the set of indecomposable representations strongly depends on the orientation of the quiver, i.e. it is not something intrinsic to the root system $\Delta^+$. Since the PBW basis $\{\mathbf{f}_{M}\;|\; M \in \mathcal{O}_{\vec{Q}}\}$ is essentially obtained by multiplying together \textit{following the order} $\preceq$ elements $\mathbf{f}_{I}$ corresponding to indecomposables, we see that the PBW basis of $\U^{res}_v(\n_+)$ depends on the quiver. Hence, strictly speaking, there are as many PBW bases of $\U^{res}_v(\n_+)$ as orientations of the Dynkin diagram of $\g$.

\vspace{.2in}

\centerline{\textbf{3.5. The cyclic quiver.}}
\addcontentsline{toc}{subsection}{\tocsubsection {}{}{\; 3.5. The cyclic quiver.}}

\vspace{.15in}

\paragraph{} We spend this Section describing the Hall algebra of the equioriented quiver of type $A_{n-1}^{(1)}$~:

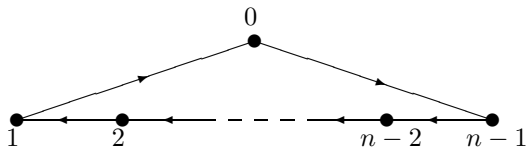
\vspace{.55in}

\centerline{
\begin{picture}(300, 10)
\put(160,30){\circle*{5}}
\put(70,0){\line(3,1){90}}
\put(160,30){\line(3,-1){90}}
\put(70,0){\circle*{5}}
\put(110,0){\circle*{5}}
\put(210,0){\circle*{5}}
\put(250,0){\circle*{5}}
\put(70,0){\line(1,0){40}}
\put(210,0){\line(1,0){40}}
\put(110,0){\line(1,0){35}}
\put(180,0){\line(1,0){40}}
\put(150,0){\line(1,0){5}}
\put(160,0){\line(1,0){5}}
\put(170,0){\line(1,0){5}}
\put(90,0){\vector(-1,0){5}}
\put(130,0){\vector(-1,0){5}}
\put(195,0){\vector(-1,0){5}}
\put(230,0){\vector(-1,0){5}}
\put(115,15){\vector(3,1){5}}
\put(205,15){\vector(3,-1){5}}
\put(156,36){$0$}
\put(66,-10){$1$}
\put(106,-10){$2$}
\put(200,-10){$n-2$}
\put(240,-10){$n-1$}
\end{picture}}

\vspace{.4in}

Though this quiver is tame, it is in some sense very close to being of finite type~: the classification of nilpotent modules is independent of the ground field $k$, and in particular there are only finitely many representations of any given dimension. This makes it possible to study not just the composition algebra $\mathbf{C}_{\vec{Q}}$ but the \textit{whole} Hall algebra ${\H}_{\vec{Q}}$. As will become clearer in Sections~3.6, 3.7 and Lecture~4, this Hall algebra appears in many places and plays an important role in the theory.

\vspace{.1in}

We start by classifying all representations of $\vec{Q}$. Fix $i \in \{0, \ldots, n-1\}$ and $l \geq 1$. Consider a $\Z$-graded vector space $V'=\bigoplus_{j=i+1-l}^{i} k e_j$ and define $x \in {End}(V')$ by $x(e_j)=e_{j-1}$ if $j \neq i+1-l$ and $x(e_{i+1-l})=0$. The induced $I=\Z/n\Z$-graded vector space
$V=\bigoplus_{h \in I} V_h$ where $V_h=\bigoplus_{j \equiv h} ke_j$, equipped with the same map $\underline{x}\in {End}(V)$ is an indecomposable representation of $\vec{Q}$, which we call $I_{[i;l]}$. As will follow from Proposition~\ref{P:classcyclic}, this is the unique indecomposable of length $l$ and socle $S_i$. Clearly $I_{[i;l]}$ only depends on the class of $i$ in $\Z/n\Z$.

\begin{prop}\label{P:classcyclic} For any field $k$, the representations $\{I_{[i;l]}\;|\; i \in \Z/n\Z,  l \geq 1\}$ form a complete collection of indecomposables.\end{prop}
\noindent
\textit{Proof.} Let $(V,x)$ be an indecomposable. As $x$ is nilpotent, $\text{Ker}\;x \neq \{0\}$. Moreover, since $x$ maps $V_j \to V_{j-1}$, the kernel splits as a direct sum 
\begin{equation}\label{E:splitx}
{Ker}\;x=\bigoplus_j (V_j \cap {Ker}\;x).
\end{equation}
By the classification of nilpotent endomorphisms of $k$-vector spaces (see Theorem~\ref{T:Jordan}), there exists a basis $u_1, \ldots, u_{{dim}\;V}$ over which $x$ is in canonical form, i.e. for which $x(u_j) \in \{u_{j-1},0\}$. By (\ref{E:splitx}), we may even choose the $u_i$'s to be homogeneous. The Proposition easily follows.\qed

\vspace{.15in}

Using the above Proposition, we see that the set of objects of $Rep^{nil}_k\vec{Q}$ is in canonical bijection with the set $\Pi^n$ of $n$-tuples of partitions $(\underline{\lambda}_1, \ldots, \underline{\lambda}_n)$, via the assignement
$$(\underline{\lambda}_1, \ldots, \underline{\lambda}_n) \mapsto \bigoplus_{i \in I} \bigoplus_{j} I_{[i;\lambda_i^j]}.$$
This set is indeed independent of $k$. The indivisible imaginary root is $\delta=\a_0+ \a_1 \cdots +\a_{n-1}$.

\vspace{.2in}

Fix $k=\mathbb{F}_q$ and set $\nu=q^{\frac{1}{2}}$. By Ringel's Theorem~\ref{T:RingelHall} there is an embedding 
\begin{align*}
\Psi: \U_{\nu}(\n_+) &\to \H_{\vec{Q}}\\
E_i &\mapsto [S_i]=[I_{[i;1]}]
\end{align*}
where $\n_+ \subset \widehat{\mathfrak{sl}}_n$ is the standard positive nilpotent subalgebra.

\vspace{.05in}

The image of $\Psi$ is the composition algebra $\mathbf{C}_{\vec{Q}} \subset \H_{\vec{Q}}$.
Let $\U'$ be the two-sided ideal of $\H_{\vec{Q}}$ generated by $\{[S_i]\;|\;i \in I\}$, and let $\mathbf{R}=(\U')^{\perp}$ be its orthogonal in $\H_{\vec{Q}}$ with respect to Green's scalar product $(\;,\;)$. As $(\;,\;)$ is a Hopf pairing, $\mathbf{R}$ is a subcoalgebra (because $\U'$ is an ideal) and a subalgebra (because $\U'$ is also a coideal).
The structure of $\H_{\vec{Q}}$ is given by the next Proposition, proved in \cite{S1}~:

\vspace{.1in}

\begin{prop}[S.]\label{P:Hallcyclic} The following hold~:
\begin{enumerate}
\item[i)] As a bialgebra, $\mathbf{R}$ is isomorphic to a polynomial ring $\mathbf{R}\simeq \C[\mathbf{x}_1, \mathbf{x}_2, \ldots]$ where $deg(\mathbf{x}_i)=i\delta$, and $\Delta'(\mathbf{x}_{i})=\mathbf{x}_i \otimes 1 + 1 \otimes \mathbf{x}_i$.
\item[ii)] $\mathbf{R}$ is a central subalgebra of $\H_{\vec{Q}}$ and the multiplication map gives an isomorphism $\mathbf{C}_{\vec{Q}} \otimes \mathbf{R} \stackrel{\sim}{\to} \H_{\vec{Q}}$.
\end{enumerate}
\end{prop}

\vspace{.1in}

A detailed analysis of $\mathbf{R}$ was made by Hubery in \cite{HuberyR}, who obtained the following remarkable result~:

\vspace{.1in}

\begin{prop}[Hubery]\label{P:HuberyR} The center of $\H_{\vec{Q}}$ coincides with $\mathbf{R}$. It is generated by the elements
$$\mathbf{c}_r=(-1)^r\nu^{-2rn}\sum_{M \in Z_r} (-1)^{\text{dim}(\text{End}(M))} a_M[M],$$
where $Z_r$ is the set of representations $M$ of dimension $r\delta$ whose socle $soc\;M$ is square-free, i.e for which $soc\;M \simeq \bigoplus_i S_i^{\oplus n_i}$ with $n_i \leq 1$. These elements
satisfy $\Delta'(\mathbf{c}_r)=\sum_{s=0}^r \mathbf{c}_{r-s} \otimes \mathbf{c}_s$.
\end{prop}

\vspace{.15in}

Observe that, by Proposition~\ref{P:Hallcyclic} i), $\mathbf{R}$ is abstractly isomorphic, as a bialgebra, to Macdonald's ring of symmetric functions $\LLambda=\C[x_1, x_2,\ldots]^{\mathfrak{S}_{\infty}}$ (see Lecture~2, Section~2.4.). Using Proposition~\ref{P:HuberyR} we can make this precise. Consider the generating function $C(z)=1+\sum_r \mathbf{c}_rT^r$. The equation $P(z)C(z)=\frac{d}{dz}C(z)$ has unique solution $P(z)=\sum_{r\geq 1} \mathbf{p}_rz^{r-1}$. The elements $\mathbf{p}_r$ are primitive (that is, $\Delta'(\mathbf{p}_r)=\mathbf{p}_r \otimes 1 + 1 \otimes \mathbf{p}_r$), and they clearly freely generate $\mathbf{R}$. Thus there is a unique isomorphism
\begin{align*}
\Phi_n: \mathbf{R} &\stackrel{\sim}{\to} \LLambda\\
\mathbf{p}_r &\mapsto (1-\nu^{-2nr})p_r
\end{align*}
where $p_r$ is the power sum symmetric function.
This generalizes the map $\Phi$ (or more precisely the map $\Phi_{(q)}$) of Section~2.4 which corresponds to the degenerate case $n=1$. We may view the rings $\mathbf{R}$ (for $n \in \N$) as another realization of $\LLambda$. Just as in Lecture~2, one may wonder what the restriction of Green's scalar product to $\mathbf{R}$, transported to $\LLambda$ via $\Phi_n$, is. This turns out to bring nothing new~:

\vspace{.1in}

\begin{prop}[Hubery] The pullback under $\Phi_n$ of Green's scalar product on $\mathbf{R}$ is the Hall-Littlewood scalar product given by
$$\{x,yz\}=\{\Delta'(x) ,y\otimes z\},$$
\begin{equation}\label{E:HLnscalarprod}
\{p_r,p_s\}=\delta_{r,s}\frac{r}{q^r-1}
\end{equation}
\end{prop}

\vspace{.15in}

To finish, let us indicate a second decomposition of $\H_{\vec{Q}}$, as a noncommutative product this time. Let $\mathcal{C}$ be the full subcategory of $Rep_k\vec{Q}$ consisting of those representatiions $(V,\underline{x})$ for which $x_i: V_i \to V_{i-1}$ is an isomorphism for all $i \neq 1$. In other words, objects of $\mathcal{C}$ are all direct sums ofindecomposables of the form $I_{[0;nl]}$ for some $l \geq 1$. It is clear that $\mathcal{C}$ is an abelian subcategory closed under extensions, and which is equivalent to $Rep_k\vec{Q}_0$ where $\vec{Q}_0$ is the Jordan quiver. Thus by Corollary~\ref{C:cathall} there is an embedding of algebras
\begin{align*}
\Theta_n: \H=\H_{\vec{Q}_0} &\to \H_{\vec{Q}}\\
[I_{\underline{\lambda}}] &\mapsto [\bigoplus_j I_{[0;n\lambda_j]}]
\end{align*}
Let $\mathbf{K} \subset \H_{\vec{Q}}$ be the image of this embedding. The subalgebras $\mathbf{C}_{\vec{Q}}$ and $\mathbf{K}$ do not commute. However,

\vspace{.1in}

\begin{prop}[S., \cite{SDuke}]\label{P:halldecomp} The multiplication map defines an isomorphism of vector spaces $\mathbf{K} \otimes \mathbf{C}_{\vec{Q}} \stackrel{\sim}{\to} \H_{\vec{Q}}$.\end{prop}

\vspace{.15in}

\addtocounter{theo}{1}
\noindent \textbf{Remarks \thetheo .} i) There exists Hall polynomials for the cyclic quiver --see \cite{Ringelcomposcyclic}.\\
ii) There is, to my knowledge, no completely elementary description of the elements of the composition subalgebra $\mathbf{C}_{\vec{Q}} \subset \H_{\vec{Q}}$.\\
iii) The decomposition $\H_{\vec{Q}} \simeq \mathbf{C}_{\vec{Q}} \otimes \mathbf{R}$ where $\mathbf{R}$ is a central polynomial ring exists in fact for \textit{all} tame quivers, as was shown (independently of \cite{S1}) by Hua and Xiao
\cite{HuaXiao}. Of course in this case the structure of $\mathbf{R}$, i.e. the number of generators per degree, depends on the ground field $k$.\\
iv) One might ask what the structure of $\H_{\vec{Q}}$ is, for an \textit{arbitrary} quiver. By a theorem of Sevenhant and Van den Bergh (see \cite{SVdB}), it is always a quantum group associated to a Borcherds algebra (usually of infinite rank). The precise determination of that Borcherds algebra is however still very much an open problem.

\vspace{.2in}

\centerline{\textbf{3.6. Structure theory for tame quivers.}}
\addcontentsline{toc}{subsection}{\tocsubsection {}{}{\; 3.6. Structure theory for tame quivers.}}

\vspace{.15in}

\paragraph{} We collect in this Section several results concerning the category $Rep_k\vec{Q}$ when $\vec{Q}$ is a tame quiver. This is necessary before we can describe the exact structure of the composition subalgebra $\mathbf{C}_{\vec{Q}} \subset \H_{\vec{Q}}$ (see Section~3.7.). This will also be important for the next Lecture.

Throughout this Section we fix a tame quiver $\vec{Q}$ \textit{which is not a cyclic quiver} and an arbitrary field $k$. Hence the categories $Rep_k \vec{Q}$ and $Rep^{nil}_k\vec{Q}$ coincide. Proofs of the assertions made below can be found in \cite{CrawleyNotes}. The fundamental tool to use here is the following

\vspace{.1in}

\begin{theo}[Auslander-Reiten] There exists a unique pair of adjoint functors $\tau,\tau^-: Rep_k\vec{Q} \to Rep_k\vec{Q}$ equipped with natural isomorphisms
\begin{align*}
{Ext}^1(M,N)^* &\simeq {Hom}(N,\tau M)\\
{Ext}^1(M,N)^* &\simeq {Hom}(\tau^-N, M).
\end{align*}
\end{theo}

\vspace{.1in}

Such pairs of adjoint functors are now called \textit{Serre functors} in general, and appear to be a very powerful tool for commutative or noncommutative algebraic geometry (see \cite{VandenBergh}). In the case of hereditary algebras they were invented by Auslander and Reiten and are known as \textit{Auslander-Reiten translations}.

\vspace{.15in}

Obviously, $\tau M=0$ if and only if  $M$ is projective while $\tau^{-} N=0$ if and only if $N$ is injective. In addition, if $M,M'$ are non zero then $(\tau M=M') \Leftrightarrow (\tau^- M'=M)$. This means that
$\tau,\tau^-$ are inverse bijections between non-projective and non-injective objects in $Rep_k\vec{Q}$.

\vspace{.1in}

An indecomposable representation $M$ is called~:
\begin{enumerate}
\item[-] \textit{preprojective} if $\tau^i M =0$ for $i \gg 0$,
\item[-] \textit{preinjective} if $\tau^{-i} M =0$ for $i \gg 0$,
\item[-] \textit{regular} if $\tau^i M \neq 0$ for $i \in \Z$.
\end{enumerate}

Let $\mathcal{P}, \mathcal{I}, \mathcal{R}$ denote the set of preprojective, preinjective or regular indecomposables respectively. Though it is not obvious from the definitions, the sets $\mathcal{P}$ and $\mathcal{I}$ are disjoint (recall that $\vec{Q}$ is tame). Clearly, the AR translations $\tau,\tau^-$ preserve each of the sets $\mathcal{P},\mathcal{I},\mathcal{R}$. Call a module preprojective if all of its indecomposable summands are preprojective, and let $\mathbb{P}$ denote the full category of $Rep_k\vec{Q}$ consisting of preprojective modules. The categories $\mathbb{I},\mathbb{R}$ are defined in a similar manner.

\vspace{.15in}

\begin{prop}\label{P:extiszero} The categories $\mathbb{P}$ and $\mathbb{I}$ are exact and stable under extensions. The category$\mathbb{R}$ is abelian and stable under extensions. In addition, if $M \in \mathbb{P}, N \in \mathbb{I}$ and $L \in \mathbb{R}$ then
\begin{align}\label{E:preprojeqs}
{Hom}(N,M) &=\mathrm{Hom}(N,L)={Hom}(L,M)=\{0\}\\
{Ext}^1(M,N)&=\mathrm{Ext}^1(L,N)={Ext}^1(M,L)=\{0\}
\end{align}
\end{prop}

\vspace{.1in}
 
Let us draw a picture of the indecomposables as follows

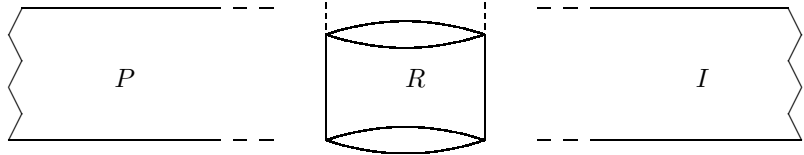
\vspace{.1in}

\centerline{
\begin{picture}(300,60)
\put(0,0){\line(1,0){80}}
\put(220,0){\line(1,0){80}}
\put(5,50){\line(1,0){75}}
\put(220,50){\line(1,0){75}}
\put(85,0){\line(1,0){5}}
\put(95,0){\line(1,0){5}}
\put(85,50){\line(1,0){5}}
\put(95,50){\line(1,0){5}}
\put(200,0){\line(1,0){5}}
\put(210,0){\line(1,0){5}}
\put(200,50){\line(1,0){5}}
\put(210,50){\line(1,0){5}}
\put(120,0){\line(0,1){40}}
\put(180,0){\line(0,1){40}}
\put(120,42){\line(0,1){2}}
\put(120,46){\line(0,1){2}}
\put(120,50){\line(0,1){2}}
\put(180,42){\line(0,1){2}}
\put(180,46){\line(0,1){2}}
\put(180,50){\line(0,1){2}}
\put(0,0){\line(1,2){5}}
\put(5,10){\line(-1,2){5}}
\put(0,20){\line(1,2){5}}
\put(5,30){\line(-1,2){5}}
\put(0,40){\line(1,2){5}}
\put(300,0){\line(-1,2){5}}
\put(295,10){\line(1,2){5}}
\put(300,20){\line(-1,2){5}}
\put(295,30){\line(1,2){5}}
\put(300,40){\line(-1,2){5}}
\put(40,20){${P}$}
\put(150,20){${R}$}
\put(260,20){${I}$}
\qbezier(120,0)(150,-10)(180,0)
\qbezier(120,40)(150,30)(180,40)
\qbezier(120,0)(150,10)(180,0)
\qbezier(120,40)(150,50)(180,40)
\end{picture}}

\vspace{.2in}

\noindent
where the projective $\{P(i)\;|\; i\in I\}$ are on the extreme left, followed by the modules $\{\tau^- P(i)\;|\; i\in I\}$, etc..; the injectives $\{I(i)\;|\; i \in I\}$ are on the extreme right, then $\{\tau I(i)\;|\; i\in I\}$, etc..; the regular modules are put in the middle. The above Proposition says that morphisms go from left to right, whereas extensions go from right to left.

\vspace{.2in}

Let us now proceed to describe completely the structure of the regular indecomposables. Call such a module \textit{simple} if it contains no nontrivial regular submodule, and call it \textit{homogeneous} if $\tau M \simeq M$. If $M$ is simple then $k_M={End}(M)$ is a field extension of $k$; the \textit{degree} of $M$ is the index $[k_M:k]$.

Any regular module $M \in \mathbb{R}$ has a composition series whose factors are regular simples. Conversely, given regular simple modules $M_1, \ldots, M_r$ we denote by $\mathcal{C}_{M_1, \ldots, M_r}$ the subcategory of $\mathbb{R}$ consisting of modules with a composition series whose factors all belong to $\{M_1, \ldots, M_r\}$.

\vspace{.15in}

\begin{prop}\label{P:noname1} The following hold~:

\vspace{.05in}

\noindent
i) If $M$ is simple homogeneous then $\underline{\text{dim}}\;M=[k_M:k]\delta$ and $\mathcal{C}_M$ is equivalent to the category of nilpotent representations of the Jordan quiver over $k_M$.

\vspace{.05in}

\noindent
ii) If $M$ is non-homogeneous then $k_M=k$ and there exists an integer $r>0$ such that $\tau^r M\simeq M$, $\tau^s M \not\simeq M$ for $0<s<r$. Moreover we have $\sum_{s=0}^{r-1} \underline{{dim}}\;\tau^sM=\delta$ and the category $\mathcal{C}_{M,\tau M, \ldots, \tau^{r-1}M}$ is equivalent to the category of nilpotent representations of the cyclic quiver of length $r$, over $k$.
\end{prop}

\vspace{.15in}

Of course, under the equivalence in statement i) above the object $M$ goes to the simple object $S$ (see Lecture~2), while under the equivalence in ii), the simple objects $M, \tau M, \ldots, \tau^{r-1} M$ are mapped to the simple objects $S_1, S_2, \ldots, S_r$ (see Section~3.5.).
\vspace{.15in}

Before stating the main structure Theorem for the category $\mathbb{R}$, let us recall that a \textit{closed
point} $x$ of $\mathbb{P}^1$ over a field $k$ is nothing but a homogeneous maximal ideal $\mathfrak{m}_x$ of the graded ring $k[X_0,X_1]$, and that the \textit{degree} of $x$ is simply the index $[k[X_1,X_2]/\mathfrak{m}_x:k]$ (hence the usual points of $\mathbb{P}^1(k)$ correspond to closed points of degree one).

\vspace{.1in}

\begin{theo}[Ringel]\label{T:Ringelstructure} Let $d$ and $p_1, \ldots, p_d$ be attached to $\vec{Q}$ as in the table (\ref{E:table}) below. Then
\begin{enumerate}
\item[i)] There is a degree-preserving bijection $M_x \leftrightarrow x$ between the set of homogeneous regular simple modules and closed points of $\mathbb{P}^1 \backslash D$ where $D$ consists of $d$ points of degree one.
\item[ii)] There are exactly $d$ $\tau$-orbits $\mathcal{O}_1, \ldots, \mathcal{O}_d$ of non-homogeneous regular simple modules, and they are of size $p_1, \ldots, p_d$ respectively.
\item[iii)] The whole category $\mathbb{R}$ decomposes as a direct sum of orthogonal blocks
$$\mathbb{R} = \prod_{x \in \mathbb{P}^1\backslash D} \hspace{-.08in}\mathcal{C}_{M_x} \times\; \prod_{l=1, \ldots, d}\hspace{-.08in} \mathcal{C}_{\mathcal{O}_l}.$$
\end{enumerate}
\end{theo}

\vspace{.1in}

\begin{equation}\label{E:table}
\begin{tabular}{|c|c|c|}
\hline
type of\;$\vec{Q}$ & $d$ & $p_1, \ldots, p_d$\\
\hline
$A_1^{(1)}$ & 0 & \\
\hline
$A_n^{(1)}, \;n>1$ & 2 & $p_1=\# \mathrm{arrows\;going\;clockwise}$\\
& & $p_2=\# \mathrm{arrows\;going\;counterclockwise}$\\
\hline
$D_n^{(1)}$ & 3 & $2,2,n-2$\\
\hline
$E_n^{(1)},\; n=6,7,8$ & 3 & $2,3,n-3$\\
\hline
\end{tabular} 
\end{equation}

\vspace{.15in}

The subcategories $\mathcal{C}_{M_x}$ and $\mathcal{C}_{\mathcal{O}_l}$ are called \textit{homogeneous} and \textit{non-homogeneous} \textit{tubes} respectively.

\vspace{.05in}

The values of $d$ and $p_1, \ldots, p_d$ can in fact be read off quite simply from the \textit{finite} Dynkin diagram of same type as $\vec{Q}$. Note that these are all shaped as ``stars'' with a central vertex out of which some branches are coming; $d$ is simply the number of such branches, and $p_1, \ldots, p_d$ are their respective lengths (there is an ambiguity in type $A$, and one has to be more careful there ).

\vspace{.15in}

\addtocounter{theo}{1}

\paragraph{\textbf{Example \thetheo.}} Let us consider again the Kronecker quiver of Example~3.11.
of which we keep the notations~:

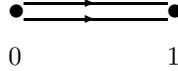
\vspace{.15in}

\centerline{
\begin{picture}(200, 10)
\put(60,0){\circle*{5}}
\put(120,0){\circle*{5}}
\put(85,-3){\vector(1,0){5}}
\put(85,3){\vector(1,0){5}}
\put(63,3){\line(1,0){54}}
\put(63,-3){\line(1,0){54}}
\put(57,-20){$0$}
\put(117,-20){$1$}
\end{picture}}

\vspace{.3in}

The indecomposable $I_{\delta}^{(\lambda,\mu)}$ is always regular, simple and homogeneous. When $k$ is not algebraically closed, there are also higher-order regular simples (see \cite{RingelTame}). Hence, all together, the category $Rep_k(\vec{Q})$ looks like this~:

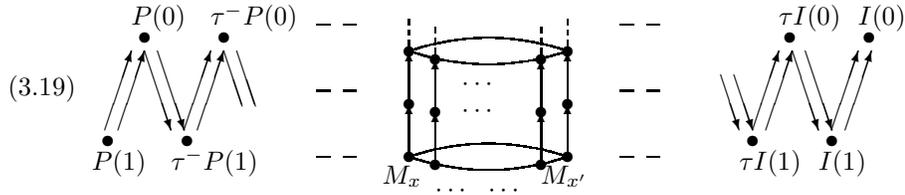
\vspace{.5in}

\begin{equation}\label{E:picturekronecker}
\centerline{
\begin{picture}(300,20)
\put(85,0){\line(1,0){5}}
\put(95,0){\line(1,0){5}}
\put(85,25){\line(1,0){5}}
\put(95,25){\line(1,0){5}}
\put(200,25){\line(1,0){5}}
\put(210,25){\line(1,0){5}}
\put(85,50){\line(1,0){5}}
\put(95,50){\line(1,0){5}}
\put(200,0){\line(1,0){5}}
\put(210,0){\line(1,0){5}}
\put(200,50){\line(1,0){5}}
\put(210,50){\line(1,0){5}}
\put(0,-5){$P(1)$}
\put(6,6){\circle*{4}}
\put(5,10){\vector(1,3){10}}
\put(10,10){\vector(1,3){10}}
\put(25,40){\vector(1,-3){10}}
\put(20,40){\vector(1,-3){10}}
\put(55,40){\line(1,-3){7}}
\put(50,40){\line(1,-3){7}}
\put(35,10){\vector(1,3){10}}
\put(40,10){\vector(1,3){10}}
\put(15,50){$P(0)$}
\put(20,45){\circle*{4}}
\put(50,45){\circle*{4}}
\put(36,6){\circle*{4}}
\put(280,6){\circle*{4}}
\put(250,6){\circle*{4}}
\put(264,45){\circle*{4}}
\put(294,45){\circle*{4}}
\put(30,-5){$\tau^-P(1)$}
\put(45,50){$\tau^-P(0)$}
\put(290,50){$I(0)$}
\put(265,40){\vector(1,-3){10}}
\put(270,40){\vector(1,-3){10}}
\put(280,10){\vector(1,3){10}}
\put(285,10){\vector(1,3){10}}
\put(238,31){\vector(1,-3){7}}
\put(243,31){\vector(1,-3){7}}
\put(255,10){\vector(1,3){10}}
\put(250,10){\vector(1,3){10}}
\put(275,-5){$I(1)$}
\put(260,50){$\tau I(0)$}
\put(245,-5){$\tau I(1)$}
\put(120,0){\line(0,1){40}}
\put(180,0){\line(0,1){40}}
\qbezier(120,0)(150,-10)(180,0)
\qbezier(120,40)(150,30)(180,40)
\qbezier(120,0)(150,10)(180,0)
\qbezier(120,40)(150,50)(180,40)
\put(120,42){\line(0,1){2}}
\put(120,46){\line(0,1){2}}
\put(120,50){\line(0,1){2}}
\put(180,42){\line(0,1){2}}
\put(180,46){\line(0,1){2}}
\put(180,50){\line(0,1){2}}
\put(130,39){\line(0,1){2}}
\put(130,43){\line(0,1){2}}
\put(130,47){\line(0,1){2}}
\put(170,39){\line(0,1){2}}
\put(170,43){\line(0,1){2}}
\put(170,47){\line(0,1){2}}
\put(120,0){\vector(0,1){20}}
\put(120,0){\circle*{4}}
\put(120,40){\circle*{4}}
\put(120,20){\vector(0,1){20}}
\put(120,20){\circle*{4}}
\put(110,-10){$M_x$}
\put(130,-15){$\cdots$}
\put(140,15){$\cdots$}
\put(140,25){$\cdots$}
\put(150,-15){$\cdots$}
\put(180,0){\vector(0,1){20}}
\put(180,0){\circle*{4}}
\put(180,40){\circle*{4}}
\put(180,20){\vector(0,1){20}}
\put(180,20){\circle*{4}}
\put(170,-10){$M_{x'}$}
\put(130,-3){\vector(0,1){20}}
\put(130,-3){\circle*{4}}
\put(130,37){\circle*{4}}
\put(130,17){\vector(0,1){20}}
\put(130,17){\circle*{4}}
\put(170,-3){\vector(0,1){20}}
\put(170,-3){\circle*{4}}
\put(170,37){\circle*{4}}
\put(170,17){\vector(0,1){20}}
\put(170,17){\circle*{4}}
\end{picture}}
\end{equation}

\vspace{.3in}

For simplicity, only the regular modules generated by simples of degree one have been drawn in the above picture, but it is important to keep in mind that when $k$ is not algebraically closed, there are infinitely many homogeneous tubes of regular indecomposables which start in degrees $l\delta$, $l >1$. Also, following custom, we drew $\text{dim}(M,M')$ arrows between neighboring indecomposables
$M$ and $M'$.

\endexample

\vspace{.25in}

\addtocounter{theo}{1}

\paragraph{\textbf{Example \thetheo.}} Let $\vec{Q}$ be the quiver of type $D_4^{(1)}$ of Example~3.12.~:

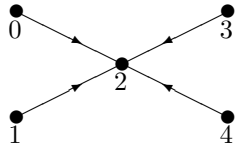
\vspace{.4in}

\centerline{
\begin{picture}(200, 10)
\put(60,-20){\circle*{5}}
\put(100,0){\circle*{5}}
\put(140,-20){\circle*{5}}
\put(140,20){\circle*{5}}
\put(60,20){\circle*{5}}
\put(57,-30){$1$}
\put(57,10){$0$}
\put(97,-10){$2$}
\put(137,-30){$4$}
\put(137,10){$3$}
\put(60,-20){\line(2,1){40}}
\put(60,20){\line(2,-1){40}}
\put(100,0){\line(2,-1){40}}
\put(100,0){\line(2,1){40}}
\put(80,-10){\vector(2,1){5}}
\put(80,10){\vector(2,-1){5}}
\put(120,-10){\vector(-2,1){5}}
\put(120,10){\vector(-2,-1){5}}
\end{picture}}

\vspace{.55in}

The indecomposable $I_{\delta}^{(\lambda_0,\lambda_1,\lambda_3,\lambda_4)}$ is always regular. It is simple if and only if $\lambda_0, \lambda_1,\lambda_3,\lambda_4$ are all distinct, in which case it is also homogeneous. The others, which are not simple, belong to the three non-homogeneous tubes generated by the $\tau$-orbits $\{I_{123},I_{024}\}, \{I_{124},I_{023}\},\{I_{012},I_{234}\}$.  Compare with the picture ~(\ref{E:ponethreepoints}).

\endexample

\vspace{.2in}

\centerline{\textbf{3.7. The composition algebra of a tame quiver.}}
\addcontentsline{toc}{subsection}{\tocsubsection {}{}{\; 3.7. The composition algebra of a tame quiver.}}

\vspace{.15in}

\paragraph{}To finish off this Lecture, we proceed to describe as precisely as possible, following Zhang \cite{Zhang} and Hubery \cite{Hubery}, the elements of the composition algebra $\mathbf{C}_{\vec{Q}}$ of a tame quiver $\vec{Q}$. We assume here again that $\vec{Q}$ is not a cyclic quiver. The ground field is $k=\mathbb{F}_q$ and as usual $\nu=q^{\frac{1}{2}}$. We will freely use the no(ta)tions of the previous Section.

\vspace{.15in}

We start by observing the following important fact. Let $\H_{\mathbb{P}}, \H_{\mathbb{R}}, \H_{\mathbb{I}}$ be the Hall algebras of the exact categories $\mathbb{P},\mathbb{R}$ and $\mathbb{I}$. As all of these categories are stable under extensions, Corollary~\ref{C:catexacthall} gives us natural algebra embeddings $\H_{\mathbb{P}} \hookrightarrow \H_{\vec{Q}}, \H_{\mathbb{R}} \hookrightarrow \H_{\vec{Q}}$ and $\H_{\mathbb{I}} \hookrightarrow \H_{\vec{Q}}$. The categories
$\mathcal{C}_{\mathcal{O}_1}, \ldots, \mathcal{C}_{\mathcal{O}_d}$ generated by non-homogeneous regular simple representations are all abelian and also closed under extensions. Thus we have embeddings $\H_i:=\H_{\mathcal{C}_{\mathcal{O}_i}} \hookrightarrow \H_{\mathbb{R}} \hookrightarrow \H_{\vec{Q}}$. By Proposition~\ref{P:noname1} ii), $\H_i$ is isomorphic to the Hall algebra of the cyclic quiver studied in Section~3.5. The composition subalgebra of $\H_i$ will be denoted $\mathbf{C}_i $. 

\vspace{.15in}

\begin{lem} The multiplication map induces an isomorphism of vector spaces $\H_{\mathbb{P}} \otimes \H_{\mathbb{R}} \otimes \mathbf{H}_{\mathbb{I}} \simeq \H_{\vec{Q}}$.
\end{lem}

\noindent
\textit{Proof.} Any object $M$ of $Rep_k\vec{Q}$ decomposes as a direct sum $M=M_P \oplus M_R \oplus M_I$ of a preprojective, regular, and preinjective module. We claim that
\begin{equation}\label{E:PRIdecomp}
[M_P]\cdot [M_R]\cdot [M_I]=\nu^e[M]
\end{equation}
where $e=\langle M_P,M_R\rangle_a + \langle M_P,M_I\rangle_a + \langle M_R,M_I\rangle_a$. First note that by Proposition~\ref{P:extiszero} any extension of $M_P,M_R,M_I$ in this order is trivial, so that the only module appearing in the product is $[M]$. Furthermore, since ${Hom}(M_I, M_P \oplus M_R)=\{0\}$ and ${Hom}(M_R,M_P)=\{0\}$, there is a unique filtration $0 \subset L_1 \subset L_2 \subset M$ satisfying $L_1 \simeq M_I, L_2/L_1 \simeq M_R$ and $M/L_2 \simeq M_P$. Thus the Hall number $P_{M_P,M_R,M_I}^M$ is equal to $1$ and (\ref{E:PRIdecomp}) follows. From this it is easy to deduce that $m: \H_{\mathbb{P}} \otimes \H_{\mathbb{R}} \otimes \H_{\mathbb{I}} \to \H_{\vec{Q}}$ is both injective and surjective.\qed

\vspace{.15in}

The composition algebra $\mathbf{C}_{\vec{Q}}$ turns out to be itself compatible with the above triangular decomposition. This will be a consequence of the fact that it is stable under the old coproduct map $\Delta'$, as the proof of the next Proposition shows. Set $\mathbf{C}_{\mathbb{P}}=\mathbf{C}_{\vec{Q}} \cap \H_{\mathbb{P}}$ and define $\mathbf{C}_{\mathbb{R}}, \mathbf{C}_{\mathbb{I}}$ in a similar fashion.

\vspace{.1in}

\begin{prop}[Zhang]\label{P:halldecompopo} The multiplication map induces an isomorphism of vector spaces $\mathbf{C}_{\mathbb{P}} \otimes \mathbf{C}_{\mathbb{R}} \otimes \mathbf{C}_{\mathbb{I}} \simeq \mathbf{C}_{\vec{Q}}$.
\end{prop}

\noindent
\textit{Proof.} Set $\mathbf{C}_{\mathbb{R}\oplus \mathbb{I}}=\mathbf{C}_{\vec{Q}} \cap \H_{\mathbb{R}\oplus \mathbb{I}}$. We start by showing that $\mathbf{C}_{\vec{Q}}=\mathbf{C}_{\mathbb{P}}\cdot \mathbf{C}_{\mathbb{R} \oplus \mathbb{I}}$.
Fix $x \in \mathbf{C}_{\vec{Q}}$ and decompose it as a sum
\begin{equation}\label{E:noname2}
x=\sum_i x_{\a_i,\beta_i}, \qquad x_{\a_i,\beta_i} \in (\H_{\mathbb{P}}[\a_i]) \cdot (\H_{\mathbb{R} \oplus \mathbb{I}}[\beta_i]).
\end{equation}
In the above, $\a_i,\beta_i$ are dimension vectors satisfying $\a_i + \beta_i=deg(x)$, and $A[\gamma]$ indicates the component of the algebra $A$ of degree $\gamma$. We will prove by descending induction on $\a_i$ that $x_{\a_i,\beta_i} \in \H_{\mathbb{P}} \cdot \mathbf{C}_{\mathbb{R}\oplus \mathbb{I}}$. Let $i$ be such that $\a_i$ is maximal, and let us consider the component $\Delta'_{\a_i,\beta_i}(x)$ of $\Delta'(x)$ which lands in $(\H_{\vec{Q}}[\a_i]) \otimes (\H_{\vec{Q}}[\beta_i])$. As $\mathbf{C}_{\vec{Q}}$ is generated by the simples $\{[S_j]\}$ and as
$\Delta'([S_j])=[S_j] \otimes 1 + 1 \otimes [S_j]$, the composition algebra $\mathbf{C}_{\vec{Q}}$ is stable under the coproduct, and thus $\Delta'_{\a_i,\beta_i}(x)$ lies in $(\mathbf{C}_{\vec{Q}}[\a_i]) \otimes (\mathbf{C}_{\vec{Q}}[\beta_i])$. Write $x_{\a_i,\beta_i}=\sum_l [P_l^i] \cdot u_{P_l^i}$ where $P_l^i$ is preprojective and $u_{P_l^i} \in \H_{\mathbb{R}\oplus \mathbb{I}}$. Let $\pi: \H_{\vec{Q}} \otimes \H_{\vec{Q}} \to \H_{\mathbb{P}} \otimes \H_{\vec{Q}}$ stand for the natural projection (with respect to the defining basis $\{[M]\;|\;M \in Obj(Rep_k\vec{Q})\}$).
We claim that 
\begin{equation}\label{E:noname3}
\pi\circ\Delta'_{\a_i,\beta_i}(x_{\a_i,\beta_i})= \sum_l  [P_l^i] \otimes u_{P_l^i}.
\end{equation}
Indeed, if $\pi\circ\Delta'_{\a_i,\beta_i}([M]) \neq 0$ then there exists a projection $M \twoheadrightarrow P$ for some preprojective module $P$ of dimension $\a_i$. By Proposition~\ref{P:extiszero}, this implies that the preprojective component of $M$ is of dimension at least $\a_i$. Since by our assumption $\a_i$ was chosen to be maximal among all dimensions appearing in (\ref{E:noname2}), we deduce the first equality of (\ref{E:noname3}). As for the second equality, it follows from the definition of $\Delta'$ and the fact that any module of the form $P\oplus T$ with $P$ preprojective and $T$ in $\mathbb{R} \oplus \mathbb{I}$ has a unique submodule isomorphic to $T$. As a consequence, we have
$$u_{P_l^i}= ([P_l^i]^* \otimes 1) \Delta'_{\a_i,\beta_i}(x) \in ([P_l^i]^*\otimes 1) \cdot \mathbf{C}_{\vec{Q}} \otimes \mathbf{C}_{\vec{Q}}\subset\mathbf{C}_{\vec{Q}}$$
and hence $u_{P_l^i} \in \mathbf{C}_{\mathbb{R}\oplus \mathbb{I}}$. This proves the first step of the induction. Next, fix $j$ and assume that $x_{\a_i,\beta_i} \in \H_{\mathbb{P}} \cdot \mathbf{C}_{\mathbb{R}\oplus \mathbb{I}}$ for any $i$ for which $\a_i > \a_j$. Let us write as before
$x_{\a_j,\beta_j}=\sum_l [P_l^j] \cdot u_{P_l^j}$ where $P_l^j$ is preprojective and $u_{P_l^j} \in \H_{\mathbb{R}\oplus \mathbb{I}}$, and consider again $\pi \circ \Delta'_{\a_j,\beta_j}(x)$. This time, we have, using Proposition~\ref{P:extiszero}
$$\pi \circ \Delta'_{\a_j,\beta_j}(x)=\sum_l [P^j_l] \otimes u_{P^j_l} +\sum_{\a_i > \a_j} \sum_h \Delta'_{\a_j,\a_i-\a_j}([P_h^i]) \cdot (1 \otimes u_{P^i_h}).$$
Hence,
\begin{equation*}
\begin{split}
([P_l^j]^* \otimes 1)&\pi \circ \Delta'_{\a_j,\beta_j}(x)\\
&= u_{P^j_l}+\sum_{\a_i > \a_j} \sum_h ([P_l^j]^* \otimes 1)\Delta'_{\a_j,\a_i-\a_j}([P_h^i]) \cdot (1 \otimes u_{P^i_h}).
\end{split}
\end{equation*}
By the induction hypothesis $u_{P^i_h} \in \mathbf{C}_{\mathbb{R} \oplus \mathbb{I}}$ hence
$$u_{P_l^j} \in (\mathbf{C}_{\vec{Q}} + \H_{\mathbb{P}} \cdot \mathbf{C}_{\mathbb{R}\oplus \mathbb{I}}) \cap \H_{\mathbb{R} \oplus \mathbb{I}}=\mathbf{C}_{\vec{Q}} \cap \H_{\mathbb{R}\oplus \mathbb{I}}=\mathbf{C}_{\mathbb{R}\oplus \mathbb{I}}$$
as desired. This concludes the induction step and proves that $\mathbf{C}_{\vec{Q}} \subset \H_{\mathbb{P}} \cdot \mathbf{C}_{\mathbb{R}\oplus \mathbb{I}}$. 

An absolutely symmetrical argument, based on a descending induction on $\beta_i$ this time, shows that 
$\mathbf{C}_{\vec{Q}} \subset \mathbf{C}_{\mathbb{P}} \cdot \mathbf{H}_{\mathbb{R}\oplus \mathbb{I}}$.
But then
$$\mathbf{C}_{\vec{Q}} \subset \big(\mathbf{C}_{\mathbb{P}} \cdot \mathbf{H}_{\mathbb{R}\oplus \mathbb{I}}\big) \cap \big(\H_{\mathbb{P}} \cdot \mathbf{C}_{\mathbb{R}\oplus \mathbb{I}}\big)=\mathbf{C}_{\mathbb{P}} \cdot \mathbf{C}_{\mathbb{R}\oplus \mathbb{I}}.$$
Finally, we may repeat the same argument again for $\mathbf{C}_{\mathbb{R}\oplus \mathbb{I}}$ in place of $\mathbf{C}_{\vec{Q}}$ to obtain the inclusion $\mathbf{C}_{\mathbb{R}\oplus\mathbb{I}}\subset \mathbf{C}_{\mathbb{R}} \cdot \mathbf{C}_{\mathbb{I}}$. This yields the inclusion $\mathbf{C}_{\vec{Q}} \subset \mathbf{C}_{\mathbb{P}} \cdot \mathbf{C}_{\mathbb{R}} \cdot \mathbf{C}_{\mathbb{I}}$. The reverse inclusion is obvious. We are done. \qed

\vspace{.15in}

The hard part is now to determine $\mathbf{C}_{\mathbb{P}}, \mathbf{C}_{\mathbb{R}}$ and $\mathbf{C}_{\mathbb{I}}$. The following Proposition takes care of a little over two-thirds of the problem.

\vspace{.15in}

\begin{prop}[Zhang] We have $\mathbf{C}_{\mathbb{P}}=\H_{\mathbb{P}}$ and $\mathbf{C}_{\mathbb{I}}=\H_{\mathbb{I}}$. Moreover, $\mathbf{C}_i \subset \mathbf{C}_{\mathbb{R}}$ for $i=1, \ldots, d$.\end{prop}

\vspace{.15in}

To complete the description of $\mathbf{C}_{\mathbb{R}}$, we need to introduce a few more notations.
Recall from Theorem~\ref{T:Ringelstructure} that $\mathbb{R}$ decomposes as a direct sum of categories $\mathcal{C}_x \simeq Rep_{k_x}\vec{Q}_0$ and $\mathcal{C}_{\mathcal{O}_i} \simeq Rep_k A_{p_i-1}^{(1)}$, where $x$ ranges over all closed points of $\mathbb{P}^1\backslash D$, $\vec{Q}_0$ is the Jordan quiver, and $i=1, \ldots, d$. As a consequence,
$$\H_{\mathbb{R}}\simeq \bigotimes_{x \in \mathbb{P}^1 \backslash D} \H_{\mathcal{C}_x} \otimes \bigotimes_i \H_{\mathcal{C}_{\mathcal{O}_i}}.$$
For $x \in \mathbb{P}^1 \backslash D$ let $\Phi_x^{-1}: \LLambda \stackrel{\sim}{\to} \H_{\mathcal{C}_x}$ be the isomorphism constructed in Section~2.4. (with ground field equal to $k_x$), and for $i=1, \ldots, d$ let $\Phi_i^{-1}$ be the composition $\LLambda \stackrel{\sim}{\to}\H \stackrel{\sim}{\to} \mathbb{K}$ where $\mathbf{K} \subset \H_{\mathcal{C}_{\mathcal{O}_i}}$ is defined in Section~3.5 (strictly speaking, the definition of $\mathbf{K}$ presupposes the choice of a distinguished simple object in $\mathcal{C}_{\mathcal{O}_i}$, corresponding to the vertex labelled $1$; the results here do not depend on such a choice).

Now, for $r \geq 1$, put
$$T_{r,x}=\begin{cases} \frac{[r]}{r}deg(x) \phi_x^{-1}(p_{\frac{r}{deg(x)}})& \text{if}\; deg(x) | r,\\
0& \text{otherwise}\end{cases}$$
$$T_{r,i}=\frac{[r]}{r}\phi_i^{-1}(p_r),$$
and $T_r=\sum_x T_{r,x}+\sum_i T_{r,i}$. Note that $deg(T_r)=r\delta$. Let $\mathbf{K}' \subset \H_{\mathbb{R}}$ be the subalgebra generated by $T_1, T_2, \ldots$. It is clear that $\mathbf{K}'\simeq \C[T_1, T_2, \ldots]$.

\vspace{.15in}

\begin{prop}[Hubery] The following hold
\begin{enumerate}
\item[i)] The algebra $\mathbf{K}'$ belongs to $\mathbf{C}_{\mathbb{R}}$. Moreover the multiplication map induces an isomorphism of vector spaces $\mathbf{K}' \otimes \bigotimes_i \mathbf{C}_i \stackrel{\sim}{\to} \mathbf{C}_{\mathbb{R}}$.
\item[ii)] The center $\mathbf{Z}$ of $\mathbf{C}_{\mathbb{R}}$ is isomorphic to $\mathbf{K}'$ and the multiplication map gives an isomorphism $\mathbf{Z} \otimes \bigotimes_i \mathbf{C}_i \stackrel{\sim}{\to} \mathbf{C}_{\mathbb{R}}$.
\end{enumerate}
\end{prop}

\vspace{.15in}

\addtocounter{theo}{1}
\noindent \textbf{Example \thetheo.} Let $\vec{Q}$ be the Kronecker quiver

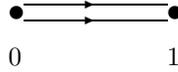
\vspace{.15in}

\centerline{
\begin{picture}(200, 10)
\put(60,0){\circle*{5}}
\put(120,0){\circle*{5}}
\put(85,-3){\vector(1,0){5}}
\put(85,3){\vector(1,0){5}}
\put(63,3){\line(1,0){54}}
\put(63,-3){\line(1,0){54}}
\put(57,-20){$0$}
\put(117,-20){$1$}
\end{picture}}

\vspace{.3in}

Then, as an easy computation shows,
$$\nu^2[S_1] \cdot [S_2]-[S_2]\cdot[S_1]=\sum_{x \in \mathbb{P}^1(k)} [I_{\delta}^x].$$
Thus we see that the average over all regular indecomposables of degree one belongs to $\mathbf{C}_{\mathbb{R}}$. This is simply $T_1$ in the notation above. The $T_r$'s for $r \geq 2$ may also be obtained as (more complicated) commutators. \endexample

\vspace{.15in}

\addtocounter{theo}{1}
\noindent \textbf{Remarks \thetheo .}  i) As was recently showed by Hubery in \cite{Hubery5}, Hall polynomials exist for any tame quiver.\\
ii)  The proof of Proposition~\ref{P:halldecompopo} may be applied \textit{verbatim} to obtain the following general statement~:

\vspace{.1in}

\begin{prop}\label{P:lasttame} Let $\mathcal{A}$ be a finitary hereditary category satisfying the finite subobject condition, and let $\mathcal{I}=\mathcal{I}_1 \sqcup \cdots \sqcup \mathcal{I}_r$ be a partition of the set of indecomposables of $\mathcal{A}$ satisfying $\text{Hom}(\mathcal{I}_i,\mathcal{I}_j)=\text{Ext}^1(\mathcal{I}_j,\mathcal{I}_i)=\{0\}$ if $i > j$.
Let $\mathbf{C}$ be a subalgebra of $\H_{\mathcal{A}}$ stable under the map $\Delta$. Then the multiplication induces an isomorphism $\mathbf{C}_1 \otimes \cdots \otimes \mathbf{C}_r \stackrel{\sim}{\to} \mathbf{C}$, where $\mathbf{C}_i=\mathbf{C} \cap \H_{\mathcal{A}_i}$ and $\mathcal{A}_i$ is the additive subcategory of $\mathcal{A}$ generated by $\mathcal{I}_i$.
\end{prop}

\newpage

\centerline{\large{\textbf{Lecture~4.}}}
\addcontentsline{toc}{section}{\tocsection {}{}{Lecture~4.}}

\setcounter{section}{4}
\setcounter{theo}{0}
\setcounter{equation}{0}

\vspace{.15in}

In this Lecture we momentarily forget about quivers to turn our attention to another very natural class of hereditary categories, namely the categories of coherent sheaves on smooth projective curves $X$.
The category $Coh(X)$ is never of finite type (i.e. there always are infinitely many indecomposables
objects) and is tame if and only if $X$ is of genus zero (i.e. $X \simeq \mathbb{P}^1$), or of genus one (i.e. $X$ is an elliptic curve). Inspired by some problems in number theory, Kapranov studied in the deep paper \cite{Kap1} the Hall algebra $\H_{X}:=\H_{Coh(X)}$ (or more precisely a certain ``composition subalgebra'' $\mathbf{C}_X \subset \H_{X}$) and observed some strong analogies with quantum \textit{loop} algebras. This becomes a precise theorem when $X=\mathbb{P}^1$, in which case the composition algebra $\mathbf{C}_X$ turns out to be a certain subalgebra of $\U_{\nu}(\widehat{\mathfrak{sl}}_2)$.

In \cite{SDuke}, Kapranov's Theorem was generalized to the case of Lenzing's \textit{weighted projective lines} $\mathbb{X}_{\underline{p},\underline{\lambda}}$ which may be viewed as noncommutative, or orbifold, projective lines, and which form a very natural and important class of hereditary categories. The Hall algebras are now related to quantum groups associated to \textit{loop algebras of Kac-Moody algebras}. The categories $Coh(\mathbb{X}_{\underline{p},\underline{\lambda}})$ themselves are still quite mysterious. In fact, the classification of (classes) of indecomposables (i.e. the analogue of Kac's Theorem~\ref{T:Kac}) was only proved recently by Crawley-Boevey \cite{CB2}.

\vspace{.1in}

After describing Kapranov's fundamental results in the case of $\mathbb{P}^1$ and the notion of a weighted projective line $\mathbb{X}_{\underline{p},\underline{\lambda}}$, we state in this Lecture the main Theorems of \cite{CB2} and \cite{SDuke} concerning the structure of $Coh(\mathbb{X}_{\underline{p},\underline{\lambda}})$ and its Hall algebra, paying special attention to the cases when $\mathbb{X}_{\underline{p},\underline{\lambda}}$ is of genus at most one. Motivated by this, we then describe, following the joint work of I. Burban and the author \cite{BS}, the Hall algebra of an elliptic curve. In the last section we compile what is known and expected of Hall algebras in higher genus. Appendices~A.5. and A.6. contain all the necessary properties of loop algebras and their quantum groups.

\vspace{.2in}

\centerline{\textbf{4.1. Generalities on coherent sheaves.}}
\addcontentsline{toc}{subsection}{\tocsubsection {}{}{\; 4.1. Generalities on coherent sheaves.}}

\vspace{.15in}

\paragraph{} For the notion of a coherent sheaf on 
an algebraic variety, we refer to \cite{Hart}. We just recall here a few basic facts; this will also be useful for comparing with weighted projective lines. Let $X$ be a smooth projective variety defined over a field $k$ which will be either $\C$ or a finite field.

A \textit{torsion sheaf} on $X$ is a sheaf whose support is a finite set of points. The notion of a \textit{locally free sheaf} on $X$ is equivalent to the notion of a vector bundle $E \to X$. As usual, $\mathcal{O}$ stands for the structure sheaf on $X$ (the trivial line bundle).
The category $Coh(X)$ of coherent sheaves on $X$, contrary to the category $Vec(X)$ of vector bundles, is abelian
(in fact, $Vec(X)$ is an exact subcategory of $Coh(X)$). By a theorem of Serre, when $X$ is smooth
we have ${gldim}(Coh(X))={dim}(X)$, and when $X$ is projective
$$\text{dim}(\text{Ext}^i(\mathcal{F},\mathcal{G})) < \infty,\qquad {\forall}\; \mathcal{F},\mathcal{G} \in Coh(X).$$
Hence $Coh(X)$ is a finitary hereditary category if and only if $X$ is a smooth projective curve defined over a finite field $k=\mathbb{F}_q$. Under this hypothesis \textit{which we henceforth make}, any sheaf $\mathcal{F}$ has a canonical (maximal) torsion subsheaf $\mathcal{T} \subset \mathcal{F}$ and canonical quotient vector bundle
$\mathcal{V}=\mathcal{F}/\mathcal{T}$. Moreover the exact sequence
$$\xymatrix{ 0 \ar[r] & \mathcal{T} \ar[r] & \mathcal{F} \ar[r]& \mathcal{V} \ar[r] & 0}$$
splits; that is, any sheaf can be decomposed as a direct sum $\mathcal{F}=\mathcal{V} \oplus \mathcal{T}$.

\vspace{.1in}

Let ${T}or(X)$ stand for the (abelian) full subcategory of $Coh(X)$ consisting of torsion sheaves. It decomposes as a direct product of blocks 
$${T}or(X)=\prod_{x \in X} {T}or_x$$
where $x$ ranges over the set of closed points of $X$ and ${T}or_x$ is the category of torsion sheaves supported at $x$. This last category is equivalent to the category of finite-dimensional modules over the local ring $\mathcal{O}^{(x)}$ at $x$ which, by the smoothness of $X$, is a discrete valuation ring. Thus, in the end, ${T}or_x$ is equivalent to the category $Rep^{nil}_{k_x}\vec{Q}_0$ of nilpotent representations of the Jordan quiver \textit{over the residue field} $k_x=\mathcal{O}^{(x)}/\mathfrak{m}^{(x)}$. In particular, there is a unique simple sheaf $\mathcal{O}_x$ supported at $x$.
  
\vspace{.1in}

The \textit{rank} of a coherent sheaf $\mathcal{F}$ is the rank of its canonical quotient vector bundle $\mathcal{V}$. The \textit{degree} of a sheaf is the only invariant satisfying $deg(\mathcal{O})=0$, $deg(\mathcal{O}_x)=deg(x)=[k_x:k]$ and which is additive on short exact sequences (i.e. which factors through the Grothendieck group $K(Coh(X))$).

\vspace{.1in}

To finish, let us state some important homological properties of $Coh(X)$. Let $\omega_X$ stand for the line bundle of differential forms. The exact functor $Coh(X) \to Coh(X)$,\; $\cdot \mapsto \cdot \otimes \omega_X$ is a \textit{Serre functor}\footnote{actually to be precise the Serre functor is the autoequivalence of the derived category given by $\cdot \mapsto \cdot \otimes \omega_X[1]$. Without the shift $[1]$, the functor is called the Auslander-Reiten translate.}
, that is there are natural isomorphisms
\begin{equation}\label{E:Serreduality}
{Ext}^1(\mathcal{F},\mathcal{G})^* \simeq {Hom}(\mathcal{G}, \mathcal{F}\otimes \omega_X).
\end{equation}
One consequence of this is that for any $\mathcal{T} \in {T}or(X)$ and $\mathcal{V} \in Vec(X)$ one has
\begin{equation}\label{E:homtorvec}
{Hom}(\mathcal{T},\mathcal{V}) ={Ext}^1(\mathcal{V},\mathcal{T})=\{0\}.
\end{equation}
In particular, ${Ext}^1(\mathcal{O}_x,\mathcal{O})\simeq k^{deg(x)}$ for any closed point $x$ of $X$, and there is a unique line bundle extension of $\mathcal{O}_x$ by $\mathcal{O}$, which is denoted by $\mathcal{O}(x)$. 

\vspace{.2in}

\centerline{\textbf{4.2. The category $Coh(\mathbb{P}^1)$.}}
\addcontentsline{toc}{subsection}{\tocsubsection {}{}{\; 4.2. The category of coherent sheaves over $\mathbb{P}^1$.}}

\vspace{.15in}

\paragraph{} Everything concerning $Coh(X)$ becomes very explicit when $X=\mathbb{P}^1$. 
Recall that the homogeneous coordinate ring of $\mathbb{P}^1$ is the $\Z$-graded ring $\mathbf{S}=k[X_1,X_2]$ where $deg(X_1)=deg(X_2)=1$. Hence a closed point $x$ of $\mathbb{P}^1$ 
is simply a maximal (homogeneous) ideal $\mathfrak{m}_x $ of $\mathbf{S}$. To any finitely generated graded $\mathbf{S}$-module $M=\bigoplus_i M_i$ is associated a coherent sheaf $M\;\widetilde{}$. The resulting functor $\widetilde{}\;:\mathbf{S}\mathrm{-}Modgr \to Coh(\mathbb{P}^1)$ is, however, not an equivalence. Let $\mathbf{S}\mathrm{-}Fin$ be the full subcategory of $\mathbf{S}\mathrm{-}Modgr$ consisting of finite-dimensional graded modules. This is a Serre subcategory, i.e. it is abelian, stable under extensions, subobjects and quotients.

\vspace{.1in}

\begin{theo}[Serre]\label{T:Serre} The functor $\widetilde{}$ induces an equivalence of categories
$$\mathbf{S}\mathrm{-}Modgr/\mathbf{S}\mathrm{-}Fin \stackrel{\sim}{\to} Coh(\mathbb{P}^1)$$
(where the quotient category is defined in the sense of Serre, see e.g. \cite{Serrequotientcategory}).
\end{theo}

\vspace{.1in}

In more plain terms, Serre's Theorem states that two graded $\mathbf{S}$-modules $M$ and $N$ give rise to isomorphic coherent sheaves if there is an isomorphism (of graded $\mathbf{S}$-modules) $M_{\geq n}:=\bigoplus_{i \geq n} M_i \simeq N_{\geq n}:=\bigoplus_{i \geq n} N_i$ for $n \gg 0$. Similarly, we have
\begin{equation}\label{E:defhomquotient}
{Hom}_{Coh(\mathbb{P}^1)}(M\;\widetilde{}, N\;\widetilde{}\;)=\underset{\longrightarrow} {{Lim}}\;{Hom}_{\mathbf{S}}(M_{\geq n},N).
\end{equation}

\vspace{.15in}

The first important result concerning the category $Coh(\mathbb{P}^1)$ is the following one~:

\vspace{.1in}

\begin{theo}[Grothendieck, \cite{Grot}]\label{P:GrothP1} Any indecomposable vector bundle is a line bundle.
Two line bundles are isomorphic if and only if they have the same degree.\end{theo}

\vspace{.1in}

For $d \in \Z$ we will denote by $\mathcal{O}(d)$ the (unique) line bundle of degree $d$. We have
$\mathcal{O}(d)=\mathbf{S}[d]\;\widetilde{}$, where $\mathbf{S}[d]$ is the $d$th shift of the trivial module $\mathbf{S}$, i.e. $\mathbf{S}[d]_i:=\mathbf{S}_{i+d}$. If $\mathcal{F}$ is any coherent sheaf we set $\mathcal{F}(d)=\mathcal{F} \otimes \mathcal{O}(d)$. An easy consequence of Grothendieck's Theorem is

\vspace{.1in}

\begin{cor} We have $K(Coh(\mathbb{P}^1))=\Z^2$. The class of a sheaf $\mathcal{F}$ is given by the pair $(rank(\mathcal{F}),deg(\mathcal{F}))$.\end{cor}

\noindent
\textit{Proof.} Let $x,y$ be any two closed points of degree one, the extensions
$$\xymatrix{ 0 \ar[r] & \mathcal{O} \ar[r] & \mathcal{O}(x) \ar[r]& \mathcal{O}_x \ar[r]& 0}$$ 
$$\xymatrix{ 0 \ar[r] & \mathcal{O} \ar[r] & \mathcal{O}(y) \ar[r]& \mathcal{O}_y \ar[r]& 0}$$ 
show that, in $K(Coh(\mathbb{P}^1))$, we have $\overline{\O}+\overline{\O_x} = \overline{\O(x)}=\overline{\O(1)}=\overline{\O(y)}=\overline{\O}+\overline{\O_y}$, so that 
$\overline{\O_y}=\overline{\O_x}$. The same argument shows that $\overline{\mathcal{O}_z}=deg(x)\overline{O_{x}}$, and finally that $\overline{\mathcal{T}}=deg(\mathcal{T}) \overline{O_x}$ for any torsion sheaf. But then
$\overline{\mathcal{F}}=rank(\mathcal{F}) \overline{\O} + deg(\mathcal{F})\overline{\O_x}$, and we are done.\qed

\vspace{.15in}

The Euler form on $K(Coh(\mathbb{P}^1))$ is now also easy to compute. Since
$${Hom}(\O(n),\O(m))=\begin{cases} \{0\} & \text{if}\; n > m\\
k^{m-n+1} & \text{if}\; n \leq m\end{cases}$$
and since $\omega_X=\O(2)$, we deduce (using (\ref{E:Serreduality}) that
$${Hom}(\O,\O)={Hom}(\O_x,\O_x)={Hom}(\O,\O_x)=k,\qquad {Hom}(\O_x,\O)=\{0\}$$
$${Ext}^1(\O,\O)={Ext}^1(\O,\O_x)=\{0\}, \qquad {Ext}^1(\O_x,\O)= {Ext}(\O_x,\O_x)=k.$$
Therefore
\begin{equation}\label{E:EulerP1}
\langle \mathcal{F},\mathcal{G} \rangle_a=rank(\mathcal{F})\big( rank(\mathcal{G}) + deg(\mathcal{G})\big) -deg(\mathcal{F})rank(\mathcal{G}),
\end{equation}
and $( \mathcal{F},\mathcal{G} )_a=2 rank(\mathcal{F}) rank(\mathcal{G})$. 

\vspace{.1in}

As first observed by Kapranov, we may summarize all this in the following manner. Let $\g=\widehat{\mathfrak{sl}}_2$ be the affine Lie algebra associated to $\mathfrak{sl}_2$, and let $\widehat{Q}=\Z \a \oplus \Z\delta$ be its root lattice (see Appendix~A.3.).

\vspace{.1in}

\begin{cor}\label{C:KacP1} The map 
\begin{align*}
\rho:\;K(Coh(\mathbb{P}^1)) &\to \widehat{Q}\\
\overline{\mathcal{F}} &\mapsto rank(\mathcal{F}) \a + deg(\mathcal{F})\delta
\end{align*}
is an isomorphism of $\Z$-modules. It sends the symmetrized Euler form $(\;,\;)_a$ to the Cartan-Killing form $(\;,\;)$ on $\widehat{Q}$. Under the identification by $\rho$, the set of classes of indecomposable sheaves is the $\mathrm{nonstandard}$ set of positive roots
\begin{equation}\label{E:nonstandardphi}
\Phi_+=\{\a + n\delta\;|\; n \in \Z\} \sqcup \{n \delta\;|\; n >0\}.
\end{equation}
Moreover, 
\begin{enumerate}
\item[i)] if $\beta \in \Phi_+$ is real then there exists a unique indecomposable sheaf $\mathcal{F}$ of class $\beta$,
\item[ii)] If $\beta \in \Phi_+$ is imaginary then there exists a $\mathbb{P}^1$-family of indecomposable sheaves of class $\beta$.
\end{enumerate}
\end{cor}

\vspace{.1in}

To finish with the presentations, let us draw a picture of $Coh(\mathbb{P}^1)$, similar to the ones appearing in Section~3.6. for tame quivers. We plot one point for each indecomposable and draw $\text{dim}(\mathcal{F},\mathcal{F}')$ arrows between neighboring indecomposables. We use the homological ordering, that is, morphisms go from left to right while extensions go from right to left.

\vspace{.2in}

\begin{equation}\label{E:picturecohp1}
\centerline{
\begin{picture}(300,60)
\put(5,0){\line(1,0){5}}
\put(15,0){\line(1,0){5}}
\put(5,25){\line(1,0){5}}
\put(15,25){\line(1,0){5}}
\put(5,50){\line(1,0){5}}
\put(15,50){\line(1,0){5}}
\put(190,25){\line(1,0){5}}
\put(200,25){\line(1,0){5}}
\put(190,0){\line(1,0){5}}
\put(200,0){\line(1,0){5}}
\put(190,50){\line(1,0){5}}
\put(200,50){\line(1,0){5}}
\put(106,6){\circle*{4}}
\put(120,45){\circle*{4}}
\put(150,45){\circle*{4}}
\put(136,6){\circle*{4}}
\put(90,45){\circle*{4}}
\put(60,45){\circle*{4}}
\put(76,6){\circle*{4}}
\put(105,10){\vector(1,3){10}}
\put(110,10){\vector(1,3){10}}
\put(125,40){\vector(1,-3){10}}
\put(120,40){\vector(1,-3){10}}
\put(155,40){\line(1,-3){7}}
\put(150,40){\line(1,-3){7}}
\put(135,10){\vector(1,3){10}}
\put(140,10){\vector(1,3){10}}
\put(75,10){\vector(1,3){10}}
\put(80,10){\vector(1,3){10}}
\put(90,40){\vector(1,-3){10}}
\put(95,40){\vector(1,-3){10}}
\put(48,19){\vector(1,3){7}}
\put(53,19){\vector(1,3){7}}
\put(65,40){\vector(1,-3){10}}
\put(60,40){\vector(1,-3){10}}
\put(80,50){$\mathcal{O}(-1)$}
\put(60,-5){$\mathcal{O}(-2)$}
\put(45,50){$\mathcal{O}(-3)$}
\put(115,50){$\mathcal{O}(1)$}
\put(130,-5){$\small{\mathcal{O}(2)}$}
\put(145,50){$\mathcal{O}(3)$}
\put(102,-5){$\mathcal{O}$}
\put(240,42){\line(0,1){2}}
\put(240,46){\line(0,1){2}}
\put(240,50){\line(0,1){2}}
\put(300,42){\line(0,1){2}}
\put(300,46){\line(0,1){2}}
\put(300,50){\line(0,1){2}}
\put(250,39){\line(0,1){2}}
\put(250,43){\line(0,1){2}}
\put(250,47){\line(0,1){2}}
\put(290,39){\line(0,1){2}}
\put(290,43){\line(0,1){2}}
\put(290,47){\line(0,1){2}}
\put(240,0){\vector(0,1){20}}
\put(240,0){\circle*{4}}
\put(240,40){\circle*{4}}
\put(240,20){\vector(0,1){20}}
\put(240,20){\circle*{4}}
\put(230,-14){$\mathcal{O}_x$}
\put(250,-15){$\cdots$}
\put(260,15){$\cdots$}
\put(260,25){$\cdots$}
\put(270,-15){$\cdots$}
\put(300,0){\vector(0,1){20}}
\put(300,0){\circle*{4}}
\put(300,40){\circle*{4}}
\put(300,20){\vector(0,1){20}}
\put(300,20){\circle*{4}}
\put(290,-14){$\mathcal{O}_{x'}$}
\put(250,-3){\vector(0,1){20}}
\put(250,-3){\circle*{4}}
\put(250,37){\circle*{4}}
\put(250,17){\vector(0,1){20}}
\put(250,17){\circle*{4}}
\put(290,-3){\vector(0,1){20}}
\put(290,-3){\circle*{4}}
\put(290,37){\circle*{4}}
\put(290,17){\vector(0,1){20}}
\put(290,17){\circle*{4}}
\put(240,0){\line(0,1){40}}
\put(300,0){\line(0,1){40}}
\qbezier(240,0)(270,-10)(300,0)
\qbezier(240,40)(270,30)(300,40)
\qbezier(240,0)(270,10)(300,0)
\qbezier(240,40)(270,50)(300,40)
\end{picture}}
\end{equation}

\vspace{.3in}

For clarity, only the subcategories $Tor_x$ with $deg(x)=1$ have been drawn.

\vspace{.2in}

\centerline{\textbf{4.3. The Hall algebra of $Coh(\mathbb{P}^1)$.}}
\addcontentsline{toc}{subsection}{\tocsubsection {}{}{\; 4.3. The Hall algebra of $\mathbb{P}^1$.}}

\vspace{.15in}

\paragraph{}The ressemblance between Corollary~\ref{C:KacP1} and Kac's Theorem~\ref{T:Kac} will certainly not have escaped the reader's sagacity. At this point it seems natural to expect that the Hall algebra $\H_{\mathbb{P}^1}$ is closely related to the quantum group $\U_v(\widehat{\mathfrak{sl}}_2)$. This is the content of Kapranov's Theorem, which we will soon state. Before this, let us work as usual through some sample computations of Hall numbers. Set $k=\mathbb{F}_q$ and put $\nu=q^{\frac{1}{2}}$. The Hall algebra of $Coh(\mathbb{P}^1)$ over $k$ will simply be denoted $\H_{\mathbb{P}^1}$.

\vspace{.15in}

\addtocounter{theo}{1}
\paragraph{\textbf{Example~\thetheo.}} Since $Tor({\mathbb{P}^1})$ is an abelian subcategory of $Coh(\mathbb{P}^1)$ which is closed under extensions, there is an natural embedding of Hall algebras $\H_{Tor({\mathbb{P}^1})} \hookrightarrow \H_{\mathbb{P}^1}$ (see Corollary~\ref{C:cathall}). Moreover, because $Tor({\mathbb{P}^1)}$ decomposes as a direct product $\prod_{x \in \mathbb{P}^1} Tor_x$, there is an isomorphism $\H_{Tor({\mathbb{P}^1})} \simeq \bigotimes_{x \in \mathbb{P}^1} \H_{Tor_x}$. Finally, each $\H_{Tor_x}$ is isomorphic to a classical Hall algebra (with respect to the ground field $k_x$).

\endexample

\vspace{.15in}

\addtocounter{theo}{1}
\paragraph{\textbf{Example~\thetheo.}} Assume that $\mathcal{F}, \mathcal{G}$ are two sheaves such that ${Hom}(\mathcal{G},\mathcal{F})=\{0\}$. Then by Serre duality
$${Ext}^1(\mathcal{F},\mathcal{G})^* \simeq {Hom}(\mathcal{G}(2),\mathcal{F}) \simeq {Hom}(\mathcal{G}, \mathcal{F}(-2))=\{0\}$$
since there is a canonical embedding $\mathcal{F}(-2) \hookrightarrow \mathcal{F}$. We deduce that
$$[\mathcal{F}] \cdot [\mathcal{G}]=\nu^{{dim\;Hom}(\mathcal{F},\mathcal{G})}[\mathcal{F} \oplus \mathcal{G}]$$
(indeed, any extension is trivial since ${Ext}^1(\mathcal{F},\mathcal{G})=\{0\}$ and there is a unique subsheaf of $\mathcal{F}\oplus \mathcal{G}$ isomorphic to $\mathcal{G}$ since ${Hom}(\mathcal{G},\mathcal{F})=\{0\}$). We may apply this to the following situations~: we have ${Hom}(\O(n),\O(m))=\{0\}$ if $n > m$, and ${Hom}(\mathcal{T},\mathcal{V})=\{0\}$ if $\mathcal{T}$ is a torsion sheaf and $\mathcal{V}$ is a vector bundle. Hence
\begin{equation}\label{E:productofos}
[\O(n_1)]  \cdots [\O(n_r)]=\nu^{\sum_{i < j} {dim\;Hom}(\O(n_i),\O(n_j))}[\mathcal{O}(n_1) \oplus \cdots \oplus \mathcal{O}(n_r)]
\end{equation}
if $n_1 < n_2 \cdots < n_r$; and
\begin{equation}\label{E:productvecttor}
[\mathcal{V}] \cdot [\mathcal{T}]=\nu^{{dim\;Hom}(\mathcal{V},\mathcal{T})}[\mathcal{V}\oplus \mathcal{T}]=\nu^{rd}[\mathcal{V}\oplus\mathcal{T}]
\end{equation}
if $\mathcal{V}$ is a vector bundle of rank $r$ and $\mathcal{T}$ is a torsion sheaf of degree $d$.

\endexample

\vspace{.15in}

\addtocounter{theo}{1}
\paragraph{\textbf{Example~\thetheo.}} Let $d \in \N$ and let us set
$$\mathbf{1}_{d\delta}=\sum_{\overline{\mathcal{T}}=(0,1)} [\mathcal{T}]$$
where the sum ranges over all torsion sheaves of degree $d$. To compute the product $\mathbf{1}_{d\delta}\cdot [\O(n)]$ we first note that any extension between a torsion sheaf of degree $d$ and $\O(n)$ is of the form $\O(n+s) \oplus \mathcal{T}'$ for some torsion sheaf $\mathcal{T}'$ of degree $d-s$. Next, a map $\phi: \mathcal{O}(n) \to \mathcal{O}(n+d-s) \oplus \mathcal{T}'$ is injective if and only if ${Im}\;\phi \not\subset \mathcal{T}'$, and such a map has a cokernel which is a torsion sheaf of degree $d$. There are, up to a scalar, $(q^{(s+1)+(d-s)}-q^{d-s})/(q-1)=q^{d-s}[s+1]_+$ such maps. Thus, all in all, we obtain
\begin{equation}\label{E:product1rdo}
\begin{split}
\mathbf{1}_{d\delta}\cdot [\O(n)]&=\nu^{-d} \sum_{s=0}^d \sum_{\overline{\mathcal{T}'}=(0,d-s)}\nu^{2(d-s)}[s+1]_+[\mathcal{O}(n+s)\oplus \mathcal{T}']\\
&=\sum_{s=0}^d \nu^{d-s} [s+1]\sum_{\overline{\mathcal{T}'}=(0,d-s)}[\mathcal{O}(n+s)\oplus \mathcal{T}']\\
&=\sum_{s=0}^d [s+1] [\mathcal{O}(n+s)] \cdot \mathbf{1}_{(d-s)\delta}
\end{split}
\end{equation}
\endexample

\vspace{.15in}

\addtocounter{theo}{1}
\paragraph{\textbf{Example~\thetheo.}} As a final example, let us compute the product $[\O(n)]\cdot[\O(m)]$ when $n > m$. Since any extension of $\O(m)$ by $\O(n)$ is of the form $\O(m+s) \oplus \O(n-s)$ for some $0 \leq s \leq (n-m)/2$, we have
$$[\O(n)] \cdot [\O(m)]=\nu^{m+1-n}\sum_{s=0}^{\lfloor \frac{n-m}{2}\rfloor} P_{\O(n),\O(m)}^{\O(m+s) \oplus\O(n-s)} [\O(m+s) \oplus \O(n-s)].$$
In order to calculate the Hall number $P_{\O(n),\O(m)}^{\O(m+s) \oplus\O(n-s)}$ we need only count the number of nonzero maps $\phi: \mathcal{O}(m) \to \O(n-s) \oplus \O(m+s)$ whose cokernel is locally free; indeed any such map is injective and by Grothendieck's Theorem~\ref{P:GrothP1} any line bundle of degree $n$ is isomorphic to $\O(n)$. Let us write $\phi=\phi_1\oplus\phi_2$ with
$$\phi_1\in {Hom}(\O(m),\O(m+s)) \simeq \underset{\longrightarrow}{{Lim}} \;{Hom}_{\mathbf{S}}(\mathbf{S}[m]_{\geq l}, \mathbf{S}[m+s]) \simeq k[X_1,X_2]_s,$$
$$\phi_2\in {Hom}(\O(m),\O(n-s)) \simeq \underset{\longrightarrow}{{Lim}} \;{Hom}_{\mathbf{S}}(\mathbf{S}[m]_{\geq l}, \mathbf{S}[n-s]) \simeq k[X_1,X_2]_{n-s-m}.$$
Thus we may view the maps $\phi_1, \phi_2$ as homogeneous polynomials in $k[X_1,X_2]$ of degrees $s$ and $n-s-m$ respectively. For ${Im}\;\phi$ to be a vector subbundle, i.e. for ${Coker}\;\phi$ to be a vector bundle it is necessary and sufficient that $\phi_1$ and $\phi_2$ be relatively prime. Indeed, the support of the torsion sheaf ${Coker}\;\phi_i$ for $i=1, 2$ is the set of homogeneous prime factors of $\phi_i$ and ${Im}\;\phi$ is a subbundle if and only if these two supports do not intersect.
\vspace{.05in}

\textbf{Claim} (\cite{BK}). Let $a_{t,w}$ stand for the number of pairs of coprime homogeneous polynomials $(P,Q)$ in $k[X_1,X_2]$ of respective degrees $t$ and $w$. Then
$$a_{t,w}=\begin{cases} (q-1)(q^{t+w+1}-1)&\text{if}\;t=0\;\text{or}\;w=0,\\
(q-1)(q^2-1)q^{t+w-1}&\text{if}\;t\geq 1\;\text{and}\;w \geq 1.
\end{cases}$$
\noindent
\textit{Proof.} Let $b_{w,t}$ stand for the number of pairs of homogeneous polynomials $(R,T)$ in $k[X_1,X_2]$ of degrees $t$ and $w$. We have $b_{t,w}=(q^{t+1}-1)(q^{w+1}-1)$ while on the other hand we may write any such pair as $(R,T)=(DP,DQ)$ where $D$ is a polynomial of degree, say, $s$, and $(P,Q)$ are coprime of degrees $t-s$ and $w-s$. As $D$ is well-defined up to a scalar, we deduce that
$$b_{w,t}=\sum_{s \leq \text{min}(t,w)} \frac{q^{s+1}-1}{q-1} a_{t-s,w-s}.$$
The claim follows by an induction on $\text{min}(t,w)$.\qed

\vspace{.05in}

Thus, putting everything together we get
\begin{equation}\label{E:prodoo}
[\O(n)]\cdot[\O(m)]=\nu^{n-m+3}[\O(m) \oplus \O(n)] + \sum_{s=1}^{\lfloor \frac{n-m}{2}\rfloor}\nu^{n-m-1}(\nu^2-1)[\O(m+s) \oplus \O(n-s)].
\end{equation}
\endexample

\vspace{.2in}

We now turn to Kapranov's Theorem. Let us regard the Lie algebra $\widehat{\mathfrak{sl}}_2$ not as a Kac-Moody algebra but rather as the loop algebra $\mathcal{L}\mathfrak{sl}_2$ of $\mathfrak{sl}_2$ (see Appendix~A.5. for the definition of the loop algebra of a Kac-Moody algebra). Let $\mathcal{L}\bo_+ \subset \mathcal{L}\mathfrak{sl}_2$ be the positive Borel subalgebra associated to the set of nonstandard weights $\Phi_+$ defined in (\ref{E:nonstandardphi}). 
Let $\U_v(\mathcal{L}{\mathfrak{sl}}_2)$ be the quantum loop algebra of  ${\mathfrak{sl}}_2$, and let
$\U_v(\mathcal{L}\bo_+)$ be the positive Borel subalgebra of $\U_v(\mathcal{L}\mathfrak{sl}_2)$
(see Appendix~A.6.). Thus $\U_v(\mathcal{L}\bo_+)$ is generated by elements $E_l$ for $l \in \Z$, $H_n,$ for $n \in \Z^*$ and $K^{\pm 1}, C^{\pm 1/2}$ modulo the relations
\begin{equation}\label{E:Eee1}
C^{\pm 1/2} \;\text{is\;central}
\end{equation}
\begin{equation}\label{E:Eee2}
[K,H_{n}]=[H_n,H_l]=0,
\end{equation}
\begin{equation}\label{E:Eee3}
K E_{k}K^{-1}=v^{2}E_{k},
\end{equation}
\begin{equation}\label{E:Eee4}
[H_{l},E_{k}]= \frac{1}{l}[2l]C^{-|l|/2}E_{k+l}, 
\end{equation}
\begin{equation}\label{E:Eee5}
E_{k+1}E_{l}-v^{2}E_{l}E_{k+1}=
v^{2}E_{k}E_{l+1}-E_{l+1}E_{k}.
\end{equation}

\vspace{.15in}

 We let $\U_{\nu}(\mathcal{L}\bo_+)$ be the specialization at $v=\nu$ of $\U_v(\mathcal{L}\bo_+)$. Recall that for all closed points $x \in \mathbb{P}^1$ there is an isomorphism $\Phi_x^{-1}: \LLambda \stackrel{\sim}{\to} \H_{Tor_x}$ where $\LLambda$ is Macdonald's ring of symmetric functions (see Lecture~2). Here, if $x$ is of degree $deg(x)>1$, we view $\H_{Tor_x}$ as the classical Hall algebra \textit{over the ground field} $k=k_x$.
For $r \geq 1$, put
$$T_{r,x}=\begin{cases} \frac{[r]}{r}deg(x) \phi_x^{-1}(p_{\frac{r}{deg(x)}})& \text{if}\; deg(x) | r,\\
0& \text{otherwise}\end{cases}$$
and $T_r=\sum_x T_{r,x}$. Note that $deg(T_r)=r\delta$. 

\vspace{.1in}

The following Theorem was stated (and a sketch of a proof was given) in \cite{Kap1}. A detailed proof first appeared in \cite{BK}. Let us slightly extend the Hall
algebra $\widetilde{\H}_{\mathbb{P}^1}$ by adding a square root $\mathbf{k}_{\delta}^{1/2}$ of $\mathbf{k}_{\delta}$ and let us  denote by $\widetilde{\H}'_{\mathbb{P}^1}=\widetilde{\H}_{\mathbb{P}^1} \otimes_{\C[\mathbf{k}_{\delta}^{\pm 1}]} \C[\mathbf{k}_{\delta}^{\pm 1/2}]$ the resulting bialgebra.

\vspace{.1in}

\begin{theo}[Kapranov, \cite{Kap1}]\label{T:Kap} The assignement $E_l \mapsto [\O(l)]$ for $l \in \Z$, $H_r \mapsto T_r\mathbf{k}_{\delta}^{-|r|/2}$ for $r \geq 1$, $K \mapsto \mathbf{k}_{\O}$ and $C^{1/2} \mapsto \mathbf{k}_{\delta}^{1/2}$ extends to an embedding of algebras $$\Psi:\U_{\nu}(\mathcal{L}\bo_+) \to \widetilde{\H}'_{\mathbb{P}^1}.$$
\end{theo}
\noindent
\textit{Proof.} Relation~(\ref{E:Eee1}) is obvious since $\delta$ is an isotropic root (i.e. $(\delta,\a)_a=0$ for any weight $\a$). Relation~{(\ref{E:Eee2}) comes from the fact that $H_n$ is of weight $n\delta$, that the classical Hall algebra is commutative and that of course the Hall algebras $\H_{Tor_x}$ and $\H_{Tor_y}$ associated to distinct points $x$ and $y$ commute. Relation~(\ref{E:Eee3}) follows from the definitions. Thus (\ref{E:Eee4}) and (\ref{E:Eee5}) are the only relations for which one has to work. We begin with the following Lemma~:

\vspace{.1in}

\begin{lem}\label{L:relation4} We have $1+\sum_{r \geq 1} \mathbf{1}_{r\delta}s^r=exp(\sum_{r \geq 1} \frac{T_r}{[r]}s^r)$.
\end{lem}
\noindent
\textit{Proof.} Since $\H_{Tor_{\mathbb{P}^1}}$ is commutative, we have
$$exp(\sum_{r\geq 1}\frac{T_r}{[r]}s^r)=\prod_{x \in \mathbb{P}^1}exp(\sum_{r\geq 1}\frac{T_{r,x}}{[r]}s^r).$$
Since on the other hand $1 + \sum_{r\geq 1} \mathbf{1}_{r\delta}s^r=\prod_{x \in \mathbb{P}^1} (1 + \sum_{r\geq 1}\mathbf{1}_{r\delta,x}s^r)$, where
$$\mathbf{1}_{r\delta,x}=\sum_{\underset{\overline{\mathcal{F}}=(0,r)}{\mathcal{F}\in Tor_x}} [\mathcal{F}],$$
it suffices to show that $exp(\sum_{r\geq 1} \frac{T_{r,x}}{[r]}s^r)=1 + \sum_{r \geq 1} \mathbf{1}_{r\delta,x}s^r$. This is in turn a consequence of the following identity in the ring of symmetric functions $\LLambda$~:
$$1+\sum_{r \geq 1} h_r s^r=exp(\sum_{r \geq 1}\frac{p_r}{r}s^r),$$
where $h_r$ is the complete symmetric function.\qed

\vspace{.1in}

With a little algebraic manipulation, we may rewite the formula (\ref{E:product1rdo}) as
\begin{equation}\label{E:arg1}
\mathbf{1}(s)\mathbf{O}(t)=\mathbf{O}(t)\mathbf{1}(s)\frac{1}{(1-\frac{s}{t\nu})(1-\frac{s\nu}{t})}
\end{equation}
where $\mathbf{1}(s)=1+\sum_{r \geq 1}\mathbf{1}_{r\delta}s^r$ and $\mathbf{O}(t)=\sum_{r \geq 0} [\O(r)]t^r$. We have, by Lemma~\ref{L:relation4}, $log(\mathbf{1}(s))=\sum_{r \geq 1} \frac{T_r}{[r]}s^r=: \mathbf{T}(s)$. Thus, by (\ref{E:arg1}),
\begin{equation}\label{E:arg2}
[\mathbf{T}(s),\mathbf{O}(t)]=-\mathbf{O}(t)log((1-\frac{s}{t\nu})(1-\frac{s\nu}{t})),
\end{equation}
which with a little effort may be rewritten as 
$$[\frac{T_r}{[r]},[\O(n)]]=\frac{\nu^r+\nu^{-r}}{r}[\O(n+r)],$$
or
$$[T_r,[\O(n)]]=\frac{[2r]}{r}[\O(n+r)].$$
Relation~(\ref{E:Eee4}) follows. The last relation (\ref{E:Eee5}) is a consequence of formula (\ref{E:prodoo}). 

We have proved that there indeed exists an algebra homomorphism 
$$\Psi:\U_{\nu}(\mathcal{L}\bo_+) \to \widetilde{\H}'_{\mathbb{P}^1}.$$
To show that this map is injective, we will use a graded dimension argument.
Let $\U_{\nu}(\mathfrak{c})$ be the subalgebra of $\U_{\nu}(\mathcal{L}\bo_+)$ generated by $E_l, l \geq 0$ and $H_n, n \geq 1$. Thus $\U_{\nu}(\mathfrak{c})$ is a deformation of the enveloping algebra $\U(\mathfrak{c})$ of the subalgebra $\mathfrak{c}:=(\mathfrak{n} \otimes \C[t]) \oplus (\h \otimes t \C[t]) \subset \mathcal{L}\bo_+$. Note that $\U(\mathfrak{c})$ has finite-dimensional weight spaces. By the PBW Theorem, we have
\begin{equation}\label{E:dimu}
\begin{split}
{dim}&\;\U(\mathfrak{c})[k\a+l\delta]\\
&=\#\left\{(n_0, n_1, \ldots ), (m_1, m_2, \ldots)\;|\; n_i, m_i \in \N, \; \sum n_i=k,\; \sum (i n_i +im_i)=l\right\}.
\end{split}
\end{equation}
The same holds for $\U_{\nu}(\mathfrak{c})$. By construction, $\Psi(\U_{\nu}(\mathfrak{c}))$ contains the elements $[\O(l)]$ for $l \geq 0$ and the elements $\mathbf{1}_{\n\delta}$ for $n \geq 1$. Hence it also contains all the (finite) products
$$[\O(0)]^{n_0} \cdot [\O(1)]^{n_1} \cdots \cdot \mathbf{1}_{\delta}^{m_1} \cdot \mathbf{1}_{2\delta}^{m_2} \cdots.$$
By Example~4.6., all these products are linearly independent, and of weight $\lambda=(\sum n_i)\a+(\sum in_i+im_i)\delta$. Comparing with (\ref{E:dimu}), we deduce that 
$${dim\;}\Psi(\U_{\nu}(\mathfrak{c}))[k\a+l\delta] \geq {dim}\;\U_{\nu}(\mathfrak{c})[k\a+l\delta],$$
whence the two dimensions are equal and the restriction of $\Psi$ to $\U_{\nu}(\mathfrak{c})$ is injective. To finish the proof, we use the following fact. There is an automorphism $\tau$ of $\U_{\nu}(\mathcal{L}\bo_+)$ such that $\tau(E_{l})=E_{l+1}$, $\tau(H_n)=H_{n+1}$ and $\tau(C)=C, \tau(K)=K$. There is a similar automorphism $\tau'$ of $\widetilde{\H}'_{\mathbb{P}^1}$ induced by the self-equivalence of categories $\mathcal{F} \mapsto \mathcal{F}(1)$. It is easy to check from the definition of $\Psi$ that $\Psi \circ \tau=\tau' \circ \Psi$. Now, any element $x$ of $\U_{\nu}(\mathcal{L}\bo_+)$ may be written as $x=y K^aC^b$, with $y$ belonging to the subalgebra $\U_{\nu}(\mathcal{L}\n_+)$ generated by $E_l$ for $l \in \Z$ and $H_n$ for $n \in \N$. In particular, there exists $m \gg 0$ such that $\tau^m(y) \in U_{\nu}(\mathfrak{c})$. But then $(\tau')^m\Psi(x)=\Psi(\tau^m(x))=\Psi(\tau^m(y))\mathbf{k}_{a\a+b\delta} \neq 0$. Thus $\Psi$ is injective, and Theorem~\ref{T:Kap} is proved.\qed
 
\vspace{.15in}
 
As $Coh(\mathbb{P}^1)$ is hereditary, its Hall algebra $\widetilde{\H}_{\mathbb{P}^1}$ has a bialgebra structure. However, because the finite subobjects condition~(\ref{E:FS}) obviously fails, this is only a topological bialgebra structure. Let us see what this coproduct looks like.

\vspace{.15in}

\addtocounter{theo}{1}
\paragraph{\textbf{Example~\thetheo.}} The abelian subcategory $Tor({\mathbb{P}^1})$ is stable under taking subobjects and quotients. Thus, by Corollary~\ref{C:cathall}, $\widetilde{\H}_{Tor({\mathbb{P}^1})}$ is a subbialgebra of $\widetilde{\H}_{\mathbb{P}^1}$. Moreover, the isomorphisms $\H_{Tor({\mathbb{P}^1})}\simeq
\bigotimes_{x \in \mathbb{P}^1} \H_{Tor_x}$ and $\Phi_x: \H_{Tor_x} \simeq \LLambda$ are compatible with the map $\Delta$. This completely determines the coproduct for $\H_{Tor(\mathbb{P}^1)}$. As an example (which is also a corollary of Lemma~\ref{L:coprodun}), we have
\begin{equation}\label{E:coprod1rd}
{\Delta}(\mathbf{1}_{r\delta})=\sum_{s+t=r} \mathbf{1}_{s\delta}\mathbf{k}_{t\delta}\otimes \mathbf{1}_{t\delta}.
\end{equation}

\endexample

\vspace{.15in}

\addtocounter{theo}{1}
\paragraph{\textbf{Example \thetheo.}} Let us now compute $\Delta([\O])$. As any subsheaf of $\O$ is a line bundle, it is clear that 
$\Delta_{(1,-l),(0,l)}([\O])=0$ if $l \geq 1$ and $\Delta_{(1,0),0}=[\O]$. The only difficult part is to calculate the component of $\Delta([\O])$ of weight $((0,l),(1,-l))$.
First observe that any quotient of $\O$ which is not equal to $\O$ is of the form
$$\O_{x_1}^{(n_1)} \oplus \cdots \oplus \O_{x_r}^{(n_r)}$$
where $x_i$ are \textit{distinct} closed points of $\mathbb{P}^1$ and $\O_x^{(n)}$ stands for the (unique) indecomposable torsion sheaf supported at $x$ and of length $n$. This comes from the fact that there are no surjective maps $\O \twoheadrightarrow \O_x^{(n_1)} \oplus \O_x^{(n_2)}$ for any $x$, because there are no surjective maps $\O \twoheadrightarrow \O_x \oplus \O_x$. Next, we claim that a morphism
$$\phi=\bigoplus_i \phi_i: \O \to \O_{x_1}^{(n_1)} \oplus \cdots \oplus \O_{x_r}^{(n_r)}$$
is surjective if and only if each of the maps $\phi_i: \O \to \O_{x_i}^{(n_i)}$ is. Indeed, we may assume that no $x_i$ is equal to $\infty$, and then restrict everything to $\mathbb{P}^1\backslash \{\infty\}=\mathbb{A}^1$. The statement is then easily seen to be a consequence of the Chinese Remainder Theorem. Finally, if $\phi: \O \to \bigoplus_i \O_{x_i}^{(n_i)}$ is surjective then by Grothendieck's Theorem, $Ker(\phi)$ is necessarily isomorphic to $\O(-\sum_i n_i deg(x_i))$. Since there are $ |Hom(\O,\O_x^{(n_i)})|-|Hom(\O,\O_{x_i}^{(n_i-1)})|=\nu^{2n_ideg(x_i)}(1-\nu^{-2deg(x_i)})$ choices for the map $\phi_i$, we see that all together,
\begin{equation}
\Delta([\O])= [\O]\otimes 1 +\sum_{l \geq 0} \nu^{l} u_l\mathbf{k}_{\O(-l)} \otimes [\O(-l)]
\end{equation}
with
\begin{equation}\label{E:formulathetal}
u_l=\sum_{x_i,n_i} \prod(1-\nu^{-2deg(x_i)}) [\O_{x_i}^{(n_i)}]
\end{equation}
where the sum ranges over the set of tuples of distinct points $x_1, \ldots, x_r$ and multiplicities $n_1, \ldots, n_r$ such that $\sum_i n_i deg(x_i)=l$.

\vspace{.1in}

Just for the pleasure, we give another, more formal, derivation of the coproduct of $[\O]$, which has the advantage of showing directly that $\Delta([\O]) \in \U_{\nu}(\mathcal{L}\bo_+) \widehat{\otimes} \U_{\nu}(\mathcal{L}\bo_+)$.
For this we set, for any class $(r,d)$ in the Grothendieck group, $\mathbf{1}_{(r,d)}=\sum_{\overline{\mathcal{F}}=(r,d)} [\mathcal{F}].$
It is easy to see that
\begin{equation}\label{E:un10}
\mathbf{1}_{(1,s)}=\sum_{l \geq 0} \nu^{-l} [\O(s-l)]\mathbf{1}_{l\delta}.
\end{equation}

Inverting (\ref{E:un10}) we obtain $[\O]=\sum_{n \geq 0} \nu^{-n} \mathbf{1}_{(1,-n)}\chi_n$, where
$$\chi_n=\sum_{r >0}(-1)^r \sum_{l_1+\cdots + l_r=n} \mathbf{1}_{(0,l_1)} \cdots \mathbf{1}_{(0,l_r)}.$$ 
Recall that $\mathbf{1}(s)=1+\sum_{l \geq 1} \mathbf{1}_{(0,l)}s^l$, and set $\chi(s)=1+\sum_{l \geq 1} \chi_l s^l$. It is easy to see  that the elements  $\{\chi_n\}$ are completely determined by the relations  
$\sum_{i+j=l} \mathbf{1}_{(0,i)} \chi_j=\delta_{l,0}$,  which can  be rewritten in the form  $\mathbf{1}(s) \chi(s)=1$. Next, by Lemma~\ref{L:coprodun}, 
\begin{equation}\label{E:un101}
\Delta(\mathbf{1}_{(1,0)})=\sum_{\a+\beta=(1,0)} \nu^{-\langle \a,\beta\rangle_a}\mathbf{1}_{\a}\mathbf{k}_{\beta} \otimes \mathbf{1}_{\beta}.
\end{equation}
As a consequence, $\Delta(\mathbf{1}(s))=\mathbf{1}(s) \otimes \mathbf{1}(s)$ and hence
$\Delta(\chi(s))=\chi(s) \otimes \chi(s)$, i.e
$\Delta(\chi_{n})=\sum_{k=0}^n\chi_k \otimes \chi_{n-k}$. This implies that
\begin{equation*}
\begin{split}
\Delta_{(0,l), (1,-l)}&([\O])\\
=&\sum_{k \geq 0}(\nu^{l-k} \mathbf{1}_{(0,l)} + 
\nu^{l-2-k} \mathbf{1}_{(0,l-1)}\chi_1 + \cdots + \nu^{-k-l} \chi_l)\mathbf{k}_{\O(-l)} \otimes \mathbf{1}_{(1,-l-k)} 
\chi_{k}.
\end{split}
\end{equation*}
Hence, setting $\theta_l=\sum_{k=0}^l \nu^{l-2k} \mathbf{1}_{(0,l-k)}\chi_k$ we obtain
$$\Delta([\O]) =[\O] \otimes 1 + \sum_{l \geq 0} \theta_l\mathbf{k}_{\O(-l)} \otimes [\O(-l)].$$
The above definition of $\theta_l$ is not very plesant or useful.
We claim that these can also be characterized through
the relation $\sum_l \theta_l s^l=
exp((\nu-\nu^{-1})\sum_{r \geq 1} T_{(0,r)}s^r)$. To see this, note that by Lemma~\ref{L:relation4} it holds
$1+\sum_{l \geq 1} \mathbf{1}_{l\delta}s^l =exp(\sum_{r \geq 1} \frac{T_{(0,r)}}{[r]} s^r)$, hence
$1+\sum_{l \geq 1} \chi_ls^l=exp(-\sum_{r \geq 1} \frac{T_{(0,r)}}{[r]} s^r)$. But then
\begin{equation*}\begin{split}
\sum_{l \geq 0}  \theta_l s^l=&\mathbf{1}(\nu s)\chi(\nu^{-1}s) = exp(\sum_{r \geq 1} 
\nu^{r} \frac{T_{(0,r)}}{[r]} s^r -\sum_{r \geq 1} \nu^{-r} \frac{T_{(0,r)}}{[r]} s^r)\\
= & exp((\nu-\nu^{-1})\sum_{r \geq 1} T_{(0,r)}s^r)
\end{split}
\end{equation*}
as desired. Of course, as a corollary, we obtain that $\theta_l=\nu^{l}u_l$ is given by the formula (\ref{E:formulathetal}).

\endexample

\vspace{.15in}

\begin{prop} Let us equip $\U_{\nu}(\mathcal{L}\bo_+)$ with Drinfeld's new comultiplication (see Appendix~A.6). Then the morphism $\Psi: \U_{\nu}(\mathcal{L}\bo_+) \to \widetilde{\H}'_{\mathbb{P}^1}$ is a morphism of bialgebras.\end{prop}

\noindent
\textit{Proof.} It suffices to compare the definitions of Drinfeld's new coproduct and the computations made in Examples~4.11. and 4.12..\qed

\vspace{.15in}

The final piece of structure which we may investigate here is the antipode map $S$. As $Coh(\mathbb{P}^1)$ does not satisfy the finite subobjects condition, $S$ is \textit{not} well-defined. However, we have

\vspace{.1in}

\begin{lem} For any two sheaves $\mathcal{F},\mathcal{G}$ in $Coh(\mathbb{P}^1)$, the set of pairs of strict filtrations $\big( \mathcal{R}_r \subsetneq \cdots \subsetneq \mathcal{R}_2 \subsetneq \mathcal{R}_1=\mathcal{F}; \;\; \mathcal{P}_r \subsetneq \cdots \subsetneq \mathcal{P}_2 \subsetneq \mathcal{P}_1=\mathcal{G} \big)$ satisfying $\mathcal{R}_i/\mathcal{R}_{i+1}=\mathcal{P}_{r+1-i}/\mathcal{P}_{r+2-i}$ for $i=1, \ldots, r$, is finite.\end{lem}
\noindent
\textit{Proof.} Recall that any sheaf $\mathcal{H}$ decomposes as $\mathcal{H}=\mathcal{V}\oplus \mathcal{T}$ where $\mathcal{V}$ is a vector bundle and $\mathcal{T}$ is a torsion sheaf. The degree of $\mathcal{T}$ will be denoted by $deg_{T}(\mathcal{H})$. To prove the Lemma, we will argue by induction on $rank(\mathcal{F})$ (which is equal to $rank(\mathcal{G})$), and then on $deg_T(\mathcal{F})$.

If $\mathcal{F}$ is a torsion sheaf then evidently there are only finite many strict filtrations of $\mathcal{F}$ and the statement of the Lemma is clear. Now let us fix a pair of sheaves $(\mathcal{F},\mathcal{G})$ and let us assume that the Lemma is proved for all pairs $(\mathcal{F}',\mathcal{G}')$ for which $rank(\mathcal{F}')<rank(\mathcal{F})$ or $rank(\mathcal{F}')=rank(\mathcal{F})$ but $deg_T(\mathcal{F}')<deg_T(\mathcal{F})$. For any pair of admissible filtrations $\big( \mathcal{R}_r \subsetneq \cdots \subsetneq \mathcal{R}_2 \subsetneq \mathcal{R}_1=\mathcal{F}; \;\; \mathcal{P}_r \subsetneq \cdots \subsetneq \mathcal{P}_2 \subsetneq \mathcal{P}_1=\mathcal{G} \big)$ as in the Lemma, there are maps $\mathcal{R}_r \hookrightarrow \mathcal{F}$ and $\mathcal{G} \twoheadrightarrow \mathcal{G}/\mathcal{P}_2\simeq\mathcal{R}_r$. Therefore $\mathcal{R}_r$ is the image of a certain morphism $\mathcal{G} \to \mathcal{F}$. As ${Hom}(\mathcal{G}, \mathcal{F})$ is finite, there are only finitely many choices for $\mathcal{R}_r$, and hence also finitely many choices for $\mathcal{P}_2$. Moreover, there is a bijection between the set of admissible pairs for $(\mathcal{F},\mathcal{G})$ which contain $\mathcal{R}_r$ and $\mathcal{P}_2$ as extreme terms, and the set of admissible pairs for $(\mathcal{F}/\mathcal{R}_r,\mathcal{P}_2)$. Note that either $rank(\mathcal{F}/\mathcal{R}_r)<rank(\mathcal{F})$ or $rank(\mathcal{F}/\mathcal{R}_r)=rank(\mathcal{F})$ but $deg_T(\mathcal{F}/\mathcal{R}_r)<deg_T(\mathcal{F})$. Hence by the induction hypothesis the number of admissible filtrations of the latter type is finite. But then, all together, the number of admissible pairs for $(\mathcal{F},\mathcal{G})$ is also finite. The Lemma is proved.\qed

\vspace{.1in}

Thanks to the above Lemma, it is possible to define an inverse antipode map $S^{-1}: \widetilde{\H}_{\mathbb{P}^1} \to \widetilde{\H}_{\mathbb{P}^1}^c$ which takes values in the formal completion of $\widetilde{\H}_{\mathbb{P}^1}$ (see Remark 1.15.). One can show that (the formal completion of) the image of $\Psi$ is stable under $S^{-1}$. Unfortunately there does not seem to be any nice known formulas for $S^{-1}$. Fortunately, it is often the existence of $S^{-1}$ which is important rather than the actual form.

\vspace{.2in}

\centerline{\textbf{4.4. Weighted projective lines.}}
\addcontentsline{toc}{subsection}{\tocsubsection {}{}{\; 4.4. Weighted projective lines.}}

\vspace{.15in}

\paragraph{}We will return later (in Section~4.9 and 4.10) to the problem of describing the Hall algebras of a smooth projective curve $X$. For the time being, we focus on a very interesting and important variant of the category $Coh(\mathbb{P}^1)$ --the so-called \textit{weighted projective lines}-- which were invented by Lenzing (see \cite{Lenzing1}, or \cite{Lenzingnotes}). There are many ways of introducing and interpreting these categories (see Remark~4.20.); we will follow here the original treatment \cite{GL} to which the reader is referred for the proofs of all the results of this Section.

\vspace{.15in}

Let $k$ be a field which is again either $\C$ or finite, and let $\mathbf{p}=(p_1, \ldots, p_N)$ be a collection of $N \geq 3$ positive integers. Let $L(\mathbf{p})$ be the quotient of the free $\Z$-module $\Z\vec{x}_1 \oplus \cdots \oplus \Z \vec{x}_N$ by the submodule generated by the elements $p_i \vec{x}_i-p_j\vec{x}_j$ for all $i,j$. Let $\vec{c} \in L(\mathbf{p})$ be the common value $p_i\vec{x}_i$. Then $\vec{c}$ generates a free $\Z$-submodule of $L(\mathbf{p})$ and we have $L(\mathbf{p})/\Z\vec{c} \simeq \prod_{i=1}^N \Z/p_i\Z$. Hence $L(\mathbf{p})$ is an abelian group of rank one. Let $k[L(\mathbf{p})]$ be the group algebra of $L(\mathbf{p})$ and put $G(\mathbf{p})=Spec\;k[L(\mathbf{p})] \subset (k^*)^N$. This is an affine algebraic group whose $k$-points are
$$G(\mathbf{p})(k)=\{(t_1, \ldots, t_N) \in (k^*)^N\;|\; t_1^{p_1}=\cdots=t_N^{p_N}\}.$$
The \textit{weighted projective space} $\mathbb{P}_{\mathbf{p}}$ is by definition the quotient of $\mathbb{A}^N\backslash\{0\}$ by the natural action of $G(\mathbf{p})$~: $(t_1, \ldots, t_N) \cdot (x_1, \ldots, x_N)=(t_1x_1, \ldots, t_Nx_N)$. Denote by $S(\mathbf{p})=k[X_1, \ldots, X_N]$ the coordinate ring of $\mathbb{A}^N$. A polynomial $P(X_1, \ldots, X_N)$ in $S(\mathbf{p})$ is projectively $G(\mathbf{p})$-invariant if and only if it is homogeneous with respect to the $L(\mathbf{p})$-grading on $S(\mathbf{p})$ obtained by assigning to $X_i$ the degree $\vec{x}_i$. Thus, morally, $\mathbb{P}_{\mathbf{p}}$ is the set $Spec_{L(\mathbf{p})}S(\mathbf{p})$ of nonzero prime homogeneous ideals  in $S(\mathbf{p})$. The \textit{weighted projective line} $\mathbb{X}_{\mathbf{p},\boldsymbol{\lambda}}$ is the curve in $\mathbb{P}_{\mathbf{p}}$ defined by the set of homogeneous equations
\begin{equation}\label{E:mot1}
X_s^{p_s}=X_2^{p_2}-\lambda_s X_1^{p_1}, \qquad (s \geq 3),
\end{equation}
where $\lambda_1, \ldots, \lambda_N$ are distinct points of $\mathbb{P}^1$, normalized in such a way that $\lambda_1=\infty, \lambda_2=0$ and $\lambda_3=1$. In other terms, if $I(\llambda)$ is the ideal of $S(\p)$ generated by (\ref{E:mot1}) and $S(\p,\llambda)=S(\p)/I(\llambda)$ then, morally, $\xpl=Spec_{L(\mathbf{p})}S(\p,\llambda)$. Observe that if $p_i=1$ for all $i$ then $L(\p)=\Z$ and $S(\p,\llambda)=k[X_1,X_2]$ so that $\xpl$ is indeed a generalization of $\mathbb{P}^1$. 

\vspace{.1in}

Following Serre's Theorem~\ref{T:Serre} we now \textit{define} the category of coherent sheaves on $\xpl$ to be
$$Coh(\xpl)=S(\p,\llambda)\mathrm{-}Modgr/S(\p,\llambda)\mathrm{-}Fin,$$
where $S(\p,\llambda)\mathrm{-}Fin$ is the Serre subcategory of $S(\p,\llambda)\mathrm{-}Modgr$ consisiting of finite-dimensional graded $S(\p,\llambda)$-modules. In plain terms, this means that, just as for $Coh(\mathbb{P}^1)$, to any finitely-generated graded $S(\p,\llambda)$-module $M=\bigoplus_{\vec{x} \in L(\p)} M_{\vec{x}}$ is associated a sheaf $M\;\widetilde{}$, and two such modules $M,N$ correspond to the same sheaf if $M_{\geq n}:=\bigoplus_{\vec{x} \geq n\vec{c}} M_{\vec{x}} \simeq N_{\geq n}:=\bigoplus_{\vec{x} \geq n\vec{c}} N_{\vec{x}}$. Moreover, the spaces of homomorphisms are computed as follows~:
\begin{equation}\label{E:defhomquotient2}
{Hom}_{Coh(\mathbb{P}^1)}(M\;\widetilde{}, N\;\widetilde{})=\underset{\longrightarrow} {{Lim}}\;{Hom}_{\mathbf{S}}(M_{\geq n},N).
\end{equation}
We will refer to objects of $Coh(\xpl)$ as \textit{sheaves} on $\xpl$. It is easy to see that adding (or removing) to $\p$ some $p_i$ with $p_i=1$ does not change $L(\p)$, or $S(\p,\llambda)$. Hence we may as well assume that $p_i >1$ for all $i$. Also, adding some extra $p_i$ with $p_i=1$ allows us to lift the restriction that $N \geq 3$, and we may hence define $\xpl$ for \text{any} sequence $\p$.
Let us look at a few examples of sheaves on $\xpl$.

\vspace{.15in}

For each $\vec{x} \in L(\p)$ there is a sheaf $\O(\vec{x})=S(\p,\llambda)[\vec{x}]\;\widetilde{}$. These are analogues of the line bundles $\O(n)$ on $\mathbb{P}^1$. If $\mathcal{F}$ is any sheaf corresponding to a module, say $M$, then we set $\mathcal{F}(\vec{x})=\mathcal{F}\otimes \O[\vec{x}]=M[\vec{x}]\;\widetilde{}$.

\vspace{.15in}

The set of nonzero homogeneous primes in $S(\p,\llambda)$ may be partitioned into two sets~: the \textit{ordinary} primes $F(X_1^{p_1},X_2^{p_2})$, where $F(T_1,T_2)$ is a prime homogeneous polynomial distinct from $T_1,T_2$, and $T_2-\lambda_s T_1$ for $s \geq 3$; and the \textit{exceptional} primes $X_1, X_2, \ldots, X_N$.

The ordinary primes correspond to the closed points of $\mathbb{P}^1 \backslash\{\lambda_1, \ldots, \lambda_N\}$--which we will call \textit{ordinary closed points}-- while the exceptional points correspond to the points $\lambda_1, \ldots, \lambda_N$.
To each of these primes is associated a simple object of $Coh(\xpl)$~: if $F(X_1^{p_1}, X_2^{p_2})$ is an ordinary prime defining a closed point $x \in \mathbb{P}^1\backslash\{\lambda_1, \ldots, \lambda_N\}$ then we put
$$\O_x=S(\p,\llambda)/F(X_1^{p_1},X_2^{p_2}) S(\p,\llambda)\;\widetilde{},$$
and for $i=1, \ldots, N$ we put
$$\O_{\lambda_i}=S(\p,\llambda)/X_i S(\p,\llambda)\;\widetilde{}.$$

\vspace{.15in} 

As the next series of Lemmas will show, the sheaves $\O(\vec{x})$, $\O_x$ and $\O_{\lambda_i}$ are the building blocks of $Coh(\xpl)$.

\vspace{.1in}

\begin{lem}[\cite{GL}]The following hold~:
\begin{enumerate}
\item[i)] Let $x$ be a closed point of $\mathbb{P}^1\backslash \{\lambda_1, \ldots, \lambda_N\}$. Then for any $\vec{y} \in L(\p)$ we have $\O_x(\vec{y})\simeq \O_x$. \\
\item[ii)] For $i=1, \ldots, N$ and $\vec{y},\vec{z} \in L(\p)$, we have $\O_{\lambda_i}(\vec{y}) \simeq \O_{\lambda_i}(\vec{z})$ if and only if $\vec{z}-\vec{y} \in \sum_{j \neq i} \Z \vec{x}_j$.
\end{enumerate}
Moreover, any simple sheaf on $\xpl$ is isomorphic to one of the above.
\end{lem}

\vspace{.15in}

A sheaf $\mathcal{F}$ will be called \textit{torsion} if it is of finite length in $Coh(\xpl)$. The full subcategory $Tor(\xpl)$ of torsion sheaves is a Serre subcategory (i.e. it is an abelian subcategory, stable under subobject, quotient and extensions). For $x$ an ordinary closed point we let $Tor_x$ be the Serre subcategory of $Tor(\xpl)$ generated by the simple sheaf $\O_x$, and for $i=1, \ldots, N$, we let $Tor_{\lambda_i}$ be the Serre subcategory generated by the simple sheaves $\{\O_{\lambda_i}, \O_{\lambda_i}(\vec{x}_i), \ldots, \O_{\lambda_i}((p_i-1)\vec{x}_i)\}$. The structure of $Tor(\xpl)$ is explained by the following result.

\vspace{.1in}

\begin{lem}[\cite{GL}]\label{L:classtorxpl} The category $Tor(\xpl)$ decomposes as a direct product of orthogonal blocks
$$Tor(\xpl)=\prod_{x \in \mathbb{P}^1\backslash\{\lambda_1, \ldots, \lambda_N\}} \hspace{-.27in}Tor_x \times \prod_{i=1}^N Tor_{\lambda_i}.$$
Moreover, $Tor_x$ is equivalent to the category $Rep^{nil}_{k_x}\vec{Q}_0$ of nilpotent representations of the Jordan quiver over the residue field $k_x$, and $Tor_{\lambda_i}$ is equivalent to the category $Rep^{nil}_kA_{p_i-1}^{(1)}$ of nilpotent representations over $k$ of the cyclic quiver of length $p_i$.
\end{lem}

\vspace{.1in}

For $i=1, \ldots, N$ the above equivalence $Tor_{\lambda_i} \stackrel{\sim}{\to} Rep^{nil}_k A_{p_i-1}^{(1)}$ sends the simple sheaf $S_{i;l}:=\O_{\lambda_i}(l\vec{x}_i)$ to the simple module $S_l$. Hence the above Lemma claims in particular that 
$${dim}\;{Ext}^1(S_{i;l},S_{i;n})=\delta_{l,n+1\;(mod\;p_i)}.$$
We may thus intuitively think of $\xpl$ as a projective line $\mathbb{P}^1$ of which the $N$ points $\{\lambda_1, \ldots, \lambda_N\}$ have been ``thickened''. In other terms, it is the result of ``inserting'' at each of the special points $\lambda_i$ a category of representations of a cyclic quiver of length $p_i$~:

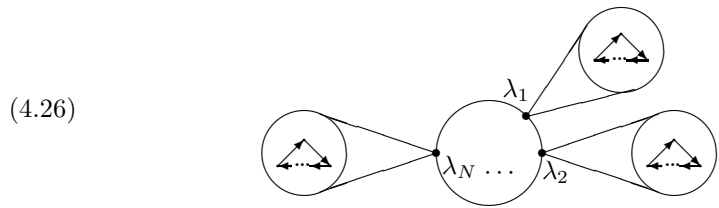
\vspace{.55in}

\begin{equation}\label{E:poneNpoints}
\centerline{
\begin{picture}(300,10)
\put(165,14){\line(4,1){40}}
\put(165,15){\line(2,3){22}}
\put(150,0){\circle{40}}
\put(164,14){\circle*{3}}
\put(155,21){$\lambda_1$}
\put(200,39){\circle{30}}
\put(200,45){\circle*{1}}
\put(190,35){\circle*{1}}
\put(195,35){\circle*{1}}
\put(205,35){\circle*{1}}
\put(210,35){\circle*{1}}
\put(190,35){\vector(1,1){10}}
\put(195,35){\vector(-1,0){5}}
\put(210,35){\vector(-1,0){8}}
\put(196,35){$.$}
\put(198,35){$.$}
\put(200,35){$.$}
\put(200,45){\vector(1,-1){10}}
\put(170,-10){$\lambda_2$}
\put(220,0){\circle{30}}
\put(220,5){\circle*{1}}
\put(210,-5){\circle*{1}}
\put(215,-5){\circle*{1}}
\put(225,-5){\circle*{1}}
\put(230,-5){\circle*{1}}
\put(210,-5){\vector(1,1){10}}
\put(215,-5){\vector(-1,0){5}}
\put(230,-5){\vector(-1,0){8}}
\put(216,-5){$.$}
\put(218,-5){$.$}
\put(220,-5){$.$}
\put(220,5){\vector(1,-1){10}}
\put(148,-10){$\cdots$}
\put(132,-8){$\lambda_N$}
\put(170,0){\line(3,1){44}}
\put(170,0){\line(3,-1){44}}
\put(130,0){\line(-3,-1){44}}
\put(130,0){\line(-3,1){44}}
\put(130,0){\circle*{3}}
\put(170,0){\circle*{3}}
\put(80,0){\circle{30}}
\put(80,5){\circle*{1}}
\put(70,-5){\circle*{1}}
\put(75,-5){\circle*{1}}
\put(85,-5){\circle*{1}}
\put(90,-5){\circle*{1}}
\put(70,-5){\vector(1,1){10}}
\put(75,-5){\vector(-1,0){5}}
\put(90,-5){\vector(-1,0){8}}
\put(76,-5){$.$}
\put(78,-5){$.$}
\put(80,-5){$.$}
\put(80,5){\vector(1,-1){10}}
\end{picture}}
\end{equation}

\vspace{.4in}

Next, let us say that a sheaf $\mathcal{F}$ corresponding to a module $M$ is \textit{locally free}, or is a \textit{vector bundle} if $M_{\geq n}$ is torsion free for $n \gg 0$. The subcategory $Vec(\xpl)$ consisting of vector bundles is an exact subcategory of $Coh(\xpl)$.

\vspace{.1in}

\begin{lem}[\cite{GL}]\label{L:vbdecomp}Any vector bundle $\mathcal{V}$ has a composition series
$$0 \subset \mathcal{V}_r \subset \mathcal{V}_{r-1} \subset \cdots \subset \mathcal{V}_1=\mathcal{V}$$ whose factors $\mathcal{V}_i/\mathcal{V}_{i+1}$ are of the form $\O(\vec{x})$ for some $\vec{x} \in L(\p)$.
\end{lem}

\vspace{.1in}

The number $r$ of terms in such a filtration is an invariant of $\mathcal{V}$ and is called the \textit{rank} of $\mathcal{V}$. Thus the above Lemma claims in particular that any line bundle (i.e. vector bundle of rank one) is of the form $\O(\vec{x})$ for some $\vec{x} \in L(\p)$.

\vspace{.15in}

The above results are already enough to give us some idea as to what the Grothendieck group $K(Coh(\xpl))$ might look like. There is an exact sequence
$$\xymatrix{ 0 \ar[r] & \O \ar[r] & \O(\vec{c}) \ar[r] & \O_x\ar[r] & 0}$$
for any closed point $x \in \mathbb{P}^1\backslash \{\lambda_1, \ldots, \lambda_N\}$ of degree one. Hence $\overline{\O_x}=\overline{\O(\vec{c})}-\overline{\O}$. Hence, just like for $\mathbb{P}^1$, the classes of the sheaves $\O_x$ for ordinary points of degree one are all equal in $K(Coh(\xpl))$. Let us fix one such closed point $x_0$.

On the other hand, from the exact sequence
$$\xymatrix{ 0\ar[r] & \O(l\vec{x}_i)\ar[r] & \O((l+1)\vec{x}_i) \ar[r] & S_{i;l+1} \ar[r] & 0}$$
we get $\overline{S_{i;l+1}}=\overline{\O((l+1)\vec{x}_i)}-\overline{\O(l\vec{x}_i)}$, and in particular
\begin{equation}\label{E:relgrothxpl}
\overline{\O(\vec{c})}-\overline{\O}=\sum_{l=0}^{p_i-1} \overline{S_{i;l}}
\end{equation}
for $i=1, \ldots, N$. 

Combining the above two remarks, we easily see that the class of any line bundle belongs to 
$\Z \overline{\O} +\Z \overline{\O_{x_0}} +\sum_{i,l} \overline{S_{i;l}}.$
By Lemma~\ref{L:vbdecomp} any vector bundle is an extension of line bundles and hence the same is true for any vector bundle. Finally, by Lemma~\ref{L:classtorxpl} the same is again true for any torsion sheaf, and thus in the end for any sheaf. It turns out, though it is not so easy to prove, that (\ref{E:relgrothxpl}) is the the \textit{only} relation between the classes $\overline{\O},\overline{\O_x},$ and $\overline{S_{i;l}}$~:

\vspace{.1in}

\begin{lem}[\cite{GL}]\label{L:grothxpl} We have
$$K(Coh(\xpl))=\left(\Z\overline{\O} \oplus \Z \overline{\O_{x_0}} \oplus \bigoplus_{i,l} \Z \overline{S_{i;l}}\right)/J,$$
where $J$ is the $\Z$-module generated by the elements $\overline{\O_{x_0}}-
\sum_{l=1}^{p_i-1}\overline{S_{i;l}}$.
\end{lem}

\vspace{.1in}

As a corollary of the above Lemma, 
\begin{equation}\label{E:Grthth}
K(Coh(\xpl))=\Z\overline{\O} \oplus \Z \overline{\O_{x_0}} \oplus \bigoplus_{i}\bigoplus_{i=1}^{p_i-1} \Z \overline{S_{i;l}}
\end{equation}

\vspace{.1in}

We have already defined the rank of a coherent sheaf on $\xpl$, we may now define the \textit{degree}. Let $d: L(\p) \to \Z$ be the group morphism for which $d(\vec{x}_i)=p/p_i$, where $p=l.c.m (p_1, \ldots, p_N)$. There is a unique group morphism $deg: K(Coh(\xpl)) \to \Z$ for which $deg(\overline{\O})=0$, $deg(\overline{\O(\vec{x})})=d(\vec{x})$ and $deg(\overline{S_{i;l}})=p/p_i$. The degree of a sheaf $\mathcal{F}$ is simply the integer $deg(\overline{\mathcal{F}})$.

\vspace{.15in}

An important tool in the study of the categories $Coh(\xpl)$ is the existence, as for any smooth projective curve, of a Serre functor.
Set $\vec{\omega}= (N-2)\vec{c}-\sum_i \vec{x}_i$. The sheaf $\O(\vec{\omega})$ plays the role of the line bundle of differential forms, and indeed when $p_i=1$ for all $i$ then $\xpl=\mathbb{P}^1$ and $\vec{\omega}=deg(\omega_{\mathbb{P}^1})=-2$.

\vspace{.1in}

\begin{lem}[\cite{GL}]\label{L:Serredualxpl} For any two coherent sheaves $\mathcal{F},\mathcal{G}$ on $\xpl$ we have
\begin{equation}\label{E:Serredualityxpl}
{Ext}^1(\mathcal{F},\mathcal{G})^* \simeq {Hom}(\mathcal{G}, \mathcal{F}(\vec{\omega})).
\end{equation}
\end{lem}

\vspace{.1in}

With this in hands, it becomes an easy matter to compute the dimensions of spaces of homomorphisms and extensions~:

\begin{equation*}
{Hom}(\O(\vec{x}),\O(\vec{y}))={S}(\p,\llambda)_{\vec{y}-\vec{x}}, \qquad {Ext}^1(\O(\vec{x}),\O(\vec{y}))={S}(\p,\llambda)_{\vec{\omega}+\vec{x}-\vec{y}},
\end{equation*}
\begin{equation*}
{Hom}(\O,\O_{x})\simeq {Ext}^1(\O_{x},\O)=k_x,\qquad {Hom}(\O_{x},\O)\simeq {Ext}^1(\O,\O_{x})=\{0\},
\end{equation*}
\begin{equation*}
{dim\;Hom}(\O,S_{i;l})=\delta_{l,0\;(mod\;p_i)}, \qquad 
{dim\;Ext}^1(S_{i;l},\O)=\delta_{l,1\;(mod\;p_i)},
\end{equation*}
\begin{equation*}
{Hom}(S_{i;l},\O)={Ext}^1(\O,S_{i;l})=\{0\}, \qquad {dim}\;{Ext}^1(S_{i;l},
S_{i;n})=\delta_{l,n+1\;(mod\;p_i)}.
\end{equation*}

\vspace{.15in}

As an immediate corollary, we get 
\begin{equation}\label{E:Eulerform}
\begin{split}
&\langle \O, \O \rangle_a =1, \qquad \langle \O, \O_{x} \rangle_a =deg(x),
\qquad \langle \O_x, \O \rangle_a =-deg(x),\\
&\langle \O_x,\O_x \rangle_a =0, \qquad \langle \O_x, S_{i;l} 
\rangle_a =0, \qquad \langle S_{i;l},\O_x \rangle_a =0,\\
&\langle \O, S_{i;l} \rangle_a =\begin{cases} 1 & \mathrm{if}\; l=0\;(\text{mod}\;p_i)\\
0 & \mathrm{if}\; l \neq 0\;(\text{mod}\;p_i)\end{cases}\\
&\langle S_{i;l}, \O \rangle_a =\begin{cases} -1 & \mathrm{if}\; l=1\;(\text{mod}\;p_i)\\
0 & \mathrm{if}\; l \neq 1 \;(\text{mod}\;p_i)\end{cases}\\
&\langle S_{i;l}, S_{i',l'} \rangle_a =\begin{cases} 1 & \mathrm{if}
\; i=i', l=l'\\
-1 & \mathrm{if}\; i=i',\;l= l'+1\;(\mathrm{mod}\;p_i)\\
0 & \mathrm{otherwise} \end{cases}.
\end{split}
\end{equation}

\vspace{.15in}

By analogy with the case of smooth projective curves, we may expect that the behavior of $Coh(\xpl)$ depends strongly on the degree $d(\vec{\omega})=(N-2)p-\sum_i \frac{p}{p_i}$ of its ``canonical bundle'' $\O(\vec{\omega})$. Accordingly, we will say that a weighted projective line is~:\\
$\bullet$ \textit{parabolic} if $d(\vec{\omega})<0$; this happens when $\p$ is of one of the forms below~:
\begin{equation}\label{E:parabolictype}
(p_1,p_2), \qquad (p_1,2,2), \qquad (2,3,3), \qquad (2,3,4), \qquad (2,3,5).
\end{equation}
$\bullet$ \textit{elliptic} if $d(\vec{\omega})=0$. Then $\p$ is one of the four sequences
\begin{equation}\label{E:elliptictype}
(2,2,2,2), \qquad (3,3,3), \qquad (2,4,4), \qquad (2,3,6).
\end{equation}
$\bullet$ \textit{hyperbolic} if $d(\vec{\omega})>0$; this is the case for all other values of $\p$.

\vspace{.1in}

For reasons which will become clear in the next Lecture, a weighted projective line for which $d(\vec{\omega})=0$ is also often called \textit{tubular}.

\vspace{.15in}

We finish with an important observation. Let ${S}_0(\p,\llambda)$ be the image of the embedding $k[T_1,T_2] \hookrightarrow {S}(\p,\llambda)$ defined by $T_1 \mapsto X_1^{p_1}, T_2 \mapsto X_2^{p_2}$. Thus ${S}_0(\p,\llambda)$ is a $\Z\vec{c}\simeq \Z$-graded subring of ${S}(\p,\llambda)$, which is easily seen to be of finite codimension. The functor
\begin{equation*}
\begin{split}
Ind~: {S}^0(\p,\llambda)\mathrm{-}Modgr &\to {S}(\p,\llambda)\mathrm{-}Modgr\\
M &\mapsto M \otimes_{{S}^0(\p,\llambda)} {S}(\p,\llambda)
\end{split}
\end{equation*}
sends finite-dimensional modules to finite-dimensional modules and gives rise to a functor
$$Ind~: {S}^0(\p,\llambda)\mathrm{-}Modgr/  {S}^0(\p,\llambda)\mathrm{-}Fin\;\to\; {S}(\p,\llambda)\mathrm{-}Modgr/{S}(\p,\llambda)\mathrm{-}Fin$$
which turns out to be a fully faithful embedding. Since 
\begin{equation*}
\begin{split}
{S}^0(\p,\llambda)\mathrm{-}Modgr/  {S}^0(\p,\llambda)\mathrm{-}Fin &\simeq k[T_1,T_2]\mathrm{-}Modgr/k[T_1,T_2]\mathrm{-}Fin \\
 &\simeq Coh(\mathbb{P}^1)
\end{split}
\end{equation*}
we have obtained in this way a canonical fully faithful functor $Coh(\mathbb{P}^1) \to Coh(\xpl)$.
Under this functor, the line bundle $\O(n)$ is mapped to the line bundle $\O(n\vec{c})$, and the simple torsion sheaf $\O_x$ associated to an ordinary closed point is mapped to the simple torsion sheaf in $Coh(\xpl)$ associated to the same ordinary point $x$. As for the simple sheaf $\O_{\lambda_i}$, it is mapped to the indecomposable sheaf $({S}(\p,\llambda)/X_i^{p_i})\;\widetilde{}$.

\vspace{.15in}

\addtocounter{theo}{1}
\paragraph{\textbf{Remark \thetheo.}} As mentionned at the beginning of this Section, there are other ways of defining the weighted projective lines~:\\
i) \textit{Parabolic coherent sheaves on} $\mathbb{P}^1$. Let $X$ be any curve, let $\lambda$ be a point of $X$ of degree one, and let $p$ be a positive integer. A $(p,\lambda)$-parabolic coherent sheaf on $X$ is a sequence of coherent sheaves $(\mathcal{F}_l)_{l \in \Z}$ on $X$ equipped with maps $\phi: \mathcal{F}_i \to \mathcal{F}_{i+1}$ and isomorphisms $\epsilon_i: \mathcal{F}_{i+l} \simeq \mathcal{F}_i \otimes \O(\lambda)$ such that the composition $\epsilon_i \circ \phi_{i+p-1} \circ \cdots \circ \phi_i : \mathcal{F}_i \to \mathcal{F}_i \otimes \O(\lambda)$ coincides with the canonical embedding. Starting from $\mathbb{P}^1$, we may iterate this procedure for each of the points $\lambda_1, \ldots, \lambda_n$ and the integers $p_1, \ldots, p_N$. The resulting category of parabolic coherent sheaves in equivalent to $Coh(\xpl)$ (see \cite{Lenzing1}).\\
ii) \textit{Equivariant coherent sheaves.} Let $X$ be a smooth projective curve and let $G \subset {Aut}(X)$ be a (finite) group of automorphisms such that $X/G \simeq \mathbb{P}^1$. When $char({k})$ and the order of $G$ are relatively prime, the category $Coh_G(X)$ of $G$-equivariant coherent sheaves on $X$ is equivalent to the category $Coh(\xpl)$ where $\llambda=\{\lambda_1, \ldots, \lambda_N\}$ are the ramification points of the quotient map $X \twoheadrightarrow \mathbb{P}^1$ and $\p=\{p_1, \ldots, p_N\}$ are the corresponding ramification indices. Then $\xpl$ is parabolic, elliptic or hyperbolic according to whether $X$ is of genus zero, one, or higher. In fact, even though $Coh(\xpl)$ has some important features of $Coh(\mathbb{P}^1)$ (like a Picard group $L(\p)$ of rank one), it behaves much more like $Coh(X)$ than like $Coh(\mathbb{P}^1)$.

\vspace{.2in}

\centerline{\textbf{4.5. Crawley-Boevey's Theorem. }}
\addcontentsline{toc}{subsection}{\tocsubsection {}{}{\; 4.5. Crawley-Boevey's theorem.}}

\vspace{.15in}

\paragraph{}Let $\xpl$ be a weighted projective line. To the weight sequence $\p=(p_1, \ldots, p_N)$ we can attach a star-shaped Dynkin diagram $\mathbb{T}_{p_1, \ldots, p_N}$, which corresponds to some Kac-Moody Lie algebra $\g$. For instance, if $\p=\emptyset$ (i.e. if $\xpl=\mathbb{P}^1$) then $\g=\mathfrak{sl}_2$ while if $\p=(p_1,p_2)$ then $\g=\mathfrak{sl}_{p_1+p_2}$. 

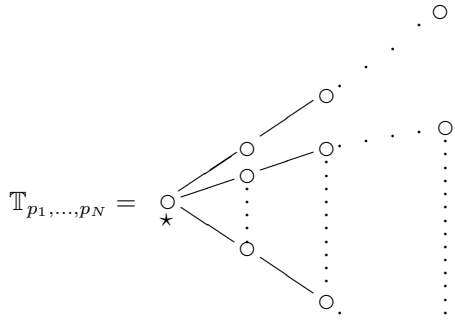
\vspace{.3in}

\centerline{
\begin{picture}(160,60)
\put(-10,-3){$\mathbb{T}_{p_1, \ldots, p_N}=$}
\put(50,0){\circle{5}}
\put(80,20){\circle{5}}
\put(110,40){\circle{5}}
\put(153,72){\circle{5}}
\put(80,10){\circle{5}}
\put(110,20){\circle{5}}
\put(155,28){\circle{5}}
\put(80,-18){\circle{5}}
\put(110,-38){\circle{5}}
\put(155,-73){\circle{5}}
\put(54,4){\line(3,2){22}}
\put(85,25){\line(3,2){20}}
\put(115,44){\circle*{1}}
\put(125,51){\circle*{1}}
\put(135,58){\circle*{1}}
\put(145,67){\circle*{1}}
\put(55,2){\line(3,1){20}}
\put(85,12){\line(3,1){20}}
\put(115,21){\circle*{1}}
\put(125,23){\circle*{1}}
\put(135,25){\circle*{1}}
\put(145,26.5){\circle*{1}}
\put(55,-2){\line(3,-2){20}}
\put(85,-21){\line(3,-2){20}}
\put(115,-42){\circle*{1}}
\put(125,-50){\circle*{1}}
\put(135,-58){\circle*{1}}
\put(145,-66){\circle*{1}}
\multiput(80,5)(0,-5){5}{\circle*{1}}
\multiput(110,15)(0,-5){10}{\circle*{1}}
\multiput(155,23)(0,-5){18}{\circle*{1}}
\put(47,-9){$\star$}
\end{picture}}

\vspace{1.2in}

Let $\mathcal{L}\g$ be the loop algebra of $\g$, as defined in Appendix~A.5. Recall that the root lattice of $\mathcal{L}\g$ is $\widehat{Q}=Q \oplus \delta$, where $Q$ is the root lattice of $\g$. Let us call $\a_{\star}$ the simple root of $\g$ associated to the central vertex of $\mathbb{T}_{p_1, \ldots, p_N}$, and $\a_{(i,l)}$ the simple root corresponding to the $l$th vertex of the $i$th branch.

\vspace{.15in}

\begin{prop}\label{P:isogo} There is an isomorphism of $\Z$-modules $\phi: K(Coh(\xpl)) \stackrel{\sim}{\to} \widehat{Q}$ satisfying
$$\phi(\overline{\O})=\a_{\star}, \qquad \phi(\overline{\O_{x}})=deg(x)\delta, \qquad \phi(\overline{S_{i;l}})=\a_{(i,l)}.$$
This isomorphism sends the symmetrized Euler form $(\;,\;)_a$ on $K(Coh(\xpl))$ to the Cartan form $(\;,\;)$ on $\widehat{Q}$.
\end{prop}
\noindent
\textit{Proof.} This follows directly from (\ref{E:Grthth}) and the formulas (\ref{E:Eulerform}).\qed

\vspace{.15in}

By analogy with the case of categories of representations of quivers, we naturally expect after the above Proposition some
deeper link between, on the one hand, the classes of indecomposable sheaves on $\xpl$, and the root system $\widehat{\Delta}$ of $\mathcal{L}\g$ on the other hand. Recall that an element
$\lambda=c_{\star}\a_{\star}+c_{\delta}\delta +\sum_{i,l} c_{(i,l)}\a_{(i,l)}$ belongs to $\widehat{Q}^+$ if $c_{\star} >0$, or $c_{\star}=0$ and $c_\delta >0$, or $c_{\star}=c_{\delta}=0$ and $c_{(i,l)} \geq 0$ for all $i,l$. Let $\widehat{\Delta}^+=\widehat{\Delta} \cap \widehat{Q}^+$ be the set of positive roots of $\mathcal{L}\g$. Observe that the class of any sheaf on $\xpl$ is automatically positive. 

\vspace{.1in}

The classes of indecomposable torsion sheaves are easily found. Let $\g_i \simeq \mathfrak{sl}_{p_i}$ be the simple Lie subalgebra of $\g$ corresponding to the $i$th branch (without the central vertex). By Lemma~\ref{L:classtorxpl} , the class of an indecomposable torsion sheaf supported at an ordinary point belongs to $\N\delta$, whereas (by Theorem~\ref{T:Godknows} applied to the cyclic quiver of length $p_i$), the class of an indecomposable torsion sheaf supported at $\lambda_i$ belongs to the positive root system $\widehat{\Delta}^+_i \subset \widehat{\Delta}^+$ corresponding to the subalgebra $\widehat{\mathfrak{sl}}_{p_i} =\mathcal{L}\g_i \subset \mathcal{L}\g$. Hence, from the definition of $\widehat{\Delta}^+$ we see that the set $\{\overline{\mathcal{T}}\;|\; \mathcal{T}\;\text{indecomposable\;torsion\;sheaf}\}$ exactly coincides with the set of roots $\a=c_{\star}\a_{\star}+c_{\delta}\delta +\sum_{i,l} c_{(i,l)}\a_{(i,l)} \in \widehat{\Delta}^+$ for which $c_{\star}=0$. 

Thus the difficulty lies in determining the classes of indecomposable vector bundles on $\xpl$. We distinguish the three types of weighted projective lines~: parabolic, elliptic, or hyperbolic.

\vspace{.15in}

Let us assume first that $\xpl$ is of parabolic type. Then, from (\ref{E:parabolictype}) we see that $\g$ is a simple Lie algebra (and moreover, any simple Lie algebra arises in this fashion). Thus $\mathcal{L}\g=\widehat{\g}$ is an affine Kac-Moody algebra.

\vspace{.1in}

\begin{theo}[Lenzing, \cite{Lenzing1}]\label{T:LL} Let $\xpl$ be a parabolic weighted projective line. Then for any $\a \in \widehat{Q}$, there exists an indecomposable coherent sheaf of class $\a$ if and only if $\a \in \widehat{\Delta}^+$. In addition, this sheaf is unique if and only if $\a$ is real.
\end{theo}

\vspace{.1in}

This is the simplest case. Since $\widehat{\Delta}^+=(\Delta \oplus \Z\delta)\backslash \{0\}$ where $\Delta$ is the root system of $\g$, and because $\Delta$ has no imaginary root, there is a bijection between the set of indecomposable vector bundles and the set of roots $\a=\a_0+l\delta$ for $l \in \Z$, where $\a_0=c_{\star}\a_{\star} +\sum_{i,l} c_{(i,l)}\a_{(i,l)}$ satisfies $c_{\star}>0$. In particular, there are, up to twisting by $\O(\vec{c})$, only finitely many indecomposable vector bundles.

\vspace{.15in}

\addtocounter{theo}{1}
\paragraph{\textbf{Example \thetheo.}} Let $\p=(p_1,p_2)$. Here $\g=\mathfrak{sl}_{p_1+p_2}$ and the roots $\a_0= c_{\star}\a_{\star} +\sum_{i,l} c_{(i,l)}\a_{(i,l)}\in \Delta$ for which $c_{\star}>0$ are
$$\a_{l_1,l_2}=\a_{\star}+\sum_{l \leq l_1} \a_{(1,l)} + \sum_{l \leq l_2} \a_{(2,l)}.$$
These correspond to the line bundles $\O(l_1\vec{x}_1+l_2\vec{x}_2)$. In particular, any indecomposable vector bundle is a line bundle in this case.
\endexample

\vspace{.15in}

\addtocounter{theo}{1}
\paragraph{\textbf{Example \thetheo.}} Let $\p=(2,2,2)$, so that $\g=\mathfrak{so}_8$. The roots
$\a_0= c_{\star}\a_{\star} +\sum_{i,l} c_{(i,l)}\a_{(i,l)}\in \Delta$ for which $c_{\star}>0$ are
\begin{equation*}
\begin{split}
&\{\a_{\star}\} \cup \{\a_{\star}+\a_{(1,1)}, \a_{\star}+\a_{(2,1)},\a_{\star}+\a_{(3,1)}\} \\
&\cup \{\a_{\star}+\a_{(1,1)}+\a_{(2,1)}, \a_{\star}+\a_{(2,1)}+\a_{(3,1)},\a_{\star}+\a_{(1,1)}+\a_{(3,1)}\}\\ 
&\cup \{\a_{\star}+\a_{(1,1)}+\a_{(2,1)}+\a_{(3,1)}\} \cup \{2\a_{\star}+\a_{(1,1)}+\a_{(2,1)}+\a_{(3,1)}\}
\end{split}
\end{equation*}
All except for the last root are again line bundles. The indecomposable vector bundle associated to 
$2\a_{\star}+\a_{(1,1)}+\a_{(2,1)}+\a_{(3,1)}$ is the unique nonsplit extension of $\O(\vec{x}_1+\vec{x}_2)$ by $\O(\vec{x}_3)$. Note that
\begin{equation*}
\begin{split}
{Ext}^1(\O(\vec{x}_1+\vec{x}_2),\O(\vec{x}_3))^* &\simeq {Hom}(\O(\vec{x}_3),\O(\vec{x}_1+\vec{x}_2 +\vec{\omega}))\\
&={Hom}(\O(\vec{x}_3),\O(\vec{c}-\vec{x}_3))\\
&={Hom}(\O(\vec{x}_3),\O(\vec{x}_3))=k
\end{split}
\end{equation*}
so there is a unique such nonsplit extension. One gets the same rank $2$ bundle by considering the extension of $\O(\vec{x}_2+\vec{x}_3)$ by $\O(\vec{x}_1)$ or of $\O(\vec{x}_1+\vec{x}_3)$ by $\O(\vec{x_2})$.
\endexample

\vspace{.15in}

Let us now turn to the elliptic case. Thus $\p$ is of one of the four types (\ref{E:elliptictype}), and $\g$ is an affine Lie algebra of type $D_4^{(1)},E_6^{(1)},E_7^{(1)}$ or $E_8^{(1)}$. Note that this is a complete list of star-shaped affine Dynkin diagram. As explained in Appendix~A.5., the loop algebra $\mathcal{L}\g$ is an elliptic Lie algebra $\mathcal{E}_{\g_0}=\g_0[t^{\pm 1},s^{\pm 1}] \oplus \mathbb{K}$ for $\g_0$ of type $D_4, E_6,E_7$ or $E_8$ respectively.

\vspace{.15in}

\begin{theo}[Lenzing, Meltzer \cite{LM}] For any $\a \in \widehat{Q}$ there exists an indecomposable sheaf on $\xpl$ of class $\a$ if and only if $\a \in \widehat{\Delta}^+$. Moreover, this sheaf is unique if $\a$ is a real root, and there is a one-parameter family of such sheaves if $\a$ is imaginary.
\end{theo}

\vspace{.1in}

Since $\g$ is an affine Lie algebra, there are imaginary roots $\a=c_{\star}\a_{\star}+\sum_{i,l} c_{(i,l)}\a_{(i,l)}$ for which $c_{\star}>0$. Hence there are now some indecomposable vector bundles which arise in families, and there are infinitely many classes of indecomposable vector bundles, even modulo twisting by line bundles. In particular, there are indecomposable vector bundles of arbitrarily high rank.

\vspace{.15in}

\addtocounter{theo}{1}
\paragraph{\textbf{Example \thetheo.}} Assume that $\p=(3,3,3)$, so that $\g=\widehat{\mathfrak{e}}_6$. We have
$\vec{\omega}=\vec{c}-\vec{x}_1-\vec{x}_2-\vec{x}_3$ and thus $3 \vec{\omega}=0$. Since ${S}(\p,\llambda)_{\vec{\omega}}=\{0\}$, the line bundles $\O,\O(\vec{\omega})$ and $\O(2\vec{\omega})$ satisfy
$${dim\;Hom}(\O(i\vec{\omega}),\O(j\vec{\omega}))=\delta_{i,j\;(\text{mod}\;3)}$$
$${dim\;Ext}^1(\O(i\vec{\omega}),\O(j\vec{\omega}))=\delta_{i+1,j\;(\text{mod}\;3)}$$
and generate an abelian subcategory closed under extensions which is equivalent to the category of (nilpotent) representations of the cyclic quiver of length $3$. Hence in particular, we may construct in this way three indecomposable vector bundles of class $n(\delta'-3\delta)$, where $\delta'=3\a_{\star}+2\sum \a_{(i,1)} + \sum \a_{(i,2)}$ is the indivisible imaginary root of $\widehat{\mathfrak{e}}_6$, and $n$ is an arbitrary positive integer. Twisting by $\O(\vec{c})$, we get more generally three indecomposable vector bundles of class $n(\delta' +3m\delta)$ for any $m \in \Z$.

These are, however, far from being all indecomposables in these classes. For instance, let $x$ be an arbitrary ordinary closed point of degree one of $\xpl$. Then $Ext^1(\O_x,\O(j\omega))=k$ for all $j$. One can show that the universal extension $\mathcal{V}_x$
$$\xymatrix{
0\ar[r] & \O \oplus \O(\vec{\omega}) \oplus \O(2\vec{\omega}) \ar[r] & \mathcal{V}_x \ar[r] & \O_x \ar[r] &0}
$$
is an indecomposable vector bundle of class $\delta'-2\delta$. But then $\mathcal{V}_x(\vec{x}_1+\vec{x}_2)$ is an indecomposable of class $\delta'$. Varying the point $x$, we obtain a whole one-parameter family of such indecomposables.
\endexample

\vspace{.15in}

The structure of $Coh(\xpl)$ in the elliptic case is in fact very well understood thanks to the work of Lenzing and Meltzer \cite{LM} (see also Section~4.8.).

\vspace{.15in}

The hyperbolic case covers all the weight sequences $\p$ not mentioned above. It turns out, though it is much harder to prove, that the analogue of Kac's theorem also holds here.

\vspace{.1in}

\begin{theo}[Crawley-Boevey, \cite{CB2}]\label{T:CBKac} For any $\a \in \widehat{Q}$ there exists an indecomposable coherent sheaf on $\xpl$ of class $\a$ if and only if $\a\in \widehat{\Delta}^+$. Moreover, this sheaf is unique if and only if $\a$ is a real root.\end{theo}

\vspace{.2in}

\centerline{\textbf{4.6. The Hall algebra of $Coh(\mathbb{X}_{\underline{p},\underline{\lambda}})$.}}
\addcontentsline{toc}{subsection}{\tocsubsection {}{}{\; 4.6. The Hall algebra of a weighted projective line.}}

\vspace{.15in}

\paragraph{}Let $\xpl$ be an arbitrary weighted projective line. As above, we associate to $\xpl$ a star-shaped Dynkin diagram $\mathbb{T}_{p_1, \ldots, p_N}$, a Kac-Moody Lie algebra $\g$, and finally a loop algebra $\mathcal{L}\g$. We have also defined a positive Borel subalgebra $\mathcal{L}\bo_+$ corresponding to the set of positive roots $\widehat{\Delta}^+$.
Let $\U_v(\mathcal{L}\g)$ and $\U_v(\mathcal{L}\bo_+)$ be the quantized enveloping algebras of $\mathcal{L}\g$ and $\mathcal{L}\bo_+$ (see Appendix~A.6.). The main result of this section relates $\U_v(\mathcal{L}\bo_+)$ to the Hall algebra $\H_{\xpl}:=\H_{Coh(\xpl)}$. Before we state this more precisely, we need to introduce a little more notation.

\vspace{.15in}

Let $k=\mathbb{F}_q$ and set $\nu=q^{\frac{1}{2}}$. We denote by $\U_{\nu}(\mathcal{L}\bo_+)$ the specialization at $v=\nu$ of $\U_v(\mathcal{L}\bo_+)$. Recall that $\U_{\nu}(\mathcal{L}\g)$ contains the subalgebras $\U_{\nu}(\widehat{\mathfrak{sl}}_{p_i})$ corresponding to each branch of $\mathbb{T}_{p_1, \ldots, p_N}$, and that $\U_{\nu}(\mathcal{L}\bo_+)$ contains the subalgebras $\U_{\nu}(\widehat{\bo}_i^+)$ for $i=1, \ldots, N$ where $\widehat{\bo}_i^+$ is the (standard) positive Borel subalgebra of the $\widehat{\mathfrak{sl}}_{p_i}$. We denote by $\{E_{(i,l)}, K_{(i,l)}\;|\; l=1, \ldots, p_{i}-1\}$ the standard Chevalley generators of $\U_{\nu}(\widehat{\bo}_i^+)$.
On the other hand, $\U_{\nu}(\mathcal{L}\bo_+)$ also contains by definition elements $E_{\star,n}$ for $n \in \Z$ and $H_{\star,r}$ for $r \geq 1$, associated to the central node $\star$ of $\mathbb{T}_{p_1, \ldots, p_N}$. These generate a subalgebra isomorphic to the positive part $\U_{\nu}(\mathcal{L}\bo_{\star}^+)$ of the quantum loop algebra $\U_{\nu}(\mathcal{L}{\mathfrak{sl}}_2)$.

As for the Hall algebra $\H_{\xpl}$, recall from Section~4.4 that the simple torsion sheaves $\{S_{i;l}\;|\; l=1, \ldots, p_i\}$ generate, for any fixed $i$, a Serre subcategory $Tor_{\lambda_i}$ equivalent to the category $Rep^{nil}_kA_{p_i-1}^{(1)}$ of nilpotent representations of the cyclic quiver. Hence, by Corollary~\ref{C:cathall}, there is an embedding $\H_{Tor_{\lambda_i}} \hookrightarrow \H_{\xpl}$ and hence, by Ringel's Theorem~\ref{T:RingelHall}, an embedding $\U_{\nu}(\widehat{\bo}_I^+) \hookrightarrow \H_{\xpl}$. In a similar manner, the fully faithful functor $Coh(\mathbb{P}^1) \to Coh(\xpl)$ discussed at the end of Section~4.4. produces an inclusion $\H_{Coh(\mathbb{P}^1)} \hookrightarrow \H_{\xpl}$ and therefore, by Kapranov's Theorem~\ref{T:Kap}, an inclusion $\U_{\nu}(\mathcal{L}\bo_\star^+) \hookrightarrow \H_{\xpl}$. By construction, this map sends the generator $H_{\star,r}$ to the element $T_r \in \H_{\xpl}$ defined as follows; for any ordinary closed point $x$ of $\xpl$, set
$$T_{r,x}=\begin{cases} \frac{[r]}{r}deg(x) \phi_x^{-1}(\frac{p_r}{deg(x)})& \text{if}\; deg(x) | r,\\
0& \text{otherwise}\end{cases}$$
and for $i=1, \ldots, N$, put
$$T_{r,{\lambda_i}}=[r]\Theta_{p_i} \circ \Phi^{-1}\bigg(\frac{p_r}{r}\bigg) \in \H_{Tor_{\lambda_i}},$$
where $\Theta_n: \H_{\vec{Q}_0} \to \H_{Tor_{\lambda_i}}$ is defined at the end of Section~3.5; now set
$T_r=\sum_{z \in \xpl} T_{r,z}$.

To get a grasp of the structure of $\H_{\xpl}$ it remains to understand how the various subalgebras identified above interact in $\H_{\xpl}$. This is given by the following Theorem, which the reader will find by now hardly surprising. As in Section~4.3, we slightly extend $\widetilde{\mathbf{H}}_{\xpl}$ by setting $\widetilde{\mathbf{H}}'_{\xpl}=\widetilde{\mathbf{H}}_{\xpl}\otimes_{\C[\mathbf{k}_{\delta}^{\pm 1}]} \C[\mathbf{k}^{\pm \delta/2}]$.

\vspace{.1in}

\begin{theo}[S. ,\cite{SDuke}]\label{T:SDuke} The assignement
$E_{\star,n}  \mapsto [\O(n\vec{c})]$ for $n \in \Z$,
$H_{\star,r} \mapsto T_r\mathbf{k}_{\delta}^{-|r|/2}$ for $r \geq 1$,
$E_{(i,l)} \mapsto [S_{i;l}]$ for $i=1, \ldots, N; 1 \leq l \leq p_i-1$,
$K_{(i,l)} \mapsto \mathbf{k}_{S_{i;l}}$,
$K_{\star} \mapsto \mathbf{k}_{\O}$, and
$C^{1/2} \mapsto \mathbf{k}_{\delta/2}$
extends to an algebra morphism 
$$\Psi: \U_{\nu}(\mathcal{L}\bo_+) \to \widetilde{\H}'_{\xpl}.$$
 This map is an embedding if $\xpl$ is either parabolic or elliptic.
\end{theo}

\vspace{.05in}

\begin{conj} The map $\Psi$ is an embedding for any weighted projective line~$\xpl$.
\end{conj}

\vspace{.15in}

The proof of the existence of the algebra morphism $\Psi$ is essentially one big, long computation, which we won't give here (see \cite{SDuke}, Section~6.). However, the fact that $\Psi$ is injective for parabolic or elliptic weighted projective lines will follow from the explicit description of the image of $\Psi$ given in the Sections~4.8. and 4.9. For reasons not entirely obvious\footnote{even to the author}, we decide to call
$Im(\Psi)$ the \textit{spherical Hall algebra} of $\xpl$. This subalgebra will be denoted $\widetilde{\mathbf{C}}_{\xpl}$. We define $\mathbf{C}_{\xpl}$ in a similar way.

\vspace{.15in}

Let us now turn to the coproduct. Because the finite subobjects condition clearly fails for $Coh(\xpl)$, $\widetilde{\H}_{\xpl}$ is a topological bialgebra only. As the categories $Tor_{\lambda_i}$ are Serre subcategories, the embedding $\widetilde{\H}_{Tor_{\lambda_i}} \to \widetilde{\H}_{\xpl}$ are compatible with the coproduct. In particular, the map $\U_{\nu}(\widehat{\bo}_i^+) \to \widetilde{\H}_{\xpl}$ is a bialgebra map, where $\U_{\nu}(\widehat{\bo}_i^+)$ is equipped with the \textit{standard} coalgebra structure (given by
(\ref{E:standardcoproduct}).)

On the other hand, the subcategory $\mathcal{C}'$ of $Coh(\xpl)$ image of the functor $Coh(\mathbb{P}^1) \to Coh(\xpl)$ is \textit{not} stable under taking subobjects or quotient (for instance, there is a map $\O \twoheadrightarrow S_{i;0}$ but $S_{i;0}$ is not an object of $\mathcal{C}'$). For a line bundle, the coproduct is determined by the following Lemma~:

\vspace{.1in}

\begin{lem}[Burban-S. \cite{BS2}] We have
$$\Delta([\O])=[\O] \otimes 1 + \sum_{\vec{x} \geq 0} \theta_{\vec{x}} \mathbf{k}_{\O(-\vec{x})}\otimes [\O(-\vec{x})],$$
where if $\vec{x}=l\vec{c} + \sum_i l_i \vec{x}_i$ with $0 \leq l_i <p_i$ then
\begin{equation*}
\begin{split}
\theta_{\vec{x}}=\nu^{2(l+\#\{i\;|\; l_i \neq 0\})}\sum_{x_j,n_j,m_i}\bigg\{ \prod_{j} (1-\nu^{-2deg(x_j)})&\hspace{-.15in}\prod_{\underset{(m_i,l_i) \neq (0,0)}{i=1}}^N \hspace{-.1in}(1-\nu^{-2}) \times \\
\times &[\bigoplus_j\O_{x_j}^{(n_j)} \oplus \bigoplus_{i=1}^N \O_{\lambda_i}^{(m_i,l_i)}]\bigg\}
\end{split}
\end{equation*}
where $\O_{\lambda_i}^{(m)i,l_i)}$ is the unique indecomposable  torsion sheaf supported at $\lambda_i$ of length $m_ip_i+l_i$ and socle $S_{i;0}$; and where the first sum ranges over all tuples of distinct ordinary points $x_j$ and multplicities $n_j, m_i$ for which $\sum_j n_jdeg(x_j)+\sum_i m_i=l$.
\end{lem}
\noindent
\textit{Proof.} It is in all points identical to the computation conducted in Example~4.12.\qed

\vspace{.15in}

Since the spherical Hall algebra $\widetilde{\mathbf{C}}_{\xpl}$ is generated by torsion sheaves and the classes of line bundles, the above is enough to completely determine the restriction of $\Delta$ and $\widetilde{\Delta}$ on $\widetilde{\mathbf{C}}_{\xpl}$. We have $\widetilde{\Delta}(\widetilde{\mathbf{C}}_{\xpl})\subset \widetilde{\mathbf{C}}_{\xpl} \widehat{\otimes} \widetilde{\mathbf{C}}_{\xpl}$. Concerning the antipode, we are in the same situation as for $\mathbb{P}^1$~: we cannot define directly the map $S$, but we may define $S^{-1}$.

\vspace{.15in}

The correspondence between the category $Coh(\xpl)$ and the quantum loop algebras may be summarized in the following table~:

\vspace{.15in}

$$
\begin{tabular}{|c|c|}
\hline
&\\
Category $\mathcal{C}=Coh(\xpl)$ & Loop algebra $\mathcal{L}\g$ of star-shaped \\
&Kac-Moody Lie algebra $\g$\\
&\\
\hline
&\\
Indices $p_1, \ldots, p_N$ of exceptional points & Length of the branches of the\\
& Dynkin diagram $\mathbb{T}_{P_1, \ldots, p_N}$ of $\g$\\
&\\
\hline
&\\
Grothendieck group $K(\mathcal{C})$ & Root lattice $\widehat{Q}=\bigoplus_{i,l} \Z \a_{(i,l)}\oplus \Z\a_{\star} \oplus\Z\delta$\\
&\\
\hline
&\\
Symmetrized additive Euler form $(\;,\;)$ & Cartan-Killing form $(\;,\;)$\\
&\\
\hline
&\\
(Classes of) simple torsion sheaves $\{S_{i;l}\}$ & Simple roots $\{\a_{(i,l)}\}$\\
&\\
\hline
&\\
(Class of) simple torsion sheaves $\O_x$ & Imaginary root $\delta$\\
&\\
\hline
&\\
(Class of) structure sheaf $\O$ & Simple root $\a_{\star}$\\
&\\
\hline
&\\
(Classes of) indecomposable objects & Positive root system $\widehat{\Delta}^+$\\
&\\
\hline
&\\
Hall algebra $\mathbf{H}_{\vec{Q}}$ & Quantum group $\U_v(\mathcal{L}\n_+)$\\
&\\
\hline
&\\
Group algebra of $K(\mathcal{C})$ & Cartan  $\U_v(\h)$\\ 
&\\
\hline
&\\
Extended Hall algebra $\widetilde{\H}_{\vec{Q}}$ & Quantum group $\U_v(\mathcal{L}\bo_+)$\\
&\\
\hline
&\\
Parabolic, elliptic or hyperbolic type & Affine Lie algebra, elliptic Lie algebra, \\
& \qquad the rest  (!)\\
&\\
\hline
\end{tabular}
$$

\newpage

\centerline{\textbf{4.7. Semistability and the Harder-Narasimhan filtration.}}
\addcontentsline{toc}{subsection}{\tocsubsection {}{}{\; 4.7. Semistability and the Harder-Narasimhan filtration.}}

\vspace{.15in}

The purpose of this section is to introduce an important technical notion for studying the categories of coherent sheaves on smooth projective curves. Throuhgout this section, $X$ will denote either a smooth projective curve or an (arbitrary) weighted projective line. A good reference here is \cite{Seshadri}.

\vspace{.1in}

Recall that there are two group morphisms $deg: K(Coh(X)) \to \Z$ and $rank: K(Coh(X)) \to \Z$, and that if $rank(\mathcal{F})=0$ then $deg(\mathcal{F})>0$. Define the \textit{slope} of a coherent sheaf $\mathcal{F}$ to be the quotient 
$$\mu(\mathcal{F})=\frac{deg(\mathcal{F})}{rank(\mathcal{F})} \in \Q \cup \{\infty\}.$$
A sheaf $\mathcal{F}$ is called \textit{semistable of slope $\eta$} if $\mu(\mathcal{F})=\eta$ and if $\mu(\mathcal{G}) \leq \eta$ for all subsheaves $\mathcal{G}\subset \mathcal{F}$. It is called \textit{stable} if in addition $\mu(\mathcal{G})<\eta$ for all proper subsheaves $\mathcal{G}\subset \mathcal{F}$.

For example, a line bundle $\mathcal{L}$ is stable of slope $\mu(\mathcal{L})=deg(\mathcal{L})$. Similarly, a torsion sheaf is semistable of slope $\infty$, and is stable if and only if it is simple.
The following simple property of the slope is very useful~: if 
$$\xymatrix{ 0 \ar[r]& \mathcal{F} \ar[r] & \mathcal{G}\ar[r] & \mathcal{H} \ar[r] & 0}$$
is exact then $min(\mu(\mathcal{F}),\mu(\mathcal{H})) \leq \mu(\mathcal{G}) \leq max(\mu(\mathcal{F}),\mu(\mathcal{H}))$. In particular, if $\mathcal{G}$ is semistable then
$\mu(\mathcal{H}) \geq \mu(\mathcal{G})$.

\vspace{.15in}

The following important Lemma justifies the notions defined above~:

\vspace{.1in}

\begin{lem}[Harder-Narasimhan] For any coherent sheaf $\mathcal{F}$ there exists a unique filtration
$$0 \subset \mathcal{F}_1 \subset \cdots \subset \mathcal{F}_l=\mathcal{F}$$
for which $\mathcal{F}_i/\mathcal{F}_{i-1}$ is semistable and $\mu(\mathcal{F}_1) > \mu(\mathcal{F}_2/\mathcal{F}_1) > \cdots > \mu(\mathcal{F}/\mathcal{F}_{l-1})$.
\end{lem}
\noindent
\textit{Proof.} Since any torsion sheaf is semistable of slope $\infty$, a semistable sheaf of slope $\eta <\infty$ is automatically a vector bundle. Hence if $\mathcal{F}=\mathcal{V} \oplus \mathcal{T}$ is a decomposition of $\mathcal{F}$ as a direct sume of a vector bundle and a torsion sheaf then we put $\mathcal{F}_1=\mathcal{T}$ and we are reduced to proving the Lemma for $\mathcal{V}$. Thus we might as well assume from the beggining that $\mathcal{F}$ is a vector bundle.

Arguing by induction on the rank, we see that it is enough to prove the following assertion~:

\vspace{.05in}

a) Any vector bundle $\mathcal{F}$ has a unique subsheaf of maximal slope and of maximal rank among\ subsheaves of that slope.

\vspace{.05in}

Of course, such a subsheaf is by construction automatically semistable. To prove a), first recall that any vector bundle has a filtration by line bundles. Let $d$ be the maximal degree of all the line bundles appearing in one such filtration. As $Hom(\mathcal{L},\mathcal{L}')=\{0\}$ if $deg(\mathcal{L}')>deg(\mathcal{L})$, we deduce from the long exact sequence in homology that $Hom(\mathcal{L},\mathcal{F}) =\{0\}$ for all line bundles $\mathcal{L}$ of degree more than $d$. A similar argument for $\Lambda^n\mathcal{F}$ in place of $\mathcal{F}$ shows that there exists $d_n$ such that $Hom(\mathcal{L},\Lambda^n\mathcal{F})=\{0\}$ as soon as $deg(\mathcal{L}) >d_n$. From this it easily follows that the degree (and hence the slope) of subsheaves of $\mathcal{F}$ is bounded from above. Therefore there exists a least one maximal semistable subsheaf $\mathcal{F}_1$ satisfying the requirements of a).

To prove that $\mathcal{F}_1$ is unique, assume that there exists on the contrary another such subsheaf $\mathcal{F}_1'$, and let $\pi: \mathcal{F} \to \mathcal{F}/\mathcal{F}_1$ be the projection. Thus $\pi(\mathcal{F}_1') \neq \{0\}$ and $\mu(\mathcal{F}_1)=\mu(\mathcal{F}'_1)$. Since $\mathcal{F}'_1$ is semistable, $\mu(\pi(\mathcal{F}'_1)) \geq \mu(\mathcal{F}'_1)$. But then, from the exact sequence
$$\xymatrix{ 0 \ar[r] & \mathcal{F}_1 \ar[r] & \mathcal{F}_1 + \mathcal{F}'_1 \ar[r] & \pi(\mathcal{F}'_1) \ar[r] & 0}$$
we get $\mu(\mathcal{F}_1+\mathcal{F}'_1)\geq \mu(\mathcal{F}_1)$. In addition, if $\mu(\mathcal{F}_1+\mathcal{F}'_1)=\mu(\mathcal{F}_1)$ then $rank(\mathcal{F}_1+\mathcal{F}'_1)>rank(\mathcal{F}_1)$, contradicting the definition of $\mathcal{F}_1$. Thus $\mathcal{F}_1$ is indeed unique, and the Lemma is proved.\qed

\vspace{.15in}

The above filtration is called the \textit{Harder-Narasimhan filtration} of the sheaf $\mathcal{F}$. It is a refinement of the canonical filtration $0 \subset \mathcal{T} \subset \mathcal{F}$. Now fix $\eta \in \Q \cup \{\infty\}$ and let $\mathcal{C}_{\eta}$ stand for the full subcategory of $Coh(X)$ whose objects are the semistable sheaves of slope $\eta$. For example, we have $\mathcal{C}_{\infty}=Tor(X)$.

\vspace{.1in}

\begin{lem}\label{L:extension=0} If $\eta>\eta'$ then $Hom(\mathcal{C}_{\eta},\mathcal{C}_{\eta'})=\{0\}$.\end{lem}
\noindent
\textit{Proof.} Let $f: \mathcal{F}\to\mathcal{G}$ be a morphism from a semistable sheaf of slope $\eta$ to another semistable sheaf of slope $\eta'$. If $f \neq 0$ then $\mu(Im(f))\geq \mu(\mathcal{F})=\eta > \eta'=\mu(\mathcal{G})$, a contradiction.\qed

\vspace{.15in}

\begin{lem} For any $\eta \in \Q \cup \{\infty\}$, the category $\mathcal{C}_{\eta}$ is abelian and closed under extension. Moreover $\mathcal{C}_{\eta}$ is a length category and the simple objects are the stable sheaves. \end{lem}
\noindent
\textit{Proof.} Let us first show that $\mathcal{C}_{\eta}$ is abelian. Let $\mathcal{F},\mathcal{G}$ be two objects of $\mathcal{C}_{\eta}$ and let $f: \mathcal{F}\to \mathcal{G}$ be a morphism. As $\mathcal{G}$ is semistable, we have $\mu(Im(f)) \leq \eta$, but as $\mathcal{F}$ is semistable, we also have $\mu(Im(f)) \geq \eta$. Therefore $\mu(Im(f))=\eta$, and it easily follows that $Im(f)$ is also semistable. The same applies to $Ker(f)$.

Next, let 
$$\xymatrix{ 0 \ar[r] & \mathcal{F}\ar[r] & \mathcal{H} \ar[r] & \mathcal{G} \ar[r] & 0}$$
be an exact sequence where $\mathcal{F}, \mathcal{G}$ belong to $\mathcal{C}_{\eta}$. It is clear that
$\mu(\mathcal{H})=\eta$. If $\mathcal{J}$ is a subsheaf of $\mathcal{H}$, then $\mathcal{J}$ is an extension
$$\xymatrix{ 0 \ar[r] & \mathcal{F}\cap \mathcal{J}\ar[r] & \mathcal{J} \ar[r] & \mathcal{J}' \ar[r] & 0}$$
where $\mathcal{J}'$ is the image of $\mathcal{J}$ under the projection $\mathcal{H} \to \mathcal{H} / \mathcal{F} \simeq \mathcal{G}$. By the semistability of $\mathcal{F}$ and $\mathcal{G}$ we have
$\mu(\mathcal{J} \cap \mathcal{F}) \leq \eta$ and $\mu(\mathcal{J}') \leq \eta$. Hence $\mu(\mathcal{J}) \leq \eta$ and $\mathcal{H}$ is semistable. 

The last assertion of the Lemma is clear since the degree of any sheaf in $\mathcal{C}_{\eta}$ is determined by its rank. \qed

\vspace{.15in}

Lemma~\ref{L:extension=0} coupled with Serre duality gives 
$${Ext}^1(\mathcal{C}_{\eta'},\mathcal{C}_{\eta})=\{0\}\qquad \text{if}\;\eta>\eta'+deg(\omega_X).$$
This has the following important corollary

\vspace{.1in}

\begin{cor}\label{C:crucial} Assume that $X$ is either a smooth projectiv curve of genus at most one, or that $X$ is a weighted projective line of parabolic or elliptic type. Then
\begin{enumerate}
\item[i)] We have $Ext^1(\mathcal{C}_{\eta'},\mathcal{C}_{\eta})=\{0\}$ if $\eta > \eta'$,
\item[ii)] The Harder-Narasimhan filtration splits,
\item[iii)] Any indecomposable sheaf is semistable.
\end{enumerate}
\end{cor}
\noindent
\textit{Proof.} Clearly, $i) \Rightarrow ii) \Rightarrow iii)$. To prove $i)$ simply observe that under the hypothesis of the Corollary, $deg(\omega_X) \leq 0$.\qed

\vspace{.15in}

When $X$ is of genus at most one, or when $X$ is a weighted projective line of parabolic or elliptic type, we may thus draw our familiar picture of the indecomposables as follows

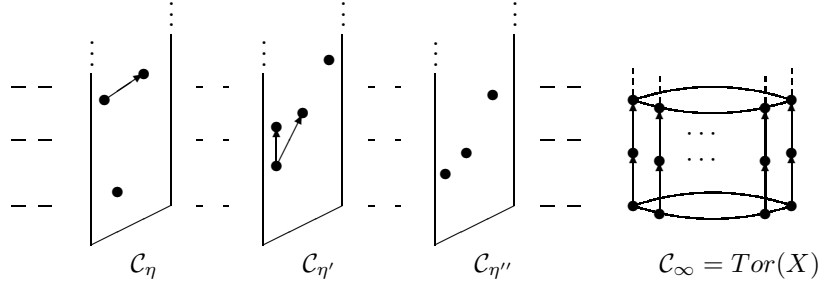
\begin{equation*}
\centerline{
\begin{picture}(300,70)
\put(5,0){\line(1,0){5}}
\put(15,0){\line(1,0){5}}
\put(5,25){\line(1,0){5}}
\put(15,25){\line(1,0){5}}
\put(5,45){\line(1,0){5}}
\put(15,45){\line(1,0){5}}
\put(215,25){\line(1,0){5}}
\put(205,25){\line(1,0){5}}
\put(205,0){\line(1,0){5}}
\put(215,0){\line(1,0){5}}
\put(205,45){\line(1,0){5}}
\put(215,45){\line(1,0){5}}
\put(75,0){\line(1,0){2}}
\put(85,0){\line(1,0){2}}
\put(75,25){\line(1,0){2}}
\put(85,25){\line(1,0){2}}
\put(75,45){\line(1,0){2}}
\put(85,45){\line(1,0){2}}
\put(140,0){\line(1,0){2}}
\put(150,0){\line(1,0){2}}
\put(140,25){\line(1,0){2}}
\put(150,25){\line(1,0){2}}
\put(140,45){\line(1,0){2}}
\put(150,45){\line(1,0){2}}
\put(50,-25){$\mathcal{C}_{\eta}$}
\put(115,-25){$\mathcal{C}_{\eta'}$}
\put(180,-25){$\mathcal{C}_{\eta''}$}
\put(34,53){$\vdots$}
\put(64,68){$\vdots$}
\put(99,53){$\vdots$}
\put(129,68){$\vdots$}
\put(35,-15){\line(0,1){65}}
\put(65,0){\line(0,1){65}}
\put(100,-15){\line(0,1){65}}
\put(130,0){\line(0,1){65}}
\put(35,-15){\line(2,1){30}}
\put(100,-15){\line(2,1){30}}
\put(165,-15){\line(0,1){65}}
\put(195,0){\line(0,1){65}}
\put(194,68){$\vdots$}
\put(164,53){$\vdots$}
\put(165,-15){\line(2,1){30}}
\put(45,5){\circle*{4}}
\put(40,40){\circle*{4}}
\put(55,50){\circle*{4}}
\put(105,15){\circle*{4}}
\put(105,30){\circle*{4}}
\put(125,55){\circle*{4}}
\put(115,35){\circle*{4}}
\put(177,20){\circle*{4}}
\put(169,12){\circle*{4}}
\put(187,42){\circle*{4}}
\put(105,15){\vector(0,1){15}}
\put(40,40){\vector(3,2){15}}
\put(105,15){\vector(1,2){10}}
\put(240,42){\line(0,1){2}}
\put(240,46){\line(0,1){2}}
\put(240,50){\line(0,1){2}}
\put(300,42){\line(0,1){2}}
\put(300,46){\line(0,1){2}}
\put(300,50){\line(0,1){2}}
\put(250,39){\line(0,1){2}}
\put(250,43){\line(0,1){2}}
\put(250,47){\line(0,1){2}}
\put(290,39){\line(0,1){2}}
\put(290,43){\line(0,1){2}}
\put(290,47){\line(0,1){2}}
\put(240,0){\vector(0,1){20}}
\put(240,0){\circle*{4}}
\put(240,40){\circle*{4}}
\put(240,20){\vector(0,1){20}}
\put(240,20){\circle*{4}}
\put(260,15){$\cdots$}
\put(260,25){$\cdots$}
\put(300,0){\vector(0,1){20}}
\put(300,0){\circle*{4}}
\put(300,40){\circle*{4}}
\put(300,20){\vector(0,1){20}}
\put(300,20){\circle*{4}}
\put(250,-3){\vector(0,1){20}}
\put(250,-3){\circle*{4}}
\put(250,37){\circle*{4}}
\put(250,17){\vector(0,1){20}}
\put(250,17){\circle*{4}}
\put(290,-3){\vector(0,1){20}}
\put(290,-3){\circle*{4}}
\put(290,37){\circle*{4}}
\put(290,17){\vector(0,1){20}}
\put(290,17){\circle*{4}}
\put(240,0){\line(0,1){40}}
\put(300,0){\line(0,1){40}}
\qbezier(240,0)(270,-10)(300,0)
\qbezier(240,40)(270,30)(300,40)
\qbezier(240,0)(270,10)(300,0)
\qbezier(240,40)(270,50)(300,40)
\put(250,-25){$\mathcal{C}_{\infty}=Tor(X)$}
\end{picture}}
\end{equation*}

\vspace{.3in}

Note that since $deg(\omega_X) <0$ when $X \simeq \mathbb{P}^1$ or when $X$ is of parabolic type,  we even have $Ext^1(\mathcal{C}_{\eta},\mathcal{C}_{\eta})=\{0\}$ in these cases for $\eta < \infty$ (i.e. each category $\mathcal{C}_{\eta}$ is semisimple).

\vspace{.15in}

We conclude this section with an easy but useful corollary of the existence (and uniqueness) of the Harder-Narasimhan filtration. For any collection of elements
$\a_1, \ldots, \a_l \in K(Coh(X))$ such that $\mu(\a_1) < \cdots < \mu(\a_l)$ we let $HN^{-1}(\a_1, \ldots, \a_l)$ be the set of all sheaves $\mathcal{F}$ with Harder-Narasimhan
filtration $\mathcal{F}_1 \subset \cdots \subset \mathcal{F}_l=\mathcal{F}$ satisfying $\overline{\mathcal{F}_i/\mathcal{F}_{i-1}}=\a_i$ for all $i$. Recall that $\mathbf{1}_{\a}
=\sum_{\overline{\mathcal{F}}=\a} [\mathcal{F}]$ for any class $\a \in K(Coh(X))$. Similarly , we put 
$$\mathbf{1}_{HN^{-1}(\a_1, \ldots, \a_l)}=\sum_{\mathcal{F} \in HN^{-1}(\a_1, \ldots, \a_l)} [\mathcal{F}].$$
For instance, $\mathbf{1}^{ss}_{\a}:=\mathbf{1}_{HN^{-1}(\a)}$ is the characteristic function of all semistable sheaves of class $\a$.

\vspace{.1in}

\begin{cor}\label{C:ssumr} For any $\a_1, \ldots, \a_l$ with $\mu(\a_1) < \cdots < \mu(\a_l)$ we have 
$$\mathbf{1}_{HN^{-1}(\a_1, \ldots, \a_l)}=\nu^{-\sum_{i<j} \langle \a_i, \a_j\rangle}\mathbf{1}^{ss}_{\a_1} \cdots \mathbf{1}^{ss}_{\a_l}.$$
\end{cor}
\noindent
\textit{Proof.} By construction, we have $\nu^{-\sum_{i<j} \langle \a_i, \a_j\rangle}\mathbf{1}^{ss}_{\a_1} \cdots \mathbf{1}^{ss}_{\a_l}=\sum_{\mathcal{F}} P^{\mathcal{F}}_{\a_1, \ldots, \a_l} [\mathcal{F}]$, where $P^\mathcal{F}_{\a_1, \ldots, \a_l}$ is the number of filtrations $\mathcal{F}_1 \subset \cdots \subset \mathcal{F}_l=\mathcal{F}$ for which $\mathcal{F}_i/\mathcal{F}_{i-1}$ is a semistable sheaf of class $\a_i$ for all $i$. But this precisely means that $\mathcal{F}_1 \subset \cdots \subset \mathcal{F}_l$ is the Harder-Narasimhan filtration of $\mathcal{F}$.\qed

\vspace{.1in}

It is clear that $\mathbf{1}_{\a}=\sum_{\a_1, \ldots, \a_l} \mathbf{1}_{HN^{-1}(\a_1, \ldots, \a_l)}$ where the sum ranges over all
collections $(\a_1, \ldots, \a_l)$ with $\mu(\a_1) < \cdots < \mu(\a_l)$ and $\sum \a_l=\a$. Therefore
\begin{equation}\label{E:ssumr}
\mathbf{1}_{\a}=\sum_{\a_1, \ldots, \a_l} \nu^{-\sum_{i<j} \langle \a_i, \a_j\rangle}\mathbf{1}^{ss}_{\a_1} \cdots \mathbf{1}^{ss}_{\a_l}
\end{equation}
(with the same range of summation as above). In other words, Corollary~\ref{C:ssumr} expresses the constant functions $\mathbf{1}_{\a}$ in terms of the
characteristic functions of semistables. Conversely, one may formally invert (\ref{E:ssumr}) to obtain the following expression of $\mathbf{1}^{s}_{\a}$ in terms of the elements
$\mathbf{1}_{\beta}$ 
\begin{equation}\label{E:ssumr2}
\mathbf{1}^{ss}_{\a}=\sum_{\beta_1, \ldots, \beta_s}(-1)^{s-1} \nu^{-\sum_{i<j} \langle \beta_i, \beta_j\rangle}\mathbf{1}_{\beta_1} \cdots \mathbf{1}_{\beta_s}
\end{equation}
where the sum now ranges over all tuples $(\beta_1, \ldots, \beta_s)$ satisfying $\mu(\sum_{l=k}^s \beta_l) > \mu(\a)$ for all $k=1, \ldots, s$.
This formula is known as \textit{Reineke's inversion formula} and a proof of it may be found in \cite{Reinekeinv}.

\vspace{.15in}

\addtocounter{theo}{1}
\paragraph{\textbf{Remarks \thetheo.}} The notions of semistability, Harder-Narasimhan filtration, etc. may be fruitfully applied to the categories of representations of quivers. However there is no \textit{a priori} canonical choice for the slope function and it may be useful to consider different such functions. See \cite{Reinekeinv} for some important applications of this circle of ideas (using a version of Corollary~\ref{C:ssumr} and equation~(\ref{E:ssumr2}) for quivers) to the computation of the cohomology of moduli spaces of representations of quivers.

\vspace{.2in}

\centerline{\textbf{4.8. The spherical Hall algebra of a parabolic weighted projective line.}}
\addcontentsline{toc}{subsection}{\tocsubsection {}{}{\; 4.8. The spherical Hall algebra of a parabolic weighted projective line.}}

\vspace{.15in}

In this section we assume that $\xpl$ is a weighted projective line of parabolic type, which is equivalent to assuming that $\g$ is a simple Lie algebra and $\mathcal{L}\g=\widehat{\g}$ is an affine Kac-Moody algebra. Our aim is to describe as precisely as possible the spherical Hall algebra $\mathbf{C}_{\xpl}$ of $\H_{\xpl}$. Set $\CC_{\mathbb{V}}=\CC_{\xpl} \cap \H_{Vec(\xpl)}$ and $\CC_{\mathbb{T}}=\CC_{\xpl} \cap \H_{Tor(\xpl)}$.

Since any extension of a vector bundle by a torsion sheaf is no longer a torsion sheaf, $\CC_{\mathbb{T}}$ is the subalgebra generated by the elements $T_r$ for $r \geq 1$ and the composition subalgebras $\CC_i$ of the subalgebras $\H_{Tor_{\lambda_i}} \subset \H_{Tor(\xpl)}$ for $i=1, \ldots, N$. Hence, in a way similar to what happens for the composition algebra of a tame quiver (see Section~3.7.), we have

\begin{prop}[Hubery]\label{P:Hubtor}
Let $\mathbf{K}'$ be the (free, commutative) algebra generated by $T_r$ for $r \geq 1$. Then 
\begin{enumerate}
\item[i)] The multiplication map is an isomorphism of vector spaces $\mathbf{K}' \otimes \bigotimes_i \mathbf{C}_i \stackrel{\sim}{\to} \mathbf{C}_{\mathbb{T}}$.
\item[ii)] The center $\mathbf{Z}$ of $\mathbf{C}_{\mathbb{T}}$ is isomorphic to $\mathbf{K}'$ and the multiplication map gives an isomorphism $\mathbf{Z} \otimes \bigotimes_i \mathbf{C}_i \stackrel{\sim}{\to} \mathbf{C}_{\mathbb{T}}$.
\end{enumerate}
\end{prop}

\vspace{.15in}

As for $\CC_{\mathbb{V}}$, things turn to be as nice as one could hope~:

\vspace{.1in}

\begin{prop}[S., \cite{SDuke}] We have $\CC_{\mathbb{V}}=\H_{Vec(\xpl)}$.\end{prop}
\noindent
\textit{Proof.} The proof ressembles that of the surjectivity of Ringel's map $\Psi: \U_{\nu}(\n_+) \to \H_{\vec{Q}}$ for finite type quivers (see Section~3.4.). First of all, it is easy to see that for $0 \leq l_i <p_i$, we have
$${Ext}^1(S_{i;l_i+1},\O(n\vec{c}+\sum_i l_i \vec{x}_i))=k$$
and the unique nonsplit extension is given by the sequence
\begin{equation*}
\begin{split}
&\xymatrix{ 0 \ar[r] & \O(n\vec{c}+\sum_i l_i \vec{x}_i) \ar[r] & \O(n\vec{c} + l_1 \vec{x}_1 + \cdots + (l_i+1) \vec{x}_i + \cdots + l_N \vec{x}_N) \ar[r] &}\\
&\hspace{3.2in}\xymatrix{& \ar[r] & S_{i;l_i+1} \ar[r] & 0}.
\end{split}
\end{equation*}
We thus obtain
\begin{equation*}
\begin{split}
\nu [S_{i;l_i+1}] \cdot [\O(n\vec{c}+\sum_i l_i \vec{x}_i)]-&\nu^2[\O(n\vec{c}+\sum_i l_i \vec{x}_i)]\cdot [S_{i;l_i+1}]\\
&=[\O(n\vec{c} +l_1 \vec{x}_1 + \cdots + (l_i+1) \vec{x}_i + \cdots + l_N \vec{x}_N)].
\end{split}
\end{equation*}
Since $[\O(n\vec{c})] \in \CC_{\mathbb{V}}$ by definition, we inductively deduce from this that $[\O(\vec{x})] \in \CC_{\mathbb{V}}$ for any line bundle $\O(\vec{x})$.

Next, we will show that $[\mathcal{V}] \in \CC_{\mathbb{V}}$ for any indecomposable vector bundle $\mathcal{V}$. Everything rests upon the following crucial Lemma~:

\vspace{.1in}

\begin{lem}\label{L:crucial} There exists a total ordering $\preceq$ on the set of all indecomposable vector bundles on $\xpl$ such that if $\mathcal{V} \prec \mathcal{V}'$ then
$${Ext}^1(\mathcal{V}, \mathcal{V}')={Hom}(\mathcal{V}',\mathcal{V})=\{0\}.$$
\end{lem}
\noindent
\textit{Proof.} By Corrolary~\ref{C:crucial} and the remark following it, it suffices to to choose any total ordering refining the partial ordering according to slope.\qed

\vspace{.1in}
So let $\mathcal{V}$ be an indecomposable vector bundle, and let us assume that $[\mathcal{V}']$ belongs to $\CC_{\mathbb{V}}$ for any indecomposable $\mathcal{V}'$ of smaller rank. Since $\mathcal{V}$ admits a composition series whose factors are line bundles, say $\mathcal{L}_1, \ldots, \mathcal{L}_r$, we have
\begin{equation}\label{E:cvect}
[\mathcal{L}_1] \cdots [\mathcal{L}_r]=c_{\mathcal{V}}[\mathcal{V}] +\sum_i c_i [\mathcal{U}_i]
\end{equation}
where $c_{\mathcal{V}} \neq 0$ and $\mathcal{U}_i$ are vector bundles different from $\mathcal{V}$ of the same class. Since by Theorem~\ref{T:LL} $\mathcal{V}$ is the unique indecomposable vector bundle of that class, each $\mathcal{U}_i$ may be nontrivially split as a direct sum
$$\mathcal{U}_i=\mathcal{W}_{i,1}^{\oplus n_{i,1}} \oplus \cdots \oplus \mathcal{W}_{i,t}^{\oplus n_{i,t}}$$
where $\mathcal{W}_{i,s}$ for $s=1, \ldots, t$ are indecomposables, which we may choose to be ordered so that $\mathcal{W}_{i,j} \prec \mathcal{W}_{i,j+1}$. Using Lemma~\ref{L:crucial}, we get
$$[\mathcal{U}_i]=\nu^{\underline{d}} [\mathcal{W}_{i,1}^{\oplus n_{i,1}}] \cdots [\mathcal{W}_{i,t}^{\oplus n_{i,t}}],$$
where
$\underline{d}=\sum_{q<s}n_{i,q}n_{i,s}{dim\;Hom}(\mathcal{W}_{i,q},\mathcal{W}_{i,s})$. Furthermore, since $Ext^1(\mathcal{W}_{i,s},\mathcal{W}_{i,s})=\{0\}$ and $End(\mathcal{W}_{i,s})=k$ for all $s$, we have
$$[\mathcal{W}_{i,s}]^{n_{i,s}}=\nu^{n_{i,s}(n_{i,s}-1)}[\mathcal{W}_{i,s}^{\oplus n_{i,s}}].$$
Now, by assumption, $[\mathcal{W}_{i,s}] \in \CC_{\mathbb{V}}$ since $rank(\mathcal{W}_{i,s}) < rank(\mathcal{V})$. Therefore $[\mathcal{W}_{i,s}^{\oplus n_{i,s}}] \in \CC_{\mathbb{V}}$ and then
$[\mathcal{U}_i] \in \CC_{\mathbb{V}}$. But then from (\ref{E:cvect}) we see that $[\mathcal{V}] \in \CC_{\mathbb{V}}$ as well. This concludes the induction step and proves that $[\mathcal{V}] \in \CC_{\mathbb{V}}$ for \textit{any} indecomposable vector bundle.

It should now be clear how to finish the proof. Any vector bundle $\mathcal{V}$ may be decomposed as a direct sum of indecomposables, and we may order these according to $\prec$. Then $[\mathcal{V}]$ is equal, up to some scalar factor, to the corresponding product in the Hall algebra. \qed

\vspace{.15in}

To conclude the description of $\CC_{\xpl}$ it only remains to add~:

\vspace{.1in}

\begin{prop}[S.] The multiplication map induces an isomorphism of vector spaces $\CC_{\mathbb{V}} \otimes \CC_{\mathbb{T}} \stackrel{\sim}{\to} \CC_{\xpl}$.\end{prop}
\noindent
\textit{Proof}. The proof, which is based on the fact that $\Delta(\CC_{\xpl}) \subset \CC_{\xpl} \widehat{\otimes} \CC_{\xpl}$, is in all points identical to the proof of Proposition~\ref{P:halldecompopo} (we replace the induction on the dimension vectors by an induction on the rank and then on the degree).
\qed

\vspace{.15in}

\addtocounter{theo}{1}
\paragraph{\textbf{Remarks. \thetheo}} i) The proof that $[\mathcal{L}]$ belongs to the spherical Hall algebra $\CC_{\xpl}$ is valid for any weighted projective line.\\
ii) The description of the structure of $\mathbf{C}_{\xpl} \cap \H_{Tor(\xpl)}$ (given by Proposition~\ref{P:Hubtor}) is also true for any weighted projective line.

\vspace{.2in}

\centerline{\textbf{4.9. The spherical Hall algebra of a tubular weighted projective line.}}
\addcontentsline{toc}{subsection}{\tocsubsection {}{}{\; 4.9. The spherical Hall algebra of a tubular weighted projective line.}}

\vspace{.15in}

Let us now move on to the case of an elliptic weighted projective line. In that situation, $\g=\widehat{\g}_0$ is an affine Lie algebra of type $D_4^{(1)}, E_6^{(1)},E_7^{(1)}$ or $E_8^{(1)}$, and $\mathcal{L}\g$ is an elliptic Lie algebra $\mathcal{E}{\g_0}$ of corresponding type.

\vspace{.1in}

The structure of $Coh(\xpl)$ is well-understood thanks to the following Theorem~:

\vspace{.1in}

\begin{theo}[Lenzing-Melzter, \cite{LM}]\label{T:LMa} For any $\eta,\eta' \in \Q \cup \{\infty\}$ there is an equivalence of categories $\epsilon_{\eta',\eta}~:\mathcal{C}_{\eta} \stackrel{\sim}{\to} \mathcal{C}_{\eta'}$. In particular, for any $\eta \in \Q$ we have $\mathcal{C}_{\eta} \simeq Tor(\xpl)$.
\end{theo}

\vspace{.15in}

As a consequence, for each $\eta \in \Q$ the category $\mathcal{C}_{\eta}$ has a block decomposition
$$\mathcal{C}_{\eta} \simeq \prod_{x \in \mathbb{P}^1 \backslash \{\lambda_1, \ldots, \lambda_N\}} \hspace{-.3in}\mathcal{C}_{\eta,x} \times \prod_{i=1}^N \mathcal{C}_{\eta,\lambda_i},$$
where $\mathcal{C}_{\eta,x} \simeq Rep^{nil}_{k_x}\vec{Q}_0$ and $\mathcal{C}_{\eta,\lambda_i} \simeq Rep^{nil}_kA_{p_i-1}^{(1)}$ are respectively equivalent to the categories of nilpotent representations of the Jordan quiver ot the cyclic quiver of length $p_i$. As a corollary of the splitting of the Harder-Narasimhan filtration (see Corrolary~\ref{C:crucial}) there is an isomorphism of vector spaces
$$\mathop{\vec{\bigotimes}}_{{\eta \in \Q \cup \{\infty\}}} \hspace{-.12in}\H_{\mathcal{C}_{\eta}} \stackrel{\sim}{\to} \H_{\xpl}$$
which is simply given by the multiplication map. Of course, the equivalence $\epsilon_{\eta',\eta}$ induces an isomorphism of algebras $\phi_{\eta',\eta}~:\H_{\mathcal{C}_{\eta}} \stackrel{\sim}{\to} \H_{\mathcal{C}_{\eta'}}$.

\vspace{.15in}

Let $\CC_{\xpl} \subset \H_{\xpl}$ denote the spherical Hall algebra. Since $\mathcal{C}_{\eta}$ is an abelian subcategory of $Coh(\xpl)$ stable under extensions, Corollary~\ref{C:cathall} provides us with an embedding of algebras $\H_{\mathcal{C}_{\eta}} \subset \H_{\xpl}$. Set $\CC_{\eta}=\CC_{\xpl} \cap \H_{\mathcal{C}_{\eta}}$. Of course $\CC_{\infty} =\CC \cap \H_{Tor(\xpl)}$ is the subalgebra of $\H_{Tor(\xpl)}$ generated by the elements $\{T_r\;|\; r \geq 1\}$ and the composition subalgebras $\CC_i$ of the exceptional tubes $\H_{Tor_{\lambda_i}}$. Its structure is hence described by Proposition~\ref{P:Hubtor}.  As far as the algebras $\CC_{\eta}$ for other values of $\eta$ are concerned, the following result shows that $\CC_{\xpl}$ possesses as much symmetry as $\H_{\xpl}$.

\vspace{.1in}

\begin{prop}[S., \cite{SDuke}] For any $\eta,\eta' \in \Q \cup \{\infty\}$ we have $\phi_{\eta',\eta}(\CC_{\eta})=\CC_{\eta'}$. Moreover, the multiplication map gives an isomorphism of vector spaces 
$$\mathop{\vec{\bigotimes}}_{\eta \in \Q \cup \{\infty\}} \hspace{-.12in}\CC_{\eta} \stackrel{\sim}{\to} \CC_{\xpl}.$$
\end{prop}

\vspace{.05in}

The proof of the second statement is again a variant of Proposition~\ref{P:halldecompopo} in the context of coherent sheaves on curves. As for the proof of the first part of the Proposition, it goes in several steps. First we use an explicit construction of the equivalences $\epsilon_{\eta',\eta}$ (the so-called ``mutations''--see \cite{LM}) to show that the elements $[S_{i;l}^{(\eta)}]=\phi_{\eta,\infty}([S_{i;l}])$ and $\phi_{\eta,\infty}(T_r)$ belong to $\CC_{\xpl}$; Then we use a graded dimension argument to conclude that $\CC_{\eta}$ is indeed generated by these elements. We refer the interested reader to \cite{SDuke}, Section~8. for more details.

\vspace{.2in}

It is natural to ask what the decomposition $\CC_{\xpl} \simeq \vec{\bigotimes} \CC_{\eta}$ means once translated to the quantum elliptic algebra $\U_{\nu}(\mathcal{L}\bo_+)$. In other words, which decomposition of the positive root system $\widehat{\Delta}^+$ of $\mathcal{E}\g_0$ is induced by the splitting of the indecomposables according to the slope ? Of course, for this it suffices to understand the rank and degree functionals on $\widehat{Q}$. 

Recall that an affine Dynkin diagram is obtained from the corresponding finite type Dynkin diagram by the adjunction of an extra node denoted $0$. The indivisible imaginary root $\delta_0$ of $\widehat{\g}_0$ is equal to $\delta_0=\a_0 + \theta$, where $\theta=\sum_i n_i\a_i$ is the maximal root of the simple Lie algebra $\g_0$.

\vspace{.1in}

\begin{lem} Let us identify $\widehat{Q}$ and $K(Coh(\xpl))$ as in Proposition~\ref{P:isogo}. Set $p=l.c.m. (p_1, \ldots, p_N)$. Then
\begin{alignat*}{2}
&deg(\a_{\star})=0, \qquad &rank(\a_{\star})=1,\\
&deg(\a_{(i,l)})=\frac{p}{p_i}, \qquad &rank(\a_{(i,l)})=0,\\
&deg(\delta)=p, \qquad &rank(\delta)=0\\
&deg(\delta_0)=p^2, \qquad &rank(\delta_0)=p
\end{alignat*}
\end{lem}
\noindent
\textit{Proof.} Everything is a straightforward verification expect for the formulas concerning $\delta_0$. For these, we use the explicit expression for the maximal root $\theta$, whose coefficients on the simple roots $\a_i$ is given in the following table~:

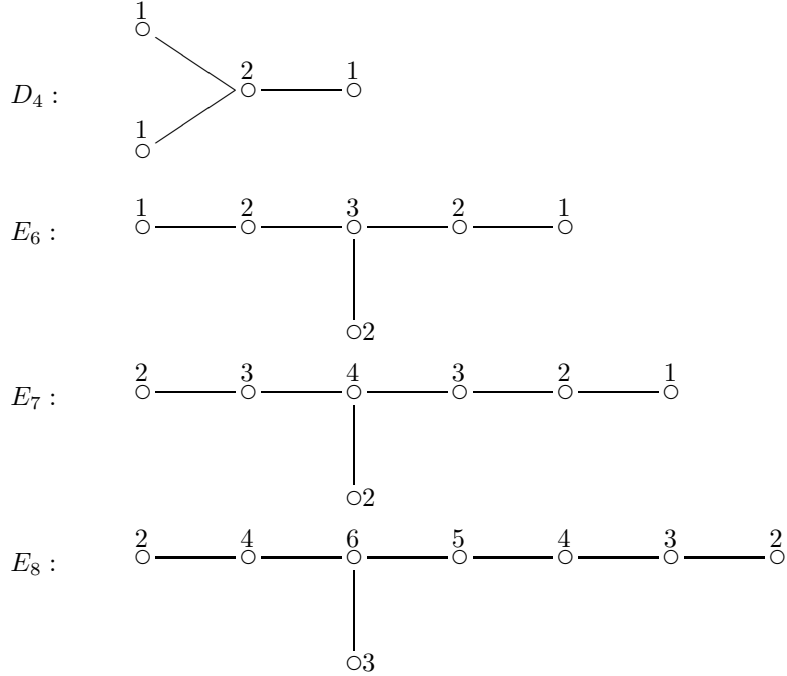
\vspace{.5in}

\centerline{
\begin{picture}(300, 10)
\put(20,-5){$D_4:$}
\put(70,23){\circle{5}}
\put(70,-23){\circle{5}}
\put(110,0){\circle{5}}
\put(150,0){\circle{5}}
\put(105,0){\line(-3,2){30}}
\put(105,0){\line(-3,-2){30}}
\put(115,0){\line(1,0){30}}
\put(67,27){$1$}
\put(67,-18){$1$}
\put(107,4){$2$}
\put(147,4){$1$}
\end{picture}}

\vspace{.55in}

\centerline{
\begin{picture}(300, 10)
\put(20,-5){$E_6:$}
\put(70,0){\circle{5}}
\put(110,0){\circle{5}}
\put(150,0){\circle{5}}
\put(190,0){\circle{5}}
\put(230,0){\circle{5}}
\put(150,-40){\circle{5}}
\put(75,0){\line(1,0){30}}
\put(115,0){\line(1,0){30}}
\put(155,0){\line(1,0){30}}
\put(150,-5){\line(0,-1){30}}
\put(195,0){\line(1,0){30}}
\put(67,4){$1$}
\put(107,4){$2$}
\put(147,4){$3$}
\put(187,4){$2$}
\put(227,4){$1$}
\put(153,-43){$2$}
\end{picture}}

\vspace{.7in}

\centerline{
\begin{picture}(300, 10)
\put(20,-5){$E_7:$}
\put(70,0){\circle{5}}
\put(110,0){\circle{5}}
\put(150,0){\circle{5}}
\put(190,0){\circle{5}}
\put(230,0){\circle{5}}
\put(270,0){\circle{5}}
\put(150,-40){\circle{5}}
\put(75,0){\line(1,0){30}}
\put(115,0){\line(1,0){30}}
\put(155,0){\line(1,0){30}}
\put(150,-5){\line(0,-1){30}}
\put(195,0){\line(1,0){30}}
\put(235,0){\line(1,0){30}}
\put(67,4){$2$}
\put(107,4){$3$}
\put(147,4){$4$}
\put(187,4){$3$}
\put(227,4){$2$}
\put(267,4){$1$}
\put(153,-43){$2$}
\end{picture}}

\vspace{.7in}

\centerline{
\begin{picture}(300, 10)
\put(20,-5){$E_8:$}
\put(70,0){\circle{5}}
\put(110,0){\circle{5}}
\put(150,0){\circle{5}}
\put(190,0){\circle{5}}
\put(230,0){\circle{5}}
\put(270,0){\circle{5}}
\put(150,-40){\circle{5}}
\put(75,0){\line(1,0){30}}
\put(115,0){\line(1,0){30}}
\put(155,0){\line(1,0){30}}
\put(150,-5){\line(0,-1){30}}
\put(195,0){\line(1,0){30}}
\put(235,0){\line(1,0){30}}
\put(310,0){\circle{5}}
\put(275,0){\line(1,0){30}}
\put(67,4){$2$}
\put(107,4){$4$}
\put(147,4){$6$}
\put(187,4){$5$}
\put(227,4){$4$}
\put(267,4){$3$}
\put(307,4){$2$}
\put(153,-43){$3$}
\end{picture}}

\vspace{.75in}

 \qed

\vspace{.2in}

Thus a root $\a=c_{\star}\a_{\star}+c_{\delta}\delta+c_{\delta_0}\delta_0+\sum_{i,l}c_{i,l}\a_{(i,l)}$ is of slope $\eta$ if 
$$\sum_{i,l} \frac{c_{(i,l)}}{p_i} +c_{\delta} +pc_{\delta_0}=\eta \left( \frac{c_{\star}}{p}+c_{\delta_0}\right).$$
Of course, each $\mathcal{E}_{\eta}:=\bigoplus_{\mu(\a)=\eta} \mathcal{E}\g_0[\a]$ is a Lie subalgebra of $\mathcal{E}{\g_0}$ and $\mathcal{E}_{\eta} \simeq \mathcal{E}_{\eta'}$ for any two $\eta,\eta'$. In particular, $\mathcal{E}_{\eta}\simeq \mathcal{E}_{\infty}=\bigoplus_i \widehat{\mathfrak{b}}_i^+ \oplus \bigoplus_{r \geq 1} \C h_{\star,r}$. Finally, the Serre functor $\cdot \otimes \omega_X$ is an automorphism of order $p$; this gives rise to a ``hidden'' action of $\Z/p\Z$ by Lie algebra automorphisms on $\mathcal{E}{\g_0}$ which preserves each $\mathcal{E}_{\eta}$. The meaning (or the usefulness) of such an action is still not clear.

All the above results naturally lift to the quantum level, though any explicit formula (for the analogues $\U_{\nu}(\mathcal{E}_{\eta})$ of the subalgebras $\mathcal{E}{\g_0}$) seems out of reach.

\vspace{.2in}

\centerline{\textbf{4.10. The Hall algebra of an elliptic curve.}}
\addcontentsline{toc}{subsection}{\tocsubsection {}{}{\; 4.10. The Hall algebra of an elliptic curve.}}

\vspace{.15in}

We now return to the question of determining the Hall algebras of smooth projective curves. Just in the same way that quivers are split according to their representation theory in the three sets: finite type, tame type and wild type, so are smooth projective curves split according to the complexity of their category of coherent sheaves into~: rational (i.e. $\mathbb{P}^1$), elliptic, and higher genus. This section deals with the second case, that of an elliptic curve $X$. The reference for everything here is the joint work \cite{BS}.

\vspace{.15in}

Let as usual $k=\mathbb{F}_q$ be a finite field and set $\nu=q^{\frac{1}{2}}$. In defining the extended Hall algebra we run into a new phenomenon: the Grohendieck group $K(Coh(X))$ is infinite-dimensional. This is not a serious issue~: by the Riemann-Roch theorem, the Euler form is
$$\langle \mathcal{F},\mathcal{G}\rangle_a=rank(\mathcal{F}) deg(\mathcal{G})-rank(\mathcal{G})deg(\mathcal{F})$$
and it makes sense to only keep track of the rank and the degree of a sheaf.
Thus we consider a ``partially extended Hall algebra'' $\widetilde{\H}_X:=\H_X \otimes
\C[\mathbf{k}_{(r,d)}]_{r,d \in \Z}$ where the imposed relations are now
$$\mathbf{k}_{(r,d)} \mathbf{k}_{(s,l)}=\mathbf{k}_{(r+s,d+l)},$$
$$\mathbf{k}_{(r,d)} [\mathcal{F}] \mathbf{k}^{-1}_{(r,d)}=[\mathcal{F}]$$
(the symmetrized Euler form vanishes). 

\vspace{.15in}

Let as usual $\mathcal{C}_{\eta}$ denote the category of semistable sheaves of slope $\eta$. The following beautiful and important result of Atiyah is the basis of all our investigations~:

\vspace{.1in}

\begin{theo}[Atiyah, \cite{A}] For any $\eta,\eta'$ there is an equivalence of categories $\epsilon_{\eta',\eta}~:\mathcal{C}_{\eta} \stackrel{\sim}{\to} \mathcal{C}_{\eta}$. In particular, $\mathcal{C}_{\eta} \simeq Tor(X)$ for any $\eta \in \Q$,
\end{theo}

\vspace{.15in}

Of course Theorem~\ref{T:LMa} is inspired by the above Theorem and can be interpreted as an equivariant  version of it. As a consequence, as in Section~4.9, there is a decomposition $\H_X \simeq \vec{\bigotimes}_{\eta} \H_{\mathcal{C}_{\eta}}$. Our aim is to understand how the various subalgebras $\H_{\mathcal{C}_{\eta}}$ commute. Since $\H_X$ and $\H_{\mathcal{C}_{\eta}}$ are of course too big, we first need to restrict ourselves to a well-chosen subalgebra (a kind of spherical subalgebra) of $\H_X$. In this end, we make the following definitions.

For $\eta \in \Q$ we set $T_{\eta,r}=\phi_{\eta,\infty}(T_r)$ where $\phi_{\eta,\infty}: \H_{\mathcal{C}_{\infty}} \stackrel{\sim}{\to} \H_{\mathcal{C}_{\eta}}$ is the isomorphism induced by the equivalence $\epsilon_{\eta,\infty}$. Thus, $T_{\eta,r}$ is a certain average of semistable sheaves of slope $\eta$. We define $\CC_X$ to be the subalgebra of $\H_X$ generated by the elements $T_{\eta,r}$ for $\eta \in \Q\cup \{\infty\}$ and $r \geq 1$.

Clearly the subalgebra of $\CC_X$ generated by $\{T_{\eta,r}\;|\; r \geq 1\}$ for any fixed $\eta$ is a free commutative polynomial algebra, canonically isomorphic to a classical Hall algebra. We put $\CC_{\eta}=\CC_X \cap \H_{\mathcal{C}_{\eta}}$.

\vspace{.1in}

\begin{theo}[Burban-S.] The following hold~:
\begin{enumerate}
\item[i)] The multiplication map induces an isomorphism of vector spaces
$$\mathop{\vec{\bigotimes}}_{\eta \in \Q \cup \{\infty\}} \hspace{-.12in}\CC_{\eta} \stackrel{\sim}{\to} \CC_X.$$
\item[ii)] For any $\eta \in \Q \cup \{\infty\}$, $\CC_{\eta}\simeq \C[T_{\eta,1}, T_{\eta,2}, \ldots]$.
\end{enumerate}
\end{theo}

\vspace{.15in}

The above Theorem gives us a precise information concerning the size of $\CC_X$. For instance, we may construct a basis of $\CC_X$ by picking any basis $\{a_{\lambda}\;|\;\lambda \in \Pi\}$ of $\CC_{\eta}$ for all $\eta \in \Q \cup \{\infty\}$ and then multiplying the elements of these bases in increasing order of their slopes. This means that in order to determine the structure of $\CC_X$ as an algebra, we need to understand how to ``reorder'' a product $T_{\eta,r}T_{\eta',s}$ when $\eta>\eta'$. This is the content of the next Theorem. Before we may properly state it, we have to introduce a bit more notation. As usual, $k=\mathbb{F}_q$ and $\nu=q^{\frac{1}{2}}$. By a Theorem of Hasse (see \cite{Hart}, App.C ) there exists conjugate algebraic numbers $\sigma, \overline{\sigma}$ satisfying $\sigma \overline{\sigma}=q$ such that $|X(q^r)|=q^r+1-(\sigma^r+\overline{\sigma}^r)$. We put $\tau=\sigma/\nu$ and $c_i(\nu,\tau)=\nu^{-i}[i]|X(q^i)|/i=(\nu^i+\nu^{-i}-\tau^i-\tau^{-i})[i]/i$. 
It will be convenient to set $T_{(p,q)}=T_{\eta,r}$ where $\eta=q/p$ and $r=g.c.d. (p,q)$. Thus, $\CC_X$ is generated by $\{T_{p,q}\;|\; (p,q) \in (\Z^{2})^+\}$ where by definition $(\Z^2)^+ =\{(p,q)\;|\; p>1\;\text{or}\;p=0, q>1\}$. Finally, for $\mathbf{x}=(p,q) \in (\Z^2)^+$ we set $deg(\mathbf{x})=g.c.d. (p,q)$ and if $\mathbf{x},\mathbf{y} \in (\Z^2)^+$ then we set $\epsilon(\mathbf{x},\mathbf{y})=sgn(det(\mathbf{x},\mathbf{y})) \in \{\pm 1\}$.

\vspace{.1in}

\begin{theo}[Burban-S.]\label{T:BS}  The assignement $T_{(p,q)} \mapsto {t}_{p,q}$ for $(p,q) \in (\Z^2)^+$ induces an isomorphism between $\CC_X$ and the algebra $\U^+_X$ generated by the elements $\{{t}_{(p,q)}\;|\; (p,q) \in (\Z^2)^+\}$ subject to the following set of relations
\begin{enumerate}
\item[i)] If $\x,\x'$ belong to the same line in $(\Z^2)^+$ then
$[{t}_\x,{t}_{\x'}]=0$,
\item[ii)] Assume that $\x,\y \in (\Z^2)^+$ are such that $deg(\x)=1$ and that 
the triangle with vertices $0, \x$ and $\x+\y$ has no interior lattice point. Then 
\begin{equation}\label{E:eqq1}
[t_\y,t_{\x}]=\epsilon_{\x,\y}c_{deg(\y)}\frac{\theta_{\x+\y}}{\nu-\nu^{-1}}
\end{equation}
where the elements $\theta_{\z}$, $\z \in (\Z^2)^+$ are obtained by equating the Fourier
 coefficients of the collection of relations
\begin{equation}\label{E:eqq2}
\sum_i \theta_{i\x_0}s^i=exp((\nu-\nu^{-1})\sum_{r \geq 1}t_{r\x_0}s^r),
\end{equation}
for any $\x_0 \in (\Z^2)^+$ such that $deg(\x_0)=1$ (note that $\theta_{\z}=(\nu-\nu^{-1})t_{\z}$ if $deg(\z)=1$). 
\end{enumerate}
\end{theo}

\vspace{.1in}

Of course, it is immediate to obtain the presentation for $\widetilde{\CC}_X:=\CC_X \otimes \C[\mathbf{k}_{(r,d)}]_{r,d \in \Z}$ from the above Theorem.

\vspace{.2in}

\addtocounter{theo}{1}
\paragraph{\textbf{Example~\thetheo.}} As an illustration of the above Theorem, let us give a sample computation of $T_{(1,2)} T_{(1,-1)}$, which we write as a path \;
\begin{picture}(20,25)
\put(0,0){\circle*{3}}
\put(20,10){\circle*{3}}
\put(0,0){\line(1,2){10}}
\put(10,20){\line(1,-1){10}}
\put(10,20){\circle*{3}}
\end{picture}.\\
We have $\left[\;
\begin{picture}(0,10)
\put(0,0){\circle*{3}}
\put(0,0){\line(0,1){10}}
\put(0,10){\circle*{3}}
\end{picture} \;,\;
\begin{picture}(10,10)
\put(0,0){\circle*{3}}
\put(0,0){\line(1,1){10}}
\put(10,10){\circle*{3}}
\end{picture}\;
\right]= c_1\;
\begin{picture}(10,25)
\put(0,0){\circle*{3}}
\put(0,0){\line(1,2){10}}
\put(10,20){\circle*{3}}
\end{picture}.$
Thus
\begin{equation*}
\begin{split}
\left[\;
\begin{picture}(10,20)
\put(0,0){\circle*{3}}
\put(0,0){\line(1,2){10}}
\put(10,20){\circle*{3}}
\end{picture} \;,\;
\begin{picture}(10,10)
\put(0,0){\circle*{3}}
\put(0,0){\line(1,-1){10}}
\put(10,-10){\circle*{3}}
\end{picture}\;
\right]&= \frac{1}{c_1}  \left\{ \left[ \left[\; 
\begin{picture}(0,10)
\put(0,0){\circle*{3}}
\put(0,0){\line(0,1){10}}
\put(0,10){\circle*{3}}
\end{picture} \;,\;
\begin{picture}(10,10)
\put(0,0){\circle*{3}}
\put(0,0){\line(1,-1){10}}
\put(10,-10){\circle*{3}}
\end{picture}\;
\right]\;,\;
\begin{picture}(10,10)
\put(0,0){\circle*{3}}
\put(0,0){\line(1,1){10}}
\put(10,10){\circle*{3}}
\end{picture}\;\right] +
\left[\;
\begin{picture}(0,10)
\put(0,0){\circle*{3}}
\put(0,0){\line(0,1){10}}
\put(0,10){\circle*{3}}
\end{picture} \;,\;
\left[\;
\begin{picture}(10,10)
\put(0,0){\circle*{3}}
\put(0,0){\line(1,1){10}}
\put(10,10){\circle*{3}}
\end{picture}\;,\;
\begin{picture}(10,10)
\put(0,0){\circle*{3}}
\put(0,0){\line(1,-1){10}}
\put(10,-10){\circle*{3}}
\end{picture}\;\right]\right]\right\}\\
&=\left[\;
\begin{picture}(10,0)
\put(0,0){\circle*{3}}
\put(0,0){\line(1,0){10}}
\put(10,0){\circle*{3}}
\end{picture} \;,\;
\begin{picture}(10,10)
\put(0,0){\circle*{3}}
\put(0,0){\line(1,1){10}}
\put(10,10){\circle*{3}}
\end{picture}\;
\right] + \left[\;
\begin{picture}(0,10)
\put(0,0){\circle*{3}}
\put(0,0){\line(0,1){10}}
\put(0,10){\circle*{3}}
\end{picture} \;,\;
\begin{picture}(20,0)
\put(0,0){\circle*{3}}
\put(0,0){\line(2,0){20}}
\put(20,0){\circle*{3}}
\end{picture} \;,
+\frac{1}{2}(\nu-\nu^{-1})\;
\begin{picture}(20,0)
\put(0,0){\circle*{3}}
\put(0,0){\line(2,0){20}}
\put(20,0){\circle*{3}}
\put(10,0){\circle*{3}}
\end{picture} \;,
\right].
\end{split}
\end{equation*}

Moreover, from
$\left[\;
\begin{picture}(10,0)
\put(0,0){\circle*{3}}
\put(0,0){\line(1,0){10}}
\put(10,0){\circle*{3}}
\end{picture} \;,\;
\begin{picture}(10,10)
\put(0,0){\circle*{3}}
\put(0,0){\line(1,1){10}}
\put(10,10){\circle*{3}}
\end{picture}\;
\right]= -c_1\;
\begin{picture}(20,10)
\put(0,0){\circle*{3}}
\put(0,0){\line(2,1){20}}
\put(20,10){\circle*{3}}
\end{picture}, \qquad
\left[\;
\begin{picture}(0,10)
\put(0,0){\circle*{3}}
\put(0,0){\line(0,1){10}}
\put(0,10){\circle*{3}}
\end{picture} \;,\;
\begin{picture}(20,0)
\put(0,0){\circle*{3}}
\put(0,0){\line(2,0){20}}
\put(20,0){\circle*{3}}
\end{picture}\;
\right]= c_2\;
\begin{picture}(20,10)
\put(0,0){\circle*{3}}
\put(0,0){\line(2,1){20}}
\put(20,10){\circle*{3}}
\end{picture} $
and $
\left[\;
\begin{picture}(0,10)
\put(0,0){\circle*{3}}
\put(0,0){\line(0,1){10}}
\put(0,10){\circle*{3}}
\end{picture} \;,\;
\begin{picture}(10,0)
\put(0,0){\circle*{3}}
\put(0,0){\line(1,0){10}}
\put(10,0){\circle*{3}}
\end{picture}\;
\right]= c_1\;
\begin{picture}(10,10)
\put(0,0){\circle*{3}}
\put(0,0){\line(1,1){10}}
\put(10,10){\circle*{3}}
\end{picture}$,  we deduce that \; 
\begin{picture}(20,20)
\put(0,0){\circle*{3}}
\put(20,10){\circle*{3}}
\put(0,0){\line(1,2){10}}
\put(10,20){\line(1,-1){10}}
\put(10,20){\circle*{3}}
\end{picture}\; is equal to~:
\begin{equation*}
\begin{split}
&\begin{picture}(20,20)
\put(0,0){\circle*{3}}
\put(10,-10){\circle*{3}}
\put(0,0){\line(1,-1){10}}
\put(10,-10){\line(1,2){10}}
\put(20,10){\circle*{3}}
\end{picture}\;
+ (c_2-c_1)\;
\begin{picture}(20,10)
\put(0,0){\circle*{3}}
\put(0,0){\line(2,1){20}}
\put(20,10){\circle*{3}}
\end{picture} + \frac{1}{2}(\nu-\nu^{-1})c_1\; \left( 
\begin{picture}(20,10)
\put(0,0){\circle*{3}}
\put(10,10){\circle*{3}}
\put(0,0){\line(1,1){10}}
\put(10,10){\line(1,0){10}}
\put(20,10){\circle*{3}}
\end{picture}\;+\;
\begin{picture}(20,10)
\put(0,0){\circle*{3}}
\put(10,0){\circle*{3}}
\put(0,0){\line(1,0){10}}
\put(10,0){\line(1,1){10}}
\put(20,10){\circle*{3}}
\end{picture}\;\right)\\
&=
\begin{picture}(20,20)
\put(0,0){\circle*{3}}
\put(10,-10){\circle*{3}}
\put(0,0){\line(1,-1){10}}
\put(10,-10){\line(1,2){10}}
\put(20,10){\circle*{3}}
\end{picture}\;
+ (c_2-c_1)\;
\begin{picture}(20,10)
\put(0,0){\circle*{3}}
\put(0,0){\line(2,1){20}}
\put(20,10){\circle*{3}}
\end{picture} + \frac{1}{2}(\nu-\nu^{-1})c_1\; \left( c_1\;
\begin{picture}(20,10)
\put(0,0){\circle*{3}}
\put(0,0){\line(2,1){20}}
\put(20,10){\circle*{3}}
\end{picture} \;
+2\;\begin{picture}(20,10)
\put(0,0){\circle*{3}}
\put(10,0){\circle*{3}}
\put(0,0){\line(1,0){10}}
\put(10,0){\line(1,1){10}}
\put(20,10){\circle*{3}}
\end{picture}\;\right)\\
&=\begin{picture}(20,20)
\put(0,0){\circle*{3}}
\put(10,-10){\circle*{3}}
\put(0,0){\line(1,-1){10}}
\put(10,-10){\line(1,2){10}}
\put(20,10){\circle*{3}}
\end{picture}\;+
(\nu-\nu^{-1})\;
\begin{picture}(20,10)
\put(0,0){\circle*{3}}
\put(10,0){\circle*{3}}
\put(0,0){\line(1,0){10}}
\put(10,0){\line(1,1){10}}
\put(20,10){\circle*{3}}
\end{picture}\;+
c_1([3]-\nu^{-1}c_1)\;
\begin{picture}(20,10)
\put(0,0){\circle*{3}}
\put(0,0){\line(2,1){20}}
\put(20,10){\circle*{3}}
\end{picture} 
\end{split}
\end{equation*}
\endexample

\vspace{.2in}

It is easy to show using Theorem~\ref{T:BS} that $\CC_X$ is generated as an algebra by the classes $\{T_{(p,q)}\;|\; p\leq 1\}$ of torsion sheaves and line bundles. Therefore, the coproduct is completely 
fixed once we have determined it for these elements. In that direction, we have~:

\vspace{.1in}

\begin{prop}[Burban-S.] The algebra $\widetilde{\CC}_X$ is stable under the coproduct map $\Delta$. In addition, we have
\begin{equation}\label{E:eqq3}
\widetilde{\Delta}(T_{(0,r)})=T_{(0,r)} \otimes 1 + \mathbf{k}_{(0,r)} \otimes T_{(0,r)},
\end{equation}
\begin{equation}\label{E:eqq4}
\widetilde{\Delta}(T_{(1,n)})=T_{(1,n)} \otimes 1 + \sum_{l \geq 0} \theta_l \mathbf{k}_{(1,n-l)} \otimes T_{(1,n-l)}
\end{equation}
where $\theta_l$ is determined by the following set of equations~:
$$\sum_{l \geq 0} \theta_l s^l=exp((\nu-\nu^{-1})\sum_{r \geq 1} T_{(0,r)}s^r).$$
\end{prop}

\noindent
\textit{Proof.} The proof of this Proposition follows \textit{verbatim} that in Example~4.12.\qed

\vspace{.2in}

Observe that all the formulas (\ref{E:eqq1}), (\ref{E:eqq2}), (\ref{E:eqq3}) and (\ref{E:eqq4}) are given by some Laurent polynomials in $\nu$ and $\tau$. Just like for the classical Hall algebra, this means that we may consider a \textit{generic} algebra $\underline{\U}^+$ which is defined over the ring $\C[\nu^{\pm 1}, \tau^{\pm 1}]$, and which specializes to $\U^+_X$ for any particular choice of an elliptic curve $X$. See \cite{BS}, Section~6 for more on this algebra.

\vspace{.15in}

\addtocounter{theo}{1}
\paragraph{\textbf{Remarks \thetheo.}} i) There is an evident similarity between the spherical Hall algebra $\CC_X$ of an elliptic curve and the spherical Hall algebra $\CC_{\xpl}$ of a weighted projective line of elliptic type. In fact, $\CC_X$ looks like a ``degenerate'' version of $\CC_{\xpl}$. Recall that these are isomorphic to the positive parts $\U_{\nu}(\mathcal{L}\bo_+)$ of quantum elliptic algebras $\U_{\nu}(\mathcal{E}{\g_0})$ when $\g_0$ is of type $D_4, E_6, E_7$ or $E_8$. As explained in \cite{SV}, $\CC_X$ may be interpreted as the positive part of a quantum elliptic algebra $\U_{\nu, \tau}(\mathcal{E}{\mathfrak{gl}(1)})$, though such an algebra has not formally been defined before. One noteworthy difference between $\CC_X$ and $\CC_{\xpl}$ is that $\CC_X$ contains \textit{two} deformation parameters-- one for the size of the finite field and one to account for the choice of the elliptic curve-- whereas $\CC_{\xpl}$ only contains one such parameter.\\
ii) The classical Hall algebra $\underline{\H}$ of Lecture~2 provides a canonical one-parameter deformation of Macdonald's ring of symmetric functions $\LLambda=\C[x_1,x_2, \ldots]^{\mathfrak{S}_{\infty}}$. Something similar happens here~: $\underline{\U}^+$ is a two-parameter  deformation of the algebra of diagonal invariants
$$\mathbf{R}=\C[x_1, x_2, \ldots, y_1^{\pm 1}, y_2^{\pm 1}, \ldots ]^{\mathfrak{S}_{\infty}},$$
where the symmetric group $\mathfrak{S}_{\infty}$ acts by permutation on the variables $x_i$ and $y_i$ simultaneously. In fact, $\underline{\U}^+$ may be interpreted as some sort of stable limit of (spherical) Cherednik Hecke algebras of type $A_n$ as $n$ tends to infinity (see \cite{SV}).

\vspace{.2in}

\centerline{\textbf{4.11. The Hall algebra of an arbitrary curve.}}
\addcontentsline{toc}{subsection}{\tocsubsection {}{}{\; 4.11. The Hall algebra of an arbitrary curve.}}

\vspace{.15in}

To conclude this Lecture, we collect the little that is known concerning the Hall algebras of (smooth) projective curves of higher genus. Let $X$ be such a curve of genus $g>1$, defined over the field $k=\mathbb{F}_q$. Let $\H_{X}$ be the Hall algebra of $Coh(X)$. As in the case of an elliptic curve, the Grothendieck group $K(Coh(X))$ is infinite dimensional. It will again be enough to take into account only the degree and the rank of a sheaf. The Riemann-Roch theorem reads
$$\langle \mathcal{F},\mathcal{G}\rangle_a=(1-g)rank(\mathcal{F})rank(\mathcal{G})+rank(\mathcal{F}) deg(\mathcal{G})-rank(\mathcal{G})deg(\mathcal{F})$$
and we may hence define a ``partially extended'' Hall algebra as $\widetilde{\H}_X:=\H_{X} \otimes \C[\mathbf{k}_{(r,d)}]$, where
$$\mathbf{k}_{(r,d)} \mathbf{k}_{(s,l)}=\mathbf{k}_{(r+s,d+l)},$$
$$\mathbf{k}_{(r,d)} [\mathcal{F}] \mathbf{k}^{-1}_{(r,d)}=\nu^{2r(1-g)rank(\mathcal{F})}[\mathcal{F}].$$

\vspace{.1in}

The first problem that we face here is to find a meaningful candidate for
$\widetilde{\CC}_{X} \subset \widetilde{\H}_X$. Though this may not be the only choice, let us define $\CC_X$ to be the subalgebra generated by the elements
$$\mathbf{1}_{(0,d)}=\sum_{\underset{\mathcal{F} \in Tor(X)}{deg(\mathcal{F})=d}}[\mathcal{F}],\qquad (d \geq 1),$$
$$\mathbf{1}^{ss}_{(1,n)}=\sum_{\underset{\mathcal{F} \;\text{line\;bundle}}{deg(\mathcal{F})=n}}[\mathcal{F}],\qquad (n \in \Z)$$
and let us set $\widetilde{\CC}_X=\CC_X \otimes \C[\mathbf{k}_{(r,d)}]_{r,d \in \Z}$.
When $X$ is of genus zero or one, this leads to the spherical Hall algebras defined in the previous Sections (see \cite{SV}, Section~6). The commutation relations satisfied by the elements $\mathbf{1}_{(0,d)}$ and $\mathbf{1}^{ss}_{(1,n)}$ are not very easy to describe. However, the coproduct may be explicitly computed~:

\vspace{.1in}

\begin{lem}\label{L:copgenus} Define elements $T_{(0,d)}$, $d \geq 1$ and $\theta_l$, $l \geq 0$ by the formulas
$$1+\sum_{r \geq 1} \mathbf{1}_{(0,d)}s^d=exp\bigg(\sum_{d \geq 1} \frac{T_{(0,d)}}{[d]}s^d\bigg),$$
$$\sum_{l \geq 0} \theta_l s^l=exp\bigg((\nu-\nu^{-1})\sum_{ r \geq 1} T_{(0,r)}s^r\bigg).$$
Then we have
\begin{equation}\label{E:coprodgenus1}
{\Delta}(T_{(0,d)})=T_{(0,d)} \otimes 1 + \mathbf{k}_{(0,d)}\otimes T_{(0,d)},
\end{equation}
\begin{equation}\label{E:coprodgenus2}
{\Delta}(\mathbf{1}^{ss}_{(1,n)})=\mathbf{1}^{ss}_{(1,n)} \otimes 1 + \sum_{l \geq 0} \theta_l \mathbf{k}_{(1,n-l)} \otimes \mathbf{1}^{ss}_{(1,n-l)}.
\end{equation}
\end{lem}
\noindent
\textit{Proof.} The proof of Example~4.12. works here as well. \qed

\vspace{.15in}

A consequence of the above Lemma is that $\widetilde{\CC}_{X}$ is stable under the coproduct. 
Next, let $(\,,\,)$ be the restriction to $\widetilde{\CC}_X$ of Green's scalar product.

\vspace{.1in}

\begin{lem}\label{L:Greengenus} We have
\begin{equation}\label{E:scalargenus}
\begin{split}
({T}_{(0,d)},\mathbf{1}^{ss}_{(1,n)})&=0\\
({T}_{(0,d)},T_{(0,l)})&=\delta_{d,l} \frac{|X(q^d)|[d]^2}{d(q^d-1)},\\
(\mathbf{1}^{ss}_{(1,n)},\mathbf{1}^{ss}_{(1,m)})&=\delta_{n,m} \frac{|Jac_X(q)|}{q-1}
\end{split}
\end{equation}
where $Jac_X$ is the Jacobian of $X$.
\end{lem}
\noindent
\textit{Proof.} The first and the third formulas are obvious. We prove the second. By definition,
$$(T_{(0,d)},T_{(0,d)})=\sum_{l |d} \sum_{\underset{deg(x)=l}{x \in X}} (T_{x,d},T_{x,d}),$$
where $T_{x,d}=\frac{[d]}{d}deg(x) \phi_x^{-1}\big(p_{\frac{d}{deg(x)}}\big)$ and $\phi^{-1}_x: \LLambda
\stackrel{\sim}{\to} \H_{Tor_x}$ is the isomorphism constructed in Section~2.4.
Using (\ref{E:HLscalarprod}), we deduce that
$$(T_{x,d},T_{x,d})=\frac{[d]^2 l}{d(q^d-1)}.$$
The relation 
$$\sum_{l|d} \sum_{\underset{deg(x)=l}{x \in X}} l=|X(q^d)|$$
finishes the proof. \qed

\vspace{.15in}

Let $Z_X(t)$ be the zeta function of $X$, which is defined as
$$Z_X(t)=exp\left(\sum_{ r \geq 1} |X(q^r)| \frac{t^r}{r}\right).$$
The Weil conjectures (proved by Deligne) assert that there exists algebraic integers $\a_1, \ldots, \a_{2g}$ all of norm $|\a_i|=q^{\frac{1}{2}}$ and satisfying $\a_i \a_{2g+1-i}=q$ such that
$$Z_X(t)=\frac{\prod_{i} (1-\a_it)}{(1-t)(1-qt)}.$$
Moreover (see e.g. \cite{Milne}, 7. Section 14.), we have $|Jac_X(q)|=\prod_i (1-\a_i)$. Thus we see that the scalar products (\ref{E:scalargenus}) are all expressed as rational functions of $\a_1, \ldots, \a_{2g}$, or even $\a_1, \ldots, \a_g$.

\vspace{.15in}

It turns out that the mere knowledge of the coproduct (Lemma~\ref{L:copgenus}) and the value of Green's scalar product (Lemma~\ref{L:Greengenus}) completely determines $\widetilde{\CC}_X$ !
To see this, let us introduce a formal topological bialgebra ${\U}'_X$ generated by elements
$\underline{T}_{(0,d)}$ for $d \geq 1$, $\underline{\mathbf{1}}^{ss}_{(1,n)}$ for $n \in \Z$ and $\underline{\mathbf{k}}_{(r,d)}$ for $r,d \in \Z$ subject to the relations
\begin{align*}
&\underline{\mathbf{k}}_{(r,d)} \underline{T}_{(0,l)}= \underline{T}_{(0,l)}\underline{\mathbf{k}}_{(r,d)},\\
&\underline{\mathbf{k}}_{(r,d)} \underline{\mathbf{1}}^{ss}_{(1,n)}\underline{\mathbf{k}}^{-1}_{(r,d)}= 
\nu^{2r(1-g)} \underline{\mathbf{1}}^{ss}_{(1,n)}
\end{align*}
and
\begin{align*}
&{\Delta}(\underline{\mathbf{k}}_{(r,d)})=\underline{\mathbf{k}}_{(r,d)} \otimes \underline{\mathbf{k}}_{(r,d)}, \\
&{\Delta}(\underline{T}_{(0,d)})=\underline{T}_{(0,d)} \otimes 1 + \underline{\mathbf{k}}_{(0,d)}\otimes \underline{T}_{(0,d)},\\
&{\Delta}(\underline{\mathbf{1}}^{ss}_{(1,n)})=\underline{\mathbf{1}}^{ss}_{(1,n)} \otimes 1 + \sum_{l \geq 0} \underline{\theta}_l \underline{\mathbf{k}}_{(1,n-l)} \otimes \underline{\mathbf{1}}^{ss}_{(1,n-l)},
\end{align*}
where $\underline{\theta}_l$ is again defined by the relation
$$\sum_{l \geq 0} \underline{\theta}_l s^l=exp\bigg((\nu-\nu^{-1})\sum_{ r \geq 1} \underline{T}_{(0,r)}s^r\bigg).$$

\vspace{.15in}

\begin{theo}[Kapranov, \cite{Kap1}]\label{T:Kappo2} There exists a unique Hopf scalar product $(\,,\,)$ on ${\U}'_X$ satisfying
\begin{equation}\label{E:scalargenus2}
\begin{split}
(\underline{T}_{(0,d)},\underline{\mathbf{1}}^{ss}_{(1,n)})&=0\\
(\underline{T}_{(0,d)},\underline{T}_{(0,l)})&=\delta_{d,l} \frac{|X(q^d)|[d]^2}{d(q^d-1)},\\
(\underline{\mathbf{1}}^{ss}_{(1,n)},\underline{\mathbf{1}}^{ss}_{(1,m)})&=\delta_{n,m} \frac{|Jac_X(q)|}{q-1}\\
(\underline{\mathbf{k}}_{(r,d)},\underline{\mathbf{k}}_{(s,l)})&=\nu^{2(1-g)rs}.
\end{split}
\end{equation}
Let $I' \subset {\U}_X'$ be the intersection of the radical of $(\,,\,)$ with the subalgebra generated by $\{\underline{T}_{(0,d)}\;|\; d \geq 1\}$ and $\{\underline{\mathbf{1}}^{ss}_{(1,n)}\;|\; n \in \Z\}$. Then $\widetilde{\CC}_X \simeq {\U}_X'/I'$.
\end{theo}
\noindent
\textit{Proof.} The first statement follows from the fact that ${\U}_X'$ is essentially a free algebra over the subalgebra generated by $\{\mathbf{k}_{(r,d)}\;|\; r,d \in \Z\}$ (of course this is not strictly true because
the elements $\mathbf{k}_{(r,d)}$ are not central). In other words we may, using the Hopf property of $(\,,\,)$, express any scalar product $(a,b)$ in terms of some fundamental scalar products 
(\ref{E:scalargenus2}) in a unique fashion.

As for the second statement, observe that by construction there is a natural surjective homomorphism
${\U}_X' \tto \widetilde{\CC}_X$, and that this surjection preserves the scalar products. Hence it only remains to show that the kernel of the restriction of Green's scalar product to $\widetilde{\CC}_X$ lies entirely in the subalgebra generated by $\{\mathbf{k}_{(r,d)},\;|\; r,d \in \Z\}$, or equivalently that the restriction of Green's scalar product to $\CC_X$ is nondegenerate. Note that because all the structure constants for the multiplication in $\H_{X}$ are real, $\CC_X$ admits an obvious real form $\CC_{X,\mathbb{R}}$-- the $\mathbb{R}$-subalgebra of $\H_{X}$ generated
by the elements $\{T_{(0,d)}\}$ and $\{\mathbf{1}^{ss}_{(1,n)}\}$. But from the definition it is equally clear that the restriction of Green's scalar product to the real form $\H_{X,\mathbb{R}}$ of $\H_X$ is positive definite. Therefore the same scalar product is positive definite on $\CC_{X,\mathbb{R}}$ and nondegenerate on $\CC_X$ as wanted. We are done.\qed

\vspace{.15in}
 
Of course, the above description is rather implicit and not so useful in practice for actual calculations. Using some classical results in the theory of automorphic functions for function fields, Kapranov managed to find certain (but \textit{a priori} not all) elements lying in the kernel $I'$~:

\vspace{.1in}

\begin{prop}[Kapranov, \cite{Kap1}, Thm. 3.3.] Set 
$$E(t)=\sum_{n \in \Z} \underline{\mathbf{1}}^{ss}_{(1,n)}t^n, \qquad \psi(s)=exp\left(\sum_{d \geq 1} 
\frac{\underline{T}_{(0,d)}}{[d]} s^d\right).$$
Then the Fourier coefficients of the following series belong to $I'$~:
$$Z_X(t_2/t_1){E}(t_1){E}(t_2)-Z_{X}(t_1/t_2){E}(t_2){E}(t_1),$$
$$\psi(s)E(t)-Z_X(v^{-1}st^{-1}){E}(t)\psi(s),$$
$$\psi(s_1)\psi(s_2)-\psi(s_2)\psi(s_1),$$
where as before $Z_X(t)$ is the zeta function of the curve $X$.
\end{prop}

\vspace{.1in}

The last relation simply says that the subalgebra of $\CC_X$ generated by $\{T_{(0,d)}\}$ is commutative, and the second says that this subalgebra acts in some kind of ``semisimple'' manner on the whole algebra. Finding more relations (most likely some higher order relations involving the series
$E(t)$) seems to be an important open problem.

\vspace{.15in}

\addtocounter{theo}{1}
\paragraph{\textbf{Remarks \thetheo.}} i) As we mentioned earlier, the values of the scalar products 
(\ref{E:scalargenus}) are all expressed as rational functions of the Weil numbers $\a_1, \ldots, \a_g$.
From the universal construction given in Theorem~\ref{T:Kappo2} we see that $\CC_X$ and $\widetilde{\CC}_X$ depend rationally on these $\g$ parameters as well, in addition to the parameter $q$. Hence, whatever the algebra ${\U}_X={\U}_X'/I'$ is, it should depend on $g+1$ deformation parameters. Indeed, we have seen that when $X$ is $\mathbb{P}^1$, or when $X$ is a weighted projective line then there is a single deformation parameter whereas when $X$ is an elliptic curve there are two such deformation parameters.\\
ii)It is also possible to define the quantum enveloping algebras $\U_v(\bo_+)$ of Kac-Moody algebras as quotients of some (essentially) free algebra equipped with a particular coproduct and a Hopf pairing, by the radical of that Hopf pairing (this is the approach taken by Lusztig \cite{Lu1}). It would be interesting to single out and study in a unified manner a class of quantum algebras obtained in this fashion, encompassing all examples appearing in these lectures.\\
iii) There is a dual construction of the spherical Hall algebra $\CC_X$ as a \textit{subalgebra} of a certain `shuffle' algebra, whose structure constants again depend rationally on the Weil numbers $\a_1, \ldots, \a_g$. We refer the interested reader to \cite{SV3} for this.\\
iv) The definition of the spherical Hall algebra given here is inspired by the interpretation of the (whole) Hall algebra of a curve $X$ as an algebra of automorphic forms for the general linear groups $GL(n, \mathbb{A}_X)$ over the adele ring of $X$. We refer the reader to \cite[Lecture~5]{SLectures2} for more.

\newpage

\centerline{\large{\textbf{Lecture~5.}}}
\addcontentsline{toc}{section}{\tocsection {}{}{Lecture~5.}}

\setcounter{section}{5}
\setcounter{equation}{0}
\setcounter{theo}{0}

\vspace{.2in}

This final Lecture brings together some examples of Hall algebras studied in Lectures~3 and 4
under the new perspective of derived equivalences. As observed roughly at the same time by Xiao and Peng (see \cite{Xiao}, \cite{PX}, see also \cite{P?}), and by Kapranov (\cite{Kap1}, \cite{Kap2}) the Hall algebras $\H_{\A}$ and $\H_{\mathcal{B}}$ of two derived equivalent \textit{hereditary} categories $\A$ and $\mathcal{B}$, though in general not isomorphic, are two different ``positive halves'' of the same algebra. Algebraically, the relevant process of ``doubling up'' a Hopf algebra is provided by the \textit{Drinfeld double} construction, which associates to any Hopf algebra $\H$ another Hopf algebra $\mathbf{DH}$ isomorphic, as a vector space, to $\H \otimes \H^*$. 

\vspace{.05in}

Thus, heuristically, one should consider the Drinfeld double $\mathbf{D}\H_{\mathcal{A}}$ as some kind of invariant (a ``Hall algebra'') of the derived category $D^b(\A)$, or more precisely of the 2-periodic category $D^b(\mathcal{A})/T^2$. Motivated by this guiding principle, Peng and Xiao constructed a Lie algebra $\mathcal{H}_\A$ attached to $D^b(\A)/T^2$ for any $\A$ of finite global dimension, and Kapranov defined an associative algebra $L(\A)$ attached to $D^b(\mathcal{A})$ (but not $D^b(\A)/T^2$) for any hereditary $\A$.

\vspace{.05in}

This line of thought naturally leads one to the more general problem of associating directly a ``Hall algebra'' $\H_{\mathcal{C}}$ to an arbitrary derived, or more generally triangulated category $\mathcal{C}$. This was only recently settled by To\"en \cite{Toen} who constructed such a Hall algebra for any dg-category under the extra assumption that $\mathcal{C}$ be \textit{of finite global dimension} (so that the case of $D^b(\A)/T^2$ is unfortunately again excluded).

\vspace{.1in}

After giving some motivating examples, we recall in this Lecture the definition of the Drinfeld double of a Hopf algebra and state various conjectures relating the Drinfeld double $\mathbf{D}\H_{\mathcal{A}}$ of a hereditary category $\mathcal{A}$ and the 2-periodic category $D^b(\mathcal{A})/T^2$. We then (succintly) give the definition of Xiao and Peng's  Hall Lie algebra of $D^b(\mathcal{A})/T^2$ and describe (even more succintly) the constructions of Kapranov and To\"en. A vey last Section contains a few open problems and directions for further investigations. The article \cite{Keller} can serve as a good reference for all the notions of homological algebra which will be used here.

\vspace{.2in}

\centerline{\textbf{5.1. Motivation.}}
\addcontentsline{toc}{subsection}{\tocsubsection {}{}{\; 5.1. Motivation.}}

\vspace{.15in}

 Although the categories of representations of quivers (discussed in Lecture~3) and the categories of coherent sheaves on smooth projective curves (which form the content of Lecture~4) are of very different origin there are, as the reader has beyond doubt already noticed, some strinking similarities. 
 
 The simplest example is that of the categories of representations of the Kronecker quiver $\vec{Q}$ (Examples~3.11 and 3.34.) and of coherent sheaves on the projective line $\mathbb{P}^1$ (see Section~4.2.)~: the picture (\ref{E:picturekronecker}) of the category $Rep_k\vec{Q}$ may be obtained from the picture (\ref{E:picturecohp1}) of $Coh(\mathbb{P}^1)$ by slicing in the middle and moving the rightmost piece to the left. This is in fact a consequence of the existence of an equivalence of \textit{derived categories} $F~:D^b(Coh(\mathbb{P}^1)) \stackrel{\sim}{\to}D^b(Rep_k\vec{Q})$. The equivalence $F$ sends the subcategory $Tor(\mathbb{P}^1)$ of torsion sheaves to the subcategory $\mathbb{R}$ of regular modules, sends the line bundles $\O,\O(1), \O(2),\ldots$ to the preprojective indecomposables $P(1),P(0),\tau^-P(1),\ldots$ and sends the line bundles $\O(-1),\O(-2),\O(-3), \ldots$ to the preinjective indecomposables $I(0),I(1), \tau I(0),\ldots$. Equivalently, we may view $Rep_k\vec{Q}$ and $Coh(\mathbb{P}^1)$ as the hearts of two distinct t-structures on the same triangulated category. This is best understood by drawing a picture ~:

\vspace{.2in}

\begin{equation}\label{E:pictureequiv}
\centerline{
\begin{picture}(330,60)
\put(2.5,0){\line(1,0){2.5}}
\put(7.5,0){\line(1,0){2.5}}
\put(2.5,12.5){\line(1,0){2.5}}
\put(7.5,12.5){\line(1,0){2.5}}
\put(2.5,25){\line(1,0){2.5}}
\put(7.5,25){\line(1,0){2.5}}
\put(95,12.5){\line(1,0){2.5}}
\put(100,12.5){\line(1,0){2.5}}
\put(95,0){\line(1,0){2.5}}
\put(100,0){\line(1,0){2.5}}
\put(95,25){\line(1,0){2.5}}
\put(100,25){\line(1,0){2.5}}
\put(53,3){\circle*{2}}
\put(60,22.5){\circle*{2}}
\put(75,22.5){\circle*{2}}
\put(68,3){\circle*{2}}
\put(45,22.5){\circle*{2}}
\put(30,22.5){\circle*{2}}
\put(38,3){\circle*{2}}
\put(52.5,5){\vector(1,3){5}}
\put(55,5){\vector(1,3){5}}
\put(62.5,20){\vector(1,-3){5}}
\put(60,20){\vector(1,-3){5}}
\put(77.5,20){\line(1,-3){4}}
\put(75,20){\line(1,-3){4}}
\put(67.5,5){\vector(1,3){5}}
\put(70,5){\vector(1,3){5}}
\put(37.5,5){\vector(1,3){5}}
\put(40,5){\vector(1,3){5}}
\put(45,20){\vector(1,-3){5}}
\put(47.5,20){\vector(1,-3){5}}
\put(24,9.5){\line(1,3){4}}
\put(26.5,9.5){\line(1,3){4}}
\put(32.5,20){\vector(1,-3){5}}
\put(30,20){\vector(1,-3){5}}
\put(120,21){\line(0,1){2}}
\put(120,23){\line(0,1){2}}
\put(120,25){\line(0,1){2}}
\put(150,21){\line(0,1){2}}
\put(150,23){\line(0,1){2}}
\put(150,25){\line(0,1){2}}
\put(125,19.5){\line(0,1){2}}
\put(125,21.5){\line(0,1){2}}
\put(125,23.5){\line(0,1){2}}
\put(145,19.5){\line(0,1){2}}
\put(145,21.5){\line(0,1){2}}
\put(145,23.5){\line(0,1){2}}
\put(120,0){\vector(0,1){10}}
\put(120,0){\circle*{2}}
\put(120,20){\circle*{2}}
\put(120,10){\vector(0,1){10}}
\put(120,10){\circle*{2}}
\put(128,5){$\cdots$}
\put(128,12){$\cdots$}
\put(150,0){\vector(0,1){10}}
\put(150,0){\circle*{2}}
\put(150,20){\circle*{2}}
\put(150,10){\vector(0,1){10}}
\put(150,10){\circle*{2}}
\put(125,-1.5){\vector(0,1){10}}
\put(125,-1.5){\circle*{2}}
\put(125,18.5){\circle*{2}}
\put(125,8.5){\vector(0,1){10}}
\put(125,8.5){\circle*{2}}
\put(145,-1.5){\vector(0,1){10}}
\put(145,-1.5){\circle*{2}}
\put(145,18.5){\circle*{2}}
\put(145,8.5){\vector(0,1){10}}
\put(145,8.5){\circle*{2}}
\put(120,0){\line(0,1){20}}
\put(150,0){\line(0,1){20}}
\qbezier(120,0)(135,-5)(150,0)
\qbezier(120,20)(135,15)(150,20)
\qbezier(120,0)(135,5)(150,0)
\qbezier(120,20)(135,25)(150,20)
\put(172.5,0){\line(1,0){2.5}}
\put(177.5,0){\line(1,0){2.5}}
\put(172.5,12.5){\line(1,0){2.5}}
\put(177.5,12.5){\line(1,0){2.5}}
\put(172.5,25){\line(1,0){2.5}}
\put(177.5,25){\line(1,0){2.5}}
\put(265,12.5){\line(1,0){2.5}}
\put(270,12.5){\line(1,0){2.5}}
\put(265,0){\line(1,0){2.5}}
\put(270,0){\line(1,0){2.5}}
\put(265,25){\line(1,0){2.5}}
\put(270,25){\line(1,0){2.5}}
\put(223,3){\circle*{2}}
\put(230,22.5){\circle*{2}}
\put(245,22.5){\circle*{2}}
\put(238,3){\circle*{2}}
\put(215,22.5){\circle*{2}}
\put(200,22.5){\circle*{2}}
\put(208,3){\circle*{2}}
\put(222.5,5){\vector(1,3){5}}
\put(225,5){\vector(1,3){5}}
\put(232.5,20){\vector(1,-3){5}}
\put(230,20){\vector(1,-3){5}}
\put(247.5,20){\line(1,-3){4}}
\put(245,20){\line(1,-3){4}}
\put(237.5,5){\vector(1,3){5}}
\put(240,5){\vector(1,3){5}}
\put(207.5,5){\vector(1,3){5}}
\put(210,5){\vector(1,3){5}}
\put(215,20){\vector(1,-3){5}}
\put(217.5,20){\vector(1,-3){5}}
\put(194,9.5){\line(1,3){4}}
\put(196.5,9.5){\line(1,3){4}}
\put(202.5,20){\vector(1,-3){5}}
\put(200,20){\vector(1,-3){5}}
\put(290,21){\line(0,1){2}}
\put(290,23){\line(0,1){2}}
\put(290,25){\line(0,1){2}}
\put(320,21){\line(0,1){2}}
\put(320,23){\line(0,1){2}}
\put(320,25){\line(0,1){2}}
\put(295,19.5){\line(0,1){2}}
\put(295,21.5){\line(0,1){2}}
\put(295,23.5){\line(0,1){2}}
\put(315,19.5){\line(0,1){2}}
\put(315,21.5){\line(0,1){2}}
\put(315,23.5){\line(0,1){2}}
\put(290,0){\vector(0,1){10}}
\put(290,0){\circle*{2}}
\put(290,20){\circle*{2}}
\put(290,10){\vector(0,1){10}}
\put(290,10){\circle*{2}}
\put(298,5){$\cdots$}
\put(298,12){$\cdots$}
\put(320,0){\vector(0,1){10}}
\put(320,0){\circle*{2}}
\put(320,20){\circle*{2}}
\put(320,10){\vector(0,1){10}}
\put(320,10){\circle*{2}}
\put(295,-1.5){\vector(0,1){10}}
\put(295,-1.5){\circle*{2}}
\put(295,18.5){\circle*{2}}
\put(295,8.5){\vector(0,1){10}}
\put(295,8.5){\circle*{2}}
\put(315,-1.5){\vector(0,1){10}}
\put(315,-1.5){\circle*{2}}
\put(315,18.5){\circle*{2}}
\put(315,8.5){\vector(0,1){10}}
\put(315,8.5){\circle*{2}}
\put(290,0){\line(0,1){20}}
\put(320,0){\line(0,1){20}}
\qbezier(290,0)(305,-5)(320,0)
\qbezier(290,20)(305,15)(320,20)
\qbezier(290,0)(305,5)(320,0)
\qbezier(290,20)(305,25)(320,20)
\put(0,45){\line(0,-1){5}}
\put(155,45){\line(0,-1){5}}
\put(160,45){\line(0,-1){5}}
\put(325,45){\line(0,-1){5}}
\put(45,-15){\line(0,1){5}}
\put(50,-15){\line(0,1){5}}
\put(215,-15){\line(0,1){5}}
\put(220,-15){\line(0,1){5}}
\put(0,45){\line(1,0){155}}
\put(160,45){\line(1,0){165}}
\put(0,-15){\line(1,0){45}}
\put(50,-15){\line(1,0){165}}
\put(220,-15){\line(1,0){105}}
\put(-20,15){$\cdots$}
\put(330,15){$\cdots$}
\put(0,-25){$Rep_k\vec{Q}[1]$}
\put(120,-25){$Rep_k\vec{Q}$}
\put(260,-25){$Rep_k\vec{Q}[-1]$}
\put(70,50){$Coh(\mathbb{P}^1)$}
\put(220,50){$Coh(\mathbb{P}^1)[-1]$}
\end{picture}}
\end{equation}

\vspace{.5in}

In (\ref{E:pictureequiv}), we have used the fact that since $Coh(\mathbb{P}^1)$ (or $Rep_k\vec{Q}$) is of global dimension one, any complex $M$ of $D^b(Coh(\mathbb{P}^1))$ (or of $D^b(Rep_k\vec{Q})$) is quasi-isomorphic to a direct sum of stalk complexes $M \simeq \bigoplus_i H^i(M)[-i]$. Thus any indecomposable complex is of the form $I[l]$ for some indecomposable object $I$ in $Coh(\mathbb{P}^1)$ (or $D^b(Rep_k\vec{Q})$). In other words, we get a picture of these derived categories simply by putting along a line $\Z$ copies of the picture of the corresponding abelian category.

\vspace{.1in}

The explicit form of the above equivalence is as follows. Put $T=\O\oplus \O(1)$. Then $T$ generates $D^b(Coh(\mathbb{P}^1))$ as a triangulated category and ${Ext}^1(T,T)=\{0\}$. It is easy to check that $\Lambda={End}(T)$ is the path algebra of $\vec{Q}$.
 The main Theorem of tilting theory (see \cite{Happelbook}) thus gives an equivalence
$$RHom(T,\cdot)~: D^b(Coh(\mathbb{P}^1)) \stackrel{\sim}{\to} D^b({Mod}\text{-}\Lambda) = D^b(Rep_k(\vec{Q})).$$

\vspace{.1in}

The existence of this derived equivalence has a some immediate important consequences~: we obtain an isomorphism of Grothendieck groups $K(Coh(\mathbb{P}^1)) \stackrel{\sim}{\to} K(Rep_k\vec{Q})$ compatible with Euler forms; the indecomposables in $D^b(Coh(\mathbb{P}^1))$ are sent bijectively to the indecomposables in $D^b(Rep_k(\vec{Q}))$. What can we deduce concerning the classes of indecomposables in the Grothendieck group ? Since the class of an indecomposable $I[l]$ is $\overline{I[l]}=(-1)^l\overline{I}$, we see that, under the above isomorphism of Grothendieck groups,
\begin{equation*}
\begin{split}
\pm \{ \overline{M}\;|\; M \;\text{indecomposable\;in}& \; D^b(Coh(\mathbb{P}^1))\}\\
&=\pm \{ \overline{M}\;|\; M \;\text{indecomposable\;in} \; D^b(Rep_k\vec{Q})\}
\end{split}
\end{equation*}

Indeed, if we plot the classes of indecomposables in $Coh(\mathbb{P}^1)$ then we obtain the following set in $K(Coh(\mathbb{P}^1))=\Z^2$~:

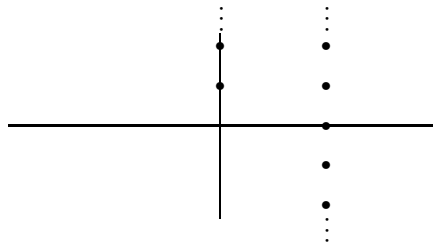
\centerline{
\begin{picture}(400, 60)
\multiput(179,34)(0,4){3}{$\cdot$}
\multiput(219,-46)(0,4){3}{$\cdot$}
\multiput(219,34)(0,4){3}{$\cdot$}
\multiput(220,-30)(0,15){5}{\circle*{3}}
\multiput(180,15)(0,15){2}{\circle*{3}}
\put(180,-35){\line(0,1){70}}
\put(100,0){\line(1,0){160}}
\end{picture}}

\vspace{.9in}

On the other hand, if we plot the indecomposables in $Rep_k\vec{Q}$ then we obtain the following set~:

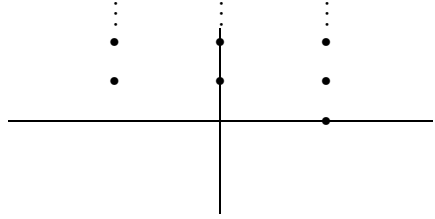
\centerline{
\begin{picture}(400, 60)
\multiput(179,34)(0,4){3}{$\cdot$}
\multiput(139,34)(0,4){3}{$\cdot$}
\multiput(219,34)(0,4){3}{$\cdot$}
\multiput(220,0)(0,15){3}{\circle*{3}}
\multiput(180,15)(0,15){2}{\circle*{3}}
\multiput(140,15)(0,15){2}{\circle*{3}}
\put(180,-35){\line(0,1){70}}
\put(100,0){\line(1,0){160}}
\end{picture}}

\vspace{.7in}

By Theorem~\ref{T:Godknows} and Corollary~\ref{C:KacP1}, we see that the sets of indecomposables
in $Coh(\mathbb{P}^1)$ and $Rep_k\vec{Q}$ correspond to two \textit{distinct} (and even \textit{nonconjugate}) sets of positive roots for the root system $\widehat{\Delta}$ of $\widehat{\mathfrak{sl}}_2$. These in turn are attached to two distinct and \textit{nonconjugate} Borel subalgebras $\widehat{\bo}_+$ and $\mathcal{L}\bo_+$ of $\widehat{\mathfrak{sl}}_2$ (see Appendix~A.2, Example A.15.). Alternatively, we can say that the appearance of the same root system $\widehat{\Delta}$ and Kac-Moody algebra $\widehat{\mathfrak{sl}}_2$ in both contexts is explained by the existence of a derived equivalence
$D^b(Coh(\mathbb{P}^1)) \stackrel{\sim}{\to} D^b(Rep_k\vec{Q})$.

\vspace{.1in}

Based on this observation it seems natural to consider the derived category $\mathcal{D}=D^b(Coh(\mathbb{P}^1))\simeq D^b(Rep_k\vec{Q})$ as ``representing'' the \textit{whole} Lie algebra $\widehat{\mathfrak{sl}}_2$, in the same way as the abelian categories $Coh(\mathbb{P}^1)$ and $Rep_k\vec{Q}$ ``represent'' some Borel subalgebras of $\widehat{\mathfrak{sl}}_2$. To be more accurate, one should consider the $2$-\textit{periodic} derived category $\mathcal{D}/T^2$ rather than the whole $\mathcal{D}$ since $\widehat{\mathfrak{sl}}_2$ is (roughly) made up of two copies of $\widehat{\bo}_+$ or $\mathcal{L}\bo_+$. Recall from \cite{PX0} that $\mathcal{D}/T^2$ is a triangulated category whose objects are (necessarily unbounded) complexes $P$ equipped with an isomorphism $P \stackrel{\sim}{\to} T^2(P)$, where $T$ stands here for the shift (suspension) functor.

\vspace{.1in}

At this point, it is very tempting to try to define a Hall algebra for the derived category $\mathcal{D}$ and its $2$-periodic version $\mathcal{D}/T^2$, in the hope of obtaining in this fashion a realization of the whole quantum group $\U_v(\widehat{\mathfrak{sl}}_2)$ rather than of its various ``positive parts''. There are in this direction several important partial results, due to Peng and Xiao, Kapranov and To\"en respectively, which we survey in Sections~5.5 and 5.6. Another approach is to recover the whole quantum group $\U_v(\widehat{\mathfrak{sl}}_2)$ from its positive parts $\U_v(\widehat{\bo}_+)$ or $\U_v(\mathcal{L}\bo_+)$ in a purely algebraic manner (by the so-called \textit{Drinfeld double} construction), and to relate properties of $\mathcal{D}$ and $\mathcal{D}/T^2$ to this algebra.
This is the line we follow in Sections~5.2--5.4.

\vspace{.2in}

Of course, the pair consisting of $Coh(\mathbb{P}^1)$ and $Rep_k\vec{Q}$ for $\vec{Q}$ the Kronecker quiver is far from being the only example around of derived equivalent hereditary categories. By a Theorem of Geigle and Lenzing \cite{GL}, any weighted projective line $\xpl$ of \textit{parabolic type} is derived equivalent to $Rep_k\vec{Q}$ for an appropriate \textit{tame} quiver. In this case, one may construct a derived equivalence just as was done above for $Coh(\mathbb{P}^1)$, but using the object
$T=\bigoplus_{0 \leq \vec{x} \leq \vec{c}} \O(\vec{x})$ instead. Pictorially, this is represented by~:

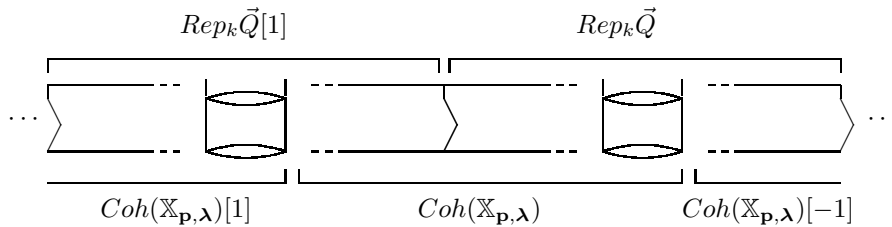
\vspace{.6in}

\centerline{
\begin{picture}(300,20)
\put(0,0){\line(1,0){40}}
\put(110,0){\line(1,0){40}}
\put(0,25){\line(1,0){40}}
\put(110,25){\line(1,0){40}}
\put(42.5,0){\line(1,0){2.5}}
\put(47.5,0){\line(1,0){2.5}}
\put(42.5,25){\line(1,0){2.5}}
\put(47.5,25){\line(1,0){2.5}}
\put(100,0){\line(1,0){2.5}}
\put(105,0){\line(1,0){2.5}}
\put(100,25){\line(1,0){2.5}}
\put(105,25){\line(1,0){2.5}}
\put(60,0){\line(0,1){20}}
\put(90,0){\line(0,1){20}}
\put(60,21){\line(0,1){2}}
\put(60,23){\line(0,1){2}}
\put(60,25){\line(0,1){2}}
\put(90,21){\line(0,1){2}}
\put(90,23){\line(0,1){2}}
\put(90,25){\line(0,1){2}}
\put(0,0){\line(1,2){5}}
\put(5,10){\line(-1,2){5}}
\put(0,20){\line(0,1){5}}
\put(150,0){\line(1,2){5}}
\put(155,10){\line(-1,2){5}}
\put(150,20){\line(0,1){5}}
\qbezier(60,0)(75,-5)(90,0)
\qbezier(60,20)(75,15)(90,20)
\qbezier(60,0)(75,5)(90,0)
\qbezier(60,20)(75,25)(90,20)
\put(150,0){\line(1,0){40}}
\put(260,0){\line(1,0){40}}
\put(150,25){\line(1,0){40}}
\put(260,25){\line(1,0){40}}
\put(192.5,0){\line(1,0){2.5}}
\put(197.5,0){\line(1,0){2.5}}
\put(192.5,25){\line(1,0){2.5}}
\put(197.5,25){\line(1,0){2.5}}
\put(250,0){\line(1,0){2.5}}
\put(255,0){\line(1,0){2.5}}
\put(250,25){\line(1,0){2.5}}
\put(255,25){\line(1,0){2.5}}
\put(210,0){\line(0,1){20}}
\put(240,0){\line(0,1){20}}
\put(210,21){\line(0,1){2}}
\put(210,23){\line(0,1){2}}
\put(210,25){\line(0,1){2}}
\put(240,21){\line(0,1){2}}
\put(240,23){\line(0,1){2}}
\put(240,25){\line(0,1){2}}
\put(150,0){\line(1,2){5}}
\put(155,10){\line(-1,2){5}}
\put(150,20){\line(0,1){5}}
\put(300,0){\line(1,2){5}}
\put(305,10){\line(-1,2){5}}
\put(300,20){\line(0,1){5}}
\qbezier(210,0)(225,-5)(240,0)
\qbezier(210,20)(225,15)(240,20)
\qbezier(210,0)(225,5)(240,0)
\qbezier(210,20)(225,25)(240,20)
\put(-15,10){$\cdots$}
\put(310,10){$\cdots$}
\put(0,35){\line(1,0){148}}
\put(0,35){\line(0,-1){5}}
\put(148,35){\line(0,-1){5}}
\put(152,35){\line(1,0){148}}
\put(152,35){\line(0,-1){5}}
\put(300,35){\line(0,-1){5}}
\put(0,-12){\line(1,0){90}}
\put(90,-12){\line(0,1){5}}
\put(95,-12){\line(1,0){145}}
\put(95,-12){\line(0,1){5}}
\put(240,-12){\line(0,1){5}}
\put(245,-12){\line(0,1){5}}
\put(245,-12){\line(1,0){55}}
\put(50,45){$Rep_k\vec{Q}[1]$}
\put(200,45){$Rep_k\vec{Q}$}
\put(20,-25){$Coh(\xpl)[1]$}
\put(140,-25){$Coh(\xpl)$}
\put(240,-25){$Coh(\xpl)[-1]$}
\end{picture}}

\vspace{.4in}

This explains the occurence of affine Lie algebras $\widehat{\g}$ in both contexts, and matches the picture of regular indecomposables for tame quivers and the picture for torsion sheaves in weighted projective lines. It also provides a conceptual explanation for the existence of a total order on the set of indecomposables of a quiver of finite type satisfying the conditions stated in Lemma~\ref{L:orderfinite} (it suffices to choose any order refining the order given by the slope function).
Other examples are given in Section~5.4.

\vspace{.2in}

\centerline{\textbf{5.2. The Drinfeld double.}} 
\addcontentsline{toc}{subsection}{\tocsubsection {}{}{\; 5.2. The Drinfeld double.}}

\vspace{.15in}

In this Section we will describe a fundamental algebraic construction due to Drinfeld which allows one to construct many new Hopf algebras from old ones.
Let $\H$ be a finite-dimensional Hopf algebra, and let $\H_{coop}^*$ be its dual Hopf algebra, with opposite \textit{coproduct}. 

\vspace{.1in}

\begin{theo}[Drinfeld] There exists a unique Hopf algebra structure on the vector space $\mathbf{D}\H=\H \otimes \H^*_{coop}$ such that
\begin{enumerate}
\item[i)] $\H$ and $\H_{coop}^*$ are Hopf subalgebras of $\mathbf{D}\H$,
\item[ii)] If $h \in \H$ and $h' \in \H^*_{coop}$ then $h \cdot h'=h \otimes h' \in \mathbf{D}\H$,
\item[iii)] The natural pairing $\{\;,\;\}$ on $\mathbf{D}\H$ is a Hopf pairing.
\end{enumerate}
Moreover, for $h \in \H$ and $h' \in \H^*_{coop}$ we have
$$h' \cdot h=\sum \{h'_{(1)},h_{(3)}\}\{S^{-1}(h'_{(3)}),h_{(1)}\}h_{(2)}h'_{(2)}.$$
\end{theo}

In the above, we have used Sweedler's notation for the coproduct of Hopf algebras, namely
$\Delta(x)=\sum x_{(1)} \otimes x_{(2)}$ and $\Delta^2(x)=\sum x_{(1)} \otimes x_{(2)} \otimes x_{(3)}$.

\vspace{.1in}

Now if in addition $\H$ is self-dual, i.e. that $\H \simeq \H^*$ then we may identify $\mathbf{D}\H$ with $\H \otimes \H$ and the above construction endows $\H \otimes \H$ with the structure of a Hopf algebra. This happens for instance when $\H$ is itself equipped with a nondegenerate Hopf pairing
$(\;,\;)$ (not to be confused with the natural pairing $\{\;,\;\}$ on $\mathbf{D}\H$ ). In that situation we may rephrase Drinfeld's construction without mentioning $\H^*$ at all~:

\vspace{.1in}

\begin{cor}\label{C:Dridri} Let $\H$ be a finite-dimensional Hopf algebra with a nondegenerate Hopf pairing. Let $\H^+, \H^-$ be two isomorphic copies of $\H$, and equip $\H^-$ with the opposite coproduct. Let $\mathbf{D}$ be the Hopf algebra generated by $\H^+$ and $\H^-$ modulo the relations
\begin{enumerate}
\item[i)] $\H^{\pm}$ are Hopf subalgebras,
 \item[ii)] For $h^\pm \in \H^{\pm}$ we have
 $$h^- \cdot h^+=\sum  (h^-_{(1)},h^+_{(3)})(S^{-1}(h^-_{(3)}),h^+_{(1)})h^+_{(2)}h^-_{(2)}.$$
\end{enumerate}
Then the multiplication map $\H^+ \otimes \H^- \to \mathbf{D}$ is an isomorphism of vector spaces.
\end{cor}

\vspace{.1in}

The advantage of this reformulation is twofold~: we may apply it to an infinite-dimensional Hopf algebra (even a topological Hopf algebra), and we may also apply it when $\H$ is equipped with a Hopf pairing which is allowed to be degenerate. This produces in all cases a topological Hopf algebra structure on the vector space $\mathbf{D}\H=\H \otimes \H$. 

\vspace{.15in}

One of the motivations behind Drinfeld's construction is its implication for quantum groups, as the following Proposition demonstrates. Let $\g$ be a Kac-Moody algebra and let $\U_v(\bo'_+)$ be the quantized enveloping algebra of a Borel subalgebra $\bo' \subset \g'$ (see Appendix~A.2). This is a Hopf subalgebra equipped with a (usually degenerate) Hopf pairing.

\vspace{.1in}

\begin{prop}[Drinfeld]  Let $\mathbf{K}=\C(v)[K_i^{\pm 1}]_{i \in I} \subset \U_v(\bo_+)$ be the quantized enveloping algebra of the Cartan subalgebra $\h \subset \bo'_+$.
Then $K_i \otimes K_i^{-1} \in \mathbf{D}\U_v(\bo'_+)$ are central elements for all $i$, and there is an isomorphism of Hopf algebras
$$\mathbf{D}\U_v(\bo'_+)/\langle K_i \otimes K_i^{-1} -1 \rangle_i \simeq \U_v(\g').$$
\end{prop}

\vspace{.1in}

In other words, the Drinfeld double is exactly the gadget needed to pass from the positive part
$\U_v(\bo'_+)$ to the whole quantum group $\U_v(\g)$. It is important to note that this is a purely \textit{quantum} phenomenon-- the classical enveloping algebras $\U(\g)$ do not possess natural Hopf pairings with good properties, and are not self-dual. 

\vspace{.1in}

Extended Hall algebra $\widetilde{\H}_{\A}$ of some hereditary finitary category $\A$ form another wide class of Hopf algebras which we may force through Corollary~\ref{C:Dridri}. Given the similarility between quantum groups and Hall algebras, we define, following Xiao \cite{Xiao1} the \textit{reduced} Drinfeld double $\widetilde{\mathbf{D}}{\H}_{\A}$ to be the quotient of $\mathbf{D}\widetilde{\H}_{\A}$ by the ideal generated by the elements $(\mathbf{k}_{\mathcal{F}} \otimes \mathbf{k}_{\mathcal{F}}^{-1}-1) \in \widetilde{\H}_\A \otimes \widetilde{\H}_{\A}$ for all $\mathcal{F} \in \A$. It is easy to check that $\widetilde{\mathbf{D}}\H_{\A}$ admits a triangular decomposition
$$\widetilde{\mathbf{D}}\H_{\A}=\H_{\A} \otimes \mathbf{K}  \otimes \H_{\A}$$
where $\mathbf{K}=\C[K(\A)]$.

\vspace{.2in}

\centerline{\textbf{5.3. Conjectures and Cramer's Theorem.}}
\addcontentsline{toc}{subsection}{\tocsubsection {}{}{\; 5.3. Conjectures and Cramer's theorem.}}

\vspace{.2in}

We are now ready to formulate certain conjectures relating the structure of the (reduced) Drinfeld double $\widetilde{\mathbf{D}}{\H}_A$ of a finitary \textit{hereditary} category $\A$ with that of the derived category $D^b(\A)$ and of the 2-periodic derived category $D^b(\A)/T^2$. Throughout, we let $\mathcal{C}_1, \mathcal{C}_2$ be finitary hereditary connected categories, which we suppose in addition to be $\mathbb{F}_q$-linear for simplicity. An exact functor $F~:D^b(\mathcal{C}_1) \to D^b(\mathcal{C}_2)$ will be called \textit{friendly} if there exists an integer $n$ such that
$$F(\mathcal{C}_1) \subset \mathcal{C}_2[-n] \oplus \mathcal{C}_2[-1-n],$$
i.e. for any object $M$ of $\mathcal{C}_1$, the complex $F(M)$ is concentrated in degrees $n$ and $n+1$. Examples of unfriendly exact functors are known when $\mathcal{C}_1=Rep_kA_n$ for instance. However, it is believed that if $\mathcal{C}_1$ and $\mathcal{C}_2$ are not of finite type then \textit{any} derived equivalence is friendly (\cite{Lenzingoral}).

\vspace{.1in}

\begin{conj}\label{Conj:1} If $\mathcal{C}_1$ and $\mathcal{C}_2$ are derived equivalent then
$\widetilde{\mathbf{D}}{\H}_{\mathcal{C}_1}$ and $\widetilde{\mathbf{D}}{\H}_{\mathcal{C}_2}$ are isomorphic \textit{as algebras}.
\end{conj}

\vspace{.1in}

For any object $L$ of a hereditary category $\A$, we set 
$$[L]^+=[L] \otimes 1 \in \widetilde{\mathbf{D}}\H_{\A}, \qquad [L]^-=1 \otimes [L] \in \widetilde{\mathbf{D}}\H_{\A}.$$

\vspace{.1in}

\begin{conj}\label{Conj:2} Let $F: D^b(\mathcal{C}_1) \to D^b(\mathcal{C}_2)$ be an exact, fully faithful and friendly functor. For any object $M$ of $\mathcal{C}_1$, define $n \in \Z$ and objects $N, N'$ of $\mathcal{C}_2$ such that $F(M) \simeq N[-n] \oplus N'[-1-n]$. Then, setting $\epsilon(n)=-(1)^n$, the assignement
\begin{equation}\label{E:juki1}
\mathbf{k}_{\mathcal{F}} \mapsto \mathbf{k}_{N}^{\epsilon(n)}\mathbf{k}_{N'}^{-\epsilon(n)},
\end{equation}
\begin{equation}\label{E:juki2}
[M]^{\pm} \mapsto \nu^{\mathbf{d}(N,N',n)} [N]^{\pm \epsilon(n)} \mathbf{k}^{\pm n \epsilon(n)}_{N} \cdot [N']^{\mp \epsilon(n)} \mathbf{k}_{N'}^{\mp (n+1)\epsilon(n)}
\end{equation}
with $\mathbf{d}(N,N',n)=n\langle N,N\rangle_a+(n+1)\langle N',N'\rangle_a-\langle N,N'\rangle_a$,
extends to an embedding of reduced Drinfeld doubles $F_*~: \widetilde{\mathbf{D}}\H_{\mathcal{C}_1} \hookrightarrow \widetilde{\mathbf{D}}\H_{\mathcal{C}_2}$. This embedding is an isomorphism if $F$ is a derived equivalence.\end{conj}

\vspace{.1in}

Now let us assume that the symmetrized Euler form $(\;,\;)_a$ of a hereditary category $\A$ vanishes identically (the so-called \textit{Calabi-Yau} case). In that situation, the Hall algebra $\H_\A$ is already a Hopf algebra and there is no need to extend it by $\C[K(\A)]$ (see Remark~1.11.). Thus it is legal to consider the (genuine) Drinfeld double $\mathbf{D}\H_{\A}$ of the Hall algebra of such a category.
 
\vspace{.1in}
 
\begin{conj}\label{Conj:3} Assume that the symmetrized Euler forms of both $\mathcal{C}_1$ and $\mathcal{C}_2$ vanish identically.  Let $F: D^b(\mathcal{C}_1)/T^2 \to D^b(\mathcal{C}_2)/T^2$ be an exact, fully faithful functor. For any object $M$ of $\mathcal{C}_1$ define objects $N,N'$ of $\mathcal{C}_2$ such that $F(M) \simeq N \oplus N'[-1]$. Then the map
$$[M]^{\pm} \mapsto  \nu^{-\langle N,N'\rangle_a}[N]^{\pm} \cdot [N']^{\mp}$$
extends to an embedding of Drinfeld doubles $F_*~: {\mathbf{D}}\H_{\mathcal{C}_1} \hookrightarrow {\mathbf{D}}\H_{\mathcal{C}_2}$. This embedding is an isomorphism if $F$ is a derived equivalence.\end{conj}

\vspace{.1in}

To finish, we mention the following weaker conjecture, which is a direct consequence of Conjectures~\ref{Conj:2} and \ref{Conj:3}~:

\vspace{.1in}

\begin{conj}\label{Conj:4} Let $\mathcal{C}$ be any hereditary and finitary category. If any autoequivalence of $D^b(\mathcal{C})$ is friendly then the group $Aut(D^b(\mathcal{C}))$ acts on the reduced Drinfeld double $\widetilde{\mathbf{D}}\H_{\mathcal{C}}$ by algebra automorphisms. If the symmetrized Euler form $(\;,\;)_a$ of $\mathcal{C}$ vanishes then the group $Aut(D^b(\mathcal{C})/T^2)$ acts by algebra automorphisms on ${\mathbf{D}}\H_{\mathcal{C}}$.
\end{conj}

\vspace{.1in}

In a recent preprint \cite{Cramer}, T. Cramer made a big step towards resolving the above conjectures. Namely,

\vspace{.1in}

\begin{theo}[Cramer]
Conjectures~\ref{Conj:1}, \ref{Conj:2} hold under the assumption that one of $\mathcal{C}_1, \mathcal{C}_2$ satisfies the finite subobjects condition.
\end{theo}

\vspace{.2in}

\centerline{\textbf{5.4. Examples and applications.}} 
\addcontentsline{toc}{subsection}{\tocsubsection {}{}{\; 5.4. Example and applications.}}

\vspace{.15in}

In this Section we illustrate the principles explained above and give some evidence for  Conjectures~\ref{Conj:1}--\ref{Conj:4}. 

\vspace{.1in}

The first type of applications is related to derived equivalences between categories of representations of quivers with the same underlying graph but different orientations. Let $\vec{Q}$ be an arbitrary quiver and suppose that the vertex $i$ of $\vec{Q}$ is a sink, i.e. there are no oriented edges going into $i$. Let $s_i\vec{Q}$ be the quiver obtained from $\vec{Q}$ by reversing all the arrows leaving $i$ (so that $i$ is now a source of $s_i\vec{Q}$). To the quivers $\vec{Q}$ and $s_i\vec{Q}$ is attached the same Kac-Moody Lie algebra $\g$. Bernstein, Gelfand and Ponomaraev discovered in
\cite{BGP}  the so-called \textit{reflection functor}, which is a derived equivalence
$$\mathcal{S}_i~: D^b(\vec{Q}) \stackrel{\sim}{\to} D^b(s_i\vec{Q}).$$ 
Observe that in $Rep_k\vec{Q}$ we have $Ext^1(S_i,I)=\{0\}$ for any indecomposable object $I$ different from $S_i$, whereas in $Rep_ks_i\vec{Q}$, $Ext^1(I,S_i)=\{0\}$ for any indecomposable $I$ distinct from $S_i$. The equivalence $\mathcal{S}_i$ may be represented as~:

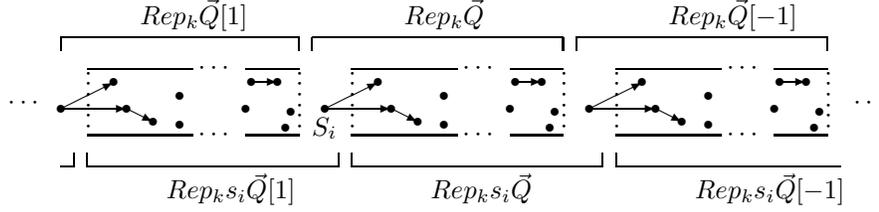
\vspace{.3in}

\centerline{
\begin{picture}(300,40)
\put(0,32){\line(0,1){5}}
\put(0,37){\line(1,0){90}}
\put(90,37){\line(0,-1){5}}
\put(95,32){\line(0,1){5}}
\put(95,37){\line(1,0){95}}
\put(190,37){\line(0,-1){5}}
\put(195,32){\line(0,1){5}}
\put(195,37){\line(1,0){95}}
\put(290,37){\line(0,-1){5}}
\put(0,-12){\line(1,0){5}}
\put(5,-12){\line(0,1){5}}
\put(10,-12){\line(0,1){5}}
\put(10,-12){\line(1,0){95}}
\put(105,-12){\line(0,1){5}}
\put(110,-12){\line(0,1){5}}
\put(110,-12){\line(1,0){95}}
\put(205,-12){\line(0,1){5}}
\put(210,-12){\line(0,1){5}}
\put(210,-12){\line(1,0){85}}
\put(30,42){$Rep_k\vec{Q}[1]$}
\put(130,42){$Rep_k\vec{Q}$}
\put(230,42){$Rep_k\vec{Q}[-1]$}
\put(40,-24){$Rep_ks_i\vec{Q}[1]$}
\put(140,-24){$Rep_ks_i\vec{Q}$}
\put(240,-24){$Rep_ks_i\vec{Q}[-1]$}
\put(10,0){\line(1,0){40}}
\put(10,25){\line(1,0){40}}
\multiput(9,1)(0,4){6}{$\cdot$}
\put(52,-2){$\cdots$}
\put(52,23){$\cdots$}
\put(70,0){\line(1,0){20}}
\put(70,25){\line(1,0){20}}
\multiput(89,1)(0,4){6}{$\cdot$}
\put(0,10){\circle*{3}}
\put(0,10){\vector(1,0){25}}
\put(0,10){\vector(2,1){20}}
\put(25,10){\circle*{3}}
\put(35,5){\circle*{3}}
\put(25,10){\vector(2,-1){10}}
\put(72,20){\circle*{3}}
\put(20,20){\circle*{3}}
\put(45,4){\circle*{3}}
\put(70,10){\circle*{3}}
\put(72,20){\vector(1,0){10}}
\put(85,3){\circle*{3}}
\put(45,15){\circle*{3}}
\put(87,9){\circle*{3}}
\put(82,20){\circle*{3}}
\put(110,0){\line(1,0){40}}
\put(110,25){\line(1,0){40}}
\multiput(109,1)(0,4){6}{$\cdot$}
\put(152,-2){$\cdots$}
\put(152,23){$\cdots$}
\put(170,0){\line(1,0){20}}
\put(170,25){\line(1,0){20}}
\multiput(189,1)(0,4){6}{$\cdot$}
\put(95,0){$S_i$}
\put(100,10){\circle*{3}}
\put(100,10){\vector(1,0){25}}
\put(100,10){\vector(2,1){20}}
\put(125,10){\circle*{3}}
\put(135,5){\circle*{3}}
\put(125,10){\vector(2,-1){10}}
\put(172,20){\circle*{3}}
\put(120,20){\circle*{3}}
\put(145,4){\circle*{3}}
\put(170,10){\circle*{3}}
\put(172,20){\vector(1,0){10}}
\put(185,3){\circle*{3}}
\put(145,15){\circle*{3}}
\put(187,9){\circle*{3}}
\put(182,20){\circle*{3}}
\put(210,0){\line(1,0){40}}
\put(210,25){\line(1,0){40}}
\multiput(209,1)(0,4){6}{$\cdot$}
\put(252,-2){$\cdots$}
\put(252,23){$\cdots$}
\put(270,0){\line(1,0){20}}
\put(270,25){\line(1,0){20}}
\multiput(289,1)(0,4){6}{$\cdot$}
\put(200,10){\circle*{3}}
\put(200,10){\vector(1,0){25}}
\put(200,10){\vector(2,1){20}}
\put(225,10){\circle*{3}}
\put(235,5){\circle*{3}}
\put(225,10){\vector(2,-1){10}}
\put(272,20){\circle*{3}}
\put(220,20){\circle*{3}}
\put(245,4){\circle*{3}}
\put(270,10){\circle*{3}}
\put(272,20){\vector(1,0){10}}
\put(285,3){\circle*{3}}
\put(245,15){\circle*{3}}
\put(287,9){\circle*{3}}
\put(282,20){\circle*{3}}
\put(-20,10){$\cdots$}
\put(300,10){$\cdots$}
\end{picture}}

\vspace{.5in}

This equivalence is friendly in the sense of the preceding Section.
The functor $\mathcal{S}_i$ descends at the level of Grothendieck groups to an isomorphism~:
$$S_i~:K(D^b(Rep_k\vec{Q}) \stackrel{\sim}{\to} K(D^b(Rep_ks_i\vec{Q}))$$
which, after the identifications of these Grothendieck groups with the root lattice $\bigoplus_i\Z\a_i$ of $\g$, coincides with the simple reflection $S_i$ in the Weyl group $W$ of $\g$.
By Ringel's Theorem~\ref{T:RingelHall}, there are embeddings 
$$\Psi~: \U_\nu(\bo'_+) \hookrightarrow \widetilde{\H}_{\vec{Q}},$$ 
$$\Psi'~: \U_\nu(\bo'_+) \hookrightarrow \widetilde{\H}_{s_i\vec{Q}}.$$
Cramer's theorem implies the existence of an isomorphism of reduced Drinfeld doubles
$$\mathcal{S}_{i,\star}~: \widetilde{\mathbf{D}}\H_{\vec{Q}} \stackrel{\sim}{\to} \widetilde{\mathbf{D}}\H_{s_i\vec{Q}},$$
given by the formulas~:
$$\mathcal{S}_{i,\star}([S_i]^+)=\nu [S_i]^-\mathbf{k}_i^{-1}, \qquad \mathcal{S}_{i,\star}([S_i]^-)=\nu [S_i]^+ \mathbf{k}_i, \qquad \mathcal{S}_{i,\star}(\mathbf{k}_i)=\mathbf{k}_i^{-1},$$
$$\mathcal{S}_{i,\star}([R]^{\pm})=[\mathcal{S}_i(R)]^{\pm}, \qquad \text{for\;any\;indecomposable}\; R.$$
This restricts to an isomorphism of quantum groups
$$\mathcal{S}_{i,\star}~: \U_\nu(\g') \stackrel{\sim}{\to} \U_\nu(\g')$$
satisfying (among other relations)
\begin{equation}\label{E:Si}
\mathcal{S}_{i,\star}(E_i)=\nu F_i{K}_i^{-1}, \qquad \mathcal{S}_{i,\star}(F_i)=\nu E_i {K}_i, \qquad \mathcal{S}_{i,\star}({K}_i)={K}_i^{-1}.
\end{equation}
Sevenhant and Van den Bergh in \cite{SVdB2}, and later Xiao and Yang in \cite{XiaoBGP}, checked this directly some time ago.
The formulas (\ref{E:Si}) were first written down by Lusztig (see \cite{Lu1}), but they acquire under this new light a very clear conceptual meaning. 

\vspace{.05in}

There are, of course, completely similar results for the reflection functors 
$$\mathcal{S}_i~: D^b(Rep_k\vec{Q}) \stackrel{\sim}{\to} D^b(Rep_ks_i\vec{Q})$$
associated to a source $i$ rather than to a sink of the quiver $\vec{Q}$. In that situation, the pictorial representation of the equivalence id as follows~: 

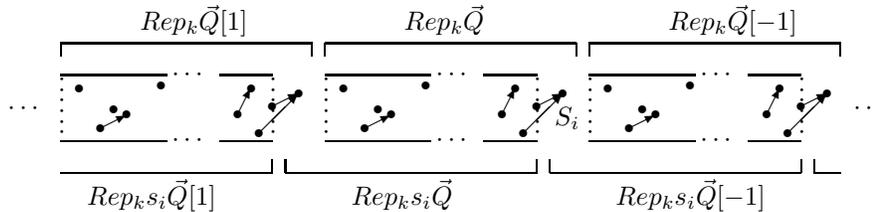
\vspace{.3in}

\centerline{
\begin{picture}(300,40)
\put(0,32){\line(0,1){5}}
\put(0,37){\line(1,0){95}}
\put(95,37){\line(0,-1){5}}
\put(100,32){\line(0,1){5}}
\put(100,37){\line(1,0){95}}
\put(195,37){\line(0,-1){5}}
\put(200,32){\line(0,1){5}}
\put(200,37){\line(1,0){95}}
\put(295,37){\line(0,-1){5}}
\put(0,-12){\line(1,0){80}}
\put(80,-12){\line(0,1){5}}
\put(85,-12){\line(0,1){5}}
\put(85,-12){\line(1,0){95}}
\put(180,-12){\line(0,1){5}}
\put(185,-12){\line(0,1){5}}
\put(185,-12){\line(1,0){95}}
\put(280,-12){\line(0,1){5}}
\put(285,-12){\line(0,1){5}}
\put(285,-12){\line(1,0){10}}
\put(30,42){$Rep_k\vec{Q}[1]$}
\put(130,42){$Rep_k\vec{Q}$}
\put(230,42){$Rep_k\vec{Q}[-1]$}
\put(10,-24){$Rep_ks_i\vec{Q}[1]$}
\put(110,-24){$Rep_ks_i\vec{Q}$}
\put(210,-24){$Rep_ks_i\vec{Q}[-1]$}
\put(0,0){\line(1,0){40}}
\put(0,25){\line(1,0){40}}
\multiput(-1,1)(0,4){6}{$\cdot$}
\put(42,-2){$\cdots$}
\put(42,23){$\cdots$}
\put(60,0){\line(1,0){20}}
\put(60,25){\line(1,0){20}}
\multiput(79,1)(0,4){6}{$\cdot$}
\put(90,18){\circle*{3}}
\put(15,5){\circle*{3}}
\put(25,10){\circle*{3}}
\put(15,5){\vector(2,1){10}}
\put(7,20){\circle*{3}}
\put(20,12){\circle*{3}}
\put(38,21){\circle*{3}}
\put(67,10){\circle*{3}}
\put(67,10){\vector(1,2){5}}
\put(75,3){\circle*{3}}
\put(75,3){\vector(1,1){15}}
\put(80,13){\circle*{3}}
\put(72,20){\circle*{3}}
\put(80,13){\vector(2,1){10}}
\put(100,0){\line(1,0){40}}
\put(100,25){\line(1,0){40}}
\multiput(99,1)(0,4){6}{$\cdot$}
\put(142,-2){$\cdots$}
\put(142,23){$\cdots$}
\put(160,0){\line(1,0){20}}
\put(160,25){\line(1,0){20}}
\multiput(179,1)(0,4){6}{$\cdot$}
\put(190,18){\circle*{3}}
\put(187,6){$S_i$}
\put(115,5){\circle*{3}}
\put(125,10){\circle*{3}}
\put(115,5){\vector(2,1){10}}
\put(107,20){\circle*{3}}
\put(120,12){\circle*{3}}
\put(138,21){\circle*{3}}
\put(167,10){\circle*{3}}
\put(167,10){\vector(1,2){5}}
\put(175,3){\circle*{3}}
\put(175,3){\vector(1,1){15}}
\put(180,13){\circle*{3}}
\put(172,20){\circle*{3}}
\put(180,13){\vector(2,1){10}}
\put(200,0){\line(1,0){40}}
\put(200,25){\line(1,0){40}}
\multiput(199,1)(0,4){6}{$\cdot$}
\put(242,-2){$\cdots$}
\put(242,23){$\cdots$}
\put(260,0){\line(1,0){20}}
\put(260,25){\line(1,0){20}}
\multiput(279,1)(0,4){6}{$\cdot$}
\put(290,18){\circle*{3}}
\put(215,5){\circle*{3}}
\put(225,10){\circle*{3}}
\put(215,5){\vector(2,1){10}}
\put(207,20){\circle*{3}}
\put(220,12){\circle*{3}}
\put(238,21){\circle*{3}}
\put(267,10){\circle*{3}}
\put(267,10){\vector(1,2){5}}
\put(275,3){\circle*{3}}
\put(275,3){\vector(1,1){15}}
\put(280,13){\circle*{3}}
\put(272,20){\circle*{3}}
\put(280,13){\vector(2,1){10}}
\put(-20,10){$\cdots$}
\put(300,10){$\cdots$}
\end{picture}}

\vspace{.5in}

Note that although all these functors $\mathcal{S}_i$ lift the action of the Weyl group on the root lattices to the level of derived categories, the corresponding automorphisms $\mathcal{S}_{i,\star}$ of the Drinfeld double $\U_\nu(\g')$ only give rise to an action of the \textit{braid group}. In other words, there is a slight loss of symmetry when one goes from the Grothendieck group to the (Drinfeld doubles of) Hall algebras.

\vspace{.2in}

We now turn to another class of examples of applications of the conjectures stated in the previous Section. Let $X$ be an elliptic curve defined over the field $k=\mathbb{F}_q$, and let $\H_X$ be its Hall algebra. Recall from Section~4.10. that the indecomposable objects of $Coh(X)$ are all semistable and that the subcategories $\mathcal{C}_\mu$ of semistables are all equivalent to one another for $\mu \in \mathbb{Q} \cup \{\infty\}$. Moreover, the symmetrized Euler form $(\;,\;)_a$ vanishes identically on the Grothendieck group. In that situation, motivated by Conjecture~\ref{Conj:4}, we consider the group of autoequivalences of the $2$-periodic category $D^b(Coh(X))/T^2$. We may represent this triangulated category as the result of gluing the right end of the abelian subcategory $Coh(X)[-1]$ to the left end of
$Coh(X)$, as depicted below~:

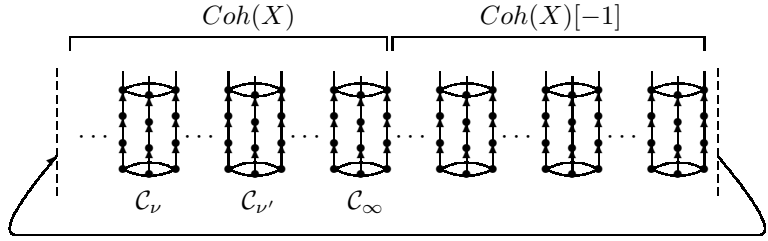
\vspace{.3in}

\centerline{
\begin{picture}(200,60)
\qbezier(0,0)(10,-5)(20,0)
\qbezier(0,0)(10,5)(20,0)
\put(5,-15){$\mathcal{C}_{\nu}$}
\qbezier(40,0)(50,-5)(60,0)
\qbezier(40,0)(50,5)(60,0)
\put(45,-15){$\mathcal{C}_{\nu'}$}
\qbezier(80,0)(90,-5)(100,0)
\qbezier(80,0)(90,5)(100,0)
\put(85,-15){$\mathcal{C}_{\infty}$}
\qbezier(0,30)(10,25)(20,30)
\qbezier(0,30)(10,35)(20,30)
\qbezier(40,30)(50,25)(60,30)
\qbezier(40,30)(50,35)(60,30)
\qbezier(80,30)(90,25)(100,30)
\qbezier(80,30)(90,35)(100,30)
\qbezier(120,0)(130,-5)(140,0)
\qbezier(120,0)(130,5)(140,0)
\qbezier(160,0)(170,-5)(180,0)
\qbezier(160,0)(170,5)(180,0)
\qbezier(200,0)(210,-5)(220,0)
\qbezier(200,0)(210,5)(220,0)
\qbezier(120,30)(130,25)(140,30)
\qbezier(120,30)(130,35)(140,30)
\qbezier(160,30)(170,25)(180,30)
\qbezier(160,30)(170,35)(180,30)
\qbezier(200,30)(210,25)(220,30)
\qbezier(200,30)(210,35)(220,30)
\qbezier(225,5)(250,-25)(240,-25)
\put(240,-25){\line(-1,0){280}}
\qbezier(-25,5)(-50,-25)(-40,-25)
\multiput(225,-10)(0,5){10}{\line(0,1){3}}
\multiput(-25,-10)(0,5){10}{\line(0,1){3}}
\put(-30,0){\vector(1,1){5}}
\multiput(0,0)(40,0){6}{\line(0,1){30}}
\multiput(20,0)(40,0){6}{\line(0,1){30}}
\multiput(-17,10)(40,0){6}{$\cdots$}
\multiput(0,0)(0,10){4}{\circle*{3}}
\multiput(0,0)(0,10){3}{\vector(0,1){10}}
\multiput(10,-2)(0,10){4}{\circle*{3}}
\multiput(10,-2)(0,10){3}{\vector(0,1){10}}
\multiput(20,0)(0,10){4}{\circle*{3}}
\multiput(20,0)(0,10){3}{\vector(0,1){10}}
\multiput(40,0)(0,10){4}{\circle*{3}}
\multiput(40,0)(0,10){3}{\vector(0,1){10}}
\multiput(50,-2)(0,10){4}{\circle*{3}}
\multiput(50,-2)(0,10){3}{\vector(0,1){10}}
\multiput(60,0)(0,10){4}{\circle*{3}}
\multiput(60,0)(0,10){3}{\vector(0,1){10}}
\multiput(80,0)(0,10){4}{\circle*{3}}
\multiput(80,0)(0,10){3}{\vector(0,1){10}}
\multiput(90,-2)(0,10){4}{\circle*{3}}
\multiput(90,-2)(0,10){3}{\vector(0,1){10}}
\multiput(100,0)(0,10){4}{\circle*{3}}
\multiput(100,0)(0,10){3}{\vector(0,1){10}}
\multiput(120,0)(0,10){4}{\circle*{3}}
\multiput(120,0)(0,10){3}{\vector(0,1){10}}
\multiput(130,-2)(0,10){4}{\circle*{3}}
\multiput(130,-2)(0,10){3}{\vector(0,1){10}}
\multiput(140,0)(0,10){4}{\circle*{3}}
\multiput(140,0)(0,10){3}{\vector(0,1){10}}
\multiput(160,0)(0,10){4}{\circle*{3}}
\multiput(160,0)(0,10){3}{\vector(0,1){10}}
\multiput(170,-2)(0,10){4}{\circle*{3}}
\multiput(170,-2)(0,10){3}{\vector(0,1){10}}
\multiput(180,0)(0,10){4}{\circle*{3}}
\multiput(180,0)(0,10){3}{\vector(0,1){10}}
\multiput(200,0)(0,10){4}{\circle*{3}}
\multiput(200,0)(0,10){3}{\vector(0,1){10}}
\multiput(210,-2)(0,10){4}{\circle*{3}}
\multiput(210,-2)(0,10){3}{\vector(0,1){10}}
\multiput(220,0)(0,10){4}{\circle*{3}}
\multiput(220,0)(0,10){3}{\vector(0,1){10}}
\multiput(0,30)(20,0){12}{\line(0,1){7}}
\multiput(10,28)(40,0){6}{\line(0,1){7}}
\put(-20,45){\line(0,1){5}}
\put(-20,50){\line(1,0){120}}
\put(100,45){\line(0,1){5}}
\put(102,45){\line(0,1){5}}
\put(102,50){\line(1,0){118}}
\put(220,45){\line(0,1){5}}
\put(30,55){$Coh(X)$}
\put(135,55){$Coh(X)[-1]$}
\end{picture}}

\vspace{.65in}

The derived category $D^b(Coh(X))$ and its $2$-periodic version posses many interesting autoequivalences, as was discovered by Mukai \cite{Mukai}. Let $\mathcal{E}$ be a \emph{spherical sheaf} in  $Coh(X)$, i.e.~an object satisfying ${Hom}(\mathcal{E}, \mathcal{E}) = {Ext}^1(\mathcal{E}, \mathcal{E}) = k$, for example the structure sheaf the curve $\O$ or 
the structure sheaf of a closed point $\O_{x_0}$ of degre one. Consider the functor 
$$
T_\mathcal{E}~: D^b(Coh(X)) \lto D^b(Coh(X))
$$
defined by 
$$ T_\mathcal{E}(\mathcal{F}) = \mathrm{cone}({{RHom}}(\mathcal{E}, \mathcal{F})\overset{k}\otimes \mathcal{E} \stackrel{ev}\lto \mathcal{F}).
$$
The functor $T_\mathcal{E}$  is exact  and induces an equivalence of derived categories for any spherical sheaf $\mathcal{E}$. Moreover, as shown by Seidel and Thomas (see \cite{ST}), there is a natural transformation of functors~:
$$T_{\O} T_{\O_{x_0}} T_{\O} \simeq T_{\O_{x_0}} T_{\O} T_{\O_{x_0}}$$
and if we set $\Phi=T_{\O} T_{\O_{x_0}} T_{\O}$ then there exists an involution $i$ of $X$ such that
$\Phi^2 \simeq i^* \circ T$, where $T$ is the translation $\mathcal{F} \mapsto \mathcal{F}[1]$.
This means that there is a group homomorphism from the braid group $B_3$ on three strands to the group $Aut(D^b(Coh(X)))$. This homomorphism is injective (see \cite{ST}). The braid group $B_3$ is isomorphic to the universal covering  $\widetilde{SL}(2,\Z)$ of $SL(2,\Z)$ given by a  central  extension of $SL(2,\Z)$ by $\Z$. Since in ${Aut}(D^b(Coh(X))/T^2)$ we have $[1]^2 \simeq
id$, the action of $B_3$ on the $2$-periodic category descends 
to an action of  $\widehat{SL}(2,\Z)$, where  $\widehat{SL}(2,\Z)$ is a two-fold covering of 
$SL(2,\Z)$. All these equivalences are friendly. Moreover, the action of $\widehat{SL}(2,\Z)$ on the Grothendieck group $K(Coh(X))$ factors through the quotient $K(Coh(X)) \tto \Z^2$ given by $\overline{\mathcal{F}} \mapsto (rank(\mathcal{F}),deg(\mathcal{F}))$, and coincides there with the tautological action of $SL(2,\Z)$ on $\Z^2$. All this may be summarized in the following commutative diagram~:
$$\xymatrix{
B_3\simeq\widetilde{SL}(2,\Z) \ar@{->>}[d] \ar@{^{(}->}[r] & \mathrm{Aut}(D^b(Coh(X))) \ar[d]\\
\widehat{SL}(2,\Z) \ar@{^{(}->}[r] \ar@{->>}[rd] & \mathrm{Aut}(D^b(Coh(X))/T^2) \ar@{->>}[d]\\
& SL(2,\Z)=\mathrm{Aut}(Z^2)}
$$

\vspace{.1in}

\begin{prop}[Burban-S., \cite{BS}]\label{P:BSauto} Formulas (\ref{E:juki1}) and (\ref{E:juki2}) define an action of $\widehat{SL}(2,\Z)$ on the reduced Drinfeld double $\widetilde{\mathbf{D}}\H_{X}$ of the Hall algebra of $X$.\end{prop}

\vspace{.15in}

Using the above Proposition, we were able in \cite{BS} to give an algebraic description of the reduced Drinfeld double $\widetilde{\mathbf{D}}\mathbf{C}_X$ of the spherical Hall algebra $\widetilde{\mathbf{C}}_X$ (see Section~4.10.). We keep the notations of Theorem~\ref{T:BS}, and in addition we set 
$$\epsilon_\x=\begin{cases} 1 & \text{if}\; \x \in (\Z^2)^+,\\ -1 & \text{if}\; -\x \in (\Z^2)^+.\end{cases}$$ 

\vspace{.1in}

\begin{theo}[Burban-S., \cite{BS}]\label{T:BSauto} The algebra $\widetilde{\mathbf{D}}\mathbf{C}_X$ is isomorphic to the the algebra $\U_X$ generated by elements $\{{t}_{(p,q)}\;|\; (p,q) \in (\Z^2)^*\}$ and $\{\mathbf{k}_{r,d}\;|\; r,d \in \Z\}$ subject to the following set of relations
\begin{enumerate}
\item[i)] $\mathbf{k}_{r,d} \mathbf{k}_{r',d'}=\mathbf{k}_{r+r',d+d'}$,
\item[i)] If $\x,\x'$ belong to the same line in $\Z^2$ then
$$[{t}_\x,{t}_{\x'}]=\delta_{\x,-\x'}c_{deg(\x)}\frac{\mathbf{k}_{\x}-\mathbf{k}_{x}^{-1}}{\nu-\nu^{-1}},$$
\item[ii)] Assume that $\x,\y \in \Z^2$ are such that $deg(\x)=1$ and that 
the triangle with vertices $0, \x$ and $\x+\y$ has no interior lattice point. Then 
\begin{equation}
[t_\y,t_{\x}]=\epsilon_{\x,\y}c_{deg(\y)}\mathbf{k}_{\a(\x,\y)}\frac{\theta_{\x+\y}}{\nu-\nu^{-1}}
\end{equation}
where $\a(\x,\y)$ is given by the cocycle~:
\begin{equation}\label{E:alphard}
\a(\x,\y)=\begin{cases}
\frac{1}{2}\epsilon_\x(\epsilon_{\x}\x+\epsilon_{\y}\y-\epsilon_{\x+\y}(\x+\y)) & \text{if}\; \epsilon_{\x,\y}=1,\\
\frac{1}{2}\epsilon_\y(\epsilon_{\x}\x+\epsilon_{\y}\y-\epsilon_{\x+\y}(\x+\y)) & \text{if}\; \epsilon_{\x,\y}=-1,
\end{cases}
\end{equation}
and the elements $\theta_{\z}$, $\z \in \Z^2$ are obtained by equating the Fourier
 coefficients of the collection of relations
\begin{equation}
\sum_i \theta_{i\x_0}s^i=exp((\nu-\nu^{-1})\sum_{r \geq 1}t_{r\x_0}s^r),
\end{equation}
for any $\x_0 \in (\Z^2)^+$ such that $deg(\x_0)=1$ (note that $\theta_{\z}=(\nu-\nu^{-1})t_{\z}$ if $deg(\z)=1$). 
\end{enumerate}
\end{theo}

The action of $\widehat{SL}(2,\Z)$ on $\widetilde{\mathbf{D}}\H_X$ preserves the Drinfeld double of the spherical subalgebra $\widetilde{\mathbf{D}}\mathbf{C}_X$, and is nothing else than the (obvious) action
$$\Phi(\boldsymbol{\kappa}_{\x}) = \boldsymbol{\kappa}_{\Phi(\x)},\qquad
\Phi(t_{\x})= t_{\Phi(\x)}\boldsymbol{\kappa}_{\Phi(\x)}^{n(\Phi,\mu(\x))},$$
for $\Phi \in SL(2,\Z)$.

\vspace{.1in}

Hence, we see that the very symmetric nature of the derived category $D^b(Coh(X))$ is reflected in the equally symmetric structure of the algebra $\U_X$.  Set $\U_X^0=\C[\mathbf{k}_{r,d}]_{r,d \in \Z}$.
Acting by the automorphisms in $\widehat{SL}(2,\Z)$, we obtain, in addition to the standard decomposition $\U_X \simeq \U_X^+ \otimes \U_X^0 \otimes \U_X^-$, a whole family of triangular decomposition, one for each line in $\Z^2$ through the origin.

\vspace{.15in}

\addtocounter{theo}{1}
\paragraph{\textbf{Remark \thetheo.}} The algebra $\U_X$, or more precisely its generic version $\underline{\U}^+$ may be interpreted as a limit when $n$ tends to infinity of the spherical Cherednik algebras of type $A_n$ (see \cite{SV}). The action of $SL(2,\Z)$ by automorphisms is a well-known and crucial feature
of these spherical Cherednik algebras. The above thus provides one conceptual interpretation for this symmetry.

\vspace{.2in}

Let us jump now from the case of an elliptic curve to the closely related case of a weighted projective line $\xpl$ of elliptic type. Here also, the group of autoequivalences of the derived category has a rich structure. Recall that we have a degree function $deg~:L(\p) \to \Z$; the kernel $L^0(\p)=deg^{-1}(0)$ is a finite group (isomorphic to $(\Z/2\Z)^3$ if $\p=(2,2,2,2)$, to $(\Z/3\Z)^2$ if $\p=(3,3,3)$, to
$\Z/2\Z \times \Z/4\Z$ if $\p=(2,4,4)$ and to $\Z/2\Z \times \Z/3\Z$ if $\p=(2,3,6)$). The group of automorphisms $Aut(\xpl)$ of $\xpl$ is also finite.

\vspace{.1in}

\begin{prop}[Lenzing-Meltzer, \cite{LM2}] Let $\xpl$ be a weighted projective line of elliptic type.There is a short exact sequence of groups
$$
\xymatrix{
1 \ar[r] & L^0(\p) \triangleleft {Aut}(\xpl)  \ar[r]   
& {Aut}(D^b(Coh(\xpl)) 
\ar[r]  & B_3 \ar[r] &  1. }
$$
\end{prop}

\vspace{.1in}

\begin{prop}[Burban-S., \cite{BS2}] Formulas (\ref{E:juki1}) and (\ref{E:juki2}) define an action of
$Aut(D^b(Coh(\xpl))$ on the reduced Drinfeld double $\widetilde{\mathbf{D}}\H_{\xpl}$ of the Hall algebra of $\xpl$. This action preserves the Drinfeld double $\widetilde{\mathbf{D}}\CC_{\xpl}$ of the spherical subalgebra of $\H_{\xpl}$.
\end{prop}

\vspace{.15in}

\addtocounter{theo}{1}
\paragraph{\textbf{Remarks \thetheo.}} i) As far as smooth projective curves are concerned, the above example of the elliptic curve in fact exhausts all possible applications of the conjectures in Section~5.3. Indeed, by a theorem of Bondal and Orlov \cite{BO}, the group of autoequivalences of $D^b(Coh(X))$ for $X$ of genus $g >1$ is essentially trivial~: it is generated by the group $Pic^0(X)$ of line bundles of degree zero, by the group of automorphisms of the curve $X$ itself, and by the suspension functor $T$. All of these clearly act on $\widetilde{\mathbf{D}}\H_{X}$, and preserve $\widetilde{\mathbf{D}}\CC_X$. Completely similar results hold for hyperbolic weighted projective lines (see \cite{LM2}).\\
ii) The case of the derived equivalence between $Coh(\mathbb{P}^1)$ and $Rep\vec{Q}$ for $\vec{Q}$ the Kronecker quiver is dealt with in many details in \cite{BS3}.\\
iii) There is one important point which we would like to stress. Although the Drinfeld double of a Hopf algebra is again a Hopf algebra, the Drinfeld doubles of the Hall algebras of two derived equivalent hereditary categories are in general \textit{not} isomorphic as Hopf algebras (i.e. the coproduct is not invariant under derived equivalence). For instance, $D^b(Coh(\mathbb{P}^1)) \simeq D^b(Rep_k\vec{Q})$ where $\vec{Q}$ is the Kronecker quiver and $\widetilde{\mathbf{D}}\CC_{\mathbb{P}^1} \simeq \widetilde{\mathbf{D}}\CC_{\vec{Q}}\;\; (\simeq \U_\nu(\widehat{\mathfrak{sl}}_2))$ as algebras but not as colagebras~: Drinfeld's coproduct is not equal to the standard coproduct (see Appendix~A.6.).

Another way of saying this is that a given quantum group (say, coming from a Hall algebra) may have several \textit{distinct} coproduct structures, all compatible with the \textit{same} multiplication --one for each polarization
$\U_{\nu} \simeq \U_{\nu}^+ \otimes \U_{\nu}^0 \otimes \U_{\nu}^-$.

\vspace{.2in}

\centerline{\textbf{5.5. The Hall Lie algebra of Peng and Xiao.}} 
\addcontentsline{toc}{subsection}{\tocsubsection {}{}{\; 5.5. The Hall Lie algebra of Peng and Xiao.}}

\vspace{.2in}

Let $\mathcal{A}$ be a finitary and hereditary $\mathbb{F}_q$-linear category and let $D^b(\mathcal{A})/T^2$ be its $2$-periodic derived category. Observe that if $M \in D^b(\mathcal{A})/T^2$ is indecomposable then the class of $M$ in the Grothendieck group $K(D^b(\A)/T^2) \simeq K(\A)$ is nonzero. 

Let $\mathcal{I}_{\A}$ stand for the set of indecomposable objects of $D^b(\A)/T^2$, up to isomorphism. For any $M \in \mathcal{I}_{\A}$ put $d(M)=dim(End\;M/rad\;End\;M)$, and let $\mathbf{h}$ be the $\Z$-submodule of $\Q \otimes K(\A)$ generated by the elements $\frac{h_N}{d(N)}$ for $N \in \mathcal{I}_{\A}$, where $h_N$ stands for the class of an object $N$. We may extend the symmetrized Euler form $(\;,\;)_a$ to $\mathbf{h}$ by setting
\begin{equation*}
\begin{split}
(M,N)_a={dim}&({Hom}(M,N))+{dim}({Hom}(N,M))\\
&-{dim}({Hom}(M,N[1]))-{dim}({Hom}(N,M[1])).
\end{split}
\end{equation*}
Consider a free $\Z$-module
$$\mathbf{n}=\bigoplus_{M \in \mathcal{I}_\A} \Z[M]$$
with basis indexed by elements of $\mathcal{I}_{\A}$, and put $\mathbf{g}=\mathbf{g}(\A)=\mathbf{n}\oplus \mathbf{h}$. 

\vspace{.1in}

For any three elements $X,Y,L$ of $\mathcal{I}_{\A}$, let $W(X,Y;L)$ stand for the (finite) set of all distinguished triangles
$$(f,g,h)~: \xymatrix{ X \ar[r]^{f} & L \ar[r]^{g} & Y \ar[r]^h & TX},$$
and let $F_{X,Y}^L$ be the cardinality of the set of orbits in $W(X,Y;L)$ under the natural action of $Aut(X) \times Aut(Y)$. Of course, the definition of $F_{X,Y}^L$ strongly ressembles that of Hall numbers of abelian categories. Unfortunately, contrary to what one might hope at first, the numbers $F_{X,Y}^L$ are \textit{not} the structure constants for an associative algebra. Nevertheless, the differences $\gamma_{X,Y}^L=F_{X,Y}^L-F_{Y,X}^L$ are, ``at $q=1$'', the structure constants for a \textit{Lie} algebra as the following first main theorem of \cite{PX} states~:

\vspace{.1in}

\begin{theo}[Peng-Xiao, \cite{PX}] Put $\mathbf{g}_{(q-1)}=\mathbf{g}/(q-1)\mathbf{g}$ (a $\Z/(q-1)\Z$-module). The following brackets
$$\left[ [X],[Y] \right]=\begin{cases} \sum_{L \in \mathcal{I}_{\A}} \gamma_{X,Y}^L [L]\;&\;\text{if}\; Y \not\simeq TX\\
\frac{\mathbf{h}_X}{d(X)}\; & \; \text{if}\; Y \simeq TX
\end{cases},$$
$$\left[\mathbf{h}_X,[Y]\right]=-\left[ [Y],\mathbf{h}_X\right]=-(\mathbf{h}_X,\mathbf{h}_Y)[Y],$$
$$[\mathbf{h}_X,\mathbf{h}_Y]=0$$
endow $\mathbf{g}_{(q-1)}$ with the structure of a Lie algebra.
\end{theo}

\vspace{.15in}

We will call this Lie algebra the \textit{Peng-Xiao Hall Lie algebra of the category $\A$}. Let us now assume that $\A=Rep_k\vec{Q}$ for $\vec{Q}$ an arbitrary quiver. As usual, as soon as $\vec{Q}$ is not of finite type, the Hall Lie algebra is too big. Define the \textit{composition Lie algebra} $\mathcal{LC}(\A)_{(q-1)}$ as the subalgebra of $\mathbf{g}(\A)_{(q-1)}$ generated by the elements $\{[S], [TS]\}$ where $S$ runs among all simple objects of $\A$. To define a generic form of $\mathcal{LC}(\A)_{(q-1)}$, defined over $\Z$ or $\C$, we consider the product over all finite extensions
$E \supseteq \mathbb{F}_q$
$$\Pi=\prod_E \mathbf{g}(\A^E)_{(|E|-1)}$$
and let $\mathcal{LC}(\A)$ be the $\Z$-Lie subalgebra of $\Pi$ generated by elements $\{[S^E],[TS^E]\}$, where $S$ again runs over all simple objects of $\A$. This is a free $\Z$-module.

\vspace{.1in}

\begin{theo}[Peng-Xiao, \cite{PX}] The Lie algebra $\mathcal{LC}(\A) \otimes_{\Z} \C$ is isomorphic to the Kac-Moody algebra $\g'$ associated to the quiver $\vec{Q}$.
\end{theo}

\vspace{.2in}

\addtocounter{theo}{1}
\paragraph{\textbf{Remarks \thetheo.}} i) One may alternatively define $\mathcal{LC}(\A)_{(q-1)}$ as the Lie subalgebra generated by the classes of all \textit{exceptional} indecomposable objects in $\A$ (recall that an object $M$ is called exceptional if $Ext^1(M,M)=\{0\}$). This definition has the advantage of being absolutely intrinsic to the $2$-periodic derived category $D^b(\A)/T^2$ (i.e. independent of the choice of an abelian subcategory $\A$).\\
ii) The Hall Lie algebra $\mathbf{g}(\A)_{(q-1)}$ and the spherical Lie algebra $\mathcal{LC}(\A)_{(q-1)}$ of the category $\A=Coh(\xpl)$ of coherent sheaves on a weighted projective line have been studied in \cite{CB2} and are strongly related to loop algebras $\mathcal{L}\g$ of star-shaped Kac-Moody algebras, as one can expect given the results of Section~4.6. The spherical Lie algebra plays a key role in Crawley-Boevey's proof of the important Theorem~\ref{T:CBKac}.\\
iii) Using a small variant of the constructions of Peng and Xiao, Lin and Peng consider in \cite{LP} the composition Lie algebra in the case $\A=Rep_k\Lambda$, where $\Lambda$ is a tubular canonical algebra, a certain finite-dimensional associative algebra of global dimension \textit{two}, which is derived equivalent to $Coh(\xpl)$ for $\xpl$ of elliptic type. They show that $\mathcal{LC}(\A)$ is isomorphic to an elliptic Lie algebra $\mathcal{E}_{\g}$ of type $D_4,E_6,E_7$ or $E_8$. This result, which appeared at the same time as and independently of \cite{SDuke}, is another illustration of the principle according to which the Hall algebras of two derived equivalent categories have isomorphic ``doubles''.

\vspace{.2in}

\centerline{\textbf{5.6. Kapranov and To\"en's derived Hall algebras.}} 
\addcontentsline{toc}{subsection}{\tocsubsection {}{}{\; 5.6. Kapranov and To\"en's derived Hall algebras.}}

\vspace{.2in}

The work of Peng and Xiao associates to the $2$-periodic derived category $D^b(\A)/T^2$ of any finitary category of finite global dimension $\A$ a Lie algebra. If instead of the $2$-periodic derived category one considers the whole derived category $D^b(\A)$ then it is possible to define an \textit{associative algebra}, very reminiscent of Hall algebras, and which is an invariant of the derived category. This is done by Kapranov when $\A$ is hereditary and by To\"en for any abelian category $\A$ of finite global dimension (the two constructions are closely related but not equivalent). We briefly present these and refer the reader to \cite{Kap2} and \cite{Toen} for details~:

\vspace{.2in}

Let $\A$ be a finitary, hereditary category, and let $\H_{\A}$ and $\widetilde{\H}_{\A}$ be the Hall algebra and extended Hall algebra of $\A$. Kapranov's idea is to associate to $D^b(\A)$ a kind of ``universal covering of the Drinfeld double'' of $\widetilde{\A}$, in the same way as one would like to associate to the $2$-periodic derived category $D^b(\A)/T^2$ the Drinfeld double of $\widetilde{\H}_{\A}$.

If $M^{\bullet},N^\bullet,R^{\bullet}$ are objects of $D^b(\A)$ we let $W(M^{\bullet},N^{\bullet}; R^{\bullet})$ stand for the set of exact triangles
$$(f,g,h)\;: \xymatrix{ N^{\bullet} \ar[r]^f & R^{\bullet} \ar[r]^g& M^{\bullet} \ar[r]^h & N^{\bullet}[1]}$$
and put
$$F_{M^{\bullet},N^{\bullet}}^{R^{\bullet}}=\frac{|W(M^{\bullet},N^{\bullet};R^{\bullet})|}{|{Aut}(M^{\bullet})| \cdot |{Aut}(N^{\bullet})|}.$$
By definition, the \textit{Lattice algebra} $L(\A)$ is the associative $\C$-algebra generated by elements
$[M]^{(i)}$ for $M \in \text{Ob}(\A)/\sim$ and $i \in \Z$, and $\mathbf{k}_{\a}$ for $\a \in K(\A)$, subject to the following set of relations~:
$$\mathbf{k}_{\a}\mathbf{k}_{\b}=\mathbf{k}_{\a+\b},$$
$$\mathbf{k}_{\a} [M]^{(i)}\mathbf{k}_{\a}^{-1} =(\a,M)^{(-1)^i}[M]^{(i)},$$
$$[M]^{(i)}[N]^{(i)}=\langle M,N\rangle \sum_{R \in \text{Ob}(\A)} F_{M,N}^R [R]^{(i)},$$
$$[M]^{(i)}[N]^{(i+1)}=\sum_{P,Q} \langle M-P,N-Q\rangle F_{M[1],N}^{P[1]\oplus Q} \mathbf{k}_{M-P}^{(-1)^i}[Q]^{(i+1)}[P]^{(i)},$$
$$[M]^{(i)}[N]^{(j)}=(M,N)^{(-1)^{i-j}(i-j+1)}[N]^{(j)}[M]^{(i)}\qquad \text{if}\; |i-j| \geq 2.$$

From the definitions, it is clear that the assignement $[M]^{(i)} \mapsto [M]^{(i+1)}, \mathbf{k}_{\a} \mapsto \mathbf{k}_{\a}^{-1}$ gives rise to an automorphism $\Sigma~: L(\A) \stackrel{\sim}{\to} L(\A)$.

\vspace{.1in}

\begin{prop}[Kapranov, \cite{Kap2}] Let $\mathbf{K}$ be the subalgebra generated by elements $\mathbf{k}_{\a}$ for $\a \in K(\A)$, and let $L^{(i)}(\A)$ be the subalgebra of $L(\A)$ generated by $[M]^{(i)}$ for $M \in \text{Ob}(\A)/~$. Then $L^{(i)}(\A)$ is isomorphic to the Hall algebra $\H_{\A}$ and the multiplication map
$$\mathbf{K} \otimes \vec{\bigotimes} L^{(i)}(\A) \to L(\A)$$
is an isomorphism of vector spaces.
\end{prop}

\vspace{.1in}

In the above, $\vec{\bigotimes} L^{(i)}(\A)$ stands for the direct sum of all finite tensor products
$L^{(n_1)}(\A) \otimes \cdots \otimes L^{(n_r)}(\A)$ for $n_1>n_2 \cdots$. Hence $L(\A)$ is made up of infinitely many copies of the Hall algebra $\H_{\A}$, one for each abelian subcategory $\A[i]$ for $i \in \Z$. Moreover, the Hall algebras associated to categories $\A[i]$ and $\A[j]$ with $|i-j| \geq 2$ (so that $\A[i]$ and $\A[j]$ are orthogonal) essentially commute with each other.

\vspace{.1in}

The important point concerning the above algebraic construction is its invariance under derived equivalences~:

\vspace{.1in}

\begin{theo}[Kapranov, \cite{Kap2}] Let $\B$ be another finitary hereditary category and let $F: D^b(\B) \stackrel{\sim}{\to} D^b(A)$ be a derived equivalence. Then the correspondence
$$[M]^{(i)} \mapsto \Sigma^i [F(M)], \qquad \mathbf{k}_{\beta} \mapsto \mathbf{k}_{F(\beta)}$$
extends to an isomorphism of algebras, where if $F(M)= \cdots \to N_i \to N_{i+1} \to \cdots$ is a complex of objects of $\A$ then $[F(M)]:=\vec{\prod}_i [H^i(F(M))]^{(-i)}$.
\end{theo}

\vspace{.2in}

Let us now turn to the constructions of To\"en. Here, we let $\A$ be an abelian category satisfying the following finiteness conditions~:
\begin{enumerate}
\item[i)] $\A$ is of finite global dimension.
\item[ii)]For any two objects $M,N \in \text{Ob}(\mathcal{A})$ and $i \geq 0$ we have $|{Ext}^i(M,N)| < \infty$.
\end{enumerate}

Let $\mathcal{X}$ denote the set of objects of $D^b(\A)$, up to quasi-isomorphism, and let $\H_{D^b(\A)}$ be the free $\C$-vector space with basis $[M], M \in \mathcal{X}$. To\"en found in \cite{Toen} the right way to ``count distinguished triangles'' in order to endow $\H_{D^b(\A)}$ with
the structure of an associative algebra.

Let as before $W(M^{\bullet},N^{\bullet};R^{\bullet})$ stand for the set of all distinguished triangles
$$(f,g,h)\;: \xymatrix{ N^{\bullet} \ar[r]^f & R^{\bullet} \ar[r]^g& M^{\bullet} \ar[r]^h & N^{\bullet}[1]}$$
and let $(N^{\bullet},R^{\bullet})_{M^{\bullet}}$ be the set of orbits of $W(M^{\bullet},N^{\bullet};R^{\bullet})$ under the action of $Aut(M^{\bullet})$. We define
$$g_{M^{\bullet},N^{\bullet}}^{R^{\bullet}}=
\frac{|(N^{\bullet},R^{\bullet})_{M^{\bullet}}|}{|Aut(N^{\bullet})|} \cdot \prod_{i >0}
\frac{|{Hom}(N^{\bullet}[i],R^{\bullet})|^{(-1)^i}}{|{Hom}(N^{\bullet}[i],N^{\bullet})|^{(-1)^i}}.$$

\vspace{.1in}

\begin{theo}[To\"en, \cite{Toen}]\label{T:TOEN} The following defines on the vector space $\H_{D^b(\A)}$ the structure of an associative algebra~:
$$[M^{\bullet}] \cdot [N^{\bullet}]=\sum_{R^{\bullet}} g_{M^{\bullet},N^{\bullet}}^{R^{\bullet}}[R^{\bullet}].$$
\end{theo}

\vspace{.1in}

Of course, since $\H_{D^b(\A)}$ is defined directly in terms of $D^b(\A)$ without any reference to
the abelian category $\A$ itself, it is perfectly invariant under derived equivalences. In addition, it is clear that $\H_{D^b(\A)}$ contains the usual Hall algebra $\H_{\A}$ (where we ignore the twist by the Euler form $\langle\,,\,\rangle$) of the abelian category $\A$. The same holds for the Hall algebra of the heart of any t-structure on $D^b(\A)$.

\vspace{.15in}

\addtocounter{theo}{1}
\paragraph{\textbf{Remarks \thetheo.}} i) The definition of the derived Hall algebra $\H_{D^b(\A)}$ makes sense and is valid for an arbitrary finitary triangulated category of left-finite homological dimension, that is such that ${Hom}(M[i],N)$ is finite for all objects $M$ and $N$ and zero (for fixed $M$ and $N$) whenever $i \gg 0$. To\"en's approach makes use of homotopy theory and $dg$-categories techniques; a proof of Theorem~\ref{T:TOEN} entirely in the language of triangulated categories can be found in \cite{XiaoXu}.\\
ii) Both Kapranov's lattice algebra and To\"en's derived Hall algebra are related not to quantum groups but rather to ``universal coverings'' of quantum groups (see \cite{Kap2} for a precise definition).

\newpage

\centerline{\textbf{Windows.}} 
\addcontentsline{toc}{section}{\tocsection {}{}{Windows.}}

\vspace{.2in}

To bring these notes to a conclusion, let us briefly mention some possible directions for future research. 

\vspace{.15in}

\noindent
\textit{Hereditary categories, up to derived equivalences.} As these notes hope to show, the Hall algebras of the categories of representations of quivers, of coherent sheaves on weighted projective lines and of coherent sheaves on elliptic curves are now rather well understood (a more honest statement would be that the composition algebras of these categories are rather well understood).
If one believes the various conjectures made in Section~5.3. then the same could be said of any hereditary category which is derived equivalent to one of the above. A deep theorem of Happel \cite{Happel} and Happel-Reiten \cite{ReitenHappel} gives a complete classification, \textit{up to derived equivalence} of hereditary $k$-linear categories containing a tilting object, for an arbitrary field $k$. For $k$ algebraically closed, these are precisely the categories $Rep_k\vec{Q}$ and $Coh(\xpl)$ studied in Lectures~3 and 4. For $k$ a finite field, one has to consider in addition species (equivalently, quivers with automorphisms), and an analogue of species for weighted projective lines. It would be interesting to compute the Hall algebras of these variants of weighted projective lines, in the hope of obtaining in this fashion certain loop algebras of non simply laced Kac-Moody algebras.

 The situation for hereditary categories which do \text{not} possess a tilting object is much more mysterious~: only the limiting case of an elliptic curve is known. Unraveling the structure of the Hall algebra of higher genus curves (possibly weighted) seems like an important problem (see Section~4.11., and Section~2.5. for the relation to number theory and automorphic forms). By Happel's theorem, we may represent the set of hereditary categories, up to derived equivalence, in the following way~:

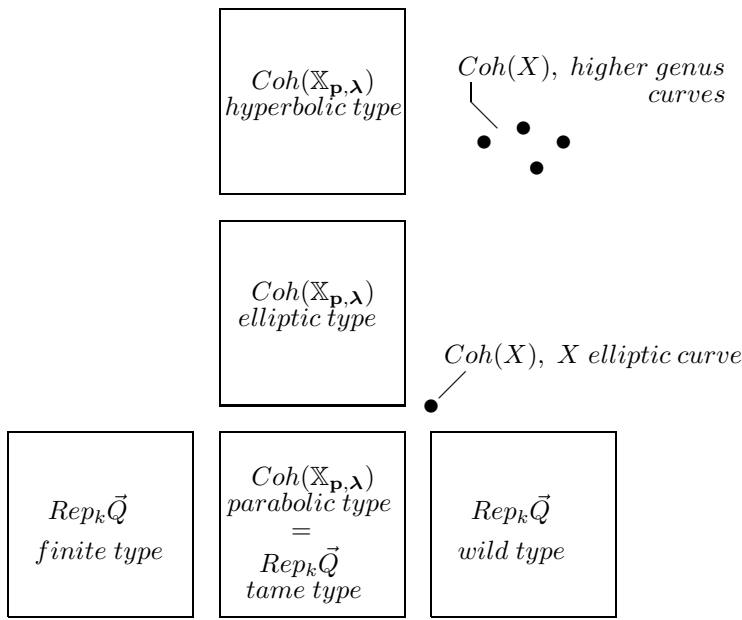
\vspace{.2in}

\centerline{
\begin{picture}(230,230)
\put(0,0){\line(0,1){70}}
\put(0,0){\line(1,0){70}}
\put(70,70){\line(-1,0){70}}
\put(70,70){\line(0,-1){70}}
\put(80,0){\line(0,1){70}}
\put(80,0){\line(1,0){70}}
\put(150,70){\line(-1,0){70}}
\put(150,70){\line(0,-1){70}}
\put(160,0){\line(0,1){70}}
\put(160,0){\line(1,0){70}}
\put(230,70){\line(-1,0){70}}
\put(230,70){\line(0,-1){70}}
\put(80,80){\line(0,1){70}}
\put(80,80){\line(1,0){70}}
\put(150,150){\line(-1,0){70}}
\put(150,150){\line(0,-1){70}}
\put(150,160){\line(0,1){70}}
\put(150,160){\line(-1,0){70}}
\put(80,230){\line(1,0){70}}
\put(80,230){\line(0,-1){70}}
\put(15,37){$Rep_k\vec{Q}$}
\put(10,22){$finite\;type$}
\put(95,17){$Rep_k\vec{Q}$}
\put(90,7){$tame\;type$}
\put(175,37){$Rep_k\vec{Q}$}
\put(170,22){$wild\;type$}
\put(92,50){$Coh(\xpl)$}
\put(83,40){$parabolic\;type$}
\put(107,30){$=$}
\put(92,120){$Coh(\xpl)$}
\put(87,110){$elliptic\;type$}
\put(92,200){$Coh(\xpl)$}
\put(82,190){$hyperbolic\;type$}
\put(160,80){\circle*{5}}
\put(163,83){\line(1,1){10}}
\put(165,96){$Coh(X),\;X\;elliptic\;curve$}
\put(180,180){\circle*{5}}
\put(200,170){\circle*{5}}
\put(195,185){\circle*{5}}
\put(210,180){\circle*{5}}
\put(185,185){\line(-1,1){10}}
\put(175,195){\line(0,1){7}}
\put(170,205){$Coh(X),\;higher\;genus$}
\put(242,195){$curves$}
\end{picture}}

\vspace{.4in}

\noindent
\textit{Higher dimensional varieties--derived setting.} These notes are almost entirely devoted to the problem of describing the Hall algebras of hereditary categories, but the definition of the Hall algebra (as an \textit{algebra}) is valid for any finitary category. Up to our knowledge, except for the work of Lin and Peng \cite{LP} which deals with certain categories of global dimension two which are derived equivalent to hereditary categories, not a single example of a Hall algebra of a category of higher global dimension has been computed. The first geometric example which comes to mind is that of $Coh(\mathbb{P}^2)$. Of course, it is necessary to restrict oneself to a well-chosen `spherical' subalgebra $\CC_{\mathbb{P}^2} \subset \H_{\mathbb{P}^2}$-- one could take for instance the subalgebra generated by elements $\mathbf{1}_{r,d,c}$, sum over all sheaves of rank $r$, degree $d$ and first/second Chern class $c$, for all values of $r,d$ and $c$-- but there may be other choices as well (taking into account the singularities of the moduli spaces in question, for instance). Moreover, since there is no more reason to expect that $\CC_{\mathbb{P}^2}$ is \textit{half} of something interesting, it is probably better to consider the Hall algebra of $D^b(Coh(\mathbb{P}^2))$, in To\"en's version (Section~5.5.) or in Peng and Xiao's version (Section~5.4.). Finally, there is another crucial difference with the one dimensional situation: the moduli spaces of objects in a category of global dimension $>1$ tend to be intrinsically singular; this has to be taken into account one way or another.

\vspace{.2in}

\noindent
\textit{Singular curves-- derived setting.} What about Hall algebras of singular curves, for example a reducible curve $xy=0$ or a nodal curve ? The category $Coh(X)$ for $X$ a singular curve is of infinite global dimension, and it seems more natural to study its Hall algebra in the derived setting than in the abelian setting. When $X$ is a rational curve with simple double points or transversal intersections, Burban and Drozd have given a classification of indecomposable objects of $D^b(Coh(X))$ (see \cite{BD}). Of course, when $X$ is reducible, the Hall algebra of $Coh(X)$ contains as a subalgebra the Hall algebra of $Coh(Z)$ for any irreducible component $Z$ of $X$, but it is the relation between these subalgebras which remains to be understood.

\vspace{.2in}

\noindent
\textit{The $\chi$-Hall algebra.} Lusztig introduced in \cite{Luseuler} a variant of the Hall algebra, in which one replaces the {number} $\mathcal{P}_{M,N}^R$ of subojects of an object $R$ of type $N$ and cotype $M$ by the \textit{Euler characteristic} $\chi(\mathcal{V}_{M,N}^R)$ of the (constructible) space of all such subobjects. Of course, one needs to prove the existence of a well-behaved space $\chi(\mathcal{V}_{M,N}^R)$ parametrizing such subojects (this is easy to do for categories of modules over a finite-dimensional algebra, or for categories of coherent sheaves on a projective curve). When applied to the category $Rep_k\vec{Q}$ for some quiver $\vec{Q}$ one obtains in this way the \textit{classical} (i.e. non quantum) enveloping algebra $\U(\bo'_+)$. The case of the categories $Coh(X)$ for $X$ a weighted projective line hasn't been treated but it is highly likely that one obtains $\U(\mathcal{L}\bo_+)$. Recently, Xiao, Xu and Zhang have given the definition of the $\chi$-Hall Lie algebra of a triangulated category of arbitrary finite global dimension (see \cite{XXZ}). For higher dimensional varieties such as $\A=Coh(\mathbb{P}^2)$ or for singular curves this Lie algebra might be easier to study than Peng and Xiao's Hall Lie algebra (Section~5.4.). The $\chi$-Hall algebra appears in recent work of Joyce \cite{Joyce} and Bridgeland-Toledano-Laredo \cite{BTL}.

\vspace{.2in}

\noindent
\textit{Kapranov's Hall category.} Let $\mathcal{A}$ be a finitary category. Let $\mathcal{A}^{iso}$ be the abelian category whose objects are the objects of $\mathcal{A}$ and whose morphisms are the \textit{isomorphisms} in $\mathcal{A}$. Following Kapranov \cite{Kaporal}, one may consider the category $\mathcal{H}_\mathcal{A}$ of \textit{functors} $F~:\mathcal{A} \to Vec_{\mathbb{F}_q}$ for which $F(M)=\{0\}$ for almost all $M$. This category is equipped with a monoidal structure, defined as follows~:
$$(F \odot G) (M)=\bigoplus_{N \subset M} F(M/N) \otimes G(N).$$
For instance, in the simplest case $\mathcal{A}=\mathbb{F}_q$, an object of $\mathcal{H}_\mathcal{A}$ is given by a collection of representations of $GL(n,\mathbb{F}_q)$ for $n \in \N$ (the images of the vectors spaces $k^n$ for $n \in \N$ ). The monoidal structure is given by the parabolic induction functor~: if $F \in Obj(\mathcal{H}_{\mathcal{A}})$ corresponds to a representation $V$ of $GL(n,\mathbb{F}_q)$ and $G \in Obj(\mathcal{H}_{\A})$ corresponds to a representation $W$ of $GL(m,\mathbb{F}_q)$ then $F \odot G$ is associated to the representation 
$Ind_{P(n,m)}^{GL(m+n)}(V \otimes W)$, where $P(n,m)$ is the maximal parabolic subgroup of $GL(m+n)$ of type $n,m$.

This definition is clearly a ``categorification'' of the notion of a Hall algebra~: the space of functions
$Obj(\mathcal{A}) \to \N$ is replaced by the space of functors $\A \to Vec_{\mathbb{F}_q}$, addition is replaced by $\oplus$ and multiplication is replaced by $\otimes$.
Note that for general categories, just as the Hall algebra is not commutative, so are $F \odot G$ and $G \odot F$ not isomorphic. Even when they are, i.e. when there exist  functorial natural transformations $\beta_{F,G}~: F \odot G \stackrel{\sim}{\to} G \odot F$ endowing $\mathcal{H}_\A$ with the structure of a \textit{braided} monoidal category, we may have $\beta_{F,G} \circ \beta_{G,F} \neq Id$. This is for instance the case for $\A=Vec_{\mathbb{F}_q}$. Describing the Hall categories of quivers or curves seems an intersesting and challenging problem which is, up to my knowledge, still untouched.

\vspace{.2in}

\noindent
\textit{Hall algebras and cluster algebras.} Fomin and Zelevinsky introduced around the year 2000 a fundamental new class of (commutative) algebras
associated to an \textit{antisymmetric} square matrix, \cite{FZ}. Such a matrix may also be constructed from a graph or from a quiver $\vec{Q}$.
When the quiver $\vec{Q}$ is a Dynkin diagram of finite type, the structure constants of the associated cluster algebra $A_{\vec{Q}}$ are shown in \cite{CalKel} to be some kind of Hall numbers (evaluated at $q=1$), not for the category of representations $Rep_k\vec{Q}$ but rather for the so-called \textit{cluster category} $\mathcal{C}_{\vec{Q}}=D^b(Rep_k\vec{Q})/\tau^{-1} \circ [1]$. Here $\tau$ is the Auslander-Reiten transform on $D^b(Rep_k\vec{Q})$ and $[1]$ is the shift functor. Note that $\mathcal{C}_{\vec{Q}}$ is not abelian, but it is still triangulated. Its indecomposable objects are in natural bijection with a set of generators for $A_{\vec{Q}}$. This makes it tempting to try to realize $A_{\vec{Q}}$ as some sort of ``Hall algebra'' (at $q=1$) of $\mathcal{C}_{\vec{Q}}$. We refer to the paper \cite{CalKeller} for more in this direction, and for some important applications to the structure theory of cluster algebras.

\vspace{.2in}

\noindent
\textit{Lusztig and Nakajima quiver varieties.} The theory of Hall algebras of quivers finds a very natural and important continuation in the work of Lusztig (see e.g. \cite{Lusztignilp}) and Nakajima (see \cite{Nak1}), among others. By considering a convolution-like product on various spaces of functions on
the \textit{cotangent bundle} to the ``moduli space'' $\mathcal{X}$ of objects in $Rep_k\vec{Q}$ Nakajima and Lusztig realize directly the whole enveloping algebra $\U(\g)$ as well as all the integrable highest weight representations of $\g$. Generalizing this theory to other abelian categories (such as $Coh(X)$ for smooth projective curves or $Coh(\xpl)$ for weighted projective lines), or even for an arbitrary abelian or triangulated category of finite global dimension seems very interesting and very promising (and most likely very difficult). The companion survey \cite{SLectures2} gives an exposition of Lusztig's theory.

\newpage

\centerline{\large{\textbf{Appendix.}}}
\addcontentsline{toc}{section}{\tocsection {}{}{Appendix.}}

\vspace{.25in}

\setcounter{section}{6}
\setcounter{equation}{0}
\setcounter{theo}{0}

\centerline{\textbf{A.1. Simple Lie algebras.}}
\addcontentsline{toc}{subsection}{\tocsubsection {}{}{\; A.1. Simple Lie algebras}}

\vspace{.15in}

\paragraph{}We will only consider Lie algebras over $\C$. An excellent reference for all the results quoted below is \cite{Serrecsla}.
Let $\g$ be a finite-dimensional simple Lie algebra, (i.e. a Lie algebra which has no nontrivial ideals). To unravel the structure of $\g$ we consider the adjoint representation~: 
\begin{equation*}
\begin{split}
\g &\to {End}(\g)\\
 x &\mapsto ad_x~: (y \mapsto [x,y])
 \end{split}
 \end{equation*}
Recall that an element $x \in \g$ is \textit{semisimple} if $ad_x$ is a semisimple (diagonalizable) operator in $\g$. A \textit{Cartan subalgebra} is a maximal commutative subalgebra $\h \subset \g$ consisting of semisimple elements. There are many such subalgebras, but they are all conjugate under the action of ${Aut}(\g)$. Any element $h \in \h$ gives rise to a decomposition of $\g$ into $ad_h$-eigenspaces, and since $\h$ is commutative there is a simultaneous decomposition
$\g= \bigoplus_{\alpha \in \h^*} \g_{\a}$ where
$$\g_{\a}=\{ x \in \g\;|\; [h,x]=\alpha(h) x \; \text{for\;all\;} h \in \h\}.$$ 
Of course since $\g$ is finite-dimensional the spaces $\g_{\a}$ are zero for all but finitely many values of $\a \in \h^*$. A nonzero element $\alpha \in \h^*$ such that $\g_{\a} \neq \{0\}$ is called a \textit{root} of $\g$ and the set of roots $\Delta=\{\a \in \h^* \backslash\{0\}\;|\; \g_{\a} \neq \{0\}\}$ is called the \textit{root system} of $\g$ (relative to $\h$). The $\Z$-span $Q=\sum \Z \a$ of the roots is the \textit{root lattice} of $\g$. Finally, $\h$ is equipped with a natural symmetric bilinear form defined by $(h,h')=Tr_{{End}\;\g}(ad_h ad_{h'})$ called the \textit{Cartan-Killing form}. This form is nondegenerate, and it therefore induces a form $(\;,\;)$ on $\h^*$, also called the Cartan-Killing form.

\vspace{.1in}

\paragraph{\textbf{Theorem~ A.1.}} \textit{The following hold~:
\begin{enumerate}
\item[i)] We have $\g_0=\h$ and ${dim}(\g_{\a})=1$ if $\a \in \Delta$,
\item[ii)] $\Delta$ spans $\h^*$ and if $x \in \Delta$ then $-x \in \Delta$,
\item[iii)] There exists a real structure $\h^*_{\mathbb{R}} \subset \h^*$ such that $\Delta \subset \h_{\mathbb{R}}^*$, and the restriction of the Cartan-Killing form to $\h^*_{\mathbb{R}}$ is positive definite.
\end{enumerate}}

\vspace{.1in}

The Jacobi identity implies that $[\g_{\a}, \g_{\beta}] \subset \g_{\a+\beta}$.  In other words, $\g$ is graded by $\Delta$ and by property i) above, the decomposition
$$\g=\h \oplus \bigoplus_{\a \in \Delta} \g_{\a}$$
gives an accurate ``skeleton'' of $\g$. 

\vspace{.15in}

\addtocounter{theo}{1}
\paragraph{\textbf{Example~A.2.}} Let us first consider the case of $\g=\mathfrak{sl}_n(\C)$. For $\h$ we may take the space of all traceless diagonal matrices
$$\h=\left\{ \begin{pmatrix} d_1 & 0 &\cdots & 0\\ 0 & d_2 & \cdots & 0\\ \vdots &\vdots & \ddots & \vdots \\ 0 & 0 & \cdots & d_n \end{pmatrix}, \; \sum_i d_i =0\right\}.$$
Let $\epsilon_i$ be the linear form on $\h$ such that $\epsilon_i(Diag(d_1, \ldots, d_n))=d_i$. We have
$\h^*=\bigoplus_i \C \epsilon_i / \C \cdot (\sum \epsilon_i)$. The root system and root spaces are given by
$$\Delta=\{ \epsilon_i-\epsilon_j\;|\; i \neq j\}, \qquad \g_{\epsilon_i-\epsilon_j}=\C E_{ij},$$
where $E_{ij}$ denotes the elementary matrix.  The Cartan form on $\h^*$ is simply given by
$(\sum u_i \epsilon_i, \sum v_i \epsilon_i)=\sum u_i v_i$ and we may take $\h_{\mathbb{R}}=\sum_{i,j} \mathbb{R} (\epsilon_i-\epsilon_j)$.

\vspace{.1in}

For instance for $n=2$, we have $\h^* \simeq \C$, and $\Delta=\{ \a, -\a\}$ where $\a=\epsilon_1-\epsilon_2$~:

\vspace{.25in}

\centerline{
\begin{picture}(300, 10)
\put(60,0){$\Delta=$}
\put(100,0){\circle*{5}}
\put(180,0){\circle*{5}}
\put(92,4){$-\alpha$}
\put(182,4){$\alpha$}
\put(100,0){\line(1,0){80}}
\end{picture}}

\vspace{.25in}

\noindent
while for $n=3$ and taking into account the Cartan form on $\h^*$, the root system may be depicted as follows~:

\vspace{.1in}

\centerline{
\begin{picture}(300, 60)
\put(30,0){$\Delta=$}
\put(80,0){\circle*{5}}
\put(200,0){\circle*{5}}
\put(110,45){\circle*{5}}
\put(110,-45){\circle*{5}}
\put(170,45){\circle*{5}}
\put(170,-45){\circle*{5}}
\put(67,4){$\epsilon_3-\epsilon_2$}
\put(187,4){$\epsilon_2-\epsilon_3$}
\put(97,49){$\epsilon_1-\epsilon_2$}
\put(97,-53){$\epsilon_3-\epsilon_1$}
\put(157,49){$\epsilon_1-\epsilon_3$}
\put(157,-53){$\epsilon_2-\epsilon_1$}
\put(80,0){\line(1,0){120}}
\put(110,45){\line(2,-3){60}}
\put(110,-45){\line(2,3){60}}
\end{picture}}

\vspace{.9in}

\endexample

\vspace{.15in}

\addtocounter{theo}{1}
\paragraph{\textbf{Example~A.3.}} Let us now consider the case of $\g=\mathfrak{so}_{2n}(\C)$ (for $n \geq 3$), which we view as the set of all matrices
$$\mathfrak{so}_{2n}=\left\{ \begin{pmatrix} X & Y \\ Z & -X^t \end{pmatrix} \;|\; X,Y,Z \in \mathfrak{gl}_n(\C), \; Y^t=-Y, \;Z^t=-Z\right\}$$
(this is equivalent to the usual definition as traceless antisymmetric matrices). We may take
$$\h=\left\{ \begin{pmatrix} d_1 & 0 &\cdots & 0\\ 0 & d_2 & \cdots & 0\\ \vdots &\vdots & \ddots & \vdots \\ 0 & 0 & \cdots & d_n \end{pmatrix}, \;  d_{i+n}=-d_i\right\}.$$
Then $\h^* =\bigoplus_{i=1}^n \C \epsilon_i$ where $\epsilon_i(Diag(d_1, \ldots, d_n))=d_i$. The root system and root spaces are now 
$$\Delta=\{ \pm \epsilon_i \pm \epsilon_j\;|\; i < j\}, $$
$$\g_{\epsilon_i -\epsilon_j}=\C(E_{ij}-E_{n+j,n+i}), \quad \g_{\epsilon_i +\epsilon_j}=\C(E_{i,n+j}-E_{j,n+i}),$$
$$\g_{-\epsilon_i -\epsilon_j}=\C(E_{n+i,j}-E_{n+j,i}).$$
The Cartan form on $\h^*$ is again given by $(\sum u_i \epsilon_i, \sum v_i \epsilon_i)=\sum u_i v_i$ and we may take $\h_{\mathbb{R}}=\sum_{i} \mathbb{R} \epsilon_i$. \\
\endexample

\vspace{.15in}

As the above examples suggest, the set $\Delta$ has a large group of symmetries.  Let $W \subset {Aut}(\h^*)$ be the group generated by the set of orthogonal reflections $s_{\alpha}$ with respect to the hyperplanes perpendicular to the roots~:
$$s_{\alpha} \in {Aut}(\h^*), \qquad s_{\a}(x)=x-2 \frac{(\alpha,x)}{(\alpha,\alpha)} \alpha.$$
This group is called the \textit{Weyl group} of $\g$. 

\vspace{.1in}

\paragraph{\textbf{Proposition~A.4.}}(Weyl).  \textit{The group $W$ is finite and preserves the root system $\Delta$. }

\vspace{.1in}

Morally, the Weyl group is the group of symmetries of $\g$ (although the action on $\Delta$ doesn't directly lift to an action on $\g$).

\vspace{.1in}

Quite remarkably, it turns out that the data of $\h^*$ equipped with the Cartan-Killing form and with the root system $\Delta \subset \h^*$ completely characterizes the Lie algebra $\g$. In fact, as we will explain below, starting from this data, one can write down a rather simple presentation of $\g$ by generators and relations.

\vspace{.1in}

To see how this is done, let us try to find a minimal system of Lie algebra generators for $\g$. Since $[\g_{\a},\g_{\beta}] \subset \g_{\a+\beta}$, the picture of the root system suggests the following method~: pick any hyperplane $H \subset \h^*_{\mathbb{R}}$ which does not contain any root and call $\h^*_+, \h^*_{-}$ the two half-spaces delimited by $H$. This divides the root system $\Delta$ into positive and negative halves $\Delta =\Delta^+ \sqcup \Delta^-$, corresponding to nilpotent subalgebras
$$\n_+ =\bigoplus_{\a \in \Delta^+} \g_{\a}, \qquad \n_-=\bigoplus_{\a \in \Delta^-} \g_{\a}$$
of $\g$.  Now choose a root $\alpha_1 \in \Delta^+$ closest to $H$, then a root $\alpha_2$ next closest to $H$, and so on, and stop as soon as $\Delta^+ \subset \N\a_1 + \N \a_2 \cdots + \N \a_r$. 

\vspace{.1in}

\paragraph{\textbf{Proposition~A.5.}}(Cartan).  \textit{The roots $\a_1, \ldots \a_r$ are linearly independent and span $\h^*$. In particular, $r={dim}\; \h$. Moreover $\g_{\a_1}, \ldots, \g_{\a_r}$ generate $\n_+$ as a Lie algebra.}

\vspace{.1in}

The roots $\{\a_1, \ldots, \a_r\}$ are the (positive) \textit{simple roots}. It is important to note that this (and hence also $\n_{+}, \n_-$) depends on the choice of $H$. Obviously, if $w \in W$ then the set of simple roots relative to $w \cdot H$ is $\{w \a_1, \ldots , w \a_r\}$. In fact,

\vspace{.1in}

\paragraph{\textbf{Proposition~A.6.}}(Weyl). \textit{The Weyl group $W$ acts simply transitively on the collection of sets of simple roots.}

\vspace{.15in}

Any root $\a \in \Delta^+$ may be written as a positive sum $\a=\sum_i n_i \a_i$ with $n_i \in \N$.
Moreover, there is a unique root $\a \in \Delta^+$ for which $\sum n_i$ is maximal, called the \textit{highest root} and denoted $\theta$.

\vspace{.15in}

\addtocounter{theo}{1}
\paragraph{\textbf{Example~A.7.}} Let us go back to Example~A.2. In that situation, it is easy to check that the reflection $s_{\epsilon_i-\epsilon_j}$ acts on $\Delta$ by permutation of the indices $i$ and $j$. The Weyl group is thus canonically identified with the symmetric group $\mathfrak{S}_n$.  

The collection $\Pi=\{\epsilon_i-\epsilon_{i+1}\;|i=1, \ldots n-1\}$ is a set of simple roots, and the corresponding nilpotent subalgebra is simply the space of upper triangular matrices. Setting $\a_i=\epsilon_i-\epsilon_{i+1}$, we see that the set of positive roots is
$$\Delta^+=\{\a_{i} + \a_{i+1} + \cdots + \alpha_{j}\;|i \leq j\}$$
and the highest root is $\theta=\a_1 + \cdots + \a_{n-1}$.

This is by no means the only choice for $\Pi$. For instance, for $n=3$, and for the choice of hyperplane $H$ below, the set of simple roots is $\Pi=\{\epsilon_1-\epsilon_3, \epsilon_2-\epsilon_1\}$, and the highest root is $\epsilon_2-\epsilon_3$.

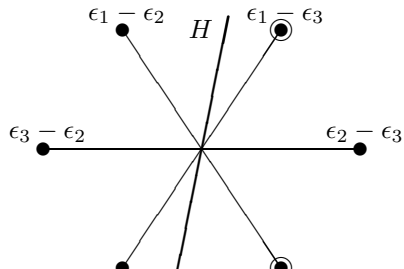
\vspace{.1in}

\centerline{
\begin{picture}(300, 60)
\put(80,0){\circle*{5}}
\put(200,0){\circle*{5}}
\put(110,45){\circle*{5}}
\put(110,-45){\circle*{5}}
\put(170,45){\circle*{5}}
\put(170,45){\circle{8}}
\put(170,-45){\circle*{5}}
\put(170,-45){\circle{8}}
\put(67,4){$\epsilon_3-\epsilon_2$}
\put(187,4){$\epsilon_2-\epsilon_3$}
\put(97,49){$\epsilon_1-\epsilon_2$}
\put(97,-53){$\epsilon_3-\epsilon_1$}
\put(157,49){$\epsilon_1-\epsilon_3$}
\put(157,-53){$\epsilon_2-\epsilon_1$}
\put(80,0){\line(1,0){120}}
\put(110,45){\line(2,-3){60}}
\put(110,-45){\line(2,3){60}}
\thicklines
\put(130,-50){\line(1,5){20}}
\thinlines
\put(135,42){$H$}
\end{picture}}

\vspace{.9in}

\noindent Here there are all together $6$ possible choices of sets of simple roots, one for each region of $\h^*$ cut out by the lines generated by the roots.
\endexample

\vspace{.15in}

\addtocounter{theo}{1}

\paragraph{\textbf{Example~A.8.}} Let us turn to Example~A.3. The collection
$$\Pi=\{\epsilon_i-\epsilon_{i+1}\;|\; i <n\} \cup \{\epsilon_{n-1} + \epsilon_n\}$$
is a set of simple roots. If we set $\a_i=\epsilon_i-\epsilon_{i+1}$ for $i <n$ and $\alpha_n=\epsilon_{n-1} + \epsilon_n$ then the set of positive roots is
\begin{equation*}
\begin{split}
\Delta^+=& \{\a_i + \a_{i+1} + \cdots + \a_j\;|\; i\leq j<n\} \cup \{\a_i + \cdots + \a_{n-2} + \a_{n}\;|\; i \leq n-2 \} \cup \{\a_n\}\\
&\cup  \{\a_i + \cdots + \a_{n-2} + \a_{n-1} +\a_{n}\;|\; i \leq n-2 \} \\
&\cup \{ \{\a_i + \cdots + 2\a_{n-2} + \a_{n}\;|\; i < n-2 \}.
\end{split}
\end{equation*}
  
As for the Weyl group is a little bit more complicated to describe. Any $w \in W$ induces a permutation of the indices $\{1, \ldots, n\}$, but also changes the signs. As an abstract group, $W$ is an extension
$$\xymatrix{
1 \ar[r] & (\Z/2\Z)^{n-1} \ar[r] & W \ar[r] & \mathfrak{S}_n \ar[r] & 1}$$
\endexample

\vspace{.15in}

Let us fix one choice of simple roots $\Pi=\{\alpha_i\;|\; i=1, \ldots, r\}$, and pick for each $i$ a nontrivial vector $e_i \in \g_{\a_i}$. Set $\h_i=[\g_{\a_i},\g_{-\a_i}] \subset \h$. It is known that $\h_i$ is nonzero (hence of dimension one since $\g_{\a_i}$ and $\g_{-\a_i}$ are both of dimension one), that $\h=\bigoplus_i \h_i$, and that there exists a unique $h_i \in \h_i$ such that $\a_i(h_i)=2$. Let $f_i $ be the unique vector of $\g_{-\a_i}$ satisfying $[e_i,f_i]=h_i$. Observe that since $\a_i-\a_j$ cannot be a root if $i \neq j$, we have $[e_i,f_j]=0$ for $i \neq j$.

Set $a_{ij}=\alpha_i(h_j)$ and put these numbers in a matrix $A=(a_{ij})_{i,j=1}^r$ (the \textit{Cartan matrix} of $\g$). By construction, we have $a_{ii}=2$ for all $i$. One shows that the coefficients of $A$ are integers, and that any off-diagonal entry $a_{ij}, i \neq j$ is nonpositive.

To sum up, we have that $\g$ is generated by $\{e_i, h_i, f_i\;|\; i=1, \ldots , r\}$ and that these satisfy the following set of relations~:
\begin{equation}\label{E:slarelationsI}
\forall\; i,j, \qquad 
\begin{cases}
[h_i,h_j]=0 & \\
[h_i, e_j]=a_{ji} e_j & \\
[h_i, f_j]=-a_{ji} f_j &\\
[e_i,f_j]= \delta_{ij} h_i
\end{cases}
\end{equation}

\vspace{.1in}

It is clear that this set of relations is not enough; for instance the $\{e_j\}$'s should generate $\n_+$ which is finite-dimensional, hence they should satisfy some relation between themselves. So the question arises as to which new relations one should add to get a presentation of $\g$. There is an abstract answer given by~:

\vspace{.15in}

\noindent
\textbf{Proposition~A.9.} \textit{Let $A$ be the Cartan matrix of some simple Lie algebra, and let $\widetilde{\g}$ be the Lie algebra generated by elements $\{e_i, h_i, f_i\;|\; i=1, \ldots , r\}$ subject to the relations (\ref{E:slarelationsI}). Then there exists a unique maximal ideal $I$ of $\widetilde{\g}$ such that $I \cap \h =\{0\}$.}

\vspace{.15in}

As $\g$ is simple, it follows that if $I$ is the ideal in the above Proposition, then $\g \simeq \widetilde{\g}/ I$. In particular, this shows that $\g$ is completely determined by its Cartan matrix $A$.
As discovered by Serre, one can even spell out the extra relations to be added~:

\vspace{.15in}

\noindent
\textbf{Theorem~A.10.}(Serre).\textit{The simple Lie algebra $\g$ is isomorphic to the Lie algebra generated by elements $\{e_i, h_i, f_i\;|\; i=1, \ldots , r\}$ subject to (\ref{E:slarelationsI}) together with
\begin{equation}\label{E:Serrerel}
\forall\; i,j, \qquad 
\begin{cases}
ad^{1-a_{ij}}(e_i)(e_j)=0 & \\
ad^{1-a_{ij}}(f_i)(f_j)=0 & 
\end{cases}
\end{equation}}

\vspace{.1in}

These new relations are called the \textit{Serre relations.} The ony question that now remains is~: which are the possible values for the Cartan matrix $A$ ? This
is answered by the following Theorem~:

\vspace{.1in}

\paragraph{\textbf{Theorem~A.11.}}(Cartan). \textit{An irreducible matrix $A=(a_{ij})_{i,j=1}^r$ with integer coefficients is the Cartan matrix of a simple Lie algebra if and only if the following hold~:}
\begin{equation}\label{E:Cartanmatrix}
a_{ii}=2, \qquad a_{ij} \leq 0\; \text{if}\; i \neq j, \qquad a_{ij}=0 \Rightarrow a_{ji}=0,
\end{equation}
\textit{and $A$ is positive definite}.

\vspace{.1in}

By irreducible matrix, we mean a matrix which cannot be decomposed as a block diagonal matrix in a nontrivial way. The set of all matrices $A$ satisfying the conditions in the above Theorem is well-known. To describe this classification, it is most convenient to ``encode'' the Cartan matrix $A$ into a \textit{Dynkin diagram}. This is done as follows~: to each simple root $\a_i$ is associated a vertex $i$, and two vertices $i$, $j$ are connected by $max\{ |a_{ij}|, |a_{ji}|\}$ edges. In case $|a_{ij}| \neq |a_{ji}|$ the integers $|a_{ij}|$ and $|a_{ji}|$ are written on top of the vertices $i$ and $j$ respectively.

\vspace{.15in}

\paragraph{\textbf{Theorem~A.12.}}(Cartan). \textit{The set of Cartan matrices of simple Lie algebras correspond to the following Dynkin diagrams~:}

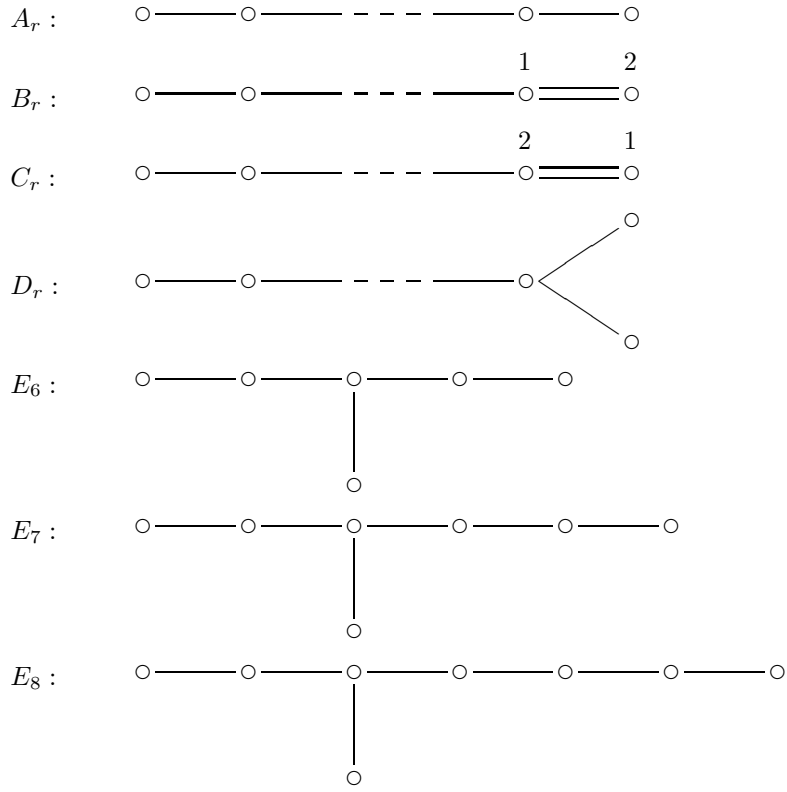
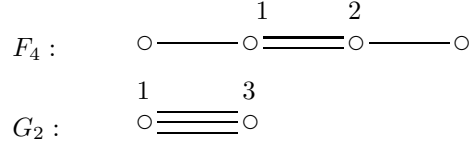
\vspace{.15in}

\centerline{
\begin{picture}(300, 10)
\put(20,-5){$A_r:$}
\put(70,0){\circle{5}}
\put(110,0){\circle{5}}
\put(215,0){\circle{5}}
\put(255,0){\circle{5}}
\put(75,0){\line(1,0){30}}
\put(220,0){\line(1,0){30}}
\put(115,0){\line(1,0){30}}
\put(180,0){\line(1,0){30}}
\put(150,0){\line(1,0){5}}
\put(160,0){\line(1,0){5}}
\put(170,0){\line(1,0){5}}
\end{picture}}

\vspace{.25in}

\centerline{
\begin{picture}(300, 10)
\put(20,-5){$B_r:$}
\put(70,0){\circle{5}}
\put(110,0){\circle{5}}
\put(215,0){\circle{5}}
\put(255,0){\circle{5}}
\put(75,0){\line(1,0){30}}
\put(220,2){\line(1,0){30}}
\put(220,-2){\line(1,0){30}}
\put(212,9){$1$}
\put(252, 9){$2$}
\put(115,0){\line(1,0){30}}
\put(180,0){\line(1,0){30}}
\put(150,0){\line(1,0){5}}
\put(160,0){\line(1,0){5}}
\put(170,0){\line(1,0){5}}
\end{picture}}

\vspace{.25in}

\centerline{
\begin{picture}(300, 10)
\put(20,-5){$C_r:$}
\put(70,0){\circle{5}}
\put(110,0){\circle{5}}
\put(215,0){\circle{5}}
\put(255,0){\circle{5}}
\put(212,9){$2$}
\put(252,9){$1$}
\put(75,0){\line(1,0){30}}
\put(220,2){\line(1,0){30}}
\put(220,-2){\line(1,0){30}}
\put(115,0){\line(1,0){30}}
\put(180,0){\line(1,0){30}}
\put(150,0){\line(1,0){5}}
\put(160,0){\line(1,0){5}}
\put(170,0){\line(1,0){5}}
\end{picture}}

\vspace{.4in}

\centerline{
\begin{picture}(300, 10)
\put(20,-5){$D_r:$}
\put(70,0){\circle{5}}
\put(110,0){\circle{5}}
\put(215,0){\circle{5}}
\put(255,23){\circle{5}}
\put(255,-23){\circle{5}}
\put(75,0){\line(1,0){30}}
\put(220,0){\line(3,2){30}}
\put(220,0){\line(3,-2){30}}
\put(115,0){\line(1,0){30}}
\put(180,0){\line(1,0){30}}
\put(150,0){\line(1,0){5}}
\put(160,0){\line(1,0){5}}
\put(170,0){\line(1,0){5}}
\end{picture}}

\vspace{.35in}

\centerline{
\begin{picture}(300, 10)
\put(20,-5){$E_6:$}
\put(70,0){\circle{5}}
\put(110,0){\circle{5}}
\put(150,0){\circle{5}}
\put(190,0){\circle{5}}
\put(230,0){\circle{5}}
\put(150,-40){\circle{5}}
\put(75,0){\line(1,0){30}}
\put(115,0){\line(1,0){30}}
\put(155,0){\line(1,0){30}}
\put(150,-5){\line(0,-1){30}}
\put(195,0){\line(1,0){30}}
\end{picture}}

\vspace{.6in}

\centerline{
\begin{picture}(300, 10)
\put(20,-5){$E_7:$}
\put(70,0){\circle{5}}
\put(110,0){\circle{5}}
\put(150,0){\circle{5}}
\put(190,0){\circle{5}}
\put(230,0){\circle{5}}
\put(270,0){\circle{5}}
\put(150,-40){\circle{5}}
\put(75,0){\line(1,0){30}}
\put(115,0){\line(1,0){30}}
\put(155,0){\line(1,0){30}}
\put(150,-5){\line(0,-1){30}}
\put(195,0){\line(1,0){30}}
\put(235,0){\line(1,0){30}}
\end{picture}}

\vspace{.6in}

\centerline{
\begin{picture}(300, 10)
\put(20,-5){$E_8:$}
\put(70,0){\circle{5}}
\put(110,0){\circle{5}}
\put(150,0){\circle{5}}
\put(190,0){\circle{5}}
\put(230,0){\circle{5}}
\put(270,0){\circle{5}}
\put(150,-40){\circle{5}}
\put(75,0){\line(1,0){30}}
\put(115,0){\line(1,0){30}}
\put(155,0){\line(1,0){30}}
\put(150,-5){\line(0,-1){30}}
\put(195,0){\line(1,0){30}}
\put(235,0){\line(1,0){30}}
\put(310,0){\circle{5}}
\put(275,0){\line(1,0){30}}
\end{picture}}

\vspace{.75in}

\centerline{
\begin{picture}(300, 10)
\put(20,-5){$F_4:$}
\put(70,0){\circle{5}}
\put(110,0){\circle{5}}
\put(150,0){\circle{5}}
\put(190,0){\circle{5}}
\put(75,0){\line(1,0){30}}
\put(115,-2){\line(1,0){30}}
\put(115,2){\line(1,0){30}}
\put(155,0){\line(1,0){30}}
\put(112,9){$1$}
\put(147,9){$2$}
\end{picture}}

\vspace{.25in}

\centerline{
\begin{picture}(300, 10)
\put(20,-5){$G_2:$}
\put(70,0){\circle{5}}
\put(110,0){\circle{5}}
\put(75,-4){\line(1,0){30}}
\put(75,0){\line(1,0){30}}
\put(75,4){\line(1,0){30}}
\put(67,9){$1$}
\put(107,9){$3$}
\end{picture}}

\vspace{.3in}

The index $r$ in the diagram $X_r$ indicates the number of vertices, (the \textit{rank}  of the corresponding Lie algebra).  The above Dynkin diagrams are sometimes refered to as \textit{finite type Dynkin diagram}. The examples $\mathfrak{sl}_n(\C)$ and $\mathfrak{so}_{2n}(\C)$ are of type $A_{n-1}$ and $D_n$ respectively. The Dynkin diagrams corresponding to symmetric Cartan matrices (i.e. those without any number on top of vertices) are called \textit{simply laced}. They are the types $A_r, D_r, E_6, E_7$ and $E_8$.

\vspace{.3in}

\centerline{\textbf{A.2. Kac-Moody algebras.}}
\addcontentsline{toc}{subsection}{\tocsubsection {}{}{\; A.2. Kac-Moody algebras.}}

\vspace{.15in}

\paragraph{}These are infinite-dimensional generalizations of the simple Lie algebras described in the previous Section, introduced independently by Kac and Moody in the late 60's. A standard reference here is \cite{Kac}. For simplification, we will only be concerned with \textit{symmetric} (simply-laced) Kac-Moody algebras, which are the ones most relevant for us.

\vspace{.1in}

The idea is to \textit{start} with an irreducible, symmetric ``generalized Cartan matrix'' $A=(a_{ij})$ with integer coefficients, which is only subject to the conditions (\ref{E:Cartanmatrix}), namely
$$a_{ii}=2, \qquad a_{ij} \leq 0\; \text{if}\; i \neq j, \qquad a_{ij}=a_{ji}.$$
A \textit{realization} of $A$ consists of a pair of vector spaces $\h$ and $\h^*$ in duality, sets of vectors
$\Pi=\{\a_1, \ldots, \a_r\} \subset \h^*$ and $\Pi^\vee=\{h_1, \ldots, h_r\} \subset \h$ satisfying the following conditions~: both $\Pi$ and $\Pi^\vee$ are linearly independent, ${dim}(\h)=2r-{rank}(A)$, and
$$\alpha_i(h_j)=a_{ji}\qquad \text{for\;all\;} i,j=1, \ldots, r.$$
Any two realizations are isomorphic.

\vspace{.1in}

Let $\widetilde{\g}$ be the Lie algebra generated by elements $\{e_i,f_i, h\;|\; i=1, \ldots ,r, \; h \in \h\}$ subject to relations~:
\begin{equation}\label{E:slarelations2}
\forall\; i,j, \text{and}\; h \in \h, \qquad 
\begin{cases}
[h,h']=0 & \\
[h, e_j]=\a_j(h) e_j & \\
[h, f_j]=-\a_{j}(h) f_j &\\
[e_i,f_j]= \delta_{ij} h_i
\end{cases}
\end{equation}

\vspace{.1in}

Proposition~A.9 is still valid in this context, that is, there exists a unique maximal ideal $I$ of $\widetilde{\g}$ such that $I \cap \h=\{0\}$. The Kac-Moody algebra associated to $A$ is defined to be $\g=\widetilde{\g}/I$. Perhaps more surprisingly, the analogue of Serre's Theorem~A.10 also holds, though the proof is rather difficult~:

\vspace{.1in}

\paragraph{\textbf{Theorem~A.13.}}(Gabber-Kac) \textit{The Kac-Moody Lie algebra $\g$ associated to $A$ is isomorphic to the Lie algebra generated by elements $\{e_i, f_i, h\;|\; i=1, \ldots , r, \; h \in \h\}$ subject to the relations
\begin{equation}\label{E:SerrerelGK}
\begin{split}
&[h,h']=0  \\
&[h, e_j]=\a_{j}(h) e_j  \\
&[h, f_j]=-\a_{j}(h) f_j \\
&[e_i,f_j]= \delta_{ij} h_i\\
&ad^{1-a_{ij}}(e_i)(e_j)=0  \\
&ad^{1-a_{ij}}(f_i)(f_j)=0  
\end{split}
\end{equation}
for all $i,j=1, \ldots, r$ and $h,h' \in \h$.}

\vspace{.15in}

\noindent
\textbf{Remark~A.14.} When $A$ is nondegenerate, $\h=\bigoplus_i \C h_i$, $\h^*=\bigoplus_i \C \a_i$ and relations (\ref{E:slarelations2}) coincide with (\ref{E:slarelationsI}). In general, it is necessary to introduce a realization of $A$ (rather than setting directly $\h=\bigoplus_i \C h_i$) in order to have a reasonable weight theory. 

\vspace{.15in}

Much of the structure of simple Lie algebras can be extended by analogy to the Kac-Moody setting.
The Cartan subalgebra $\h$ still acts semisimply on $\g$, and we define the root system to be
$$\Delta=\{\alpha \in \h^*\;|\; \g_\a \neq \{0\}\}.$$
The set $\Delta$ comes with a canonical polarization $\Delta=\Delta^+ \sqcup  
\Delta^-$ where $\Delta^{\pm}=\Delta \cap \pm \big( \bigoplus\N \a_i \big)$ and $-\Delta^+=\Delta^-$. At the level of Lie algebras, this corresponds to a splitting $\g=\n_- \oplus \h \oplus \n_+$, where $\n_{+}$ (resp. ($\n_-$) is the subalgebra generated by the $e_i$s (resp. by the $f_i$s). The \textit{root lattice} is $Q=\bigoplus_i \Z \a_i \subset \h^*$.

\vspace{.1in}

Next, let $(\;,\;)$ be a nondegenerate bilinear form on $\h$ satisfying 
$( h_i, h')=\a_{i}(h')$ for $i=1, \ldots , r$ (again, there is a unique such form up to isomorphism), and let $(\;,\;)$ be the induced form on $\h^*$. These are analogues of Cartan-Killing forms. By construction, we have
$$(\a_i,\a_j)=a_{ij}, \qquad \text{for}\; i, j=1, \ldots, r.$$
We define the Weyl group $W$ of $\g$ to be the subgroup of ${Aut}(\h^*)$ generated by the simple reflections $s_i$ where
$$s_i~: x \mapsto x-( \a_i,x ) \a_i.$$
The Weyl group preserves the set $\Delta$ and is, by construction, orthogonal with respect to $(\;,\;)$.

\vspace{.25in}

Although, as we have said, Kac-Moody algebras ressemble their finite-dimensional cousins in more than one way, there are still several important differences which we would like to insist upon. Let $\g$ be a Kac-Moody algebra which is not a finite-dimensional simple Lie algebra (i.e. $A$ is not positive definite). Then~:
\begin{enumerate}
\item[i)]  ${dim}\;\g=\infty$ and $W$ is an infinite group,
\item[ii)]  Not every root $\a \in \Delta$ is $W$-conjugate to one of the simple roots $\a_1, \ldots , \a_r$,
\item[iii)] There may be roots $\a \in \Delta$ with $\dim\;\g_{\a} >1$.
\end{enumerate}

\vspace{.1in}

A root belonging to $W \cdot \Pi$ is called \textit{real} and obviously satisfies ${dim}\;\g_{\a}=1$. A root which is not real is called \textit{imaginary}, and may have higher multiplicity. In fact, determining the imaginary roots and their multiplicities of a given Kac-Moody algebra is a central, difficult, and very much open problem in the theory. It is known that $\a \in \Delta$ is real if $( \a,\a )=2$ and imaginary if $( \a,\a ) \leq 0$.

\vspace{.1in}

There is a last fundamental difference between simple Lie algebras and Kac-Moody Lie algebras, which will is crucial for us in these notes. Call a subalgebra $\bo=\h \oplus \n$ of $\g$ a \textit{Borel subalgebra} if $\n=\bigoplus_{\a \in \Psi_+} \g_{\a}$ and $\Psi_+$ is a ``half'' of $\Delta$, i.e. $\Delta=\Psi_+ \cup -\Psi_+$ while $\Psi_+ \cap -\Psi_+ = \emptyset$.
\begin{enumerate}
\item[iv)] There may be Borel subalgebras $\mathfrak{b}$ not conjugate to either of the standard Borel subalgebras $\bo_+=\h \oplus \n_+$ or $\bo_-=\h \oplus \n_-$.
\end{enumerate}

\vspace{.15in}

\addtocounter{theo}{1}

\paragraph{\textbf{Example~A.15.}} Let us describe what is probably the most important class of (symmetric) Kac-Moody algebras besides the simple Lie algebras. These are the so-called \textit{affine Lie algebras}, or \textit{affine Kac-Moody Lie algebras}, and correspond to generalized Cartan matrices $A$ which are positive semi-definite and have corank equal to one. One shows that the associated Dynkin diagrams (the simply laced \textit{affine Dynkin diagrams}) are as follows~:

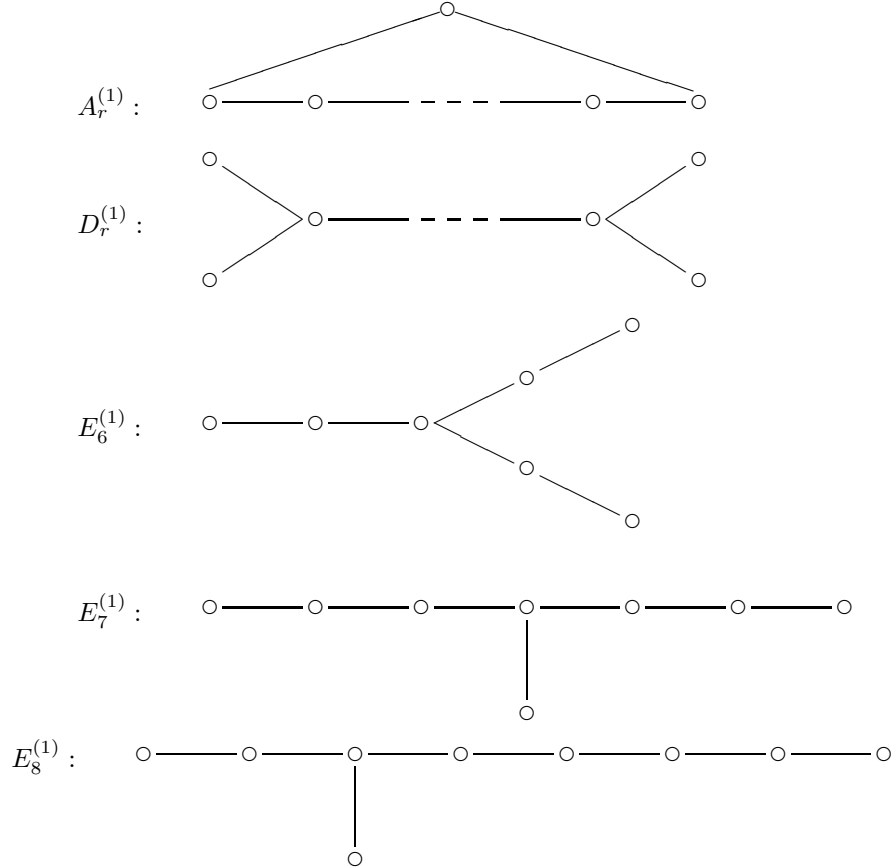
\vspace{.45in}

\centerline{
\begin{picture}(300, 10)
\put(20,-5){$A^{(1)}_r:$}
\put(160,35){\circle{5}}
\put(70,5){\line(3,1){87}}
\put(163,34){\line(3,-1){90}}
\put(70,0){\circle{5}}
\put(110,0){\circle{5}}
\put(215,0){\circle{5}}
\put(255,0){\circle{5}}
\put(75,0){\line(1,0){30}}
\put(220,0){\line(1,0){30}}
\put(115,0){\line(1,0){30}}
\put(180,0){\line(1,0){30}}
\put(150,0){\line(1,0){5}}
\put(160,0){\line(1,0){5}}
\put(170,0){\line(1,0){5}}
\end{picture}}

\vspace{.45in}

\centerline{
\begin{picture}(300, 10)
\put(20,-5){$D^{(1)}_r:$}
\put(70,23){\circle{5}}
\put(70,-23){\circle{5}}
\put(75,20){\line(3,-2){30}}
\put(75,-20){\line(3,2){30}}
\put(110,0){\circle{5}}
\put(215,0){\circle{5}}
\put(255,23){\circle{5}}
\put(255,-23){\circle{5}}
\put(220,0){\line(3,2){30}}
\put(220,0){\line(3,-2){30}}
\put(115,0){\line(1,0){30}}
\put(180,0){\line(1,0){30}}
\put(150,0){\line(1,0){5}}
\put(160,0){\line(1,0){5}}
\put(170,0){\line(1,0){5}}
\end{picture}}

\vspace{.9in}

\centerline{
\begin{picture}(300, 10)
\put(20,-5){$E^{(1)}_6:$}
\put(70,0){\circle{5}}
\put(110,0){\circle{5}}
\put(150,0){\circle{5}}
\put(190,17){\circle{5}}
\put(190,-17){\circle{5}}
\put(230,37){\circle{5}}
\put(230,-37){\circle{5}}
\put(75,0){\line(1,0){30}}
\put(115,0){\line(1,0){30}}
\put(155,0){\line(2,1){30}}
\put(155,0){\line(2,-1){30}}
\put(195,20){\line(2,1){30}}
\put(195,-20){\line(2,-1){30}}
\end{picture}}

\vspace{.8in}

\centerline{
\begin{picture}(300, 10)
\put(20,-5){$E^{(1)}_7:$}
\put(70,0){\circle{5}}
\put(110,0){\circle{5}}
\put(150,0){\circle{5}}
\put(190,0){\circle{5}}
\put(230,0){\circle{5}}
\put(270,0){\circle{5}}
\put(190,-40){\circle{5}}
\put(75,0){\line(1,0){30}}
\put(115,0){\line(1,0){30}}
\put(155,0){\line(1,0){30}}
\put(190,-5){\line(0,-1){30}}
\put(195,0){\line(1,0){30}}
\put(235,0){\line(1,0){30}}
\put(310,0){\circle{5}}
\put(275,0){\line(1,0){30}}
\end{picture}}

\vspace{.6in}

\centerline{
\begin{picture}(350, 10)
\put(20,-5){$E^{(1)}_8:$}
\put(70,0){\circle{5}}
\put(110,0){\circle{5}}
\put(150,0){\circle{5}}
\put(190,0){\circle{5}}
\put(230,0){\circle{5}}
\put(270,0){\circle{5}}
\put(150,-40){\circle{5}}
\put(75,0){\line(1,0){30}}
\put(115,0){\line(1,0){30}}
\put(155,0){\line(1,0){30}}
\put(150,-5){\line(0,-1){30}}
\put(195,0){\line(1,0){30}}
\put(235,0){\line(1,0){30}}
\put(310,0){\circle{5}}
\put(275,0){\line(1,0){30}}
\put(350,0){\circle{5}}
\put(315,0){\line(1,0){30}}
\end{picture}}

\vspace{.75in}

Each of the above extended Dynkin diagram is obtained by adding one node to a (simply laced) Dynkin diagram of finite type. Traditionally, the simple root corresponding to this new node is called $\a_0$, and the new generators are called $e_0, h_0$ and $f_0$. 
Let $X_n^{(1)}$ be one of the diagrams above. It is customary to denote the corresponding affine Lie algebra by $\widehat{\g}$, where $\g$ is the simple Lie algebra of type $X_n$. A special feature of this Kac-Moody Lie algebra is that it has a very concrete realization in terms of the \textit{loop algebra} $L\g=\g\otimes \C[t,t^{-1}]$ of $\g$. We first consider the universal central extension
$$\xymatrix{
0 \ar[r] & \C c \ar[r] & \widehat{\g}' \ar[r] & L\g \ar[r] & 0,}$$
of $L\g$. As a vector space, we have $\widehat{\g}' =(\g\otimes \C[t,t^{-1}]) \oplus \C c$ and the Lie brackets are given by
$$
\begin{cases}
&c\;\text{is\;central\;},\\
&[x\otimes t^n,y\otimes  t^m]=[x,y]\otimes t^{n+m} + \delta_{n,-m}(x,y)nc.
\end{cases}
$$
The affine Lie algebra is obtained as a further one-dimensional extension 
$$\widehat{\g}=L\g \oplus \C c \oplus \C \partial,$$
where $[\partial, c]=0$ and $\partial$ acts on a loop $x\otimes t^n$ via the derivation $t \frac{\partial}{\partial t}$, i.e.
$$[\partial, x\otimes t^n]=nx\otimes t^n.$$

\vspace{.1in}

The Cartan subalgebra of $\widehat{\g}$ is $\widehat{\h}=\h \oplus \C c \oplus \C \partial$ and its dual is $\widehat{\h}^*=\h^* \oplus \C \delta + \C \Lambda_0$, where
$$\h^{\perp}=\C \delta \oplus \C \Lambda_0, \qquad \langle \delta, \partial \rangle=\langle \Lambda_0, c \rangle =1, \qquad \langle \Lambda_0, \partial \rangle=\langle \delta, c \rangle =0.$$
With these notations, we have $\a_0=\delta-\theta$.

The expression for the generators in this presentation are 
$$e_i, f_i, h_i \mapsto e_i, f_i, h_i \in \g \subset L\g \; \text{if}\; i \geq 1, $$
$$e_0=f_{\theta} \otimes t, \qquad f_0=e_{\theta}\otimes t^{-1}, \qquad h_0= c-h_{\theta},$$
where $e_{\theta} \in \g_{\theta}$ is a vector of highest weight, $f_{\theta} \in \g_{-\theta}$ is such that $(e_{\theta},f_{\theta})=1$, and $h_{\theta}=[e_{\theta},f_{\theta}]$. 

\vspace{.1in}

Finally, let us look at the root system of $\widehat{\g}$.  Using the formula for $e_{0}, f_0$ and the expression of $\a_0$ in terms of $\delta$, it is easy to see that the weight of an element $x\otimes t^n$ with 
$x \in \g_{\a}$ is $\a + n\delta$. Hence, the root system is
$$\widehat{\Delta}=\{\Delta \oplus \Z \delta\} \cup \Z^* \delta,$$
where $\Delta$ is the root system of $\g$. We have
$$\widehat{\g}_{\a+n\delta}=\g_{\a} \otimes t^n, \qquad \widehat{\g}_{n\delta}=\h \otimes t^n.$$
In particular, the weights $n\delta$ for $n \neq 0$ are of multiplicity ${dim}\;\h={rank}(\g)$. These roots are precisely the imaginary roots of $\widehat{\g}$.

\vspace{.35in}

We illustrate this for $\g=\mathfrak{sl}_2$ by giving a picture of the root system~:

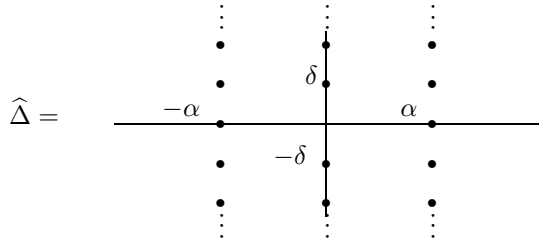
\centerline{
\begin{picture}(400, 60)
\put(60,0){$\widehat{\Delta}=$}
\multiput(140,-30)(0,15){5}{\circle*{3}}
\multiput(139,-45)(0,4){3}{$\cdot$}
\multiput(139,34)(0,4){3}{$\cdot$}
\multiput(179,-45)(0,4){3}{$\cdot$}
\multiput(179,34)(0,4){3}{$\cdot$}
\multiput(219,-45)(0,4){3}{$\cdot$}
\multiput(219,34)(0,4){3}{$\cdot$}
\multiput(220,-30)(0,15){5}{\circle*{3}}
\multiput(180,-30)(0,15){2}{\circle*{3}}
\multiput(180,15)(0,15){2}{\circle*{3}}
\put(118,3){$-\alpha$}
\put(208,3){$\alpha$}
\put(172,15){$\delta$}
\put(160,-15){$-\delta$}
\put(180,-35){\line(0,1){70}}
\put(100,0){\line(1,0){160}}
\end{picture}}

\vspace{.9in}

As an example supporting point iv) of the discussion concerning the differences between simple and Kac-Moody Lie algebras, consider the standard set of positive roots $\Delta^+$, which lie inside
 the (infinite) surrounded region below
 
 
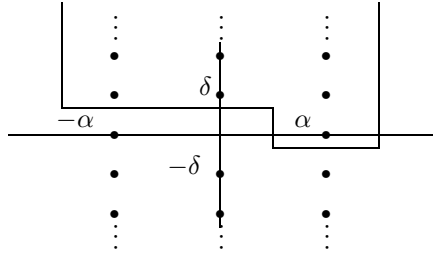
\vspace{.1in}
 
\centerline{
\begin{picture}(400, 60)
\multiput(140,-30)(0,15){5}{\circle*{3}}
\multiput(139,-45)(0,4){3}{$\cdot$}
\multiput(139,34)(0,4){3}{$\cdot$}
\multiput(179,-45)(0,4){3}{$\cdot$}
\multiput(179,34)(0,4){3}{$\cdot$}
\multiput(219,-45)(0,4){3}{$\cdot$}
\multiput(219,34)(0,4){3}{$\cdot$}
\multiput(220,-30)(0,15){5}{\circle*{3}}
\multiput(180,-30)(0,15){2}{\circle*{3}}
\multiput(180,15)(0,15){2}{\circle*{3}}
\put(118,3){$-\alpha$}
\put(208,3){$\alpha$}
\put(172,15){$\delta$}
\put(160,-15){$-\delta$}
\put(180,-35){\line(0,1){70}}
\put(100,0){\line(1,0){160}}
\put(120,10){\line(0,1){40}}
\put(120,10){\line(1,0){80}}
\put(200,10){\line(0,-1){15}}
\put(200,-5){\line(1,0){40}}
\put(240,-5){\line(0,1){55}}
\end{picture}}

\vspace{.9in}

\noindent
and which corresponds to the standard Borel subalgebra 
$$\widehat{\bo}_+=(\mathfrak{sl}_2 \otimes t \C[t]) \oplus \C e_{\a} \oplus \widehat{\h}; $$
and the ``nonstandard'' positive roots system $\Phi_+$ lying at the right of the line below

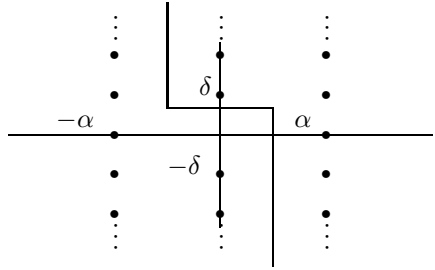
\vspace{.1in}

\centerline{
\begin{picture}(400, 60)
\multiput(140,-30)(0,15){5}{\circle*{3}}
\multiput(139,-45)(0,4){3}{$\cdot$}
\multiput(139,34)(0,4){3}{$\cdot$}
\multiput(179,-45)(0,4){3}{$\cdot$}
\multiput(179,34)(0,4){3}{$\cdot$}
\multiput(219,-45)(0,4){3}{$\cdot$}
\multiput(219,34)(0,4){3}{$\cdot$}
\multiput(220,-30)(0,15){5}{\circle*{3}}
\multiput(180,-30)(0,15){2}{\circle*{3}}
\multiput(180,15)(0,15){2}{\circle*{3}}
\put(118,3){$-\alpha$}
\put(208,3){$\alpha$}
\put(172,15){$\delta$}
\put(160,-15){$-\delta$}
\put(180,-35){\line(0,1){70}}
\put(100,0){\line(1,0){160}}
\put(160,10){\line(0,1){40}}
\put(160,10){\line(1,0){40}}
\put(200,10){\line(0,-1){60}}
\end{picture}}

\vspace{.9in}

\noindent
which corresponds to a nonstandard Borel subalgebra 
$$\widehat{\bo}=(e_{\a} \otimes \C[t,t^{-1}]) \oplus (\h \otimes t \C[t]) \oplus \widehat{\h}.$$
It is clear that $\widehat{\bo}$ is not conugate to $\widehat{\bo}_{\pm}$ under the affine Weyl group $\widehat{W}$ since $\widehat{\bo}_{\pm}$ is finitely generated whereas (by weight considerations) $\widehat{\bo}$ is not.

\vspace{.2in}
 
\centerline{\textbf{A.3. Enveloping algebras.}}
\addcontentsline{toc}{subsection}{\tocsubsection {}{}{\; A.3. Enveloping algebras.}}

\vspace{.15in}

\paragraph{}From any associative algebra $A$ we may obtain a Lie algebra $(A,[\;,\;])$ by simply setting $[x,y]=xy-yx$ for $x, y \in A$. The enveloping Lie algebra construction goes in the other direction~: it associates to any Lie algebra $\g$ an associative algebra $\U(\g)$ which has the ``same'' representation theory.

Let $\g$ be a Lie algebra. Recall that  representation of $\g$ is a pair $(V,\rho)$ where $V$ is a vector space and $\rho: \g \to ({End}(V), [\;,\;])$ is a Lie algebra map. The \textit{universal enveloping algebra} of $\g$ is the associative algebra $\U(\g)$ equipped with a Lie algebra map $i:\g \hookrightarrow (\U(\g),[\;,\;])$, solution to the following (universal) problem~: for any associative algebra $A$ and any Lie algebra map $\rho: \g \to (A,[\;,\;])$ there exists a \textit{unique}
morphism of associative algebras $ \psi: \U(\g) \to A$ making the diagram
$$\xymatrix{
\g \ar[r]^-{\forall\; \phi} \ar[d]_-{i} & A\\
\U(\g) \ar@{-->}[ru]_-{\exists \,! \,\psi} &}$$
commutative. By definition, any representation of $\g$ extends to a representation (as associative algebra) of $\U(\g)$. Vice versa, any representation of $\U(\g)$ restricts to a representation of $\g$ so that $\g$ and $\U(\g)$ indeed have the same representations.

\vspace{.15in}

In addition to the universal property, there is also an explicit construction of the enveloping algebra $\U(\g)$~: it is the quotient of the tensor algebra
$$T(\g) =\C \oplus \bigoplus_{l \geq 1} \g ^{\otimes l}$$
by the two-sided ideal generated by elements $xy-yx-[x,y]$ for $x,y \in \g$. In other words, it is the 
unital associative algebra generated by elements of $\g$ subject only to the relations 
\begin{equation}\label{E:univenvelop}
xy-yx=[x,y]\qquad \text{for\;} x,y \in \g.
\end{equation}

Let us fix a basis $x_1, x_2, \ldots$ of $\g$. Intuitively, relations (\ref{E:univenvelop}) allow to ``reorder'' any product of the $x_i$'s so that all $x_1$ appear before all $x_2$, and so on. This is made precise by the following result, known as the Poincar\'e-Birkhoff-Witt (or simply PBW) Theorem~:

\vspace{.1in}

\paragraph{\textbf{Theorem~A.16.}}(PBW Theorem). \textit{Let $\g$ be a Lie algebra and let $x_1, x_2, \ldots$ be a basis of $\g$. Then the set
$$\{x_{i_1} x_{i_2} \cdots x_{i_r}\;|\; r \geq 0, \;i_1 \leq i_2 \leq \cdots \leq i_r\}$$
forms a basis of $\U(\g)$.}

\vspace{.15in}

This shows that $\U(\g)$ is always infinite-dimensional, and in fact has the same size as the symmetric algebra on $\g$. 

\vspace{.15in}

Finally, any enveloping algebra has a canonical structure of a cocommutative Hopf algebra, completely determined by the conditions
$$\Delta(x)=x \otimes 1 + 1 \otimes x, \qquad S(x)=-x, \quad \epsilon(x)=0$$
for $x \in \g$.

\vspace{.2in}

Let us now assume that $\g$ is the Kac-Moody algebra associated to a symmetric generalized Cartan matrix $A=(a_{ij})_{i,j=1}^r$, and let $(\h,\h^*,\Pi,\Pi^{\vee}\}$ be a realization of $A$ as in Section~A.2.
By the Gabber-Kac theorem, $\g$ admits a presentation by generators $\{e_i, f_i, h\;|\; i=1, \ldots, r, h \in \h\}$ and relations (\ref{E:SerrerelGK}). Then the \textit{same} generators and relations give a presentation of $\U(\g)$, this time viewed as an associative algebra. 

\vspace{.1in}

Let $\h, \n_{+}, \n_-$ be the Cartan subalgebra and positive, resp. negative, nilpotent subalgebras, so that there is a decomposition $\g=\n_- \oplus \h \oplus \n_+$. We may consider the enveloping algebras
$\U(\h), \U(\n_+), \U(\n_-)$ as subalgebras of $\U(\g)$. Then the PBW theorem implies that the multiplication map gives an isomorphism of vector spaces $\U(\n_-) \otimes \U(\h) \otimes \U(\n_+) \simeq \U(\g)$. Of course, the same is also true for any nonstandard splitting $\g=\n'_- \oplus \h \oplus \n'_+$ corresponding to a nonstandard set $\Phi_+$ of positive roots, and in fact more generally for any splitting of $\g$ as a direct sum of subalgebras.

\vspace{.2in}

Finally, the adjoint action of $\h$ on $\g$ naturally extends to $\U(\g)$, and there is again a weight decomposition 
$$\U(\g) = \bigoplus_{\a \in Q} \U(\g)_{\a},$$
where the weights this time belong to $Q=\bigoplus_i \Z \a_i$.
The enveloping algebras of the standard Borel or nilpotent subalgebras have finite-dimensional weight spaces.

\vspace{.2in}

\centerline{\textbf{A.4. Quantum Kac-Moody algebras.}} 
\addcontentsline{toc}{subsection}{\tocsubsection {}{}{\; A.4. Quantum Kac-Moody algebras.}}

\vspace{.15in}

\paragraph{}The \textit{quantized enveloping algebra} of a Kac-Moody algebra $\g$ is a certain deformation, as a Hopf algebra, of $\U(\g)$. Let again $A=(a_{ij})_{i,j=1}^r$ be a symmetric generalized Cartan matrix with realization $(\h,\h^*,\Pi,\Pi^{\vee})$.
Let $v$ be a formal variable, and $\C(v)$ the field of rational functions in $v$. The book \cite{Jantzen} is a fine reference here. We will use the following standard notation~:
$$[n]=\frac{v^n-v^{-n}}{v-v^{-1}}= v^{-n} + v^{2-n} + \cdots + v^{n-2} + v^n, $$
$$ [n]!=[n] \cdot [n-1] \cdots [2] \cdot [1],$$
and
$$\begin{bmatrix} t\\ r \end{bmatrix}=\frac{[t]!}{[r]!\cdot [t-r]!}=\frac{[t] \cdot [t-1] \cdots [t-r+1]}{[2] \cdots [r-1][r]}.$$

\vspace{.15in}

Define $\U_v(\g)$ to be Hopf algebra over $\C(v)$ generated by elements $E_i, F_i$ for $i=1, \ldots, r$ and $v^h$, for $h \in \h$ subject to the relations
\begin{equation}\label{E:quantumSerrerel}
\begin{split}
&v^h  v^{h'}=v^{h+h'}  \\
&v^h E_j v^{-h}=v^{\a_{j}(h)} E_j  \\
&v^h F_j v^{-h}=v^{-\a_{j}(h)} F_j \\
&[E_i,F_j]= \delta_{ij} \frac{v^{h_i}-v^{-h_i}}{v -v^{-1}}\\
&\sum_{l=0}^{1-a_{ij}} (-1)^l \begin{bmatrix} 1-a_{ij} \\ l \end{bmatrix} E_i^l E_j E_i^{1-a_{ij}-l}=0  \\
&\sum_{l=0}^{1-a_{ij}} (-1)^l \begin{bmatrix} 1-a_{ij} \\ l \end{bmatrix} F_i^l F_j F_i^{1-a_{ij}-l}=0  
\end{split}
\end{equation}
for all $i,j=1, \ldots, r$ and $h, h' \in \h$.

\vspace{.05in}

\noindent
The coproduct is given by the formulas
$$\Delta(v^h)=v^h \otimes v^h, \quad \Delta(E_i)=E_i \otimes 1 + v^{h_i} \otimes E_i, \quad \Delta(F_i)=1\otimes F_i + F_i \otimes v^{-h_i} ,$$
and the antipode is
$$S(v^{h_i})=v^{-h_i}, \qquad S(E_i)=-v^{-h_i}E_i, \qquad S(F_i)=-F_iv^{h_i}.$$

\vspace{.1in}

\noindent
The relations (\ref{E:quantumSerrerel}) are a deformation of those in (\ref{E:SerrerelGK}), and these can be obtained back by setting $v=1$. The last two relations of (\ref{E:quantumSerrerel}) are known as the \textit{quantum Serre relations}. Giving the generators $E_i, F_i$ and $v^h$ the weights $\a_i, -\a_i$ and $0$ respectively defines on $\U_v(\g)$ a grading by the root lattice $Q$.

\vspace{.1in}

Let $\n_{\pm}, \bo_{\pm}$ be the standard nilpotent subalgebras of $\g$. To these correspond
the subalgebras $\U_v(\n_{\pm})$  of $\U_v(\g)$ which are generated by
$\{E_i\;|\; i =1, \ldots, r\}$ or $\{F_i\;|\; i =1, \ldots, r\}$. Similarly, to the standard Borel subalgebras $\bo_{\pm}$ correspond the quantized algebras $\U_v(\bo_{\pm})$ generated by $\{E_i, v^h\;|\; i =1, \ldots, r, \;h \in \h\}$ or $\{F_i, v^h\;|\; i =1, \ldots, r, \;h \in \h\}$. Note that $\U_v(\bo_{\pm})$ are Hopf subalgebras (i.e. stable under the comultiplication) but $\U_v(\n_{\pm})$ are not.

\vspace{.2in}

The enveloping algebra $\U(\g)$ is not commutative, but always cocommutative. Its dual Hopf algebra (in any appropriate sense) is thus commutative, but not cocommutative. The quantized enveloping algebra $\U_v(\g)$ morally stands somewhere in between these two~: it is neither commutative nor cocommutative. Moreover, as the following important Theorem of Drinfeld states, it is self-dual~:

\vspace{.1in}

\paragraph{\textbf{Theorem~A.18.}} \textit{There exists a unique Hopf pairing on $\U_v(\bo_+)$ such that}
$$ (E_i,E_j)=\frac{\delta_{ij}}{v^2-1}, \qquad (v^h,v^{h'})=v^{(h,h')}, \qquad (E_i, v^h)=0.$$
\textit{It is nondegenerate.}

\vspace{.2in}

It is clear that $\U_v(\bo_+)$ and $\U_v(\bo_-)$ are isomorphic as algebras. Hence $\U_v(\bo_+)$ and $\U_v(\bo_-)$ are put in duality with each other under the Drinfeld pairing above. This indeed means that, in the appropriate sense, $\U_v(\g)$ is self-dual. We refer to the book \cite{Jantzen} for details. 

\vspace{.1in}

Since $\U_v(\g)$ is defined over the field $\C(v)$, it is not possible to specialize directly $v$ to any complex parameter $\nu \in \C^*$. To remedy this, we introduce following Lusztig \cite{Lu1} the \textit{integral form} $\U_v^{res}(\g)$ of $\U_v(\g)$ which is defined as the $\C[v,v^{-1}]$-subalgebra
of $\U_v(\g)$ generated by elements $v^h$ for $h \in \h$ and the divided powers
$$E_i^{(n)}=\frac{E_i^n}{[n]!}, \qquad F_i^{(n)}=\frac{F_i^n}{[n]!}$$
for $i=1, \ldots, r$ and $n \geq 1$. By a theorem of Lusztig, $\U^{res}_v(\g)$ is a $\C[v,v^{-1}]$-Hopf subalgebra of $\U_v(\g)$, free as a $\C[v,v^{-1}]$-module, and we have $\U_v(\g)=\U^{res}_v(\g) \otimes \C(v)$. This allows us to define for any $\epsilon \in \C^*$ a specialization 
$$\U_\epsilon(\g):=\U^{res}_v(\g)_{|v=\epsilon}$$
which is a Hopf algebra defined over $\C$.
The integral forms $\U_v^{res}(\n_{\pm})$ and $\U_v^{res}(\bo_{\pm})$ or the specializations $\U_\epsilon(\n_{\pm})$ and $\U_\epsilon(\bo_{\pm})$ are defined in an analogous manner.

\vspace{.2in}

To finish, we need to mention a slight variant of the quantized enveloping algebra $\U_v(\g)$ of a Kac-Moody algebra which often appears, and which is relevant to the Hall algebra business. Let $\U_v(\g')$ be the Hopf subalgebra of $\U_v(\g)$ generated by $\{E_i, F_i, v^{\pm h_i}\;|\; i=1, \ldots, r\}$. Setting, as is customary, $K_i^{\pm 1}= v^{\pm h_i}$, we may write the defining relations for $\U_v(\g')$ as follows~:

\vspace{.05in}

\begin{equation}\label{E:defquantumgroups}
\begin{split}
&K_i K_j=K_j K_i  \\
&K_i E_j K_i^{-1}=v^{a_{ij}} E_j  \\
&K_iF_j K_i^{-1}=v^{-a_{ij}} F_j \\
&[E_i,F_j]= \delta_{ij} \frac{K_i-K_i^{-1}}{v -v^{-1}}\\
&\sum_{l=0}^{1-a_{ij}} (-1)^l \begin{bmatrix} 1-a_{ij} \\ l \end{bmatrix} E_i^l E_j E_i^{1-a_{ij}-l}=0  \\
&\sum_{l=0}^{1-a_{ij}} (-1)^l \begin{bmatrix} 1-a_{ij} \\ l \end{bmatrix} F_i^l F_j F_i^{1-a_{ij}-l}=0  
\end{split}
\end{equation}

The coproduct and antipode read
\begin{equation}\label{E:standardcoproduct}
\Delta(K_i)=K_i \otimes K_i, \quad \Delta(E_i)=E_i \otimes 1 + K_i \otimes E_i, \quad \Delta(F_i)=1\otimes F_i + F_i \otimes K_i^{-1} ,
\end{equation}
$$S(K_i)=K_i^{-1}, \qquad S(E_i)=-K_i^{-1}E_i, \qquad S(F_i)=-F_iK_i.$$

\vspace{.1in}

The algebra $\U_v(\g')$ is a deformation of the \textit{derived subalgebra} $\g'=[\g,\g]=
\n_- \oplus \bigoplus_i \C h_i \oplus \n_+$, and is isomorphic to $\U_v(\g)$ when $A$ is nondegenerate. It is also useful to consider the subalgebras $\U_v(\bo'_{\pm})$ of $\U_v(\g')$ generated by $\{E_i, K_i^{\pm 1}\}$ or $\{F_i, k_i^{\pm 1}\}$.

\vspace{.2in}

\centerline{\textbf{A.5. Loop algebras of Kac-Moody algebras.}}
\addcontentsline{toc}{subsection}{\tocsubsection {}{}{\; A.5. Loop algebras of Kac-Moody algebras.}}

\vspace{.15in}

\paragraph{}As we have seen in Example~A.15., the affine Lie algebras $\widehat{\g}$, i.e. the first class
of Kac-Moody algebras beyond the simple Lie algebras, may be realized as the universal central extension of the loop algebra $L\g=\g \otimes \C[t,t^{-1}]$ of a simple Lie algebra $\g$. This suggests two possible generalizations.

\vspace{.05in}

The first one is to take an arbitrary commutative $\C$-algebra $A$, a simple (simply laced) Lie algebra $\g$, and to consider the universal central extension $\overline{\g_A}$ of $\g_{A}=\g \otimes A$. According to a theorem of Kassel \cite{Kassel},
$\overline{\g_A}=\g_A \oplus \Omega^1A/dA$ where $\Omega^1A$ is the space of differentials on $A$ and $dA \subset\Omega^1A$ is the space of exact forms. The bracket is given by
$$\begin{cases}
&\Omega^1A/dA \; \text{is\;central} \\
&[x \otimes a,y\otimes b]=[x,y]\otimes ab+ (x,y) (da)b 
\end{cases}$$

When $A=\C[t_1^{\pm 1}, \ldots, t_n^{\pm 1}]$ one gets the $n$-toroidal algebra studied by Billig \cite{Billig}, Moody \cite{MRY} and Rao \cite{Rao} (see \cite{Raosurvey} for a survey).

\vspace{.05in}

The second possible generalization is actually the one that will be relevant to these lectures. We let $\g$ be an arbitrary Kac-Moody algebra and we consider its loop algebra $L\g=\g \otimes \C[t,t^{-1}]$. It turns out that it is not the universal central extension of $L\g$ which we need here, but something slightly different. We use as a motivation a Proposition of Garland~:

\vspace{.1in}

\paragraph{\textbf{Proposition~A.19.}}(Garland, \cite{Garland}). \textit{Let $\g$ be a simple Lie algebra. The affine Lie algebra $\widehat{\g}$ is isomorphic to the Lia algebra generated by $\{e_{i,k},f_{i,k},h_{i,k}\}$ for $i \in I$ and $k \in \Z$, and $c$ subject to the following relations~:
\begin{equation}\label{E:defaffinelie}
\begin{split}
&[h_{i,k}, h_{j,l}]=k \delta_{k,-l}a_{ij}{c},\\
& [e_{i,k},f_{j,l}]=\delta_{i,j}h_{i,k+l}+ k \delta_{k,-l} {c},\\
&[h_{i,k}, e_{j,l}]=a_{ij}e_{j,l+k},\\
& [h_{i,k}, f_{j,l}]=-a_{ij}f_{j,k+l},\\
&[e_{i,k+1},e_{j,l}]=[e_{i,k},e_{j,l+1}],\\
&[f_{i,k+1},f_{j,l}]=[f_{i,k},f_{j,l+1}],\\
&[e_{i,k_1},[e_{i,k_2},[\ldots [e_{i,k_n},e_{j,l}]\cdots ]=0 \qquad 
\mathrm{if\;} n=1-a_{ij},\\
&[f_{i,k_1},[f_{i,k_2},[\ldots [f_{i,k_n},f_{j,l}]\cdots ]=0 \qquad 
\mathrm{if\;} n=1-a_{ij}.
\end{split}
\end{equation}}

\vspace{.1in}

The isomorphism with $\widehat{\g}$ is induced by the assignement $e_{i,k} \mapsto e_i \otimes t^k, f_{i,k} \mapsto f_i \otimes t^k, h_{i,k} \mapsto h_i\otimes t^k, c \mapsto c$.

\vspace{.1in}

We now \textit{define} the loop algebra $\mathcal{L}\g$ of a \textit{Kac-Moody} algebra $\g$ to be the Lie algebra generated by $\{e_{i,k},f_{i,k},h_{i,k}\}$ for $i \in I$ and $k \in \Z$, and $c$ subject to relations (\ref{E:defaffinelie}). Of course, $\{a_{ij}\}$ are now the coefficients of the generalized Cartan matrix $A$ of $\g$.
The assignement $e_{i,k} \mapsto e_i \otimes t^k, f_{i,k} \mapsto f_i \otimes t^k, h_{i,k} \mapsto h_i\otimes t^k, c \mapsto c$ gives a surjective morphism of Lie algebras 
$\phi:\mathcal{L}\g \to \g \otimes \C[t,t^{-1}] \oplus \C c$, but this is usually not injective anymore.

\vspace{.1in}

The algebra $\mathcal{L}\g$ is naturally $\widehat{Q}=Q \oplus \Z\delta$-graded, where $Q$ is the root lattice of $\g$~: the elements $e_{i,k},f_{i,k},h_{i,k}$ the degrees $\a_i+k\delta, -\a_i+k\delta$ and $k\delta$ respectively. It is shown in \cite{MRY} that 
$${Ker}\;\phi \subset \bigoplus_{\underset{\a \in \Delta^{im}}{k \in \Z}} \mathcal{L}\g[\a+k\delta].$$
Hence the root system $\widehat{\Delta}$ of $\mathcal{L}\g$ is equal to that of $\g \otimes \C[t,t^{-1}]$~:
$$\widehat{\Delta}=\{\a+ l\delta\;|\; \a \in \Delta, l \in \Z\} \cup \{l\delta\;|\; l \in \Z^*\}$$
(but the multiplicities differ in general).

\vspace{.15in}

\addtocounter{theo}{1}

\paragraph{\textbf{Example~A.20.}} If $\g$ is a simple Lie algebra then of course $\mathcal{L}\g=\widehat{\g}$. \endexample

\vspace{.15in}

\addtocounter{theo}{1}

\paragraph{\textbf{Example~A.21.}} Let us now assume that $\g=\widehat{\g}_0$ is itself an affine Lie algebra. Thus we are considering the loop algebra of another loop algebra (of a simple Lie algebra)~! As shown in \cite{MRY}, one obtains in this way a \text{double-loop} or \textit{elliptic Lie algebra} with universal central extension~:
$$\mathcal{L}\g=\mathcal{E}{\g_0}:=\g_{0}\otimes \C[s^{\pm1 }, t^{\pm 1}] \oplus \mathbf{K},$$
where 
\begin{equation*}
\begin{split}
\mathbf{K}=&\Omega_1\C[s^{\pm 1}, t^{\pm 1}]/d\C[s^{\pm 1}, t^{\pm 1}]\\
=&\bigoplus_{(m,n) \neq (0,0)} \C k_{m,n} \oplus \C t^{-1}dt \oplus \C s^{-1}ds
\end{split}
\end{equation*}
with
$$k_{m,n}=\begin{cases} \frac{1}{n} s^{m-1}t^n ds& \text{if}\; n \neq 0,\\
-\frac{1}{m}s^mt^{n-1}dt & \text{if}\;m \neq 0
\end{cases}$$
Explicitly, the bracket reads
$$[x\otimes s^mt^n,y \otimes s^{m'}t^{n'}]=[x,y]s^{m+m'}t^{n+n'} + (x,y) \left| \begin{matrix} m&m'\\n&n'\end{matrix}\right| k_{m+m',n+n'}$$
if $(m,n) \neq -(m',n')$, 
$$[x\otimes s^mt^n,y \otimes s^{-m}t^{-n}]=[x,y] + (x,y)\big( nt^{-1}dt + ms^{-1}ds\big)$$
and $\mathbf{K}$ is central. The root system of $\mathcal{E}{\g_0}$ is
$$\widehat{\Delta}=\{\a+l\delta_0+l' \delta\;|\; \a \in \Delta_0, \;l, l' \in \Z\} \cup \{l \delta_0+l'\delta\;|\; (l,l') \neq (0,0)\}$$
where $\Delta_0$ is the root system of $\g_0$. Observe that the center $\mathbf{K}$ is not concentrated in degree zero~: we have $deg(s^{-1}ds)=deg(t^{-1}dt)=0$ but $deg(k_{m,n})=m\delta_0+n\delta$. It is useful to visualize the root system $\widehat{\Delta}$ as a $\Z^2$-lattice

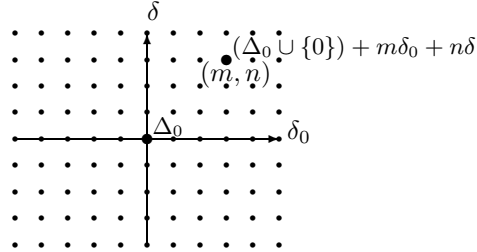
\vspace{.1in}

\centerline{
\begin{picture}(100,100)
\multiput(0,5)(10,0){11}{\circle*{2}}
\multiput(0,15)(10,0){11}{\circle*{2}}
\multiput(0,25)(10,0){11}{\circle*{2}}
\multiput(0,35)(10,0){11}{\circle*{2}}
\multiput(0,45)(10,0){11}{\circle*{2}}
\multiput(0,55)(10,0){11}{\circle*{2}}
\multiput(0,65)(10,0){11}{\circle*{2}}
\multiput(0,75)(10,0){11}{\circle*{2}}
\multiput(0,85)(10,0){11}{\circle*{2}}
\put(52,47){\small{$\Delta_0$}}
\put(50,45){\circle*{3.5}}
\put(80,75){\circle*{3.5}}
\put(82,77){\small{$(\Delta_0\cup\{0\})+m\delta_0+n\delta$}}
\put(70,67){$\small{(m,n)}$}
\put(0,45){\vector(1,0){100}}
\put(50,90){$\delta$}
\put(103,45){$\delta_0$}
\put(50,5){\vector(0,1){80}}
\end{picture}}

\vspace{.2in}

Over the vertex $(0,0)$ sits the finite root system $\Delta_0$ while over any nonzero vertex $(m,n)$ sits 
$(\Delta_0\cup \{0\})+m\delta_0+n\delta$ (to account for the imaginary weight space $\h_0 \otimes s^mt^n$).

The above picture makes the following remark obvious.

\vspace{.1in}

\paragraph{\textbf{Proposition~A.22.}}}\textit{The group $SL(2,\Z)$ acts by Lie algebra automorphisms on $\mathcal{E}{\g_0}$}.

\vspace{.15in}

Notice that for $\g=\widehat{\g}_0$ the two ``affinization'' procedures (i.e. the one leading to toroidal algebras and the one leading to loop Kac-Moody algebras) give the same result. The elliptic Lie algebras $\mathcal{E}{\g_0}$ and their root systems were first conceived by K. Saito in connection to the theory of elliptic singularities \cite{Saito}. \endexample

\vspace{.2in}

Now that we have defined the loop algebras $\mathcal{L}\g$ we will construct inside them certain suitable ``positive'' Borel subalgebras $\mathcal{L}\bo_+$ along with the corresponding sets of ``positive'' roots $\widehat{\Delta}_+ \subset \widehat{\Delta}$. Of course, since $\mathcal{L}\g$ is not a Kac-Moody algebra unless $\g$ is a simple Lie algebra, there is no \textit{a priori} notion of a standard Borel subalgebra or standard set of positive roots.
In these lectures we will only be concerned with a certain type of loop algebras (and the definition of $\mathcal{L}\bo_+$  will only make sense for these). Let us assume that $\g$ is a Kac-Moody algebra associated to a star-shaped Dynkin diagram

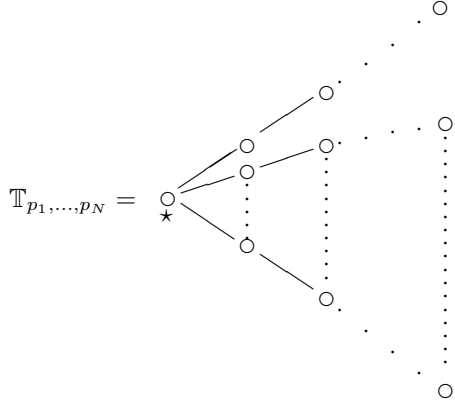
\vspace{.3in}

\centerline{
\begin{picture}(160,60)
\put(-10,-3){$\mathbb{T}_{p_1, \ldots, p_N}=$}
\put(50,0){\circle{5}}
\put(80,20){\circle{5}}
\put(110,40){\circle{5}}
\put(153,72){\circle{5}}
\put(80,10){\circle{5}}
\put(110,20){\circle{5}}
\put(155,28){\circle{5}}
\put(80,-18){\circle{5}}
\put(110,-38){\circle{5}}
\put(155,-73){\circle{5}}
\put(54,4){\line(3,2){22}}
\put(85,25){\line(3,2){20}}
\put(115,44){\circle*{1}}
\put(125,51){\circle*{1}}
\put(135,58){\circle*{1}}
\put(145,67){\circle*{1}}
\put(55,2){\line(3,1){20}}
\put(85,12){\line(3,1){20}}
\put(115,21){\circle*{1}}
\put(125,23){\circle*{1}}
\put(135,25){\circle*{1}}
\put(145,26.5){\circle*{1}}
\put(55,-2){\line(3,-2){20}}
\put(85,-21){\line(3,-2){20}}
\put(115,-42){\circle*{1}}
\put(125,-50){\circle*{1}}
\put(135,-58){\circle*{1}}
\put(145,-66){\circle*{1}}
\multiput(80,5)(0,-5){5}{\circle*{1}}
\multiput(110,15)(0,-5){10}{\circle*{1}}
\multiput(155,23)(0,-5){18}{\circle*{1}}
\put(47,-9){$\star$}
\end{picture}}

\vspace{1.2in}

\noindent
where $p_1, p_2, \ldots, p_N$ indicate the lengths of the branches (including the central vertex). A quick glance at the tables in Appendices~A.1. and A.2. shows that all such diagrams are star-shaped, while among (simply laced) affine Dynkin diagrams only $D_4^{(1)}=\mathbb{T}_{2,2,2,2},\; E_6^{(1)}=\mathbb{T}_{3,3,3}, \;E_7^{(1)}=\mathbb{T}_{4,4,2},\; E_8^{(1)}=\mathbb{T}_{6,3,2}$ are. 

Let us label the vertices of $\mathbb{T}_{p_1, \ldots, p_N}$ by calling $\star$ the central vertex and $(i,l)$ the $l-1$st vertex of the $i$th branch, so that $(i,1)$ is adjacent to $\star$ for all $i$. The corresponding simple roots of $\g$ are then $\a_{\star}$ and $\a_{(i,l)}$. We have $$\widehat{Q}=\bigoplus_{i,l}\Z\a_{(i,l)} \oplus \Z\a_{\star} \oplus \Z\delta.$$
We will say that a nonzero weight $\lambda=\sum_{i,l} c_{i,l}\a_{(i,l)} + c_{\star}\a_{\star}+c_{\delta}\delta$ is \textit{positive} if $c_{\star} >0$ or $c_{\star}=0$ and $c_{\delta}>0$, or $c_{\star}=c_{\delta}=0$ and $c_{i,l} \geq 0$ for all $(i,l)$. We simply set $\widehat{\Delta}^+=\widehat{\Delta} \cap \widehat{Q}^+$ where $\widehat{Q}^+$ is the set of positive weights. Finally, we define the \textit{positive Borel subalgebra} of $\mathcal{L}\g$ to be
$$\mathcal{L}\bo_+=\h \oplus \bigoplus_{\a \in \widehat{\Delta}^+} \mathcal{L}\g[\a],$$
and the \textit{positive nilpotent subalgebra} of $\mathcal{L}\g$ to be
$$\mathcal{L}\n_+= \bigoplus_{\a \in \widehat{\Delta}^+} \mathcal{L}\g[\a].$$

\vspace{.15in}

\addtocounter{theo}{1}

\paragraph{\textbf{Example~A.23.}} Let us first look at the simplest example $\mathbb{T}_1=A_1$, i.e. $\g=\mathfrak{sl}_2$. Thus $\mathcal{L}\g=\widehat{\mathfrak{sl}}_2$ and $\widehat{Q}^+=\Z_{\a_{\star}} \oplus \Z\delta$. In that case
$$\widehat{Q}^+=\{n\a_{\star}+m\delta\;|\; n>0\;\text{or}\; n=0, m>0\},$$
$$\widehat{\Delta}^+=\{\a+m\delta\;|\;m \in \Z\} \cup \{m\delta\;|\; m >0\}.$$
Thus $\mathcal{L}\bo_+$ is the \textit{nonstandard} Borel subalgebra considered in Example~A.15.

\endexample

\vspace{.15in}

\addtocounter{theo}{1}

\paragraph{\textbf{Example~A.24.}} Let $\g$ be an arbitrary simple Lie algebra so that $\mathcal{L}\g=\widehat{\g}$. To the branches of $\mathbb{T}_{p_1, \ldots, p_N}$ \textit{excluding the central vertex} correspond subalgebras $\g_1, \ldots, \g_N$, with $\g_i \simeq \mathfrak{sl}_{p_i}$.
There is a natural embedding $\widehat{\g}_i=\mathcal{L}\g_i \subset \mathcal{L}\g$ and we let $\widehat{\n}_+^i $ sand for the \textit{standard} positive nilpotent subalgebra of $\widehat{\g}_i$. Next, let $\n_+ \subset \g$ be the standard positive nilpotent subalgebra and let $\n^{\star}\subset \n_+$ be the Lie ideal generated by $e_{\star}$ (this is also the nilpotent radical of the maximal parabolic subalgebra associated to $\star$). Then
$$\mathcal{L}\bo_+=\h \oplus\bigoplus_i \widehat{\n}_+^i \oplus (\n^*\otimes \C[t,t^{-1}]),$$ 
and
$$\mathcal{L}\n_+=\bigoplus_i \widehat{\n}_+^i\oplus (\n^*\otimes \C[t,t^{-1}]) .$$ 
Hence in this case $\mathcal{L}\bo_+$ is some kind of hybrid between the standard Borel subalgebra and the algebra $\bo_+\otimes \C[t,t^{-1}]$ of loops in the finite Borel subalgebra. 

\endexample

\vspace{.1in}

The explicit description of $\mathcal{L}\bo_+$ when $\g$ is affine, i.e. when $\mathcal{L}\g$ is elliptic, is more complicated and we prefer to leave it to the intrepid reader.

\vspace{.2in}

\centerline{\textbf{A.6. Quantum loop algebras.}}
\addcontentsline{toc}{subsection}{\tocsubsection {}{}{\; A.6. Quantum loop algebras.}}

\vspace{.15in}

\paragraph{} These are deformations of the enveloping algebras of the loop algebras $\mathcal{L}\g$ associated to a Kac-Moody algebra $\g$. Again, the starting point is a new presentation of the quantum affine algebra $U_v(\widehat{\g}_0)$, where $\g_0$ is a simple Lie algebra. This will be a deformation of Garland's relations (\ref{E:defaffinelie}).

\vspace{.1in}

\paragraph{\textbf{Proposition~A.25.}}(Drinfeld, \cite{Drinfeld}, Beck, \cite{Beck}). \textit{ Let $\g_0$ be a simple, simply laced Lie algebra with Cartan matrix $A=(a_{ij})$. The quantum group $\U_v(\widehat{\g}_0)$ is isomorphic to the $\C(v)$-algebra generated by $E_{i,l}, F_{i,l}$, for $i \in I, l \in \Z$ and $K^{\pm 1}_{i}, H_{i,n}$ for $i \in I, n \in \Z^*$ and $C^{\pm 1/2}$ subject to the following set of relations~:
$$C^{1/2} \;\text{is\;central},$$
$$ [K_i,K_j]=[K_i,H_{j,l}]=0,$$
$$[H_{i,l},H_{j,k}]=\delta_{i,j}\delta_{l,-k}\frac{[2l]}{l} \frac{C^l-C^{-l}}{v-v^{-1}},$$
$$K_i E_{jk}K_{i}^{-1}=v^{a_{ij}}E_{jk}, \qquad K_i F_{jk}K_{i}^{-1}=v^{-a_{ij}}F_{jk},$$
$$[H_{i,l},E_{j,k}]= \frac{1}{l}[la_{ij}]C^{-|l|/2}E_{j,k+l}, $$
\begin{equation}\label{E:drinfeldpres}
 [H_{i,l},F_{j,k}]= -\frac{1}{l}[la_{ij}]C^{|l|/2}F_{j,k+l}
 \end{equation}
$$E_{i,k+1}E_{j,l}-v^{a_{ij}}E_{j,l}E_{i,k+1}=
v^{a_{ij}}E_{i,k}E_{j,l+1}-E_{j,l+1}E_{i,k},$$
$$F_{i,k+1}F_{j,l}-v^{-a_{ij}}F_{j,l}F_{i,k+1}=
v^{-a_{ij}}F_{i,k}F_{j,l+1}-F_{j,l+1}F_{i,k},$$
$$[E_{i,k},F_{j,l}]=\delta_{ij}
\frac{C^{(k-l)/2}\psi_{i,k+l}-C^{(l-k)/2}\varphi_{i,k+l}}{v-v^{-1}}$$
\begin{equation*}
\begin{split}
\mathrm{For\;}i \neq j&\;\mathrm{and\;}n=1-a_{ij},\\
&\mathrm{Sym}_{k_1, \ldots, k_n} \sum_{t=0}^{1-a_{ij}}(-1)^t 
\begin{bmatrix} n \\ t \end{bmatrix} E_{i,k_1} \cdots E_{i,k_t}E_{j,l}
E_{i,k_{t+1}} \cdots E_{i,k_n}=0\\
&\mathrm{Sym}_{k_1, \ldots, k_n} \sum_{t=0}^{1-a_{ij}}(-1)^t 
\begin{bmatrix} n \\ t \end{bmatrix} F_{i,k_1} \cdots F_{i,k_t}E_{j,l}
F_{i,k_{t+1}} \cdots F_{i,k_n}=0
\end{split}
\end{equation*}
where $\mathrm{Sym}_{k_1, \ldots, k_n}$ denotes symmetrization with respect
to the indices $k_1, \ldots , k_n$, and where $\psi_{i,k}$ and $\varphi_{i,k}$
are defined by the following equations :
$$\sum_{k \geq 0} \psi_{i,k}u^k=K_i \mathrm{exp}\bigg( (v-v^{-1})
\sum_{k=1}^\infty H_{i,k}u^k\bigg),$$ 
$$\sum_{k \geq 0} \varphi_{i,k}u^k=K_i^{-1} \mathrm{exp}\bigg( -(v-v^{-1})
\sum_{k=1}^\infty H_{i,-k}u^{-k}\bigg).$$}

\vspace{.15in}

Note that the above gives a presentation of $\U_v(\widehat{\g}_0)$ (known as \textit{Drinfeld's new realization}) as an \textit{algebra} only. To describe explicitly the coproduct in terms of these generators seems a rather difficult problem (see \cite{Hubery3}). Nevertheless, it is possible to define a \textit{new} coproduct on $\U_v(\widehat{\g}_0)$ (known as \textit{Drinfeld's new coproduct}), which is much more suited to the above presentation as the old one (see \cite{CP} or \cite{Hernandez} for applications to representation theory of quantum affine algebras).

\vspace{.1in}

\paragraph{\textbf{Proposition~A.26.}}(Drinfeld).\textit{ Introduce other elements $\zeta_{i,l}^{\pm}, \theta_{i,l}^{\pm}$ for $l \geq 1, i \in I$ via the formal relations~:
$$1+\sum_{l \geq 1} \theta_{i,l}^{\pm}s^{\pm l}= K_i^{\pm 1} exp\bigg(\pm (v^{-1}-v) \sum_{r \geq 1} H_{i,\pm r}s^{\pm r}\bigg),$$
$$1+\sum_{l \geq 1} \zeta_{i,l}^{\pm}s^{\pm l}= K_i^{\pm 1} exp\bigg( \sum_{r \geq 1} \frac{H_{i,\pm r}}{[r]}C^{\pm r/2} s^{\pm r}\bigg).$$
The sets $\{K_i^{\pm 1},\zeta_{i,l}\;|\; i \in I, l \in \Z^*\}$ or $\{K_i^{\pm 1},\theta_{i,l}\;|\; i \in I, l \in \Z^*\}$ generate the same subalgebra of $\U_v(\widehat{\g}_0)$ as $\{K_i^{\pm 1}, H_{i,l}\;|\; i \in I, l \in \Z^*\}$.
Moreover, the following formulas define on $\U_v(\widehat{\g}_0)$ the structure of a topological bialgebra~:
$$\Delta(\zeta_{i,l}^{\pm})=\sum_{s=0}^l \zeta_{i,l-s}^{\pm} C^{\pm s/2}\otimes \zeta_{i,s}^{\pm},$$
\begin{equation}\label{E:copp1}
\Delta(E_{i,l})=E_{i,l} \otimes 1 + \sum_{t \geq 0} \theta^+_{i,t}C^{(l-t)/2}\otimes E_{i,l-t},
\end{equation}
$$\Delta(F_{i,l})= 1 \otimes F_{i,l}  + \sum_{t \geq 0} F_{i,l-t} \otimes \theta^-_{i,t}C^{(t-l)/2}.$$}

\vspace{.15in}

Now let $\g$ be an arbitrary (simply laced) Kac-Moody algebra. We define the quantum loop algebra $\U_v(\mathcal{L}\g)$ to be the $\C(v)$-algebra generated by $\{E_{i,l},F_{i,l}, H_{i,k}\}$ for $ i \in I, l \in \Z, k \in \Z^*$ and $\{K_i, C\}$ for $i \in I$ subject to the relation (\ref{E:drinfeldpres}) above. If $\nu \in \C^*$ satisfies $|\nu| \neq 1$ then we also define $\U_{\nu}(\mathcal{L}\g)$ to be the $\C$-algebra with same generators and relations but with $\nu$ instead of $v$. The formulas (\ref{E:copp1}) for the coproduct define on $\U_v(\mathcal{L}\g)$ and $\U_{\nu}(\mathcal{L}\g)$ the structure of topological bialgebras.

\vspace{.2in}

Finally, let us assume that the Dynkin diagram of $\g$ is star-shaped as in the end of Appendix~A.5.
We will be interested in certain subalgebras $\U_{v}(\mathcal{L}\bo_+)$ and $\U_{v}(\mathcal{L}\n_+)$ of $\U_{v}(\mathcal{L}\g)$, which are deformations of the subalgebras $\U(\mathcal{L}\bo_+)$ and $\U(\mathcal{L}\n_+)$ which we defined previously.  For any $i=1, \ldots, N$, the elements $\{E_{(i,l),r}, F_{(i,l),s}, H_{(i,l),t}, K^{\pm 1}_{i}\;|\; l=1, \ldots, p_{i-1} \}$ corresponding to the $i$th branch generate a subalgebra isomorphic to $\U_{v}(\widehat{\mathfrak{sl}}_{p_i})$. Let $\U_{v}(\widehat{\mathfrak{b}}^+_{i})$ (resp. $\U_{v}(\widehat{\mathfrak{n}}^+_{i})$ denote the standard positive Borel subalgebra (resp. the standard positive nilpotent subalgebra) of $\U_{v}(\widehat{\mathfrak{sl}}_{p_i})$. We define $\U_{v}(\mathcal{L}\bo_+)$ to be the subalgebra generated by $E_{\star, n}$ for $n \in \Z$, $H_{\star,r}$ for $r \geq 1$, and by $\U_{v}(\widehat{\mathfrak{n}}^+_{p_i})$ for $i=1, \ldots N$. Similarly, $\U_v(\mathcal{L}\n_+)$ is the subalgebra generated by $E_{\star, n}$ for $n \in \Z$, $H_{\star,r}$ for $r \geq 1$, and by $\U_{v}(\widehat{\mathfrak{n}}^+_{p_i})$ for $i=1, \ldots N$. We do \textit{not} claim that
$\U_{v}(\mathcal{L}\bo_+)$ is stable under the coproduct $\Delta$.

\vspace{.15in}

\addtocounter{theo}{1}

\paragraph{\textbf{Example~A.27.}}  Let us spell out for the reader's convenience the simplest example, namely $\g=\mathfrak{sl}_2$. The Drinfeld generators are $F_l,E_l,H_n$ for $l \in \Z, n \in \Z^*$ and $K^{\pm 1}, C$. The subalgebra $\U_v(\mathcal{L}\bo_+)$ is generated by $E_l, H_n$ for $l \in \Z$ and $n \geq 1$, and $K^{\pm 1}, C^{\pm 1/2}$. The presentation for $\U_v(\mathcal{L}\bo_+)$ is as follows~:
\begin{equation}\label{E:Drinfeldpressl2}
\begin{cases}
&C^{1/2} \;\text{is\;central},\\
& [K_i,H_{n}]=[H_n,H_l]=0,\\
&K E_{k}K^{-1}=v^{2}E_{k},\\
&[H_{l},E_{k}]= \frac{1}{l}[2l]C^{-|l|/2}E_{k+l}, \\
&E_{k+1}E_{l}-v^{2}E_{l}E_{k+1}=
v^{2}E_{k}E_{l+1}-E_{l+1}E_{k}.
\end{cases}
\end{equation}

In this situation, $\U_v(\mathcal{L}\bo_+)$ \textit{is} stable under the coproduct and we have
\begin{equation}\label{E:copp2}
\Delta(\zeta_{l}^{\pm})=\sum_{s=0}^l \zeta_{l-s}^{\pm} C^{\pm s/2}\otimes \zeta_{s}^{\pm},
\end{equation}
\begin{equation}\label{E:copp3}
\Delta(E_{l})=E_{l} \otimes 1 + \sum_{t \geq 0} \theta^+_{t}C^{(l-t)/2}\otimes E_{l-t},
\end{equation}
where $\zeta_{l}^{\pm}, \theta_{l}^{\pm}$ are defined via the formal relations~:
$$1+\sum_{l \geq 1} \theta_{l}^{\pm}s^{\pm l}= K^{\pm 1} exp\bigg(\pm (v^{-1}-v) \sum_{r \geq 1} H_{\pm r}s^{\pm r}\bigg),$$
$$1+\sum_{l \geq 1} \zeta_{l}^{\pm}s^{\pm l}= K^{\pm 1} exp\bigg( \sum_{r \geq 1} \frac{H_{\pm r}}{[r]}C^{\pm r/2} s^{\pm r}\bigg).$$
\endexample

\vspace{.15in}

\addtocounter{theo}{1}

\paragraph{\textbf{Remark~A.28.}} When $\g$ is a simple Lie algebra, i.e. in Drinfeld's original context, Chari and Pressley found the right integral form $\U_v^{res}(\mathcal{L}\g)$ in terms of the Drinfeld generators. This allows one to specialize the parameter $v$ to any nonzero value $\nu$, including the interesting case of roots of unity. However, for us, $\nu$ will always satisfy $|\nu|>1$ and we may use the more straightforward specialization as defined above.

\newpage

\centerline{\textbf{Acknowledgements}}

\vspace{.1in}

I am grateful to the organizers of the summer school held in Grenoble in 2008, for their invitation and encouragements.
This survey was initially prepared for the proceedings of the winter school ``Representation theory of finite dimensional algebras and related topics'' held at ICTP Trieste in 2006. I would like to extend my thanks to the organizers of that school as well, and to all the participants for their helpful and fruitful remarks. Thanks are due to the Flying Dog brewery for its crucial help in overcoming the various difficulties met at various stages of this work. Finally, I am grateful to I. Burban, X.-W. Chen, F. Fauquent-Millet, D. Fratila  and V. Toledano-Laredo for pointing out several inaccuracies in an earlier version of these notes.

\vspace{.2in}

\small{}

\vspace{.2in}

\noindent
Olivier Schiffmann,\\
Institut Math\'ematique de Jussieu,\\
175 rue du Chevaleret, 75013 Paris \\
FRANCE,\\
email:\;\texttt{olive@math.jussieu.fr}
\end{document}